\RequirePackage{fix-cm}
\documentclass[11pt,a4paper,final]{article}
\usepackage[left=2.00cm, right=2.00cm, top=1.00cm, bottom=2.00cm]{geometry}
\newcommand{\keywords}[1]{\textbf{keywords}-- #1}
\usepackage{graphicx}

\usepackage{moreverb}
\usepackage{graphicx}
\usepackage{amssymb}
\usepackage{amsmath}
\usepackage{nicefrac} 
\usepackage{dsfont}

\usepackage{mathtools} 

\usepackage[colorlinks,bookmarksopen,bookmarksnumbered,citecolor=red,urlcolor=red]{hyperref}

\usepackage{color} 
\usepackage{psfrag}
\usepackage{lscape}

\usepackage{tabularx}
\usepackage{multirow}
 
\usepackage{pgf,tikz,pgfplots}
\pgfplotsset{compat=1.15}
\usepackage{mathrsfs}

\usetikzlibrary{decorations.markings}
\tikzset{->-/.style={postaction={decorate,decoration={
  markings,
  mark=at position 0.5  with {\arrow{>}}
	}}
	}
	}
	
\usetikzlibrary{arrows}
\definecolor{blublu}{rgb}{0,0,1}
\definecolor{blu2}{rgb}{0.30196078431372547,0.30196078431372547,1}
\definecolor{azzurrino}{rgb}{0.9,0.9,1}

\definecolor{verdeverde}{rgb}{0,0.6,0}
\definecolor{verdino}{rgb}{0.9,1,0.9}

\definecolor{rossorosso}{rgb}{1,0,0}
\definecolor{rossino}{rgb}{1,0.9,0.9}

\definecolor{grigiochiaro}{rgb}{0.9,0.9,0.9}
\definecolor{grigiochiaro2}{rgb}{0.85,0.85,0.85}
\definecolor{grigioscuro}{rgb}{0.65,0.65,0.65}
\tikzstyle{tratteggiato}= [dash pattern={on 5pt off 5pt}]

\newcommand{\Vfaces}{V_{\text{faces}}}
\newcommand{\Vxfaces}{V_{\text{$x$-faces}}}
\newcommand{\Vyfaces}{V_{\text{$y$-faces}}}
\newcommand{\Vzfaces}{V_{\text{$z$-faces}}}
\newcommand{\Vedges}{V_{\text{edges}}}
\newcommand{\Vxedges}{V_{\text{$x$-edges}}}
\newcommand{\Vyedges}{V_{\text{$y$-edges}}}
\newcommand{\Vzedges}{V_{\text{$z$-edges}}}

\newcommand{\Vnodes}{V_{\text{nodes}}}
\newcommand{\Vmain}{V}
\newcommand{\Vdual}{V^*}

\newcommand{\Avex}[5]{\left\{\langle   #1 \rangle_{\phantom{.}_{\!\!\!\!x}} \right\}_{#2,#3,#4}^{#5}}
\newcommand{\Avey}[5]{\left\{\langle   #1 \rangle_{\phantom{.}_{\!\!\!\!y}} \right\}_{#2,#3,#4}^{#5}}
\newcommand{\Avez}[5]{\left\{\langle   #1 \rangle_{\phantom{.}_{\!\!\!\!z}} \right\}_{#2,#3,#4}^{#5}}
\newcommand{\Avexy}[5]{\left\{\langle  #1 \rangle_{\phantom{.}_{\!\!\!\!xy}} \right\}_{#2,#3,#4}^{#5}}
\newcommand{\Avexz}[5]{\left\{\langle  #1 \rangle_{\phantom{.}_{\!\!\!\!xz}} \right\}_{#2,#3,#4}^{#5}}

\newcommand{\Aveyz}[5]{\left\{\langle  #1 \rangle_{\phantom{.}_{\!\!\!\!yz}} \right\}_{#2,#3,#4}^{#5}}

\newcommand{\Ave}[6]{\left\{\langle   #1 \rangle_{{#6}} \right\}_{#2,#3,#4}^{#5}}
\newcommand{\AveTwo}[8]{\Ave{\text{\small$\langle$}#1\text{\small$\rangle$}_{\!\text{$#7$}} #2}{#3}{#4}{#5}{#6}{#8}} 

\newcommand{\AveyxTwo}[6]{\Avex{\text{\small$\langle$} #1\text{\small$\rangle$}_{\!\text{\tiny$y$}} #2}{#3}{#4}{#5}{#6}}

\newcommand{\FSa}{\mathcal{A}}
\newcommand{\FSb}{\mathcal{B}}
\newcommand{\FSc}{\mathcal{C}}
\newcommand{\FSd}{\mathcal{D}}

\newcommand{\x}{\mathbf{x}} 
 
\newcommand{\A}{\mathbf{A}}

\newcommand{\V}{\mathbf{V}} 
 
\renewcommand{\v}{\mathbf{v}} 
\renewcommand{\u}{\mathbf{u}} 

\newcommand{\one}{\mathbb{I}} 
\newcommand{\Real}{I\!\!R}  
\newcommand{\CFL}{\texttt{CFL}} 
\newcommand{\ddelta}{\text{$\Delta$}\!\!\!\!\text{\hspace{-0.4pt}\small{$\Delta$}}}

\newcommand{\myequation}[3]{\begin{equation} \begin{array}{#1} #2 \end{array} #3 \end{equation}}

\newcommand{\halb}{{\frac{1}{2}}}
\newcommand{\viertel}{{\frac{1}{4}}}
\newcommand{\p}{{\mathbf{p}}}
\newcommand{\Q}{{\mathbf{Q}}}
\newcommand{\f}{{\mathbf{f}}}
\newcommand{\g}{{\mathbf{g}}}
\newcommand{\h}{{\mathbf{h}}}
\newcommand{\fv}{{\mathbf{f}^v}}
\newcommand{\gv}{{\mathbf{g}^v}}
\newcommand{\hv}{{\mathbf{h}^v}}
\newcommand{\fp}{{\mathbf{f}^p}}
\newcommand{\gp}{{\mathbf{g}^p}}
\newcommand{\hp}{{\mathbf{h}^p}}
\newcommand{\fb}{{\mathbf{f}^b}}
\newcommand{\gb}{{\mathbf{g}^b}}
\newcommand{\hb}{{\mathbf{h}^b}}
\newcommand{\B}{{\mathbf{B}}}
\newcommand{\E}{{\mathbf{E}}}
\newcommand{\F}{\mathbf{F}}
\newcommand{\Fv}{\mathbf{F}_{v}}
\newcommand{\Fd}{\mathbf{F}_{d}}
\newcommand{\Fdv}{\mathbf{F}_{\mu}}
\newcommand{\Fdb}{\mathbf{F}_{\eta}}

\newcommand{\Fp}{\mathbf{F}_{p}}
\newcommand{\Fb}{\mathbf{F}_{b}}
\renewcommand{\A}{{\mathbf{A}}}

\newcommand{\redP}[1]{#1} 
\newcommand{\blue}[1]{{\color{blue}#1}} 

\renewcommand{\div}{\text{div}}

\def \MyFigFolder{.}
\def \MyFigFormat{}

\usepackage{tikz} 
\usepackage{pgfplots}

\usepackage{tikz} 
\usepackage{pgfplots}

\pgfplotsset{compat=newest} 
\usepgfplotslibrary{external} 
\tikzset{external/up to date check={md5}}
\tikzset{external/mode=convert with system call} 
\tikzsetexternalprefix{build/}

\begin{document}

\title{A novel structure preserving semi-implicit finite volume method for viscous and resistive magnetohydrodynamics
}


\author{Francesco Fambri\\ Max Planck Institute for Plasma Physics, \\
              Boltzmannstr. 2, Garching, 85748, Germany \\ 
              \href{mailto:francesco.fambri@ipp.mpg.de}{\tt francesco.fambri@ipp.mpg.de} }





\date{Dec 9,  2020} 

\maketitle

\begin{abstract}

In this work we introduce a novel semi-implicit structure-preserving finite-volume/finite-difference scheme for the viscous and resistive equations of magneto-hydrodynamics (VRMHD) based on an appropriate 3-split of the governing PDE system, which is decomposed into a first convective subsystem, a second subsystem involving the coupling of the velocity field with the magnetic field and a third subsystem involving the pressure-velocity coupling.
The nonlinear convective terms are discretized explicitly, while the remaining two subsystems accounting for the Alfvén waves and the magneto-acoustic waves are treated implicitly.  Thanks to this, the final algorithm is at least formally constrained only by a mild CFL stability condition depending on the velocity field of the pure hydrodynamic convection.
To preserve the divergence-free constraint of the magnetic field exactly at the discrete level, a proper set of overlapping dual (staggered) meshes is employed. 
The resulting linear algebraic systems are shown to be symmetric and therefore can be solved by means of an efficient standard matrix-free conjugate gradient algorithm.
The final scheme can be regarded as a novel shock-capturing, conservative and structure preserving semi-implicit scheme for VRMHD. 
Several numerical tests are presented to show the main features of our novel divergence-free semi-implicit FV/FD solver: linear-stability in the sense of Lyapunov is verified at a prescribed constant equilibrium solution; a second-order of convergence is obtained for a smooth time-dependent solution; shock-capturing capabilities are proven against a standard set of stringent MHD shock-problems; accuracy and robustness are verified against a non-trivial set of two- and three-dimensional MHD problems.

\end{abstract}


\section{Introduction}
\label{intro}

While there are many single- or multi-fluid models for simulating different plasma-flow conditions, eventually by using a reduced or augmented  PDE system, this manuscript is devoted to the development of a novel structure-preserving semi-implicit finite-volume (FV) scheme for the nonlinear equations of viscous and resistive magnetohydrodynamics (MHD or VRMHD).  

Many interesting  physical systems can be described by the equations of magnetohydrodynamics, in which resistivity effects of electromagnetic fields are also important. Important applications of scientific and public interest can be found in the physics of solar flares, the plasma flow in the magnetosphere of neutron stars, as well as the study of inertial or magnetic confinement fusion for the calibration or design of fusion reactors for civil energy production, but also for the simulation and design of plasma thrusters used for the propulsion of small spacecraft or satellites, to mention a few examples \redP{\cite{Biskamp1993,Goedbloed2019,Sankaran2002}}.

The physical interactions involved in the dynamics of a plasma can operate at different timescales that can be mathematically described by the characteristics of the original PDE system. 
These are usually distinguished, independently on the fluid velocity, in terms of the so-called \emph{low}- and \emph{fast}- magnetosonic and the Alfvén speeds, $c_l$, $c_f$ and $c_a$. 
While the Alfvén and the low magnetosonic waves may propagate only in the direction parallel to the magnetic field, only the fast magnetosonics may cause a compression of the field-lines.
 Whenever solving a hyperbolic PDE system numerically, the computational time-step is severely constraint by the magnitude of the fastest characteristic speeds, e.g. the fast magnetosonics $c_f$ for MHD, together with the characteristic mesh-size. 
This constraint may become particularly severe, e.g., in low-$\beta$ Tokamak scenarios, as well as fluid-flow in the low-Mach or incompressible regime in gas-dynamics. 
Moreover, in those systems dominated by low frequency and long wavelength flows as well as in presence of important spatial gradients, the corresponding admissible computational time-step of a standard fully-explicit algorithm would make a long-time simulation economically unsustainable, and eventually also inaccurate due to the excessive accumulation of numerical dissipation and truncation error effects. 
For this reason, implicit algorithms have been for a long time  considered as an essential tool for the investigation of stabilized plasma configurations in magnetic confinement fusion reactors \redP{\cite{AYDEMIR1985,HARNED1986,NIMROD1999,SOVINEC2004,Chacon2008}.}
On the other hand, due to the nonlinearity of the MHD equations, a fully implicit solver may also lead to  
a highly nonlinear algebraic systems for the simultaneous solution of  a large 
number of coupled nonlinear equations which, to be solved, can be particularly demanding in terms of computational costs.
 For this reason, semi-implicit methods have been also investigated for long time, see e.g. \cite{SCHNACK1987,HarnedSchnack1986,Lerbinger1991,LIONELLO1999}. 
Examples of efficient implementation of nonlinear implicit Newton–Krylov methods can be found in the works by Chacon (2008) \cite{Chacon2008} and L\"utjens and Luciani (2010) \cite{LUTJENS2010}.
Because of the big interest in the topic, many different implicit-schemes have been developed, and the reader is invited to find more details in the recent review paper by Jardin (2012) \cite{Jardin2012} about implicit algorithms for magnetically confined plasma.

Regarding the implicit treatment of acoustic waves, pressure based semi-implicit schemes on staggered grids date back to Harlow and Welch (1965) \cite{markerandcell}, and became soon widely used, 
especially in the last decades, see  e.g. \cite{chorin1,chorin2,patankar,patankarspalding,BellColellaGlaz,vanKan,HirtNichols,Casulli1990,CasulliCheng1992,Casulli1999,CasulliWalters2000,Casulli2009,CasulliVOF,TavelliNS,BUSTO2018,SIDGConv,SIGPR} to mention a few applications to  incompressible fluid-dynamics. In so-called pressure-based methods,  the pressure field is obtained implicitly by solving a simple system of linear algebraic equations.  
Casulli and Greenspan (1984) soon extended this method to compressible gas-dynamics \cite{CasulliCompressible}. 
If those methods were particularly successful for the numerical solution of low-Mach number flows, some attention was needed for the treatment of shock-waves and fluid discontinuities that may appear in both low- or high-Mach regimes. 
Only very recently a novel family of \textit{conservative} semi-implicit schemes have been proposed, see e.g. \cite{MunzPark,CordierDegond,FedkiwSI,IterUpwind,DumbserCasulli2016,RussoAllMach}, making
 pressure-based method also suitable for the solution of nonlinear problems with shock-waves.   
This work is actually the extension of the conservative semi-implicit methods derived by Dumbser et al. (2019) \cite{SIMHD} for the VRMHD equations where the implicit part of the scheme was restricted to the magnetosonic waves. 

In the present work we develop a novel conservative and shock-capturing divergence-free semi-implicit finite-volume/finite-difference method for the nonlinear 
viscous and resistive MHD equations, where both, the magneto-acoustic waves as well as the Alfv\'en waves are treated implicitly, while a standard explicit discretization is only employed for the  nonlinear convective terms. 
The original properties of the sonic implicit solver of \cite{SIMHD} are preserved thanks to the use of a nested-splitting method that will be discussed later in the text.
In this work, the physical variables and the discrete operators are defined on a proper set of  overlapping and mutually dual (or staggered) meshes. The strategy of using staggered grids has to be considered as a paradigm of structure preserving schemes, but this may be considered \emph{optional}. In principle, using the same paradigm, one can decide to apply the same numerical strategy with the introduction of auxiliary variables on the dual grids, while using a standard collocation  at the cell-barycenters for all the physical variables, or even considering the entire algorithm merely on a collocated grid, see e.g. \cite{BDLTV2020}. On the other hand, one may notice that in this work the adopted discretization of the electro-magnetic fields is radically different from what commonly used in magnetohydrodynamics. Indeed, the magnetic field is defined on the \emph{edges} of the main grid, while the electric field is defined on the faces and not vice-versa, see Sec.  \ref{sec:num}. 

The rest of the paper is organized as follows: (i) in Sec.  \ref{sec:goveq} the governing PDE system is first presented in its original form and second, the chosen three-split form is presented as it has been selected for the numerical discretization in such a fashion that the magnetosonic and the Alfvén waves are separated from the  nonlinear hydrodynamic advection; a one-dimensional semi-analytic scheme is also given as representative time-discretization algorithm; (ii) then, in Sec.  \ref{sec:num}, the numerical scheme is presented, after defining the discrete space of solutions, the adopted dual (staggered) domains, the definition of the discrete differential operators, and the final algorithm; S
(iii) Sec.  \ref{sec:numval} collects all the computational results aimed to prove the main features of our novel structure-preserving semi-implicit three-split solver for the full nonlinear VRMHD equations. The paper closes with some concluding remarks and an outlook to future work in Sec.  \ref{sec:conclusion}. 

\section{Governing equations}
\label{sec:goveq} 
The viscous and resistive MHD equations can be cast in the following conservative form
\begin{equation}
\frac{\partial}{\partial t} \Q + \nabla \cdot \left(\F - \Fd\right) = 0 \label{eq:hypPDE}
\end{equation}
where $\Q=(\rho,\rho\mathbf{v},\rho E,\mathbf{B})$ is the array of the corresponding conserved variables. Here, we defined $\rho$ as the fluid matter-density, $\v=(v_x,v_y,v_z)$ is the fluid velocity, $E$ is the specific total energy that can be written in terms of the specific internal energy  $e$, kinetic energy  $k$ and magnetic energy (or pressure) $m$ as $\rho E=\rho e + \rho k + m$, where $\rho k= \frac{1}{2} \rho \v^2$, $m=\frac{1}{8\pi}\B^2$. Then,  $\B=(B_x,B_y,B_z)$ is the magnetic field, $\F$  is the tensor of fluxes for ideal MHD, while the tensor $\Fd$ contains all the purely dissipative fluxes. In particular, the fluxes $\F=\F(\Q)=(\f,\g,\h)$ can be regarded as a nonlinear function of the conserved variables $\Q$,
while the dissipative fluxes $\Fd=\Fd(\Q,\nabla \Q)=(\f^d,\g^d,\h^d)$ depend also on the corresponding gradients $\nabla \Q$.  More specifically, we have defined
\begin{equation}\begin{array}{rl}
& 	\Q:= \left( \begin{array}{c} \rho \\ \rho \mathbf{v} \\ \rho E \\ \mathbf{B} \end{array}  \right); \quad 
	 \F:= \left( \begin{array}{c} \rho \mathbf{v} \\ \rho \mathbf{v} \otimes \mathbf{v} +  \left( p + \frac{\mathbf{B}^2}{8\pi} \right) \mathbf{I} - \frac{1}{4\pi} \mathbf{B} \otimes \mathbf{B} 
	\\  \mathbf{v}^T \left( \rho E + p + \frac{1}{8\pi} \mathbf{B}^2 \right)  - \frac{1}{4\pi} \mathbf{B}^T \left( \mathbf{v} \cdot \mathbf{B} \right)   \\ \mathbf{B} \otimes \mathbf{v} - \mathbf{v} \otimes \mathbf{B} \end{array}  \right) \end{array}\label{eq:MHDQF}  
\end{equation}  
\begin{equation}\begin{array}{l}
	 \Fd := \left( \begin{array}{c} 
0 \\ 
 \mu \left( \nabla \mathbf{v} + \nabla \mathbf{v}^T - \frac{2}{3} \left( \nabla \cdot \mathbf{v} \right) \mathbf{I} \right) \\
\mu \mathbf{v}^T  \left( \nabla \mathbf{v} + \nabla \mathbf{v}^T - \frac{2}{3} \left( \nabla \cdot \mathbf{v} \right) \mathbf{I} \right) + \kappa \nabla T + \frac{\eta}{4 \pi} \mathbf{B}^T \left( \nabla \mathbf{B} - \nabla \mathbf{B}^T \right) \\
\eta \left( \nabla \mathbf{B} - \nabla \mathbf{B}^T \right)    
\end{array} \right) \end{array}  \label{eq:MHDFd}
\end{equation} 
with the identity matrix $\mathbf{I}$, $T$ is the temperature that refers to a thermal equation of state $T=T(p,\rho)$, 
$\mu$ is the kinematic viscosity, $\kappa$ is the thermal conductivity and $\eta$ is the electric resistivity of the fluid. 

It will turn out to be useful to write the momentum and the Faraday equations also in the following non-conservative form, i.e.
\begin{align} 
&\partial_t (\rho \v) + \nabla \cdot \left(   \rho \mathbf{v} \otimes \mathbf{v} +    p   \mathbf{I}\right)  - \frac{1}{4\pi} \left( \nabla \times \B \right) \times \B = 0 \label{eq:momentum}  \\
 & \partial_t \B  + \nabla \times \mathbf{E} = 0, 
\label{eq:Faraday} 
\end{align} 
with the electric field vector given by 
\begin{equation}
 \mathbf{E} = - \mathbf{v} \times \mathbf{B} + \eta \nabla \times \mathbf{B}. 
\label{eq:Evector} 
\end{equation}

One of the big issues that are widely discussed in the literature whenever dealing with the MHD equations, is the numerical preservation of the second Maxwell equation of magnetism, which is known as the Gauss law of magnetism, and widely named as the divergence-free condition of the magnetic field, i.e.  
\begin{equation}
\nabla \cdot \B = 0 \label{eq:divB}.
\end{equation}
 The physical meaning of (\ref{eq:divB}) is that magnetic monopoles do not exist in nature and/or the magnetic field lines are always closed. This is not properly a dynamic governing equation. Indeed, provided that at one initial time this condition is valid,
 the Faraday equation for the magnetic field  (\ref{eq:Faraday}) exactly preserves this condition for all future times. On the other hand, this is not true anymore for the discretized system of the PDE.
Whenever a numerical error induces a violation of the constraint (\ref{eq:divB}), due to machine error, the error of the algorithm or any other spurious effect, then the magnetic field is only approximately divergence-free, spurious numerical monopoles may be generated in the spatial domain, magnetic field lines may break into open lines, and as a consequence spurious oscillations may arise, equilibrium solutions may be compromised and the plasma may be transported in an unphysical direction, parallel to the magnetic field, see \cite{Brackbill1980}.   
These spurious and artificial effects could severely compromise any long time simulation of a Tokamak system, and for this reason, in the context of plasma physics for magnetic confinement fusion, special care should be always taken to develop numerical algorithms able to preserve the constraint (\ref{eq:divB}) up to machine precision.

Many other techniques  have been developed  in the past to preserve or at least to control the divergence-free conditions. Two of them that deserve to be noted are the staggered-formulation, widely known as '\emph{constrained transport}', see e.g. the work by Evans and Hawley (1998) and by De Vore (1991) \cite{Evans1988,DEVORE1991}, and the so-called \emph{hyperbolic Generalized Lagrangian Multiplier} approach (GLM) proposed by Dedner et al.  (2002) \cite{Dedneretal}, based on the hyperbolic cleaning technique by Munz et al. (2000)\cite{MunzCleaning}. For a comparison between these two different approaches and a detailed list of references, see \cite{BalsaraKim2004}.
Our approach is actually a staggered-formulation, and similar to the well-known schemes proposed by Balsara and Spicer (1999) \cite{BalsaraSpicer1999} and the subsequent works by Balsara  \cite{Balsara2004,balsarahlle2d,balsarahllc2d}, the magnetic and the electric fields are defined and evolved in time on two different staggered grids.

\subsection{3-split operator for the MHD system} 
\label{sec.split} 

In this work, we introduce a novel semi-implicit discretization of the nonlinear MHD equations that may become convenient especially for the numerical simulation of highly magnetized plasma flows, both in the low or high Mach number regimes, that are particularly interesting for applications of magnetic confinement fusion in a Tokamak. 

In order to proceed, the flux-tensors are split  in a numerically more convenient form
$$ \F := \Fv + \Fb + \Fp$$ 
that represents a step forward compared to the purely pressure-based splitting of \cite{SIMHD} which was based on the Toro-V\'azquez splitting of the compressible Euler equations \cite{ToroVazquez}. The new scheme presented in this paper thus results in a new \emph{triple} flux-splitting scheme, with the aim of isolating both the Alfvénic and the pressure modes within a \emph{coupled nonlinear implicit solver}. In particular, we defined the flux tensors $\Fv=(\fv,\gv,\hv)$, $\Fp=(\fp,\gp,\hp)$ and $\Fb=(\fb,\gb,\hb)$ to be the convective Euler fluxes, the $p$-related  Euler fluxes and the $\B$-related fluxes, respectively, i.e.
\begin{equation}\begin{array}{rl}
& 	\Fv:= \left( \begin{array}{c} \rho \mathbf{v} \\ \rho \mathbf{v} \otimes \mathbf{v} 
	\\  \mathbf{v}^T  \rho k   \\ 0 \end{array}  \right) 
\quad 	\Fp:= \left( \begin{array}{c}  0 \\    p \mathbf{I} 
\\    \mathbf{v}^T \rho  h   \\  0 \end{array}  \right) 
\quad 	\Fb:= \left( \begin{array}{c} 0 \\  m \mathbf{I} - \frac{1}{4\pi} \mathbf{B} \otimes \mathbf{B} 
	\\  \mathbf{v}^T 2 m  - \frac{1}{4\pi} \mathbf{B}^T \left( \mathbf{v} \cdot \mathbf{B} \right)   \\ \mathbf{B} \otimes \mathbf{v} - \mathbf{v} \otimes \mathbf{B} \end{array}  \right)\end{array}  \label{eq:FluxSplitting}
\end{equation} 
where we have used the definition of the enthalpy  $ h = e + \frac{p}{\rho}$ and $\rho E + p + \frac{1}{8\pi} \mathbf{B}^2= \rho k + \rho h + 2 m$.
In order to study the spectrum of the hyperbolic part of the system, we first rely on a simplified model, considering ideal plasma flow in one single space dimension: 
\begin{equation}
   \frac{\partial \Q}{\partial t} + \frac{\partial \f}{\partial x} = 0. 
	\label{eq:split1} 
\end{equation}
For \eqref{eq:split1} and with $B_x = const.$ one can easily compute its eight eigenvalues
\begin{equation}
\lambda_{1,8}^{\text{MHD}} = u \mp c_f, \quad \lambda_{2,7}^{\text{MHD}} = u \mp c_a, \quad \lambda_{3,6}^{\text{MHD}} = u \mp c_s, \quad \lambda_4^{\text{MHD}} = u, \quad \lambda_5^{\text{MHD}} = 0, 
\label{eq:eval.full} 
\end{equation} 
where 
\begin{equation}\begin{array}{rl}
 &c_a = B_x / \sqrt{4 \pi \rho}, \\
 &c_s^2 = \halb \left( b^2+c^2 - \sqrt{(b^2+c^2)^2-4 c_a^2 c^2} \right), \\
 &c_f^2 = \halb \left( b^2+c^2 + \sqrt{(b^2+c^2)^2-4 c_a^2 c^2} \right).  \end{array}
\label{eq:mhd.wavespeeds} 
\end{equation} 
where  $c_a$, $c_s$ and $c_f$ are the speeds of the Alfv\'en, the slow and the fast magnetosonic waves, respectively. Then,  $c$ is the adiabatic sound speed and it can be computed with respect to a general equation of state (EOS) $p=p(e,\rho)$ as $c^2 = \partial p / \partial \rho + p / \rho^2 \partial p / \partial e$, e.g. $c^2 = \gamma p/ \rho$ for the ideal gas EOS. Here, we used also the abbreviation $ b^2 = \mathbf{B}^2/ (4 \pi \rho)$.

Depending on the flow scenario, the dominant eigenvalue of (\ref{eq:eval.full}) may change. In the low Mach number limit, indeed, the sound speed becomes very high if compared to the fluid velocity and the sonic components are dominant. On the other hand, for highly magnetized plasma, all the eigenvalues become important. Whenever discretizing a hyperbolic system like (\ref{eq:split1}) through a semi-implicit or an implicit-explicit scheme, in order to take benefit of it in terms of weakening the CFL time restriction, the faster modes should be always treated implicitly in time. 

Therefore, the next step is to try to divide the fast modes from the slow ones.
From the 1D model (\ref{eq:split1}) the following 3-split form can be obtained  
\begin{equation}
   \frac{\partial \Q}{\partial t} + \frac{\partial \f^v}{\partial x} + \frac{\partial \f^p}{\partial x}+ \frac{\partial \f^b}{\partial x} = 0, 
	\label{eq:split2} 
\end{equation}
where the indexes of the fluxes refer to the notation introduced in equation (\ref{eq:FluxSplitting}).
Finally, equation (\ref{eq:split2}) leads to three subsystems that, together with their corresponding eigenspectra, read as
\begin{align}
 \left.i\right)&& \partial_t \Q +  \partial_x \f^v  = 0,  &\qquad& &
   \lambda^v_{1,2,3,4} =  0, \label{eq:1dc}\\
	&&&&& 
	 \lambda^v_{5,6,7,8} = v_x,  
	 \nonumber  \\
\left.ii\right)&&	\partial_t \Q +  \partial_x \f^p  = 0,  &&&
 \lambda^p_{1,8} = \frac{1}{2} \left( v_x \mp \sqrt{v_x^2 + 4 c^2 } \right), \label{eq:1dp}\\ &&&&&\lambda^p_{2,3,4,5,6,7} = 0,  
 \nonumber \\
\left.iii\right)&&	\partial_t \Q +  \partial_x \f^b  = 0.&& &
   \lambda^B_{1,8} = \frac{1}{2} \left( v_x \mp \sqrt{v_x^2 + 4 \left( \frac{\left | \B \right |}{ \sqrt{ 4 \pi \rho }}\right)^2 } \right), \label{eq:1db} \\ &&&&&
	 \lambda^B_{2,7} = \frac{1}{2} \left( v_x \mp \sqrt{v_x^2 + 4 \left( \frac{B_x}{\sqrt{ 4 \pi \rho }}\right)^2 } \right), \nonumber\\ &&&&&\lambda^B_{3,4,5,6} =  0. 
	 \nonumber \end{align}

 A divergence-free algorithm for the system (\ref{eq:hypPDE}-\ref{eq:MHDQF}) that uses an implicit time-discretization of the pressure terms, crudely represented by equation (\ref{eq:1dp}), has been successfully implemented in \cite{SIMHD}, and it turned out to be extremely convenient in terms of computational efficiency for the simulation of low Mach number flows and, whenever handling high Mach regimes, comparable to other fully explicit codes.

In this work, a novel algorithm for solving the viscous and resistive MHD equations (\ref{eq:hypPDE}-\ref{eq:MHDFd}) is particularly well-suited for solving highly magnetized plasma in the low Mach regime. For this reason, an implicit time-discretization in used for solving the pressure and magnetic terms, represented by the split systems (\ref{eq:1dp}) and (\ref{eq:1db}), respectively, while and explicit time-discretization is used only for the purely convective fluxes in subsystem \eqref{eq:1dc}.

It should be mentioned that 
the three subsystems (i), (ii) and (iii) need to be coupled also in time. In this work, the coupling is performed via a very simple operator-splitting algorithm. 
 Indeed, due to their nonlinearities, the sub-systems (ii) and (iii) will be solved through two nested Picard-corrector algorithms, by means of \emph{locally}-linearized steps. 
Further details will be given and the algorithm will be outlined in the following sections. It should be noticed that the time-integration algorithm of the three-splitting scheme is not a simplectic integrator
 A simplified sketch of the final algorithm is provided in equation (\ref{eq:finalalg}), and can be summarized, for the one-dimensional case as:

\begin{align} 
&\left\{
\begin{array}{l}
 	\frac{\blue{\FSa} - \Q^n}{ \Delta t } +  \partial_x \f^v(\Q^n)  = 0, \\ 
 	\frac{\blue{\FSb^{r+1}} - \FSa}{ \Delta t} + \partial_x \f^b(\blue{\FSb^{n+\theta_{\text{B}},r+1}}) + \partial_x \f^p(\FSc^{n+\theta_{\text{p}},r})    = 0, \end{array}\right. \\
	&	\left\{
	\begin{array}{l}
	\frac{\blue{\FSd}  - \Q^n }{ \Delta t} + \partial_x \f^v(\Q^n) +  \partial_x \f^b(\FSb^{n+\theta_{\text{B}},r+1}  )   = 0,\\
	 	\frac{\blue{\FSc^{r+1}} - \FSd }{ \Delta t} + \partial_x \f^p(\blue{\FSc^{n+\theta_{\text{p}},r+1}})  = 0, 
	\end{array}\right.
\label{eq:1dalgorithm} 
\end{align}
that, with a direct summation of the two coupled systems, can be cast into  the following  algorithm 
\begin{align} 
 \;\; i)\quad &   
 	\frac{\blue{\tilde{\Q}^{n+1,r+1}} - \Q^n}{ \Delta t} +  \partial_x \f^v(\Q^n) + \partial_x \f^b(\blue{\tilde{\Q}^{n+\theta_{\text{B}},r+1}}) + \partial_x \f^p(\Q^{n+\theta_{\text{p}},r})    = 0,   \\
 \;\; ii)\quad	&		 		\frac{\blue{\Q^{n+1,r+1}} - \Q^n }{ \Delta t} + \partial_x \f^v(\Q^n) +  \partial_x \f^b(\tilde{\Q}^{n+1,r+1,s+1} )+ \partial_x \f^p(\blue{\Q^{n+\theta_{\text{p}},r+1}})  = 0, \label{eq:1dalgorithm2}  
\end{align}
where $r=1$,$\ldots$, $R$ is the iteration index. Here, where the index $r$ is used, the time index $n+1$ is obviously assumed. In this way, the first system is solved with an implicit discretization of the magnetic terms obtaining an \emph{auxiliary} solution $\tilde{\Q}$, and then the second system is solved with an implicit time-discretization of the pressure terms, evaluating the other fluxes with respect to $\FSb$, obtaining a solution $\Q^{r+1}$. 
In the practice, the final numerical solution $\Q^{n+1}$ is built with the magnetic field components of $\tilde{\Q}$, and matter-density, momentum and total-energy of $\Q^{r+1}$. Here the $\theta$-variables are defined as $Q^{n+\theta}:=(1-\theta)Q^n + \theta Q^{n+1}$, with a real valued implicit-parameters $\theta_{\text{B},\text{p}} \in [\halb,1]$. 
As initial guess for the scheme, we used  $\tilde{\Q}^{r=1}=\Q^{r=1}=\Q^n$ while the pressure terms of the first iteration can be set to zero, i.e $\f^p(\Q^{n+\theta_{\text{p}},r=1})\equiv 0$. %

Finally, we mention that the diffusive parabolic terms of the PDE system can also split in two parts.
Indeed, a parabolic time step restriction can also become rather severe, in particular for large values of kinematic viscosity or magnetic resistivity, but also for highly refined meshes. The parabolic time step  constraint is actually linear in the physical dissipative parameters, magnetic resistivity $\eta$  or fluid viscosity $\mu$, and quadratic in the characteristic mesh size.

\noindent
In this manuscript,  the magnetic resistivity is also discretized  within the implicit scheme and then, the diffusive terms (\ref{eq:MHDFd}) have been split as
\begin{equation}  \begin{array}{l}
   \Fd := \Fdv + \Fdb,\hspace{1cm} \Fdb := \left( \begin{array}{c} 
0 \\ 
 0 \\
  \frac{\eta}{4 \pi} \mathbf{B}^T \left( \nabla \mathbf{B} - \nabla \mathbf{B}^T \right) \\
\eta \left( \nabla \mathbf{B} - \nabla \mathbf{B}^T \right)    
\end{array} \right),  \\
	\Fdv := \left( \begin{array}{c} 
0 \\  \mu \left( \nabla \mathbf{v} + \nabla \mathbf{v}^T - \frac{2}{3} \left( \nabla \cdot \mathbf{v} \right) \mathbf{I} \right) \\
\mu \mathbf{v}^T  \left( \nabla \mathbf{v} + \nabla \mathbf{v}^T - \frac{2}{3} \left( \nabla \cdot \mathbf{v} \right) \mathbf{I} \right) + \kappa \nabla T \\
0   
\end{array} \right), \end{array}
\end{equation} 
The summary of the chosen split systems will be
	\begin{align}
   & \partial_t \Q + \nabla \cdot \mathbf{F} = 0, & \mathbf{F}:=    & \left( \Fv - \Fdv \right)  & \text{\emph{Explicit in time}}  \label{eq:PDEex}\\
   &&  &  +\left( \Fb -\Fdb \right)     & \text{\emph{Implicit in time}}   \label{eq:PDEimB}\\
  && &  + \Fp. 										 	 & \text{\emph{Implicit in time}}   \label{eq:PDEimP}
\end{align}
It has also to be noticed that the two implicit subsystems (ii) and (iii) will be coupled within a single nested Picard-type algorithm, see Sec.  \ref{sec:num}.

\subsection{Ideal gas EOS} 
\label{sec.eos} 
 
 Although the presented numerical formulation is valid for a general nonlinear equation of state $e=e(p,\rho)$, we will test the numerics only after choosing  the thermal and caloric
equations of state for ideal gases, i.e.
\begin{equation} 
\label{eq:thermcal.ideal} 
 \frac{p}{\rho} = R T, \qquad \textnormal{ and } \qquad e = c_v T,  
\end{equation} 
where we have introduced the specific gas constant $R = c_p - c_v$, and the heat capacities $c_v$ and $c_p$ at constant volume and at constant pressure, respectively. From 
\eqref{eq:thermcal.ideal} one easily obtains 
\begin{equation}
 e = e(p,\rho) = \frac{p}{(\gamma-1) \rho} ,  
\end{equation} 
which is linear in the pressure $p$ and where $\gamma = c_p / c_v$ denotes the so-called ratio of specific heats. 

\section{Semi-implicit discretization on staggered grids} 
\label{sec:num} 

Regarding the spatial discretization, we adopt a finite-volume/finite-difference scheme on staggered grids that can be considered as an extention of the previous work \cite{SIMHD} to an implicit treatment of the magnetic forces, but it can be regarded also as a first finite-difference actuation of the theory of \emph{finite element exterior calculus} (FEEC) in the discretization of the MHD equations \cite{FEEC2010,Possanner}. 
In this work the choice of staggered grids is motivated mainly by two argument: the preservation of the divergence-free condition of the magnetic field, or the velocity field in the incompressible limit; the preservation of the duality relation between the discrete differential operators.
Dating back to Yee (1966) \cite{Yee1966}, the use of staggered grids is known to be particularly indicated for the electric and magnetic fields, and then also for the MHD equations,
 Indeed, it has been shown that an electric field $\E$ defined on a mesh that is \emph{dual} with respect to the mesh of the magnetic field $\B$ is the natural way of preserving the divergence-free condition up to machine error. This is rigorously supported by the analysis of FEEC for MHD. Indeed, depending of the topology that is induced by the differential operators of the governing PDE system, the FEEC analysis clearly indicates which variables are naturally elected for being defined on an \emph{adjoint} (or \emph{dual}) space. In this sense, FEEC becomes a fundamental mathematical tool for the definition of \emph{structure-preserving schemes} on staggered grids. Thanks to the preservation of the duality relation of the discrete differential operators, the resulting algebraic system that has to be solved for the unknown variables is \emph{symmetric} by construction and a simple \emph{matrix-free} conjugate-gradient method can be easily implemented, even without any preconditioning strategy. 

Since in this paper we mainly introduce a novel numerical scheme, for simplicity we assume the physical domain to be rectangular or cubic in two or three space-dimensions, respectively, and a segment in the one-dimensional case. Without loss of generality, we directly present the numerical discretization for the three dimensional case.

As usual in the context of staggered meshes, some variables are defined on the control volumes of the main grid, others are centered on the respective faces, others on the respective edges. The result will be a set of different staggered overlapping grids. Obviously one can always project one variable from one grid to the other by means of a suitable averaging operator, whenever needed. Indeed, it is always assumed that a product of two variables can be performed only if they are located at the same discrete spatial position. For this reason, averaging operators are used. For the implicit solver, the discrete differential and averaging operators are always defined as centered, mapping variables from one grid to the other. Because of the chosen staggering and centered discrete operators, some symmetries that are present in the analytical PDE can be preserved at the discrete level.


\subsubsection*{Discrete domain and space of solutions}

\redP{Here, $\Omega:=\Omega_x\times \Omega_y \times \Omega_z$ is the physical domain, defined as the Cartesian product of the three spatial intervals $\Omega_x = [x_L,x_R]$, $\Omega_y = [y_L,y_R]$ and $\Omega_z = [z_L,z_R]$ in the $x$, $y$ and $z$ direction, respectively. Then, a primary grid $\Omega_h$ is defined as the standard uniform Cartesian mesh obtained after choosing the three discretization numbers $N_x$, $N_y$ and $N_z$ for the three spatial directions.
In this pattern, a single control volume $T_{ijk}\in \Omega_h$ can be identified by its three spatial indexes,
so that $T_{i,j,k}=[x_{i-\halb},x_{i+\halb}]\times[y_{j-\halb},y_{j+\halb}]\times[z_{k-\halb},z_{k+\halb}]$. We define $\Delta s_i\in\Omega_s$, with $s=x$, $y$ or $z$ to be the one-dimensional spatial-interval centered in the coordinate point $s_i$, having, e.g., $\Delta s_i:= [s_{i-\halb},s_{i+\halb}]$, and one can write
\begin{equation} \begin{array}{rl}
T_{i,j,k}&= \Delta x_i \times \Delta y_j \times \Delta z_k \in \Omega_h   \\ 
&i=1,\ldots,N_x; \quad j=1,\ldots,N_y;\quad  k=1,\ldots,N_z;  \end{array}
\end{equation}}

In this work, we define also other two meshes that are dual to $\Omega_h$, i.e. the \emph{node-based} and the \emph{edge-based}  staggered grids, that correspond to the \emph{B-}  and the \emph{C-grid}, respectively, according to the nomenclature of Arakawa \& Lamb \cite{Arakawa}. 
It turns to be useful to define the space of the discrete functions where the numerical observables can be defined. In the Cartesian and structured case, we indicate with $\Vmain_x$, $\Vmain_y$ and $\Vmain_z$ the three one-dimensional spaces of piecewise-constant functions with real values that build the main space of the numerical solutions as 
$$\Vmain := \Vmain_x \otimes \Vmain_y  \otimes \Vmain_z = \left\{ f\;\big|\;f : \Omega_h \xrightarrow[\quad\quad\quad]{} \Real  \right\} $$
 Then, we indicate 
the set of the three face-based meshes with $$\Omega_{\text{faces}}:=\Omega_{\text{x-faces}}\times\Omega_{\text{y-faces}}\times\Omega_{\text{z-faces}}$$ and
\begin{equation} \begin{array}{l}
%
\Omega_{\text{$x$-faces}} \ni T_{i+\halb, j,k}=\Delta x_{i+\halb} \times \Delta y_{j} \times \Delta z_{k}   \\ 
\Omega_{\text{$y$-faces}} \ni T_{i,j+\halb, k}=\Delta x_{i} \times \Delta y_{j+\halb} \times \Delta z_{k}  \\ 
\Omega_{\text{$z$-faces}} \ni T_{i,j,k+\halb }=\Delta x_{i}\times  \Delta y_{j} \times \Delta z_{k+\halb}   \end{array}  
\end{equation} 
that are staggered to the primary grid along the $x$, $y$ and $z$ directions, respectively.  Here the notation has been extended to halfed indexes as  $\Delta s_{i+\halb}:= [s_{i},s_{i+1}]\in\Omega_{\text{$s$-faces}}$. 
The corresponding (\emph{dual}) space of discrete solutions are $\Vdual_x$, $\Vdual_y$ or $\Vdual_z$ depending on the space-direction of the staggering. In this way, the components of a three-vector $\v$ defined on the faces belong to the following discrete space-of-solutions,
\begin{equation}
 \Vfaces :=  \Vxfaces \times \Vyfaces \times \Vzfaces
\end{equation}
where
\begin{equation}\begin{array}{rl}
& \Vxfaces := (\Vdual_x \otimes \Vmain_y \otimes \Vmain_z) \\
& \Vyfaces := (\Vmain_x \otimes \Vdual_y \otimes \Vmain_z)\\
& \Vzfaces := (\Vmain_x \otimes \Vmain_y \otimes \Vdual_z) \end{array}
\end{equation}
having
\begin{equation}
 \v := (v^1_h, v^2_h, v^3_h )  \in \Vfaces, \quad\quad \v : \Omega_{\text{faces}}  \xrightarrow[\quad\quad\quad]{}  \Real^3
\end{equation}
\begin{figure}\centering
\begin{tabular}{ccc} 
   \multicolumn{3}{c}{ 
\begin{tabular}{cc} 
	\includegraphics[width=0.36\textwidth]{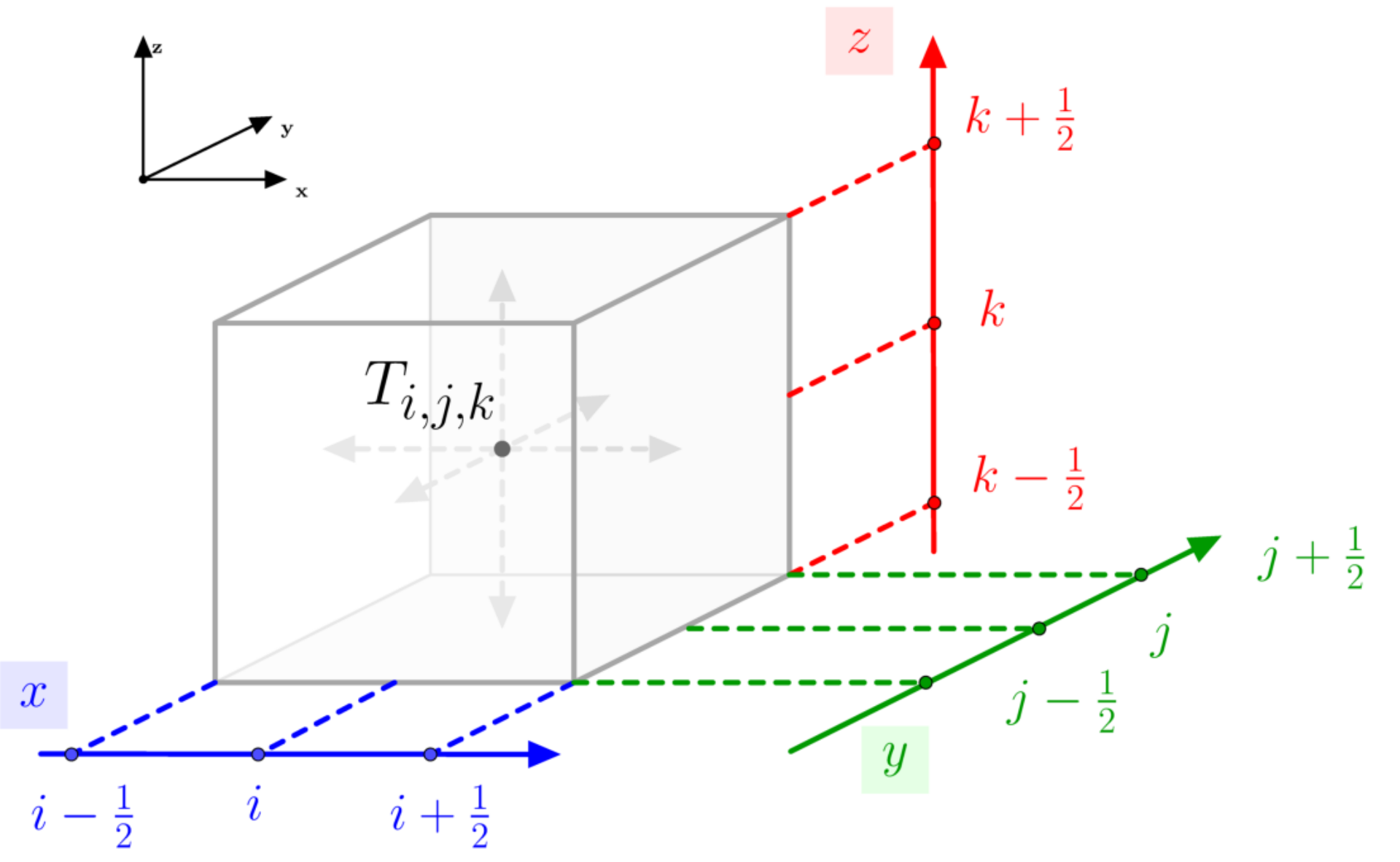}&
	\includegraphics[width=0.36\textwidth]{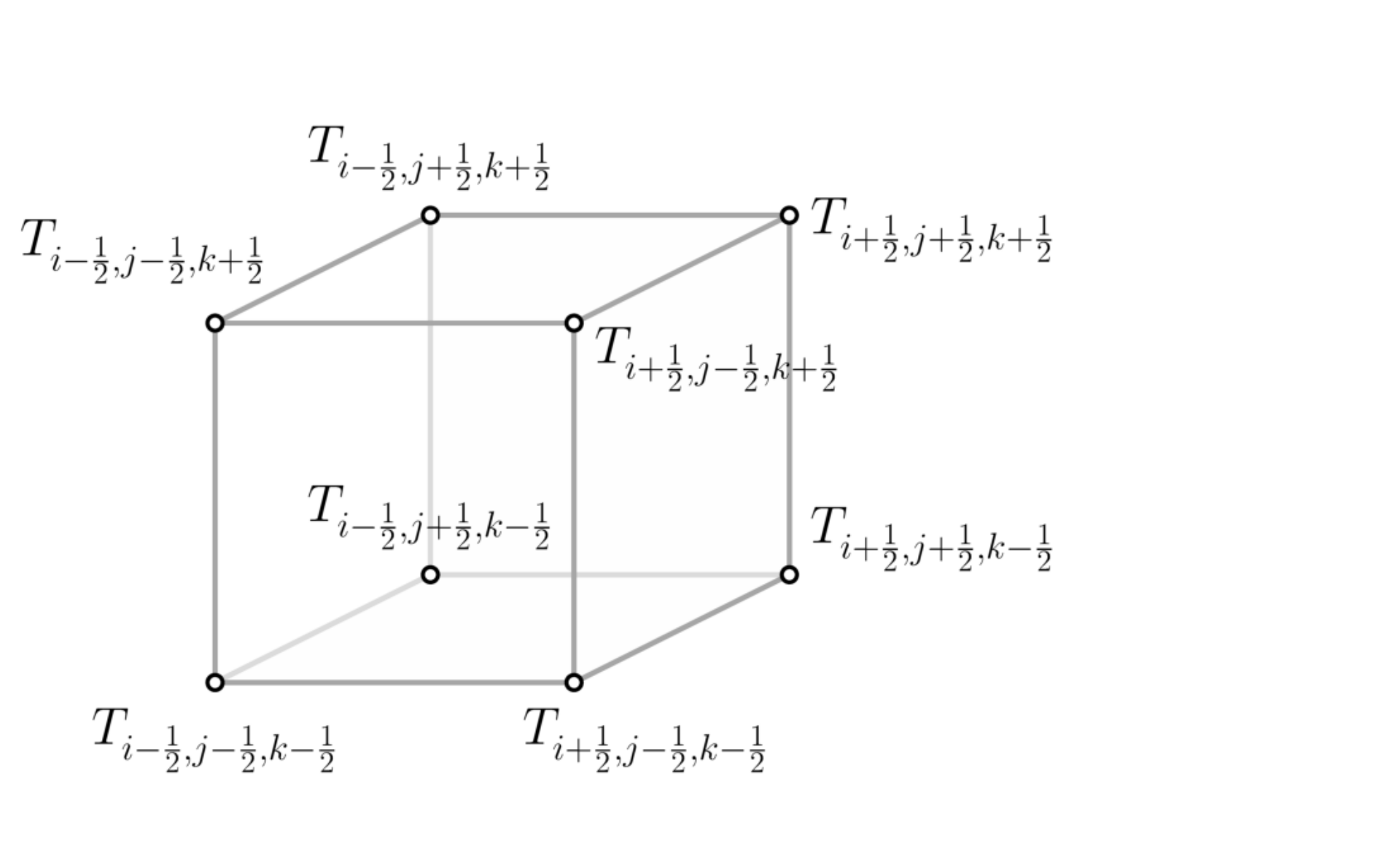} 
\end{tabular}
}  \\
\includegraphics[width=0.33\textwidth]{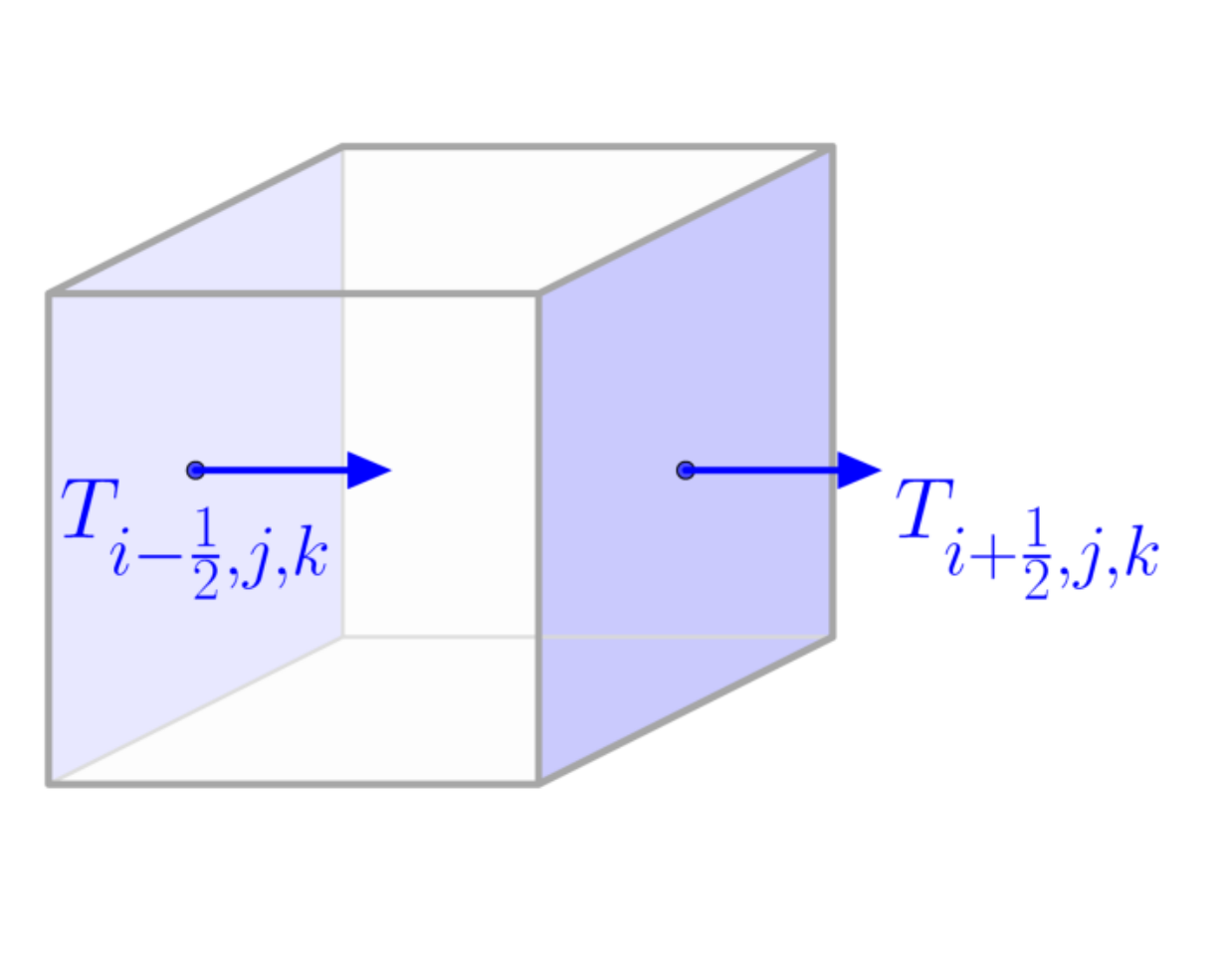} &
\includegraphics[width=0.33\textwidth]{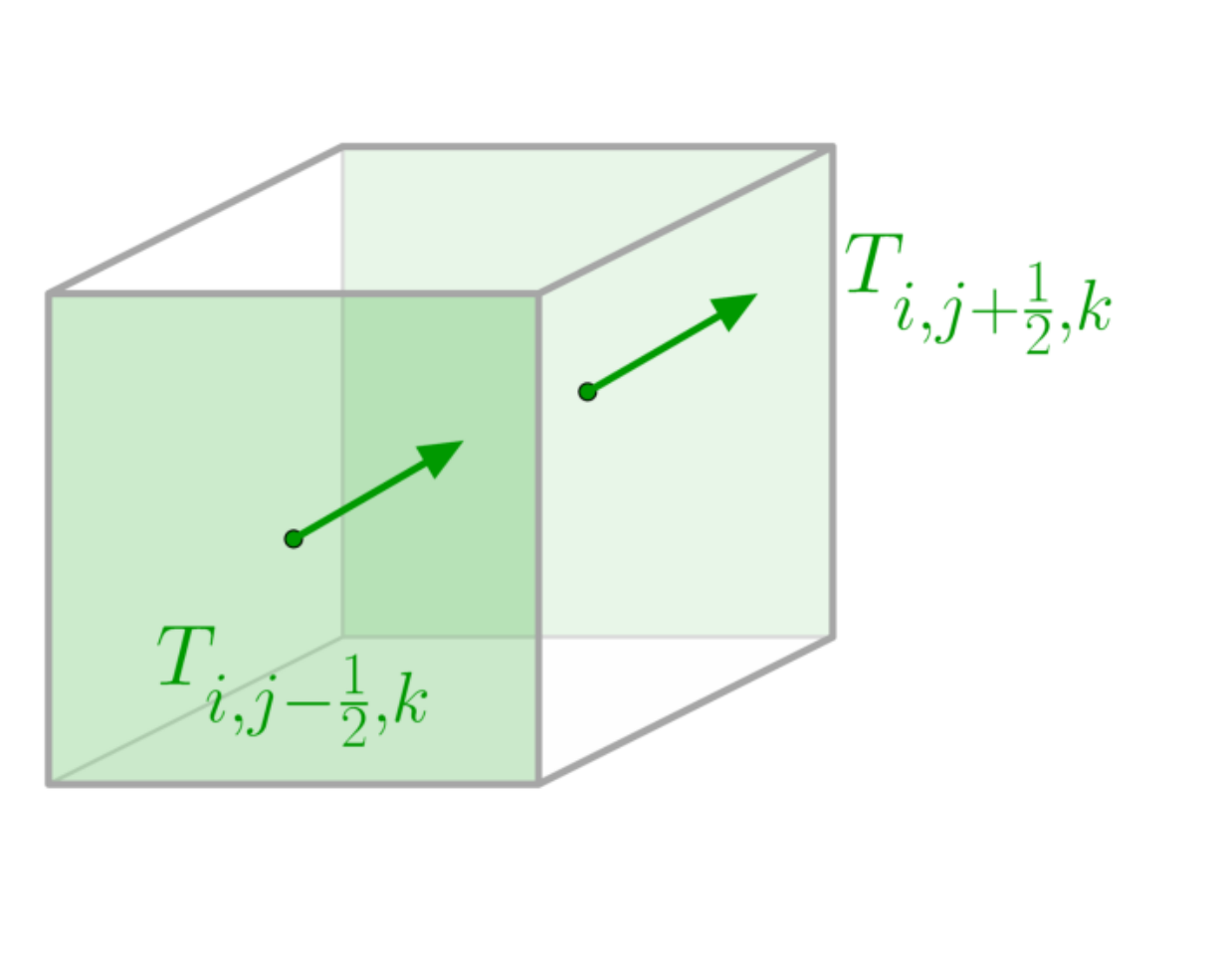} &
\includegraphics[width=0.33\textwidth]{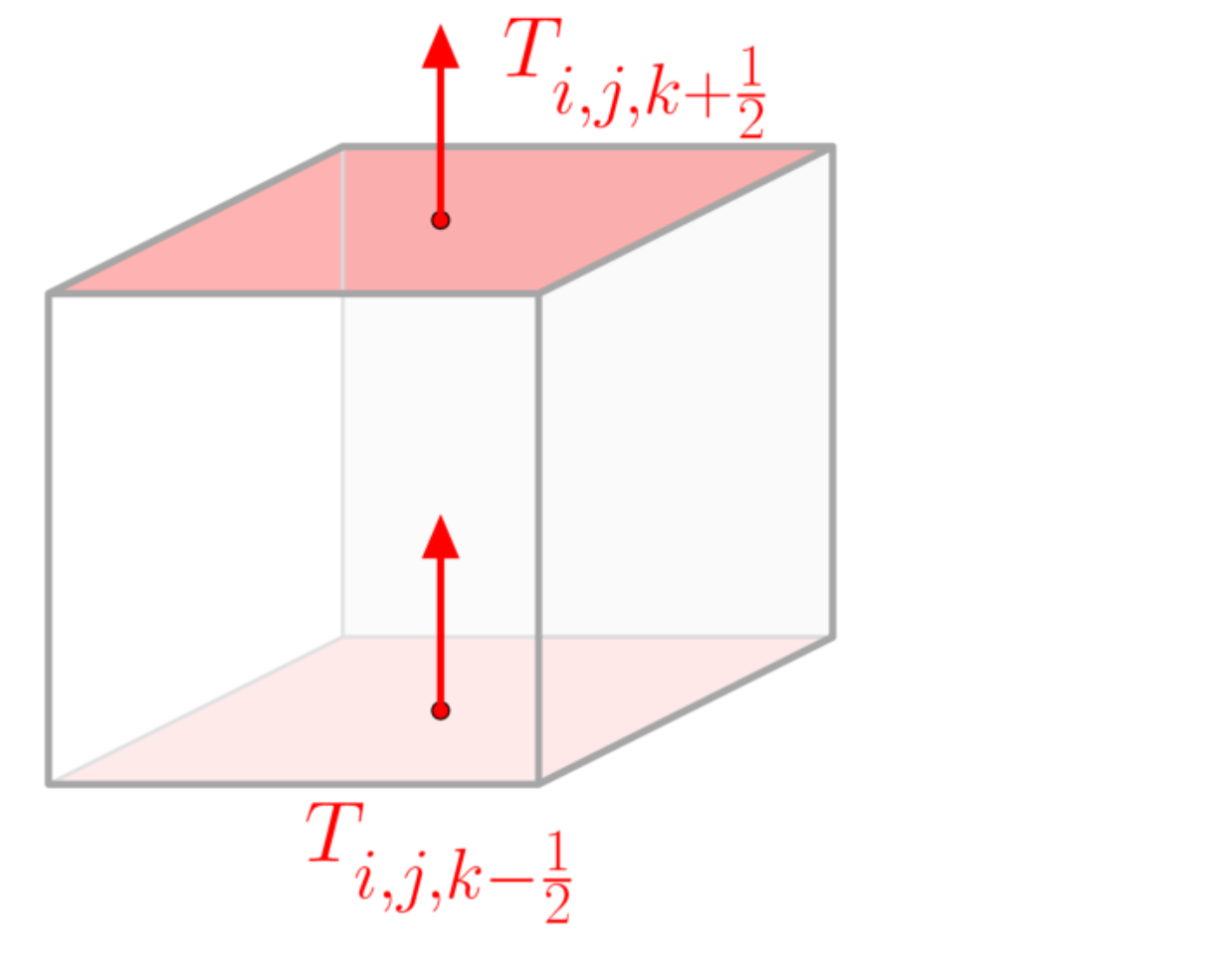} \\
\includegraphics[width=0.33\textwidth]{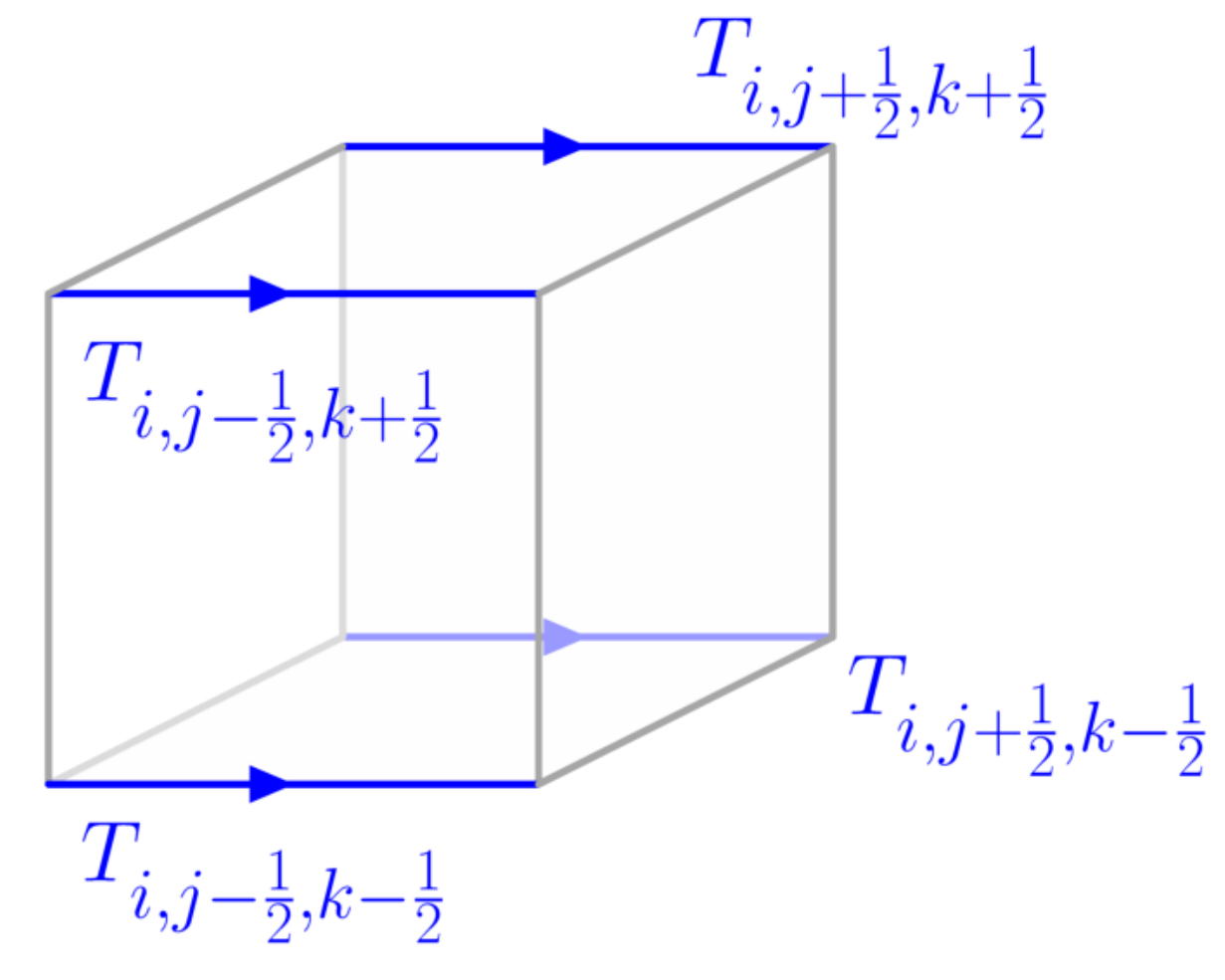} &
\includegraphics[width=0.33\textwidth]{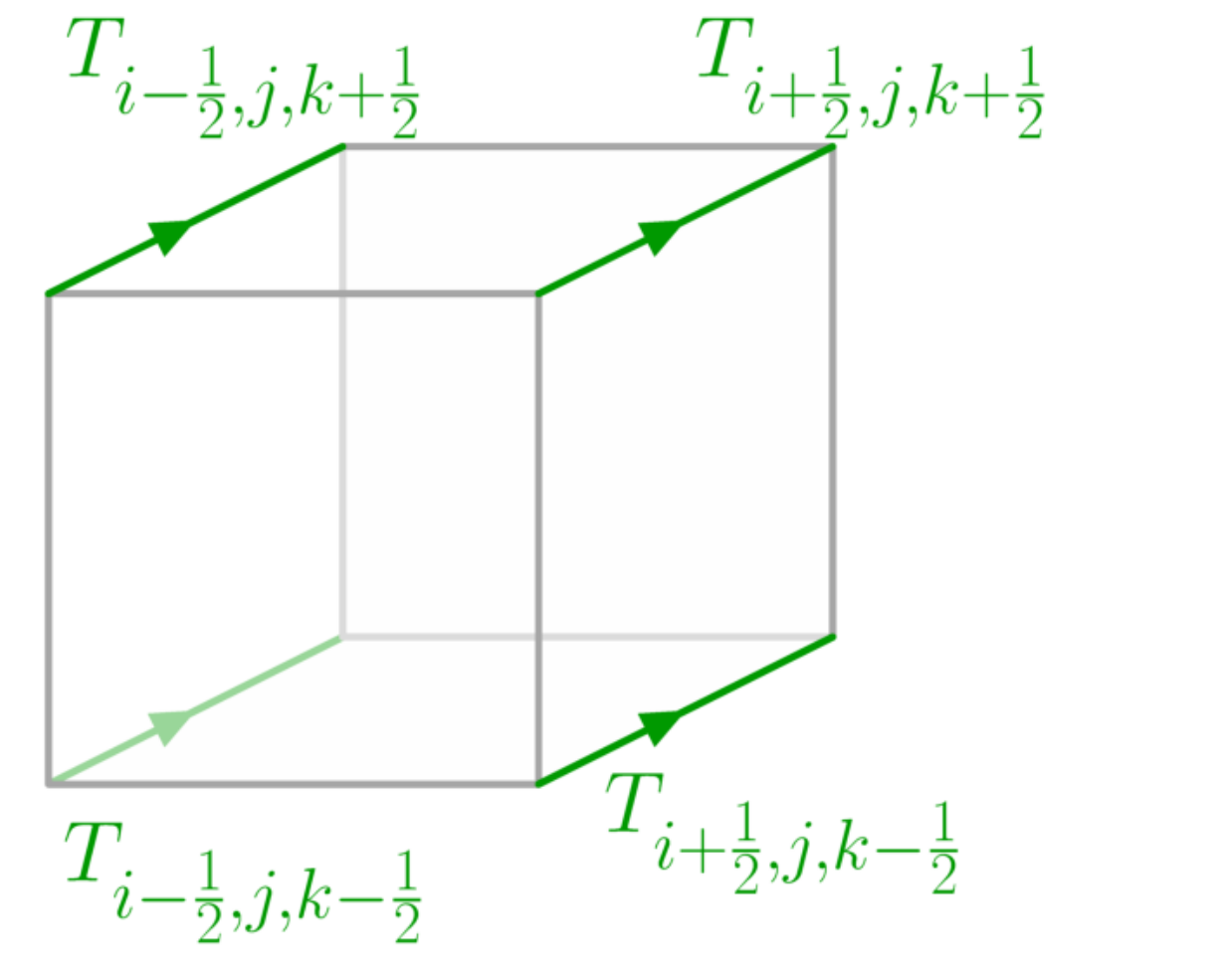} &
\includegraphics[width=0.33\textwidth]{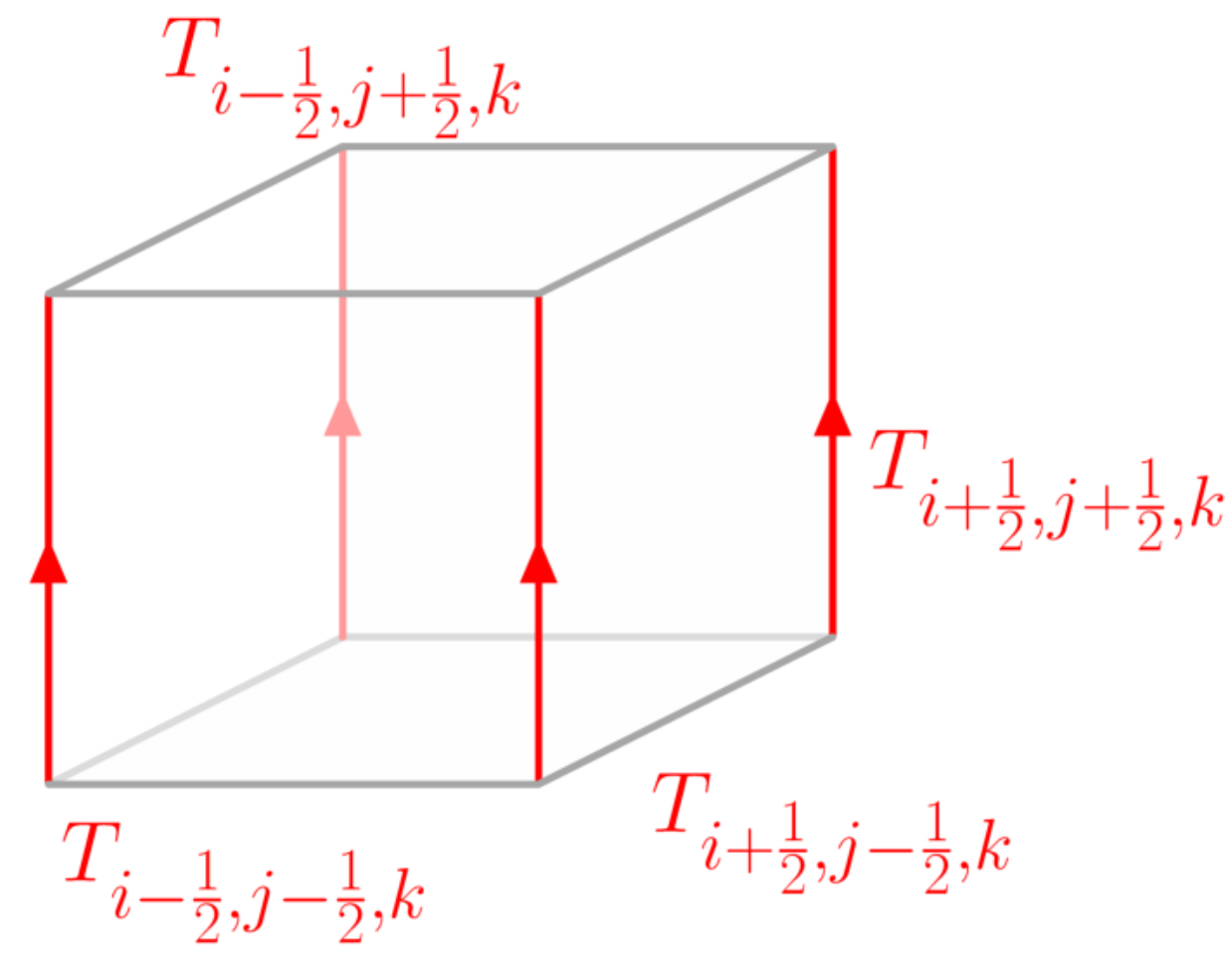}
\end{tabular}
\caption{Sketch of the location of the discrete variables. Barycenters, vertexes, faces and edges on a main three-dimensional Cartesian structured  grid. It may become useful for the reader reporting that in this work the velocity variables are defined on the \emph{faces}, the magnetic-field on the \emph{edges}.}
\label{fig:dualmeshes}  
\end{figure}
 Then, in three space-dimensions, one has also the \emph{edge-based} staggered meshes where the cells are centered at spatial discrete indexes with two halfed indexes, i.e.
  $$\Omega_{\text{edges}}:=\Omega_{\text{x-edges}}\times\Omega_{\text{y-edges}}\times\Omega_{\text{z-edges}}$$ with
\begin{equation} \begin{array}{l}
%
\Omega_{\text{$x$-edges}} \ni T_{i, j+\halb,k+\halb}=\Delta x_{i} \times \Delta y_{j+\halb} \times \Delta z_{k+\halb}   \\ 
\Omega_{\text{$y$-edges}} \ni T_{i+\halb,j, k+\halb}=\Delta x_{i+\halb} \times \Delta y_{j} \times \Delta z_{k+\halb}  \\ 
\Omega_{\text{$z$-edges}} \ni T_{i+\halb,j+\halb,k }=\Delta x_{i+\halb}\times  \Delta y_{j+\halb} \times \Delta z_{k} \end{array}    
\end{equation}
Similarly, the components of a three-vector $\u_h$ defined on the edges belong to the following discrete space-of-solutions,
\begin{equation}
 \Vedges :=  \Vxedges \times \Vyedges \times \Vzedges
\end{equation}
where
\begin{equation} \begin{array}{rl}
& \Vxedges := (\Vmain_x \otimes \Vdual_y \otimes \Vdual_z) \\
& \Vyedges := (\Vdual_x \otimes \Vmain_y \otimes \Vdual_z)\\
& \Vzedges := (\Vdual_x \otimes \Vdual_y \otimes \Vmain_z)  \end{array}
\end{equation}
having
\begin{equation}
\u_h := (u^1_h, u^2_h, u^3_h )  \in \Vedges, \quad\quad \v : \Omega_{\text{edges}}  \xrightarrow[\quad\quad\quad]{}  \Real^3
\end{equation}

Finally we define also the node-based staggered mesh  $\Omega_{\text{nodes}}$ as
\begin{equation}\begin{array}{rl}
%
\Omega_{\text{nodes}}\ni T_{i+\halb, j+\halb, k+\halb }=&\Delta x_{i+\halb}\times \Delta y_{j+\halb} \times \Delta z_{k+\halb} \end{array}
\end{equation} 
and its corresponding space of discrete solutions $\Vnodes$.
In this way, it has been defined up to 8 different and mutually overlapping meshes: one the \emph{primary} grid, three \emph{face}-based staggered meshes, three \emph{edge}-based staggered meshes and one node-based staggered mesh; see Fig.  \ref{fig:dualmeshes} for a schematic view of the all considered main and dual meshes. 
 It has to be noticed that, in principle, this choice does not necessarily affect in any way the total amount of computational memory requested for storing the respective physical variables. 
Even if one can easily project one discrete variable from its original discrete location to another, the following choice of staggering is adopted to define where the variables are located:
\begin{align}
\rho_b,\rho\v_b,\rho E_b,p_b \in \Vmain &&
\B_e \in \Vedges &&
\E_f, \redP{\v_f} \in \Vfaces
\end{align}
The subscripts are used to ease the interpretation of the variable location, i.e. $b$ for the variables at the \emph{barycenter} of the control volumes, $f$ at the \emph{faces}, $e$ at the \emph{edges}, $n$ at the \emph{nodes} of the grid.

In the next we also introduce the usual time index $n$, referring to a corresponding time $t=t^n$ that is reached after the $n$-th time-step of the simulation. Consequently, one has $$ X^n := X(t=t^n)$$ where $X_h$ is the numerical approximation of a given physical observable $X$, or  $$t^{n+1} := t^n + \Delta t^n$$ where $\Delta t^n$ is the length of the $(n+1)$-th time-step, which may be  in principle computed as a function of the conserved variables $\mathbf{Q}^n$, see eq. (\ref{eq:CFL}). The initial time of a simulation is assumed to be at $t=t^0$.  

\subsubsection*{Projection operators}

Through the text, it will become useful to define the projection operators that allow projecting the numerical information carried by a physical observable from one grid to the other. To do this, we first define the standard one-dimensional discrete \emph{averaging} operator $\mathcal{A}$ as 
\begin{equation}
\left[\langle X \rangle \right]_l^n  \equiv \mathcal{A}\left( X_{\mathcal{R}(l)}^n, X_{\mathcal{L}(l)}^n \right)  := \halb \left( X_{\mathcal{R}(l)}^n + X_{\mathcal{L}(l)}^n  \right) \label{eq:av1db} 
\end{equation}
where $\mathcal{R}(l)$ and $\mathcal{L}(l)$ are respectively the right and left halfed (or integer) indexes with respect to the integer (or halfed, correspondingly) index $l$.
Here, we have also defined an obvious alternative notation in the aim of simplifying the equations, writing the index $l$ between angular brackets to indicate in which direction the average is taken. Then, in the three dimensional and structured Cartesian case, one has
\begin{equation} \begin{array}{ll}
 &\Avex{X}{i+\halb}{j}{k}{n} := \mathcal{A}\left( X_{i,j,k}^n, X_{i+1,j,k}^n \right)  \equiv \halb \left( X_{i,j,k}^n + X_{i+1,j,k}^n \right)\\ 
 &\Avey{X}{i}{j+\halb}{k}{n} := \mathcal{A}\left( X_{i,j,k}^n, X_{i,j+1,k}^n \right)  \equiv \halb \left( X_{i,j,k}^n + X_{i,j+1,k}^n \right),  \\
 &\Avez{X}{i}{j}{k+\halb}{n} := \mathcal{A}\left( X_{i,j,k}^n, X_{i,j,k+1}^n \right)  \equiv \halb \left( X_{i,j,k}^n + X_{i,j,k+1}^n \right), \end{array}
\end{equation}
that give an estimate of the observable $X$ on the \emph{half-indexed} discrete coordinate as an average between the values of $X$ at the corresponding right and left \emph{integer-indexed} discrete coordinate. Here the subscript of $\langle \cdot \rangle$ indicates the space-direction where the average is taken, while the subscript of $\{ \cdot \}$ indicates the discrete index where the average is centered.
As it has been stated above, without loss of rigor, a similar notation can be also used for averaging the observables from the half-indexed to the integer-indexed discrete coordinates, i.e.
\begin{equation} \begin{array}{lc}
  &\Avex{X}{i}{j}{k}{n} := \mathcal{A}\left( X_{i-\halb,j,k}^n, X_{i+\halb,j,k}^n \right)  \equiv   \halb \left( X_{i-\halb,j,k}^n + X_{i+\halb,j,k}^n \right), \\
  &\Avey{X}{i}{j}{k}{n} := \mathcal{A}\left( X_{i,j-\halb,k}^n, X_{i,j+\halb,k}^n \right)  \equiv   \halb \left( X_{i,j-\halb,k}^n + X_{i,j+\halb,k}^n \right), \\
  &\Avez{X}{i}{j}{k}{n} := \mathcal{A}\left( X_{i,j,k-\halb}^n, X_{i,j,k+\halb}^n \right)  \equiv  \halb \left( X_{i,j,k-\halb}^n + X_{i,j,k+\halb}^n \right).\end{array} 
\end{equation} 
and this averaging operator can also easily extended to multi space-dimensions, e.g.,
%
\begin{equation}\begin{array}{lc}
 &\Aveyz{X}{i}{j+\halb}{k+\halb}{n}  = \viertel \left( X_{i,j,k}^n + X_{i,j+1,k}^n + X_{i,j,k+1}^n + X_{i,j+1,k+1}^n \right), \\
 &\Avexz{X}{i+\halb}{j}{k+\halb}{n}  = \viertel \left( X_{i,j,k}^n + X_{i,j,k+1}^n +  X_{i+1,j,k}^n + X_{i+1,j,k+1}^n \right), \\
 &\Avexy{X}{i+\halb}{j+\halb}{k}{n}  = \viertel \left( X_{i,j,k}^n + X_{i,j+1,k}^n +X_{i+1,j,k}^n + X_{i+1,j+1,k}^n \right),  \end{array}
\end{equation} 
but also
\begin{equation} \begin{array}{lc}
  &\Aveyz{X}{i}{j}{k}{n}     = \viertel \left( X_{i,j-\halb,k-\halb}^n + X_{i,j+\halb,k-\halb}^n + X_{i,j-\halb,k+\halb}^n + X_{i+\halb,j,k+\halb}^n \right), \\
  &\Avexz{X}{i}{j}{k}{n}      = \viertel \left( X_{i-\halb,j,k-\halb}^n + X_{i+\halb,j,k-\halb}^n + X_{i-\halb,j,k+\halb}^n + X_{i-\halb,j,k+\halb}^n \right),  \\
  &\Avexy{X}{i}{j}{k}{n}      = \viertel \left( X_{i-\halb,j-\halb,k}^n + X_{i+\halb,j-\halb,k}^n + X_{i-\halb,j+\halb,k}^n + X_{i+\halb,j+\halb,k}^n \right),   \end{array} 
\end{equation} 
and all the possible obvious combinations.
Finally, since MHD is governed by a nonlinear system of partial differential equations (PDE), it turns to be useful to define the combination of two averaging operator acting on the product of two discrete scalar observables $X$ and $Y$ that belong to different meshes, e.g. let us define this product at the discrete location $({i+\halb},{j},{k},{n})$, having $X$ on the z-edges, e.g. $\Vdual_x\times\Vdual_y\times\Vmain_z$, and $Y$ on the barycenters, e.g. $\Vmain_x\times\Vmain_y\times\Vmain_z$,
\begin{equation} \begin{array}{rl}
  \AveyxTwo{X}{Y}{i+\halb}{j}{k}{n}   := 
	  \halb &\left[  \halb \left(X_{i ,j-\halb,k} + X_{i ,j+\halb,k}\right) Y_{i ,j,k} + \right.  \\
		 & \left.+ \halb \left(X_{i+1,j-\halb,k} + X_{i+1,j+\halb,k}\right) Y_{i+1,j,k}   \right] \end{array}\label{eq:ave2}
\end{equation} 
Notice that on the discrete variable $X$, two averages are taken, while on $Y$ only one. 

As stated above, we make use of staggered grids in order to build a numerical scheme able to preserve the divergence-free condition $\nabla \cdot \B = 0 $ exactly. In this sense, we are inspired by the work by Balsara and collaborators  \cite{BalsaraSpicer1999,Balsara2004}, defining the magnetic field on a staggered mesh, but, on the contrary, instead of defining it on the three face-based staggered grid, in this work \textbf{the magnetic field components are defined on the edges} 
and, consequently, \textbf{the velocity and electric field are defined on the faces}, see Sec.\ref{sec:B}.
 \redP{An alternative choice would be to define velocities and electric fields on the edges, and magnetic fields on the faces. This is certainly a viable option, 
at the expense of introducing a new projection/average operator to recalculate velocities on the faces. Indeed, the velocities on the faces are necessary for matching the conservative scheme for the adopted implicit acoustic solver, see Sect.\ref{sec:p}. On the other hand, Faraday's equation is solved in its non-conservative form, and thus predisposes to both types of staggering.
}

\subsubsection*{Discrete operators}

At this point, we define the discrete curl operator acting on the electric field and the discrete divergence operator acting on the magnetic field 
\begin{equation} \begin{array}{rl}
&\E_f\in \Vfaces {\xrightarrow[\quad\quad\quad]{\nabla_e\times}}  \Vedges  \ni \B_h  \\
 & \B_e\in   \Vedges  {\xrightarrow[\quad\quad\quad]{\nabla_n\cdot}} \Vnodes    \end{array}
\end{equation}
%
by means of, respectively, 
\begin{equation}\begin{array}{c}
\left[\left(\nabla \times \E\right)_x\right]_{i,j+\halb,k+\halb}  :=  \frac{(E_z)_{i,j+1,k+\halb} - (E_z)_{i,j,k+\halb}}{\Delta y}  - \frac{(E_y)_{i,j+\halb,k+1} - (E_y)_{i,j+\halb,k}}{\Delta z} , \\
\left[\left(\nabla \times \E\right)_y\right]_{i+\halb,j,k+\halb}  :=  \frac{(E_x)_{i+\halb,j,k+1} - (E_x)_{i+\halb,j,k}}{\Delta z}  - \frac{(E_z)_{i+1,j,k+\halb} - (E_z)_{i,j,k+\halb}}{\Delta x} ,  \\
\left[\left(\nabla \times \E\right)_z\right]_{i+\halb,j+\halb,k} :=   \frac{(E_y)_{i+1,j+\halb,k} - (E_y)_{i,j+\halb,k}}{\Delta x} - \frac{(E_x)_{i+\halb,j+1,k} - (E_x)_{i+\halb,j,k}}{\Delta y}.   \end{array} \label{eq:dex}
\end{equation}
and
\begin{equation}\begin{array}{rcl}
 \left( \nabla \cdot \B \right)_{i+\halb,j+\halb,k+\halb} := &&\frac{(B_x)_{i+1,j+\halb,k+\halb}-(B_x)_{i,j+\halb,k+\halb}}{\Delta x} +   \\ &+& \frac{(B_y)_{i+\halb,j+1,k+\halb}-(B_y)_{i+\halb,j,k+\halb}}{\Delta y}+    \\
	   &+& \frac{(B_z)_{i+\halb,j+\halb,k+1}-(B_z)_{i+\halb,j+\halb,k}}{\Delta z} . \end{array}
	\label{eq:disc.divb} 
\end{equation} 
 Moreover, by a direct application of the above definition,   one can easily check that the discrete divergence of the discrete curl of $\E$ is \emph{exactly zero}, i.e the following discrete equation is exact up to machine precision
$$\nabla_n \cdot \left(\nabla_e \times \A_f \right) = 0 $$ for any choice of discrete vector \redP{$\A_f \in \Vfaces$}.

Then, since the velocity $\v$ is defined on the faces, while the magnetic field $\B$  on the edges,  then the averaging index-operator (\ref{eq:ave2}) has to be used to give a proper definition of the cross product $(\v\times\B)_f$ on the faces. In particular, we define here the discrete cross product that maps 
$$(\v_f,\B_e) \in (\Vfaces \times \Vedges) {\xrightarrow[\quad\quad\quad]{\times}}  \Vfaces \ni (\v\times\B)_f$$
as
\begin{equation} \begin{array}{rl}
\Big[ \left(\v\times\B \right)_x  \Big]_{i+\halb,j,k}&:=  \AveTwo{v}{B_z}{i+\halb}{j}{k}{}{x}{y}	-   \AveTwo{w}{B_y}{i+\halb}{j}{k}{}{x}{z}\\
\Big[ \left(\v\times\B \right)_y  \Big]_{i,j+\halb,k}&:=  \AveTwo{w}{B_x}{i}{j+\halb}{k}{}{y}{z}-  \AveTwo{u}{B_z}{i}{j+\halb}{k}{}{y}{x} \\
\Big[ \left(\v\times\B \right)_z  \Big]_{i,j,k+\halb}&:=  \AveTwo{u}{B_y}{i}{j}{k+\halb}{}{z}{x}-  \AveTwo{v}{B_x}{i}{j}{k+\halb}{}{z}{y}\end{array}
\label{eq:vxBx}
\end{equation}
where all the averaging operators are needed only in the aim of defining the variables in the proper spatial discrete location. For the sake of clarity, the full expression of (\ref{eq:vxBx}) reported below:
\begin{equation} \begin{array}{rl}
&\Big[ \left(\v\times\B \right)_x  \Big]_{i+\halb,j,k}= \\ 
&      \halb\left[  \frac{\left(v_{i ,j-\halb,k} + v_{i+1 ,j-\halb,k}\right)}{2} (B_z)_{i+\halb ,j-\halb,k} + \frac{\left(v_{i,j+\halb,k} + v_{i+1,j+\halb,k}\right)}{2} (B_z)_{i+\halb,j+\halb,k}    \right]   + \\
&   -  \halb\left[  \frac{\left(w_{i ,j,k-\halb} + w_{i+1 ,j,k-\halb}\right)}{2} (B_y)_{i+\halb ,j,k-\halb} + \frac{\left(w_{i,j,k+\halb} + w_{i+1,j,k+\halb}\right)}{2} (B_y)_{i+\halb,j,k+\halb}    \right]   \end{array}
\end{equation}

Finally, the resistive terms are again a discrete curl, acting on the magnetic field lying on $\Vedges$, with values on $\Vfaces$, i.e.
\myequation{l}{
\B_e\in \Vedges{\xrightarrow[\quad\quad\quad]{ \nabla^*_f\times}}  \Vfaces \ni  \nabla^*_f\times \B_e ,}
{}
and they read, together with any eventual artificial stabilizing term $\mathbf{D}_f$, as

\begin{equation} \begin{array}{rl}
\Big\{ \left[\left( \eta \nabla^* +\mathbf{D}^x\right) \times\B \right]_x  \Big\}_{i+\halb,j,k}&:=  \left[(\eta + s_{y}^x) \partial_y B_z\right]_{i+\halb,j,k}-  \left[(\eta + s_{z}^x)\partial_z B_y\right]_{i+\halb,j,k}\\
\Big\{ \left[\left( \eta \nabla^* +\mathbf{D}^y\right) \times\B \right]_y  \Big\}_{i,j+\halb,k}&:=  \left[(\eta + s_{z}^y) \partial_z B_x\right]_{i,j+\halb,k}-  \left[(\eta + s_{x}^y)\partial_x B_z\right]_{i,j+\halb,k} \\\
\Big\{ \left[\left( \eta \nabla^* +\mathbf{D}^z\right) \times\B \right]_z  \Big\}_{i,j,k+\halb}&:=  \left[(\eta + s_{x}^z) \partial_x B_y\right]_{i,j,k+\halb}-  \left[(\eta + s_{y}^z)\partial_y B_x\right]_{i,j,k+\halb}\end{array} \label{eq:dBx}
\end{equation}
where the symbol $*$ is used only indicating that the discrete operator is different from what defined in (\ref{eq:dex}).  $\mathbf{D}_f$ is a purely artificial operator that acts on the magnetic field $\B$ as stabilization term.
Then, the dissipative coefficients $s$ depends on the local convective eigenvalues and are chosen to be 
\begin{equation}
s_x:=  \alpha \frac{\max(u) \Delta x}{2}, \quad s_y:= \alpha \frac{\max(v) \Delta y}{2}, \quad s_z:= \alpha \frac{\max(w) \Delta z}{2}. \label{eq:sxyz}
\end{equation}
Upper scripts in eq. (\ref{eq:dBx}) are used only to distinguish between the dissipative terms of the three different equations. In almost all the simulations presented in the numerical validation have been performed after choosing $\alpha = 1$.
\newline
\redP{In equation (\ref{eq:sxyz}), it has to be noticed that the maximum} is evaluated \emph{locally} and only between the neighbor cells in proximity of the reference discrete spatial position, on the proper staggered grid, e.g. on the $x$-oriented faces $\Vxfaces$ for $\max(u)$, $\Vyfaces$  for $\max(v)$ and $\Vzfaces$ for $\max(w)$.
As stated at the beginning of this section, spatial derivatives are taken central and then, there is only one possible interpretation of (\ref{eq:dBx}) at the discrete level, e.g. for the $x$-component
\begin{equation}\begin{array}{rl}
\Big\{ \left[\left( \eta \nabla +\mathbf{D}^x\right) \times\B \right]_x  \Big\}_{i+\halb,j,k} 
\equiv&  (\eta + s_{y}^x)   \left[\frac{ (B_z)_{i+\halb,j+\halb,k} - (B_z)_{i+\halb,j-\halb,k}  }{\Delta y}  \right] +  \\
 & -  (\eta + s_{z}^x) \left[\frac{ (B_y)_{i+\halb,j,k+\halb} -  (B_y)_{i+\halb,j,k-\halb}}{\Delta z}\right] \end{array} \label{eq:vxBx2}  
\end{equation}
 where the magnetic resistivity $\eta$ is, in general, a space-dependent non-negative scalar field. Here, it has to be noticed that the components of the stabilizing terms are those of an approximated multi-dimensional Riemann solver (see e.g. \cite{balsarahlle2d,balsarahllc2d,balsarahlle3d,BalsaraMultiDRS,MUSIC1,MUSIC2}) but, as a direct extension of what used in \cite{SIMHD}, in this work they will be solved \emph{implicitly in time}.
 The artificial dissipating term is in this case a standard Rusanov (or local Lax-Friedrich) flux. %
It is also interesting to note that it eventually resembles a purely resistive term 
$$ \left(\mathbf{D}^x , \mathbf{D}^y, \mathbf{D}^z\right) \times\B   \longrightarrow \eta_h \nabla \times \B$$ 
after choosing $\eta_h = s_x = s_y = s_z$, with (numerical) resistivity proportional to the characteristic grid-size $h$, i.e. $\eta_h = O( h)$, and, consequently, leading to a standard hyperbolic-CFL time restriction.  
The choice $s^x_y = s^x_z$, $s^y_z=s^y_x$ and $s^z_x=s^z_y$ can be used to ensure the symmetry of the resulting algebraic system that,  in the here presented algorithm, has to be finally  solved for the $\B$ variables.

\subsection{Explicit terms: nonlinear convection and viscous subsystems} 
\label{sec:Fv}

The chosen procedure for discretizing  the nonlinear convection and the viscous terms that are summarized in the first subsystem  
\begin{equation}
\partial_t \Q + \nabla \cdot \left( \Fv - \Fdv \right) = 0,
\begin{array}{l}

\end{array}
\end{equation}
see eq. (\ref{eq:PDEex}), 
follows the one chosen in the previous paper \cite{SIMHD}. A conservative explicit finite-volume scheme of the type 
\begin{equation} 
\begin{array}{rl}
   \Q_{i,j,k}^* = \Q_{i,j,k}^n &- \frac{\Delta t}{\Delta x} \left( \mathbf{f}_{i+\halb,j,k} - \mathbf{f}_{i-\halb,j,k} \right) - \frac{\Delta t}{\Delta y} \left( \mathbf{g}_{i,j+\halb,k} - \mathbf{g}_{i,j-\halb,k} \right) +  \\
	&	- \frac{\Delta t}{\Delta z} \left( \mathbf{h}_{i,j,k+\halb} - \mathbf{h}_{i,j,k-\halb} \right),    \end{array} 	\label{eq:expfv3d} 
\end{equation}
is adopted, where the star symbol $^*$ is used to indicate that $\Q^*$ is only a local solution of a sub-system (\ref{eq:PDEex}). In particular, one has   numerical (Rusanov, or local Lax-Friedrichs) fluxes of the type
\begin{equation}\begin{array}{rl}
 \mathbf{f}_{i+\halb,j,k} :=  & \halb \left( \mathbf{f}^v( \mathbf{Q}_{i+\halb,j,k}^-) + \mathbf{f}^v( \mathbf{Q}_{i+\halb,j,k}^+) \right) -  \Aveyz{\f^{\mu}}{i+\halb}{j}{k}{}  +  \\
			 & -  \halb s_{\max}^x \left(\mathbf{Q}_{i+\halb,j,k}^+  -  \mathbf{Q}_{i+\halb,j,k}^- \right)   \end{array} \label{eq:Rusanov} %
\end{equation} 
where $\mathbf{f}^{\mu}_h := \mathbf{f}^{\mu}(\mathbf{V}_h, \nabla \mathbf{V}_h)$
and, according to the definition of the averaging operators, 
\begin{equation} \begin{array}{rl} 
		 \Aveyz{\f^{\mu}}{i+\halb}{j}{k}{} = \viertel &  \left[ \mathbf{f}^{\mu}_{i+\halb,j+\halb,k+\halb} + \mathbf{f}^{\mu}_{i+\halb,j-\halb,k+\halb}+ \right. \\ 
		 &\left. + \mathbf{f}^{\mu}_{i+\halb,j+\halb,k-\halb} +  \mathbf{f}^{\mu}_{i+\halb,j-\halb,k-\halb} \right] \end{array}
\end{equation} 
Here, the viscous terms $\f^{\mu}$ depend on the value of the primitive variables and its gradients evaluated on the nodes of the grid, while the convective terms depend on the conserved variables reconstructed on the faces.

In the finite-volume framework, the numerical solution is piecewise constant with the discontinuities lying on the faces of the control volumes. For this reason, exact or approximate Riemann solvers \cite{toro-book} are needed to define a numerical flux at the cell interfaces. Several choices can be adopted, but for simplicity, we made use of the Rusanov (or local Lax-Friedrich) flux. Then, $\Q_{i+\halb}^-$ and $\Q_{i+\halb}^+$ are the conserved variables that have been  reconstructed  at the faces, from the left and from the right, respectively. In the numerical validation we always adopted a standard first order reconstruction, i.e. $ \Q_{i+\halb}^{\pm} := \Q_{i\pm 1}$ , or, if specified in the description of the numerical results, a second-order MUSCL-Hancock method, as it has been described in \cite{toro-book,SIMHD}.

Here, $s_{\text{max}}$ is the maximum of the  signal speed  of the convective subsystem in the $x$ direction at the cell faces, defined as $s_{\max} = \max \limits_l \left( |\lambda^c_l(\Q_{i+\halb}^-) |, |\lambda^c_l(\Q_{i+\halb}^+) | \right)$, that takes account also of the viscous term. The extension to the fluxes in the remaining space-direction, i.e. $\g_{i,j+\halb,k}$ and $\h_{i,j,k+\halb}$, is obvious and omitted in the aim of brevity.

The computational time-step is chosen accordingly to the rule
\begin{equation}
 \Delta t = \CFL \frac{1}{ \frac{\max |\lambda^h_x|}{\Delta x} + \frac{ \max |\lambda^h_y|}{\Delta y}+ \frac{ \max |\lambda^h_z|}{\Delta z} + 2 \lambda^p \left( \frac{1}{\Delta x^2} + \frac{1}{\Delta y^2} + \frac{1}{\Delta z^2} \right) }, \label{eq:CFL}
\end{equation} 
with Courant number $\CFL < 1$, after choosing a proper set of 'hyperbolic'-eigenvalues $\lambda^h_s$, $s=x,y,z$  , and a proper 'parabolic'-penalty value $\lambda^p>0$. Numerical stability is usually obtained by choosing $\lambda^h_s$ 
to be the set of eigenvalues of the Jacobian of the fluxes in the $s$-th direction, and $\lambda^\mu$ the  penalty term that refers to the parabolic terms of the PDE system, both of them referred to those 
terms of the PDE that are discretized explicitly in time. Since we consider viscosity and heat-conduction, we chose $\lambda^\mu=(4/3)\mu /\rho  +\kappa/(c_v \rho)$. 
The choice of $\lambda^h_s$ may vary between the set of eigenvalues $\lambda^{\text{MHD}}$ of the full MHD system, see (\ref{eq:eval.full}), the explicit terms of the 2-splitting scheme $\lambda^{v\text{B}}$ 
\begin{equation}
 \lambda^{v\text{B}}_{1,8} = v_x \mp \sqrt{\frac{\mathbf{B}^2}{4 \pi \rho}}, \qquad  
 \lambda^{v\text{B}}_{2,7} = v_x \mp \frac{B_x}{\sqrt{ 4 \pi \rho }}, \qquad 
 \lambda^{v\text{B}}_{3,4} = 0, \qquad 
 \lambda^{v\text{B}}_{5,6} = v_x, 
	\label{eqn.eval.c} 
\end{equation}
see reference \cite{SIMHD}, or the explicit terms of the 3-splitting scheme $\lambda^v$, see (\ref{eq:1db}). 
In this way, because of the property 
$$  \max{(\lambda^v)}\leq \max{(\lambda^{v\text{B}})}\leq \max{(\lambda^{\text{MHD}})}$$
it is possible to estimate the following time-scales: 
$$\Delta t (\lambda^v) \geq \Delta t {(\lambda^{v\text{B}})}\geq \Delta t {(\lambda^{\text{MHD}})}.$$
Depending on the problem to be solved, one can change the choice of eigenvalues and $\CFL$ number in order to increase the resolution in time when needed. For some specific physical problems, the time-step chosen with respect to the split system, i.e. $\Delta t(\lambda^v)$, could be larger than the main time-scales of the problem. In such cases, if the numerical diffusion is low enough, then spurious oscillation may arise and can deteriorate the numerical solution. In these cases, one can consider the possibility of introducing some more numerical diffusion. To do this, there are essentially two main direct approaches: first, by using the intermediate set of eigenvalues $\lambda^{v\text{B}}$ or directly the full set $\lambda^{\text{MHD}}$ in eq. (\ref{eq:Rusanov}), both for the estimate of the time-step and the computation of the numerical fluxes as an \emph{explicit} stabilization term; second,  by adding an \emph{implicit} stabilization within the pressure and the magnetic-field solvers,  e.g. depending on the corresponding sub-set of eigenvalues $\lambda^p$ and $\lambda^B$ respectively, see eq. (\ref{eq:1db}-\ref{eq:1dp}).

If not specified otherwise, we set $\CFL=0.9$ in all test problems presented in this section and use the pure-convective eigenvalues $\lambda^v$, see (\ref{eq:1db}). Furthermore, for all test 
cases we have explicitly verified that up to machine precision the magnetic field is divergence-free and mass, momentum and energy are conserved. 

\subsection{Implicit and divergence-free integration of the Faraday equation.}
 \label{sec:B}

Since a splitting scheme is adopted, in order to apply an implicit time discretization to the magnetic field, a reduced system of equations can be derived by applying the splitting-scheme (\ref{eq:PDEimB}), reducing the Faraday and momentum equations (\ref{eq:momentum}-\ref{eq:Faraday}) to
\begin{equation}
\left\{\begin{array}{l}
\partial_t (\rho \v)  - \frac{1}{4\pi}\left( \nabla \times \B \right) \times \B = 0    \\
  \partial_t \B  + \nabla \times \mathbf{E} = 0,  
\end{array}\right. \label{eq:Bpde}
\end{equation}

\noindent
Then, a consistent and divergence-free discrete form of (\ref{eq:Bpde}) reads as
\begin{equation} \begin{array}{l}
   (B_x)_{i,j+\halb,k+\halb}^{n+1} =  (B_x)_{i,j+\halb,k+\halb}^n   - \Delta t    \Big[\left(\nabla \times \E\right)_x\Big] _{i,j+\halb,k+\halb}^{n+1}   ,  \\ 
   (B_y)_{i+\halb,j,k+\halb}^{n+1} = (B_y)_{i+\halb,j,k+\halb}^n     - \Delta t   \Big[\left(\nabla \times \E\right)_y\Big] _{i+\halb,j,k+\halb}^{n+1} ,    \\ 
   (B_z)_{i+\halb,j+\halb,k}^{n+1} = (B_z)_{i+\halb,j+\halb,k}^n      - \Delta t  \Big[\left(\nabla \times \E\right)_z \Big]_{i+\halb,j+\halb,k}^{n+1} . 	\end{array} \label{eq:bx}   
\end{equation}
while, for the momentum equation 
\begin{equation} \begin{array}{l}
    u _{i+\halb,j,k}^{n+1} =   u_{i+\halb,j,k}^*   + \Delta t     \Big\{ \frac{1}{4\pi\rho}\left[\left( \nabla \times \B \right) \times \B \right]_x \Big\}_{i+\halb,j,k}^{n+1},   \\ 
    v _{i,j+\halb,k}^{n+1} =   v_{i,j+\halb,k}^*     + \Delta t   \Big\{ \frac{1}{4\pi\rho}\left[\left( \nabla \times \B \right) \times \B \right]_y \Big\}_{i,j+\halb,k}^{n+1}, 	  \\ 
    w _{i,j,k+\halb}^{n+1} =   w_{i,j,k+\halb}^*      + \Delta t  \Big\{ \frac{1}{4\pi\rho}\left[\left( \nabla \times \B \right) \times \B \right]_z \Big\}_{i,j,k+\halb}^{n+1}. 	\end{array} 	\label{eq:vx} 
\end{equation}
 Here, as a reminder, the magnetic field is discretized on the edges (two halfed indexes), while  the velocity and electric field are discretized on the faces (one single halfed index). 
The upper-script $^*$ is used to remember that this is the second stage of a splitting scheme, and 
the corresponding discrete variables, namely the velocities, have been already evolved through 
the aforementioned explicit scheme for the  nonlinear convection. For the sake of brevity, one 
can introduce a more compact matrix-vector notation and rewrite (\ref{eq:bx}-\ref{eq:vx}) as 
\begin{equation}  \begin{array}{rl}
   &\B_{\text{e}}^{n+1} =  \B_{\text{e}}^n   - \Delta t    \left( \nabla \times \E_{\text{f}} \right)_{\text{e}}^{n+1}   , 	 \\ 
   & \v_{\text{f}}^{n+1} =   \v_{\text{f}}^*   + \Delta t     \left[\frac{1}{4\pi\rho}\left( \nabla^* \times \B_{\text{e}} \right)_{\text{f}} \times \B_{\text{e}} \right]_{\text{f}}^{n+1},  \\ 
		& 
 \mathbf{E}_{\text{f}} = - \left(\mathbf{v}_{\text{f}} \times \mathbf{B}_{\text{e}} + \eta \nabla^* \times \mathbf{B}_{\text{e}}\right)_{\text{f}}. \end{array}
\label{eq:Eh3dwrth} 
\end{equation}
Here, $\B_e$ and $\v_f$ are the arrays containing the full set of degrees of freedom for the three magnetic field and, correspondingly, velocity components. 

Notice that, even without specifying how the electric field components are computed, one can show that, \emph{if} the discrete divergence free condition is valid at one time step $t^n$, then, by evolving in time the magnetic field components according to (\ref{eq:bx}), then the divergence-free condition is exactly preserved by construction, (see equation \eqref{eq:dex}) at time $t^{n+1}$, i.e.
\begin{equation}
\left( \nabla \cdot \B \right)_{i+\halb,j+\halb,k+\halb}^n = 0 \quad \quad \Longrightarrow  \quad \quad \left( \nabla \cdot \B \right)_{i+\halb,j+\halb,k+\halb}^{n+1} = 0.
\end{equation} 
 
Then, for defining the electric field components, a discrete form of (\ref{eq:Evector}) is derived as 
\begin{equation} \begin{array}{rl}
  E^x_{i+\halb,j,k}  :=  &  - \Big[ \left(\v\times\B \right)_x  \Big]_{i+\halb,j,k} +  \Big\{ \left[   \left(\eta \nabla + \mathbf{D}^x \right) \times \B \right]_x \Big\}_{i+\halb,j,k}  \\
  E^y_{i,j+\halb,k}  :=  &  - \Big[ \left(\v\times\B \right)_y  \Big]_{i,j+\halb,k} +  \Big\{ \left[   \left(\eta \nabla + \mathbf{D}^y \right) \times \B \right]_y \Big\}_{i,j+\halb,k}  	  \\
  E^z_{i,j,k+\halb}  :=  &  - \Big[ \left(\v\times\B \right)_z  \Big]_{i,j,k+\halb} +  \Big\{ \left[   \left(\eta \nabla + \mathbf{D}^z \right) \times \B \right]_z \Big\}_{i,j,k+\halb}  	\end{array} 	\label{eq:ex} 
\end{equation} 


Because of the nonlinearity of the MHD equations, a fully implicit method may lead to a system that is very difficult to be solved. For this reason, the MHD equations are often linearized, dealing 
with the so-called \emph{linearized MHD}. In this work, the original \textbf{nonlinear  viscous and resistive MHD equations are solved}. 
Indeed, a fixed-point Picard procedure has been adopted in order to take care of the nonlinearity of the PDE system. For the sake of brevity, we can rewrite (\ref{eq:Eh3dwrth}) introducing a new index $r$ labeling the index of the Picard iteration, i.e.
\begin{align}  
   &\blue{\B_{\text{e}}^{n+1,{r+1}}} =  \B_{\text{e}}^{n,{r+1}}   + \Delta t    \left\{ \nabla \times \left[ \blue{\v_{\text{f}}^{r+1}} \times \B_{\text{e}}^{r} - \left( \eta \nabla^* + \mathbf{D}\right) \times \blue{\B_{\text{e}}^{r+1}} \right]_{\text{f}} \right\}_{\text{e}}^{n+1}    	\label{eq:bh3d} \\
   & \blue{\v_{\text{f}}^{n+1,{r+1}}} =   \v_{\text{f}}^*   +     \frac{\Delta t }{4\pi\rho^{n+1}_{\text{f}}}\left[\left( \nabla^* \times \blue{\B_{\text{e}}^{r+1}} \right)_{\text{f}} \times \B_{\text{e}}^r \right]_{\text{f}}^{n+1}, 	\label{eq:vh3d} 
\end{align}
for $r=1$, $\ldots$, $R$, where  $R$ is the maximum number of Picard iterations and it is chosen equal to $R=2$. 
With the aim of simplifying the notation, any time-index $n+1$ labeling a right-bracket is inherited by inner terms. 
The definition of electric field $\E$ has been already substituted in the Faraday equation. 
 Notice that, in this way, the system has been only locally linearized. Here, blue terms corresponding to iteration index $r+1$ are the real unknown of the locally linearized system.  Then, after substituting  (\ref{eq:vh3d}) into (\ref{eq:bh3d}), a decoupled system for the only variables $\B^{n+1,r+1}_{\text{e}}$ has been obtained, and it reads
\begin{equation} \begin{array}{rl}
  & \blue{\B_{\text{e}}^{n+1,{r+1}} } +  \\
	& + \Delta t^2    \left\{ \nabla \times \left[   \frac{\left[\left( \nabla^* \times \blue{\B_{\text{e}}^{r+1}} \right)_{\text{f}} \times \B_{\text{e}}^r \right]_{\text{f}}}{4\pi \rho^{n+1}_{\text{f}}} \times \B_{\text{e}}^{r}  - \left( \eta \nabla^* + \mathbf{D}\right) \times  \blue{\B_{\text{e}}^{r+1}} \right]_{\text{f}} \right\}_{\text{e}}^{n+1} = \\ 
	&=  \B_{\text{e}}^{n,{R+1}} + \Delta t^2    \left\{ \nabla \times \left[   \frac{\left[ \v_{\text{f}}^* \times \B_{\text{e}}^{n+1,r} \right]_{\text{f}}}{\rho^{n+1}_{\text{f}}} \times \B_{\text{e}}^{n+1,r}  \right]_{\text{f}} \right\}_{\text{e}}, \end{array}	\label{eq:vhryi3d} 
\end{equation}
 At the discrete level, after collecting all the degrees of freedom of the field variables, equation (\ref{eq:vhryi3d}) can be written in matrix-array form as 
\begin{equation}
\left[\hat{\one}_e + \hat{\mathbb{H}}(\mathbf{B}_e^{n+1,r}) + \ddelta \right] \cdot \blue{\B_e^{n+1,r+1} } = \mathbf{r}_e, \quad \quad r=1,2
\end{equation}
that is a linear system of $N_{e}$ equations, for $N_e$ unknowns, where $N_e$ is the total number of edges.
Here, the following definitions have been assumed:
\begin{itemize}
\item $\hat{\one}_e$ is the identity matrix;
$$\hat{\one}_e \cdot \blue{\B_e^{n+1,r+1} } := \blue{\B_e^{n+1,r+1} } $$
\item $\B_e$ is the array collecting all the spatial degrees of freedom of the discrete magnetic field, i.e. the values at the edges of the grid at the prescribed time-step and Picard-iteration; the elements of $\B_e^{n+1,r}$ are the \emph{known} components at the $(n+1)$-th time-step  and $r$-th Picard iteration , while $\B_e^{n+1,r+1}$ collects all the \emph{unknown} components at the same time-step, for the next $(r+1)$-th Picard-iteration;
\item $\hat{\mathbb{H}}(\B^{n+1,r})$ can be called a discrete 'Faraday-Lorentz' operator, arising from the coupling between the Lorentz-force appearing in the momentum equation and the Faraday equation; it represents the second term at the left-hand side of equation (\ref{eq:vhryi3d}) and, since it depends on the magnetic field at the \emph{previous} Picard iteration $r$, it can be regarded as a constant-coefficient matrix, whose coefficients may only change between the Picard iterations;
$$\hat{\mathbb{H}}(\B^{n+1,r})\cdot \blue{\B_e^{n+1,r+1} }:= \Delta t^2      \nabla \times \left[   \frac{\left[\left( \nabla^* \times \blue{\B_{\text{e}}^{r+1}} \right)_{\text{f}} \times \B_{\text{e}}^r \right]_{\text{f}}}{4\pi \rho^{n+1}_{\text{f}}} \times \B_{\text{e}}^{r}\right]^{n+1}_e $$
\item $\ddelta$ can be called a discrete 'Faraday-Ohm' operator,  arising from the coupling between the Ohm law,  appearing in the definition of the electric field in the MHD approximation, and the Faraday equation; it collects all the physical and, for compactness, also the numerical dissipative terms (i.e. the magnetic resistivity and numerically stabilizing terms);
$$\ddelta \cdot \blue{\B_e^{n+1,r+1} }:=  - \Delta t^2      \nabla \times \left[   \left( \eta \nabla^* + \mathbf{D}\right) \times  \blue{\B_{\text{e}}^{r+1}}  \right]^{n+1}_e $$
\item $\mathbf{r}_e$ is the array collecting all the known terms of the system, i.e. the right hand side of equation (\ref{eq:vhryi3d}).
\end{itemize}
This system, as a direct consequence of the chosen discretization, has been shown to be \textbf{symmetric} 
 and can thus be solved within rather few iterations of an efficient and standard \emph{matrix-free} conjugate-gradient method. In numerical experiments the system appears so well behaved that it was not necessary to employ any preconditioning strategy. 
Once the components of the magnetic field $\B^{n+1,r+1}$ are obtained, the momentum is updated following, instead of (\ref{eq:vx}), the corresponding conservative form through
\begin{equation}
\partial_t (\rho \v) + \nabla \cdot \left[\F^{(\rho \v)} - \Fd^{(\rho \v)}\right] = 0 \label{eq:CVel}
\end{equation}
using the same explicit scheme (\ref{eq:expfv3d}) that was used for the first (explicit) sub-system (\ref{eq:PDEex}), but only after evaluating the fluxes as a function of the new temporary variable
\myequation{l}{
 \tilde\Q := \left[\rho^n, (\rho \v)^n, (\rho E)^n, \B^{n+1,r+1} \right]^T,
}{\label{eq:Qtilde}}
see eq. (\ref{eq:1dalgorithm2}), i.e. using the new computed magnetic field $\B^{n+1,r+1}$. %

\subsection{Pressure subsystem} 
\label{sec:p}

In this work, the method outlined in Sec.  2.4.2 of \cite{SIMHD} has been extended to the third space-dimension. The guidelines for the discretization are: momentum $\rho \v_f$ and the enthalpy $h_f=h_f(p)$ are  discretized along the three overlapping face-based staggered meshes, the pressure $p_b$ and the energy $\rho E$ are defined on the barycenters of the main grid. Then,  one can rewrite the pressure subsystem (\ref{eq:PDEimP}) as a coupled mildly nonlinear system for the unknowns $(\rho\v)^{n+1}_f$ and $p$.
\begin{equation} 
\text{(III)}\hspace{1cm} 
\left\{ \begin{array}{rl}
   & \partial_t(\rho \v) + \nabla p =0,  \\
	& \partial_t (\rho E) +   \nabla  \cdot \left( h \rho \v \right) =0\end{array} 	 \right.
\end{equation}

Since the main used algorithm would look almost the same of \cite{SIMHD}, and this is not the main part of this research, we restrict ourselves to depict the algorithm only in its matrix-vector notation, i.e.
\begin{equation} \begin{array}{rl}
   & (\rho \v)_{\text{f}}^{n+1} =   \tilde{(\rho \v)}_{\text{f}}    -     \Delta t \left(\nabla \p_{\text{b}} \right)_{\text{f}}^{n+1},  \\
	& (\rho E)_b^{n+1} = \tilde{(\rho E)}_b -   \Delta t \left[ \nabla  \cdot \left( h \rho \v \right)_f\right]_b^{n+1} \end{array} 	\label{eq:vh3dp} 
\end{equation}
 The tilde script is used in relation to $\tilde{\Q}$ in eq. (\ref{eq:Qtilde}).
Then, since the enthalpy is in general a nonlinear function of the pressure, a Picard iteration is introduced, see \cite{SIMHD}, having
\begin{equation}\begin{array}{rl}
 &
\begin{array}{l}   \blue{(\rho \v)_{\text{f}}^{n+1,r,{s+1}}} =   \tilde{(\rho\v)}_{\text{f}}  -     \Delta t \left(\nabla \blue{\p_{\text{b}}^{r,s+1}} \right)^{n+1}_{\text{f}}, \\
	  \left[\rho^{n+1} e(\blue{\p^{n+1,r,s+1}} ,\rho^{n+1} ) \right]_b+    \Delta t \left[ \nabla  \cdot \left( h^{r,s} \blue{(\rho \v)^{r,s+1}} \right)_f\right]^{n+1}_b = \mathbf{d}^{n+1,r,s}_b \end{array}
		\\
  &\mathbf{d}^{n+1,r,s}_b := \mathbf{\tilde{(\rho E)}}_b - \mathbf{m}_b^{n+1,r+1} - (\mathbf{\rho k})^{n+1,r,s}_b  
	\end{array}\label{eq:vh3dpb}
\end{equation}
where $s$ is a new recursion index $s=1$,$\ldots$,$S=2$; one may recall the definition of the energy density $\rho E=\rho e + \rho k + m$, where $k=  \rho \v^2/2$ is the  kinetic energy,  $m=\B^2/(8\pi)$ is the magnetic energy. 
Notice that the $r$ index refers to the fact that this stage of the algorithm is \emph{nested} within the Picard method outlined in the previous section, see also the simplified algorithm in eq. (\ref{eq:1dalgorithm2}). Here, the real unknowns refers to the only $s+1$ index and are highlighted in blue (see colored version online).
Then, the first discrete equation is substituted in the second one, resulting in 
\begin{equation}
	  \left[\rho^{n+1} e(\blue{\p^{n+1,r,s+1}} ,\rho^{n+1} ) \right]_b-     \Delta t^2 \left\{ \nabla  \cdot \left[ h^{n+1,r,s}_f   \left(\nabla \blue{\p_{\text{b}}^{n+1,r,s+1} }\right)_{\text{f}} \right]\right\}_b = \mathbf{r}^{n+1,r,s}_b,  	
\end{equation}
where
$
  \mathbf{r}^{n+1,r,s}_b := \mathbf{d}^{n+1,r,s}_b   -   \Delta t [ \nabla  \cdot ( h^{n+1,r,s} \tilde{(\rho\v)}_{\text{f}} )_f]_b$
	that is a decoupled mildly nonlinear \footnote{The system becomes nonlinear only for general nonlinear equations of state, where the internal energy density is a nonlinear function of the  pressure. For the ideal gas EOS, the final pressure system is \textit{linear}.} system of $N_b$ equations, and $N_b$ unknowns, where $N_b$ is the total number of pressure points, i.e. the control volumes of the main grid. In compact form, it is
\begin{equation} \begin{array}{rl}
	&  \hat{\mathbb{I}} \rho e (\blue{\p^{n+1,r,s+1}})_b + \hat{\mathbb{L}}[\p_{\text{b}}^{n+1,r,s} ] \cdot \blue{\p_{\text{b}}^{n+1,r,s+1} } = \mathbf{r}^{n+1,r,s}_b,  	\end{array} \label{eq:psyscomp}
\end{equation}
where $\mathbb{I}$ is the identity matrix, operator $\mathbb{L}[\p_{\text{b}}^{n+1,r,s} ]$ shares the properties of a discrete Laplace operator, namely it is defined as the discrete divergence of the discrete gradient, i.e.
$$
\mathbb{L}[\p_{\text{b}}^{n+1,r,s} ] \cdot \blue{\p_{\text{b}}^{n+1,r,s+1}}:=  -     \Delta t^2 \left[ \nabla  \cdot   \left(h^{n+1,r,s}  \nabla \blue{\p_{\text{b}}^{n+1,r,s+1}}\right)_{\text{f}}\right]_b
$$ 
 it is symmetric and at least positive semi-definite. The coefficients of $\mathbb{L}$ depends on the pressure evaluated at the previous Picard iteration through the enthalpy $h=h(p)$, and then it can be considered as a constant coefficient matrix. This means that the system for the pressure \eqref{eq:psyscomp} is a mildly nonlinear with a linear part that is symmetric and as least positive semi-definite. As stated in \cite{SIMHD}, with the usual assumptions on the nonlinearity detailed in \cite{CasulliZanolli2012}, it is efficiently solved with 
the nested Newton method of Casulli and Zanolli \cite{CasulliZanolli2010,CasulliZanolli2012}, combined, in this work, with a \emph{matrix-free} conjugate gradient method, without using any preconditioner. Once system (\ref{eq:psyscomp}) is solved, the discrete momentum and energy can be directly updated through the conservative scheme (\ref{eq:vh3dp}). At the end of this \emph{inner} $s$-Picard loop, the updated variables are relayed to the \emph{external} $r$-Picard loop by setting $(\rho\v,p)^{n+1,r+1} = (\rho\v,p)^{n+1,r,S+1}$. 

In the \emph{low Mach} limit, this algebraic system (\ref{eq:psyscomp}) is known to converge to a pressure Poisson equation that is typical of almost any incompressible flow solver, see \cite{SIMHD,DumbserCasulli2016}. In this sense, the here proposed scheme is also said to be \emph{low-Mach preserving}.
 
\noindent 
The final algorithm may be summarized in the following diagram:
	\begin{equation} 
\left\{ 
				%
				%
%
\begin{array}{lll}
							\texttt{\underline{explicit scheme}:}								 &\texttt{solve for } \rho^{n+1} \texttt{, and a first estimate of } \Q &  \\
						 							 & & \\
															 \begin{array}{c}\text{\emph{ Convection-diffusion}:} \\ \text{\emph{(see Sec.\ref{sec:Fv})}}\end{array}  & \begin{array}{l}\partial_t \Q + \nabla \cdot \left( \Fv - \Fdv \right) = 0   \end{array} &  \\
						 							 & & \\
							\texttt{\underline{implcit scheme}:}& \texttt{two nested Picard loops for } (\rho\v, p,\B)^{n+1} 				  & \\
						 							 & &   \\ 
	 \quad \begin{array}{r}	\text{$\B$\emph{-subsystem: (see Sec.\ref{sec:B})}}	\\ \\ \text{$p$\emph{-subsystem: (see Sec.\ref{sec:p})}}  \end{array}				 &
											 	\left. \begin{array}{l}
			 									  \partial_t \Q + \nabla \cdot \left( \Fb -\Fdb \right) = 0       \\\\
			 									  \left.\partial_t \Q + \nabla \cdot \Fp = 0  \phantom{\int_\frac{1}{1}^\frac{1}{1}\frac{1}{1}}\color{green}\right]^{\text{inner loop:}}_{\text{$p$ and $\rho\v$}}																																	 	 
																					\end{array}	\color{red}\right]^{\text{external loop:}}_{\text{ $\B$ and $\rho\v$}		}																														& \begin{array}{l}   \end{array} 
\end{array}\right. \label{eq:finalalg}
\end{equation}

\section{Numerical validation}
\label{sec:numval}
This section is dedicated to the validation process of the novel three-split scheme presented in this paper. The new method is a structure preserving semi-implicit method that treats the Alfv\'en waves and the magneto-sonic waves implicitly, while the nonlinear convection and the diffusion terms are discretized with an explicit conservative Godunov-type FV scheme. If not stated otherwise, for all simulations the time-step is evaluated accordingly to the CFL condition (\ref{eq:CFL}), using the eigenvalues of the pure convective subsystem, i.e. $dt=dt(\lambda^v)$, and, obviously, also the penalty terms due to the explicit integration of the viscous terms.

The two- and three-dimensional tests have been selected in order to verify the robustness of the method, and the eventual gain in terms of the \emph{effective Courant number} $dt/dt_{\text{MHD}}$, where $dt_{\text{MHD}}:= dt(\lambda_{\text{MHD}})$ is the time-step evaluated with respect to the full MHD eigen-spectrum. 

We emphasize again that in this paper a completely new way of staggering the discrete (dual or primal) differential operators, as well as the electric and magnetic field quantities, has been used. The here proposed algorithm is often labeled as "`SI-p\&B"', being a Semi-Implicit method which is implicit in the pressure terms $p$ (acoustic solver) and in the magnetic field $B$ (Alfvénic solver). This may be distinguished by the "`SI-p"' method, that is implicit only in the pressure terms (acoustic solver), see \cite{SIMHD}.

\subsection{Linear stability analysis}

A mathematically rigorous analysis of the linear and nonlinear stability properties of the new scheme proposed in this paper is beyond the scope of this paper and needs to be carried out separately in a different paper.  However, in this paper we want to give at least a numerical proof of linear stability of the scheme. For this reason, in this section we present the results of a linear stability analysis that is performed numerically over the full numerical scheme. 

The final fully-discrete scheme can be cast into the following nonlinear system
\begin{align}
\Q_h^{n+1} = \mathcal{L}(\Q_h^n).
\end{align}
Then, after linearizing the right-hand-side around an equilibrium solution $\Q_h^{\text{eq.}}$, 
 assuming  a regime of small perturbations,
 one can study the following linear system for the variations of the discrete solution, i.e.
\begin{align}
\delta \Q_h^{n+1} =  \left(\frac{\delta \mathcal{L}}{\delta \Q_h^n}\right)_{|_{\Q_h^{\text{eq.}}}} \cdot \delta \Q_h^n. \label{eq:MSA}
\end{align}

\begin{figure}[!t]
\begin{center}
\begin{tabular}{cc}  
\includegraphics[width=0.45\textwidth]{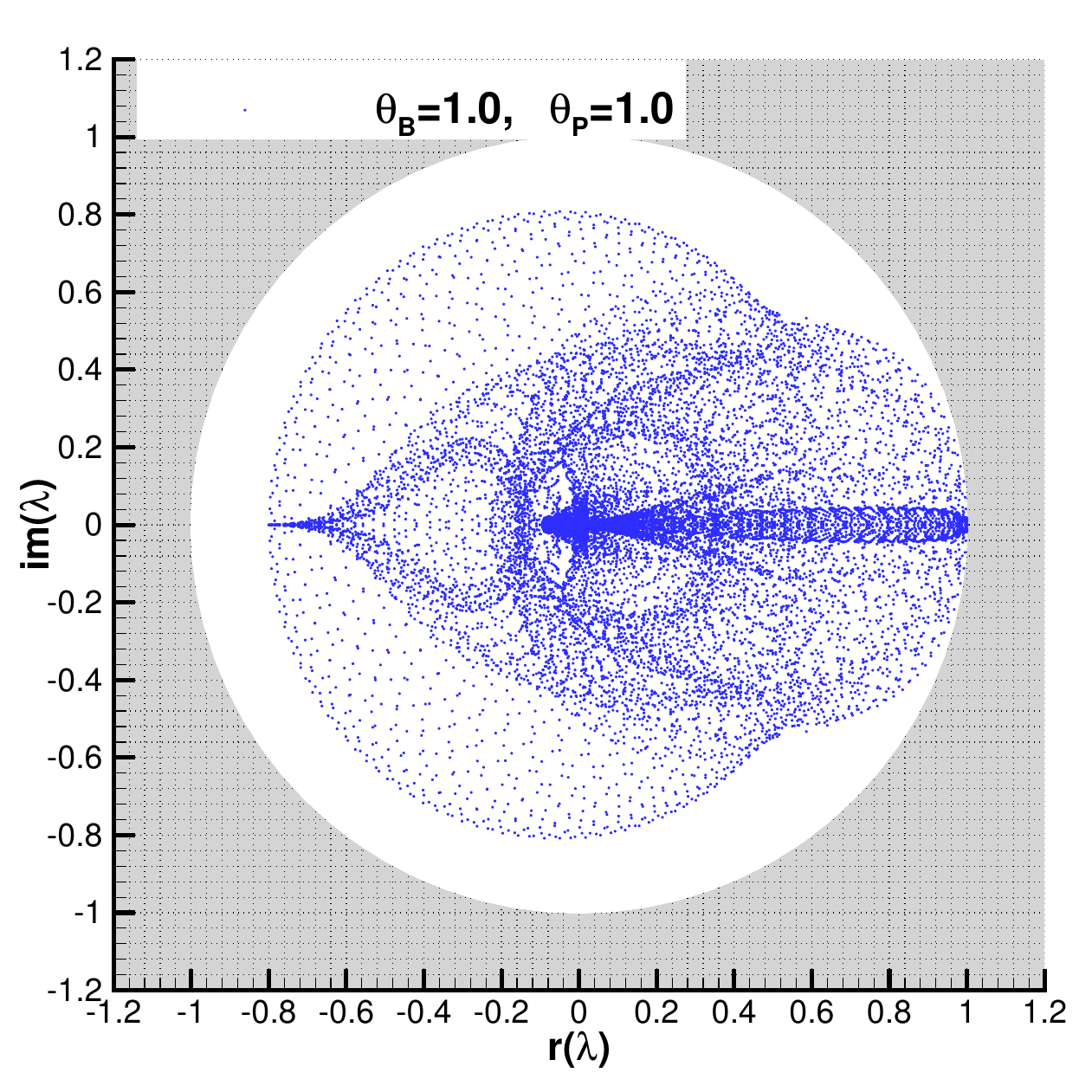} & 
\includegraphics[width=0.45\textwidth]{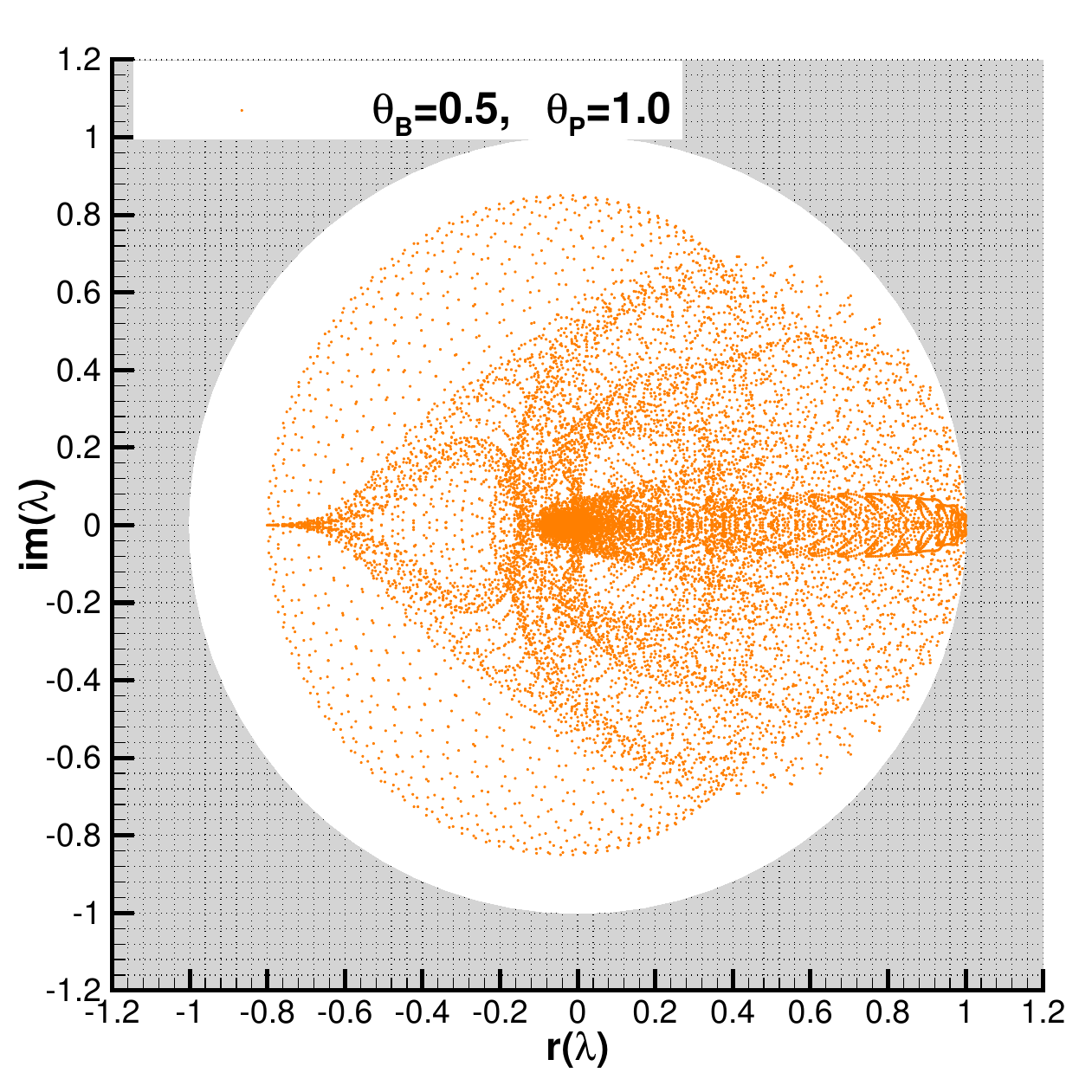} \\ 
\includegraphics[width=0.45\textwidth]{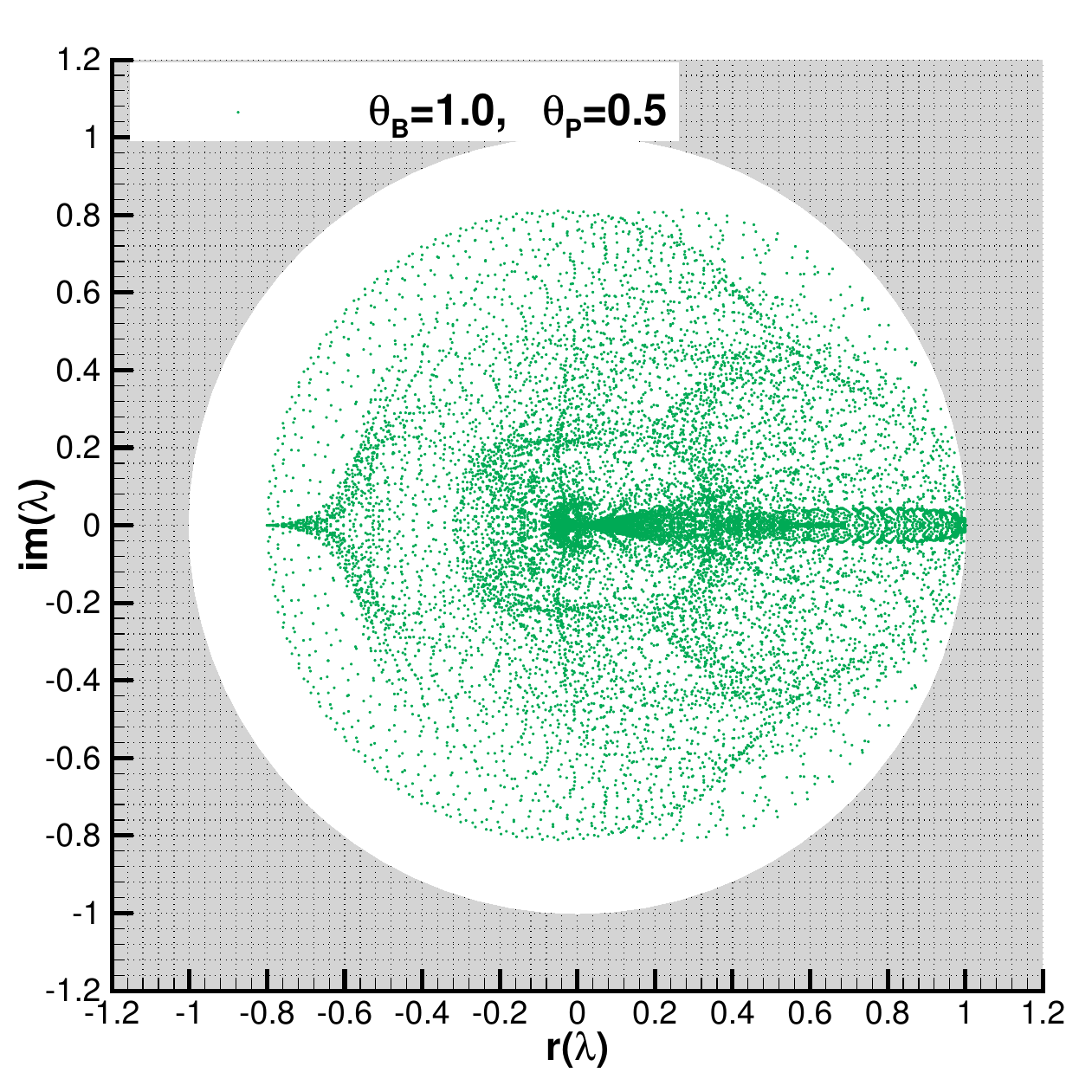} & 
\includegraphics[width=0.45\textwidth]{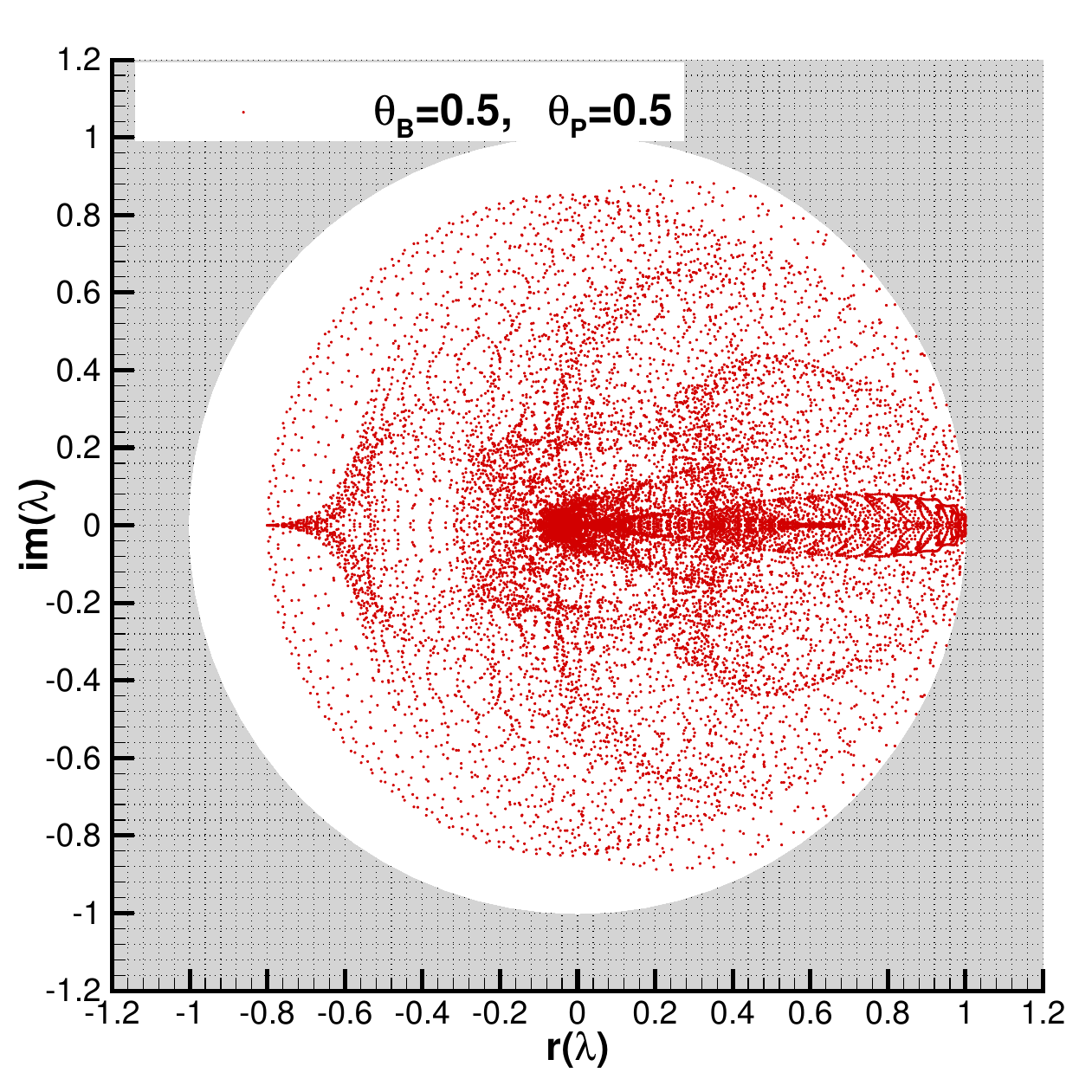}        
\end{tabular} 
\caption{Eigenvalues of the Jacobian matrix of the fully discrete scheme, evaluated around the constant equilibrium solution, varying the $\theta$-parameters of the time-discretization between $1$ and $1/2$. } 
\label{fig:MSA}
\end{center}
\end{figure}

From equation (\ref{eq:MSA}) one can argue that the perturbations do not grow whenever all the eigenvalues $\lambda_i$ of the Jacobian of the fully-discrete system  $\delta \mathcal{L}/\delta \Q_h^n$ live within the unit circle in the complex plane, i.e. $|| \lambda_i || \leq 1 $. Fig.  \ref{fig:MSA} show the results of the linear stability analysis after choosing a constant solution in two-space dimensions, in primitive variables,
\myequation{l}{
\V^{\text{eq.}} = ( \rho_0,  \v , p_0, \mathbf{b} ),\;\; \text{with}\;\; 
\v = \left(\begin{array}{c} \cos \alpha_1 \cos \alpha_2 \\ 
														\sin \alpha_1 \cos \alpha_2,\\ 
														\sin \alpha_2 		
						\end{array}								\right), \;
 \mathbf{b}  = b_0  \left(\begin{array}{c}   \cos \alpha_3 \cos \alpha_2 \\
																						\sin \alpha_3 \cos \alpha_2	 \\
																							\sin \alpha_2 
						\end{array}								\right),}
{}
$\alpha_1 = \pi/6$, 
$\alpha_2 = \pi/4$, 
$\alpha_3 = \pi/3$, 
$b_0   = \sqrt{4\pi}$,
$\rho_0 = 1             $, and
$p_0   = 1/\gamma     $. The numerical domain is chosen to be a two-dimensional periodic box $(x,y)\in [-0.5,0.5]^2$ with a mesh resolution of $\Delta x=\Delta y = 0.02$. The Jacobian matrix has been computed numerically, while the eigenspectrum has been computed using the {\tt dggev} function of the Intel Math Kernel Library.

The computed full spectrum, see Fig.  \ref{fig:MSA}, is shown to be contained within the unit circle, and for this reason the fully discrete system can be said to be \emph{linearly stable}, in the sense of \emph{Lyapunov}, at the selected equilibrium solution. Moreover, as expected, decreasing the value of the $\theta$-parameters, the eigenspectrum gets closer to the border of the stability region.

\subsection{One dimensional tests}

\paragraph{Current sheet.}

The exact solutions available for VRMHD that we can use for numerical validation are only very few. One of the more elementary exact solutions, is an analytical solution of the diffusion equation
$$ \partial_t \B + \nabla \times (\eta \nabla \times \B) = 0 $$
with initial data
\begin{align}
\left( \rho, \v, p, \B \right) = \left\{ \begin{array}{@{(}c@{,}c@{,}c@{,}c@{)\hspace{0.5cm}}r}  
  \rho_0 & (0,0,0)  & p_0 & (0,+B_y^0,B_z^0)    &  x\leq 0\\
  \rho_0 & (0,0,0)  & p_0 & (0,-B_y^0,B_z^0)       &  x> 0
\end{array}
\right.
\end{align}
and it takes the form (see \cite{boundaryLayer,Komissarov2007})
\begin{equation}
\label{eq:stokes} 
  B_y(x,t) = - B_y^0 \textnormal{erf}\left( \halb \frac{x}{\sqrt{\eta t}} \right).
\end{equation}
	The solution (\ref{eq:stokes}) is still a perturbative solution for the full set of the nonlinear VRMHD equations. In particular, this is valid in the low Mach regime and for small values of $B_y^0$.
 
Our new semi-implicit finite volume scheme is particularly well-suited for all Mach number flows, but also for highly magnetized flows. In this test the density
and the fluid pressure are set to $\rho_0=1$ and $p_0=10^5$, respectively, with a magnetic field of $B_z^0=10^4$, and in the $y$-direction a step function with amplitude $B_y^0=10^{-3}$. 
The resulting plasma beta is then $\beta = 8 \pi p/\B^2 \sim  2.5 \times 10^{-4}$.
The magnetic resistivity is $\eta=10^{-1}$, the other fluid parameters are $\mu=0$, $Pr=1$ and $c_v=1$. In order to test the capabilities of the presented implicit solver, the resolution of the mesh is set to $\Delta x= 0,02$ with a one-dimensional domain $\x \in [-50,50]$. 
The simulation is carried out from $t=0$ up to $t=1000$ with a time-step $\Delta t=10$.
With this spatial resolution, this simple test-problem would be very difficult to be solved by an explicit solver because of the large values of the magnetosonic and Alfvén wave speeds, and the corresponding $\CFL$ restriction (\ref{eq:CFL}).  
More specifically, the time-step chosen according to the CFL condition of an explicit scheme based on the $\lambda^{v\B}$  or $\lambda$ eigenvalues, see (\ref{eqn.eval.c}) and (\ref{eq:eval.full}), would be of the order of  $$\Delta t(\lambda^{v\B}) \sim \Delta t(\lambda) \sim  10^{-6},$$ while the chosen time-step in our simulation was $\Delta t=10$ and thus seven orders of magnitude larger.

\begin{figure}[!htbp]
\begin{center}
\begin{tabular}{cc}  
\includegraphics[width=0.55\textwidth]{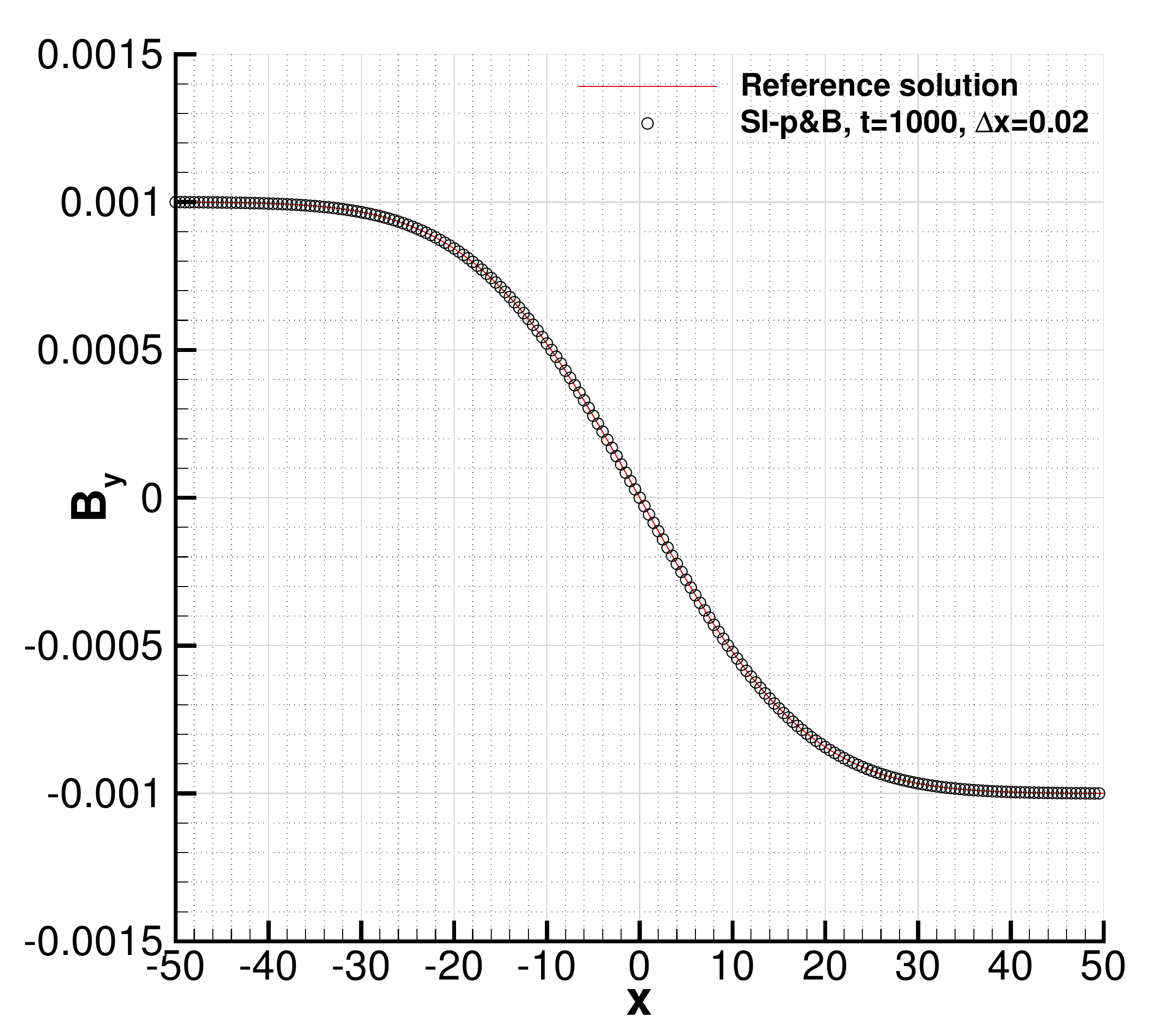}    
\end{tabular} 
\caption{Exact and numerical solution for the current sheet at time $t=1000$ solving the VRMHD equations with $\eta=10^{-1}$. The computational time-step has been set constant to $\Delta t=10$. The chosen background pressure is $p=10^6$ (low Mach regime), magnetic field in the $z$-direction $B_z=10^4$, and mesh resolution $\Delta x=0.02$  make this simple test almost unresolvable without an implicit treatment of pressure and magnetic field. } 
\label{fig.slcs}
\end{center}
\end{figure}

In Fig. \ref{fig.slcs} we report a scatter plot of the discrete solution obtained  with our new divergence-free  three-split semi-implicit FV scheme. An excellent agreement between the numerical solution and the reference solution can be observed.  
We emphasize that the full nonlinear VRMHD equations have been solved for this test problem, although the reference solution is an analytical solution only for the pure diffusion equation. 

\paragraph{1D Riemann problems.} 

Even if for fusion-oriented applications, shock waves and fluid discontinuities are not the first matter of interest, there are also many applications for which it is important to give an accurate approximation also of such nonlinear phenomena, e.g. violent explosion  or nonlinear phase-transition of plasma in astrophysics or solar physics. Moreover, since the MHD system is solved in its nonlinear form, it is important to show that the semi-implicit scheme presented here is able to face any eventual flow discontinuity that may arise during the simulation.
In this section, a set of standard nontrivial Riemann problems has been solved to test the robustness and the shock-capturing capabilities of our new scheme even at high CFL numbers.  The respective initial conditions are provided in Tab. \ref{tab.ic.mhd}. These configurations have been partially introduced and analyzed by \cite{BrioWu,RyuJones,DaiWoodward,OsherUniversal}. For a detailed analysis of the spectral properties of the MHD equations, see \cite{roebalsara}.

\begin{figure}[!t]
\begin{center}
\begin{tabular}{cc} 
\includegraphics[width=0.45\textwidth]{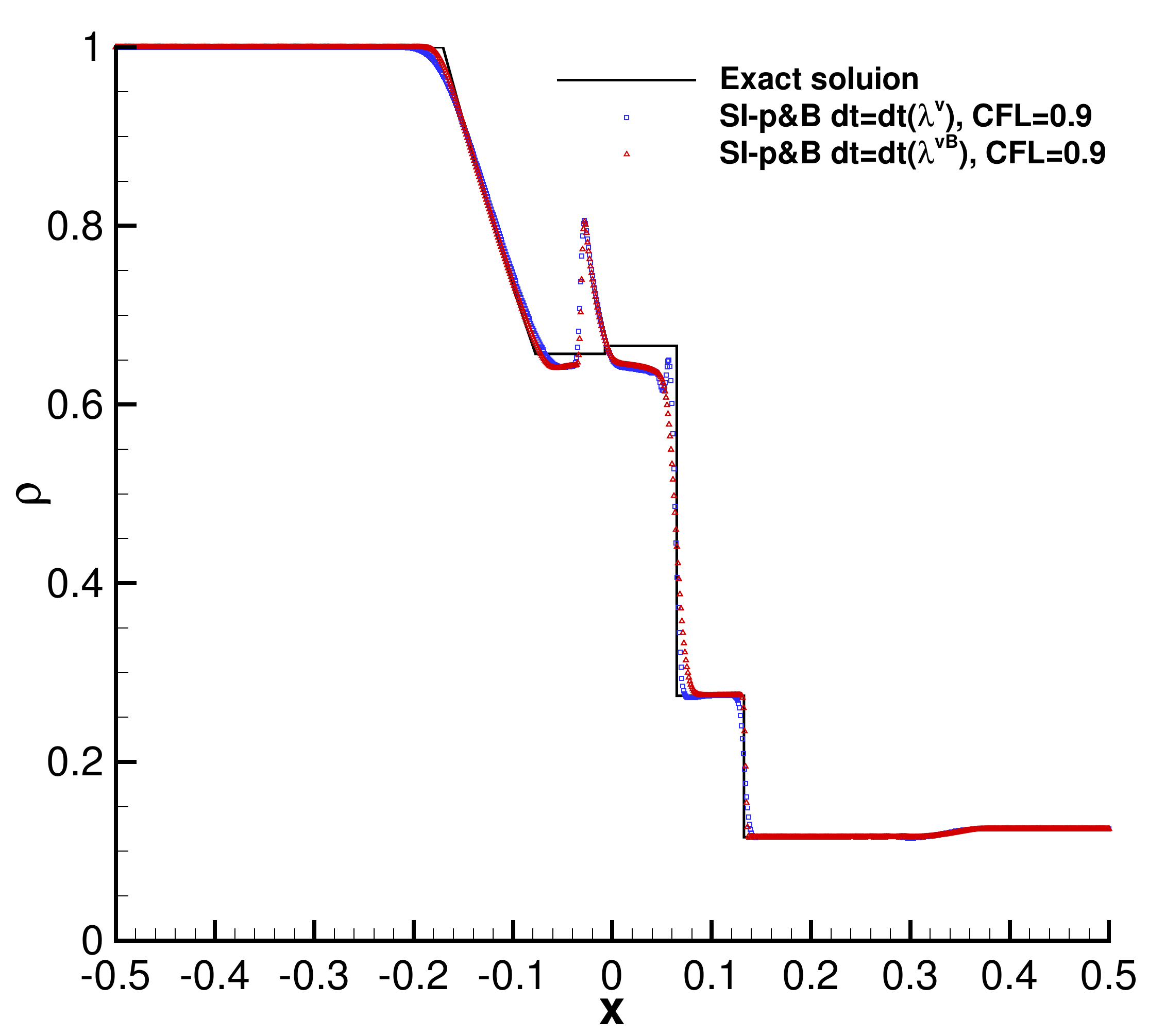}  & 
\includegraphics[width=0.45\textwidth]{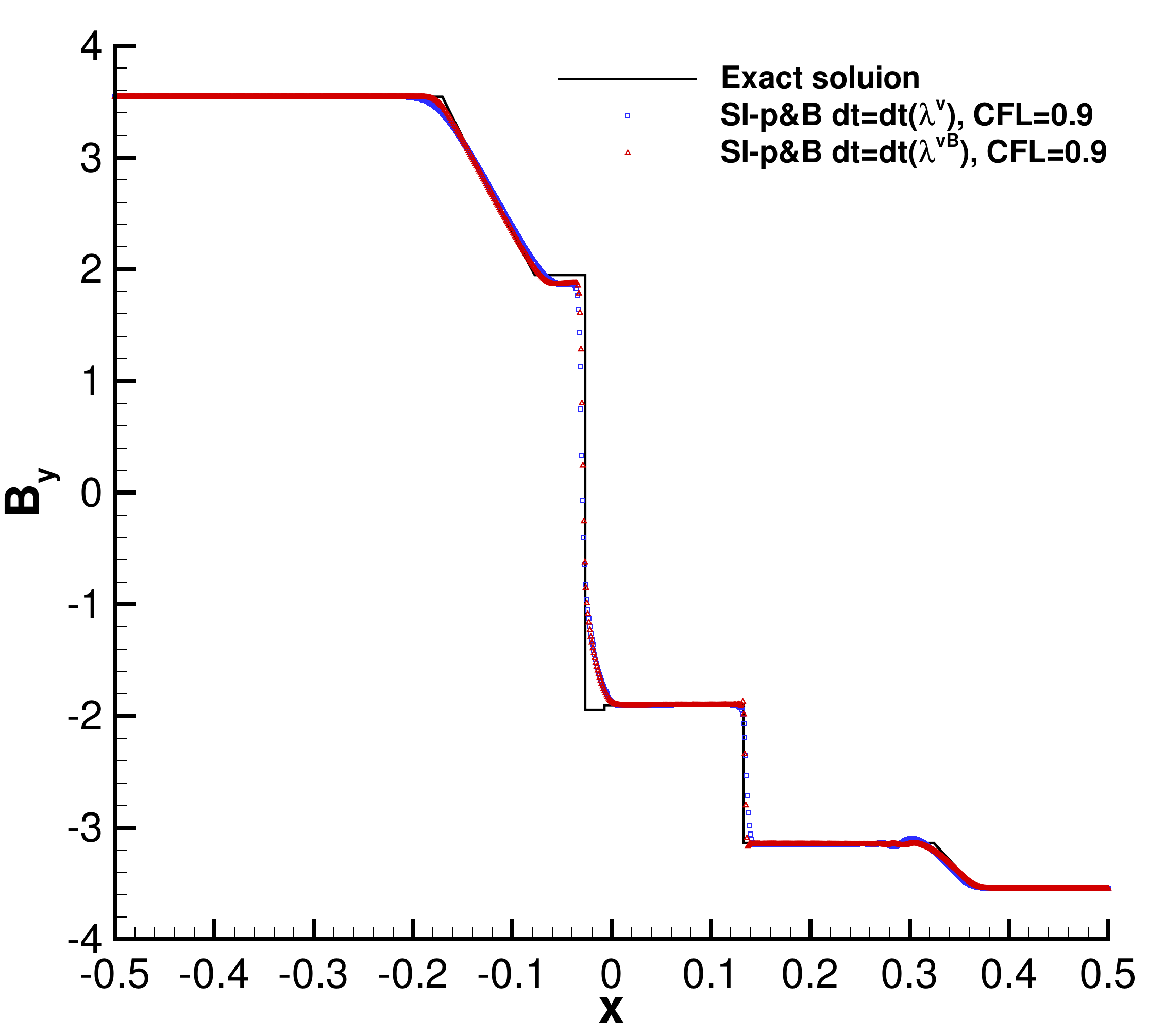}   \\ 
\includegraphics[width=0.45\textwidth]{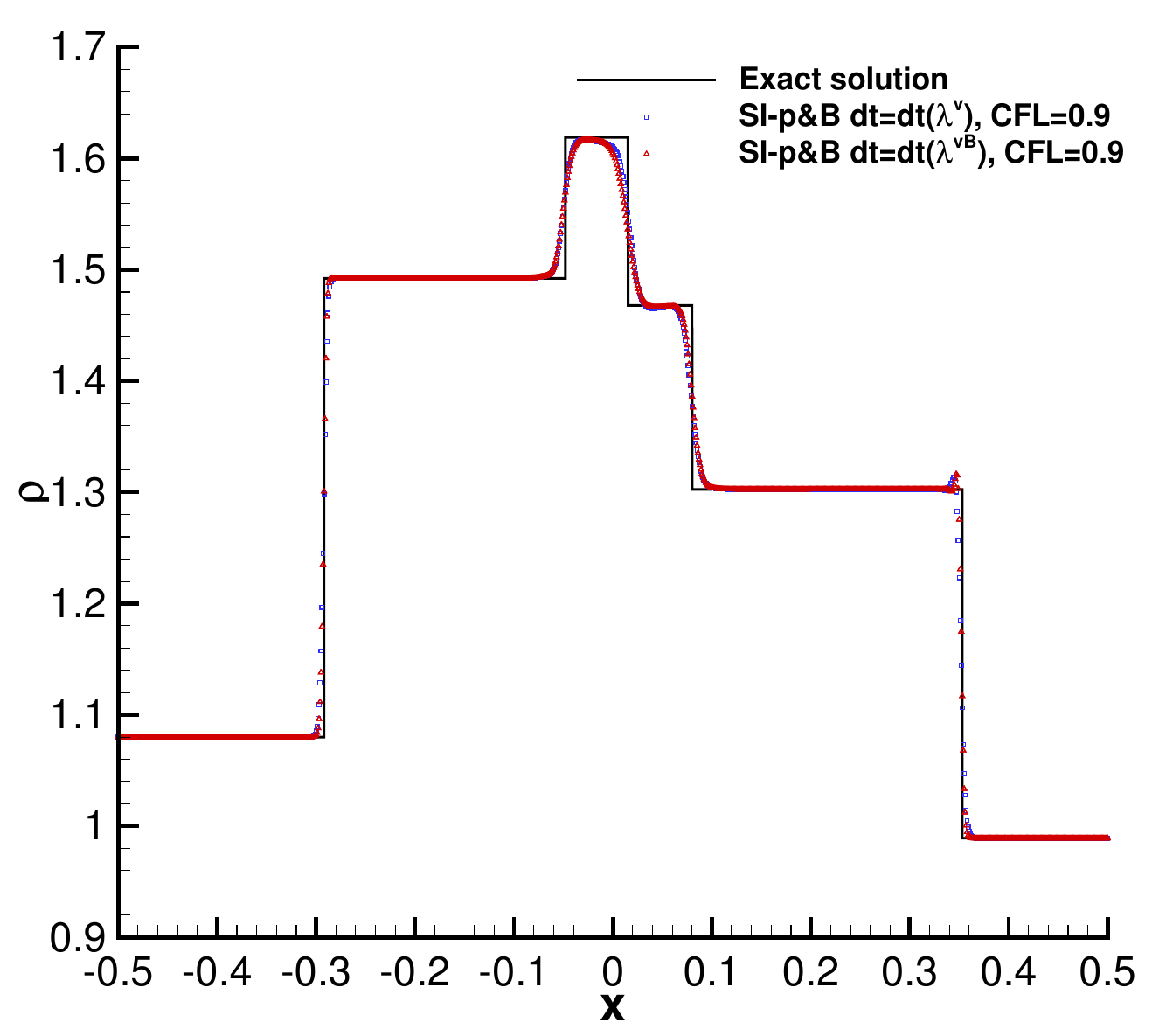}  & 
\includegraphics[width=0.45\textwidth]{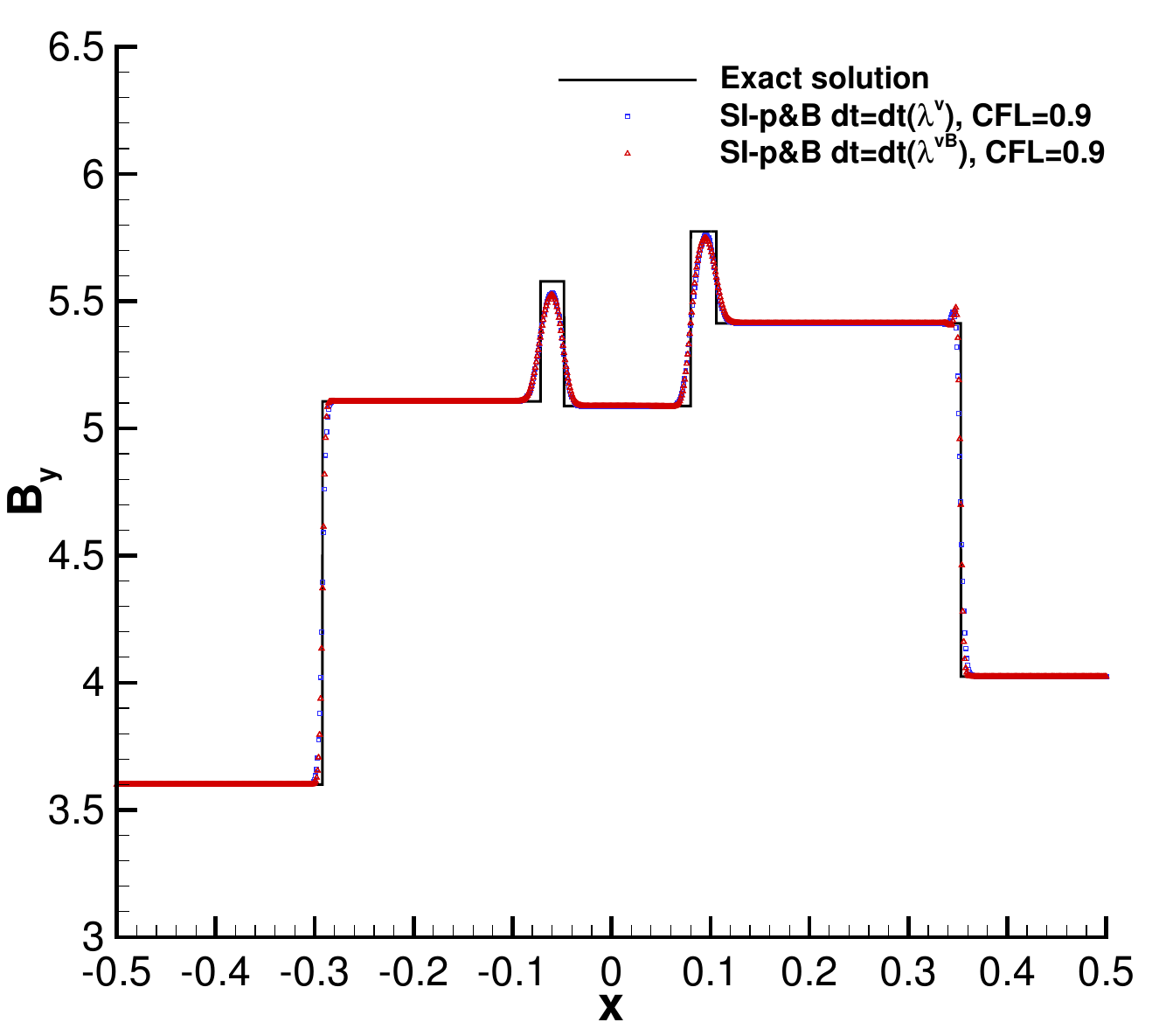}   
\end{tabular} 
\caption{Exact and numerical solution for the Riemann problem RP1 at $t=0.1$ (first row) and RP2 at $t=0.2$ (second row)
solving the ideal MHD equations with the here presented semi-implicit scheme. Density (left column) and magnetic field component $B_y$ (right column) are plotted, comparing the numerical solutions obtained after choosing a CFL time condition depending on the eigenvalues of the pure (convective) explicit terms, i.e. $\Delta t(\lambda^{v})$ (blue squares), and the eigenvalues $\Delta t(\lambda^{v\B})$ (red triangles). 
} 
\label{fig:rp1}
\end{center}
\end{figure}

\begin{figure}[!t]
\begin{center}
\begin{tabular}{cc} 
\includegraphics[width=0.45\textwidth]{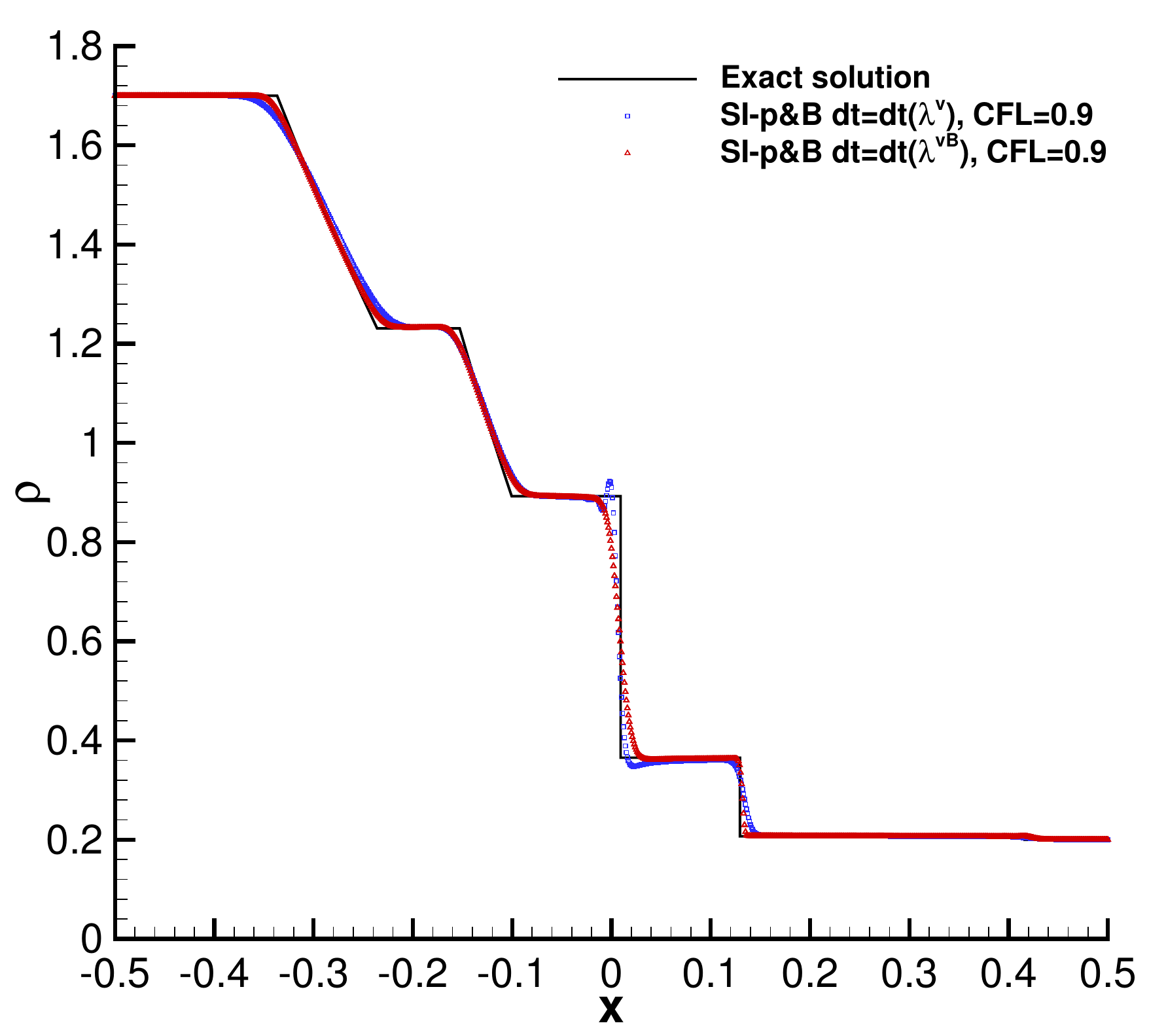}  & 
\includegraphics[width=0.45\textwidth]{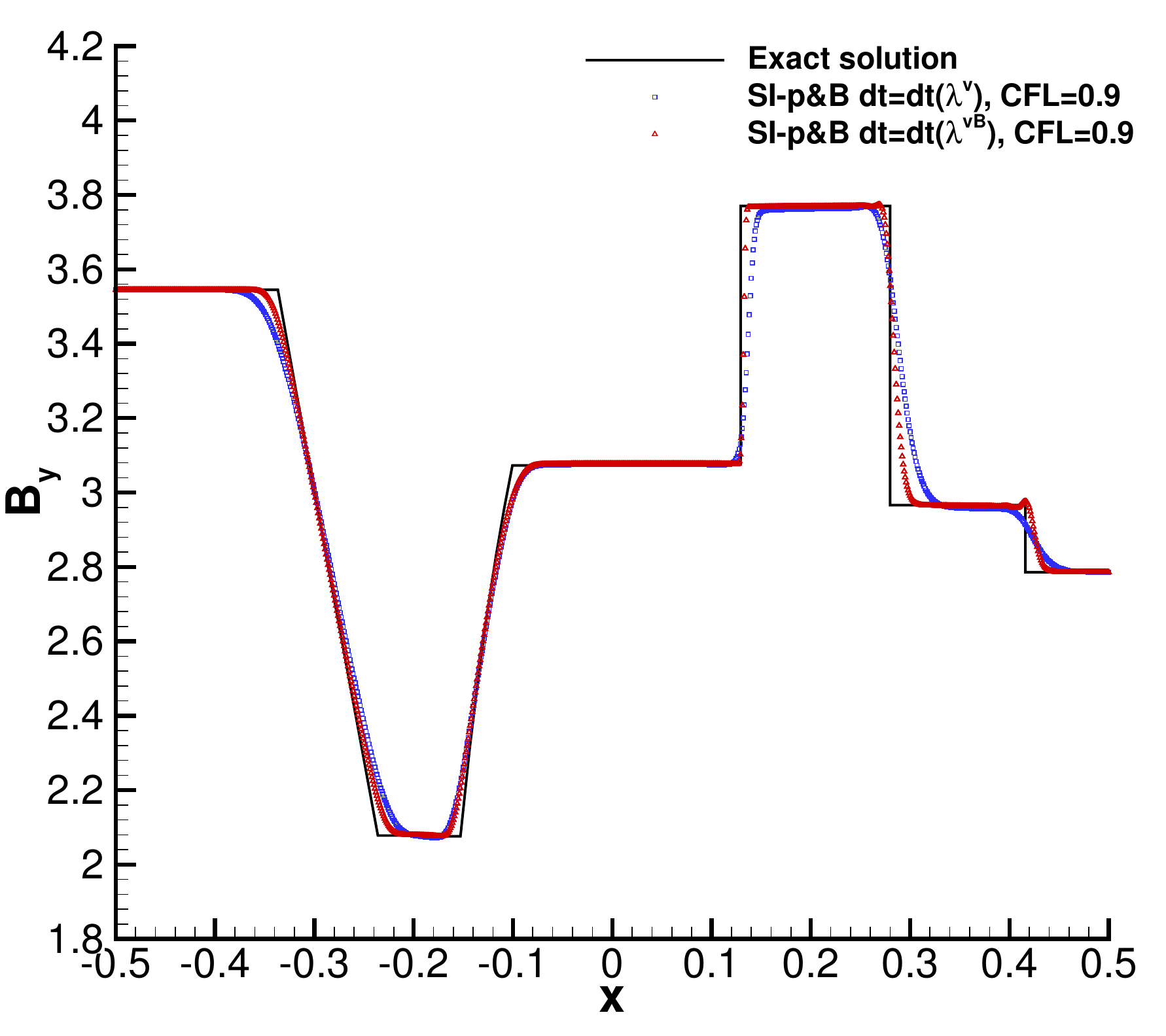}  \\ 
\includegraphics[width=0.45\textwidth]{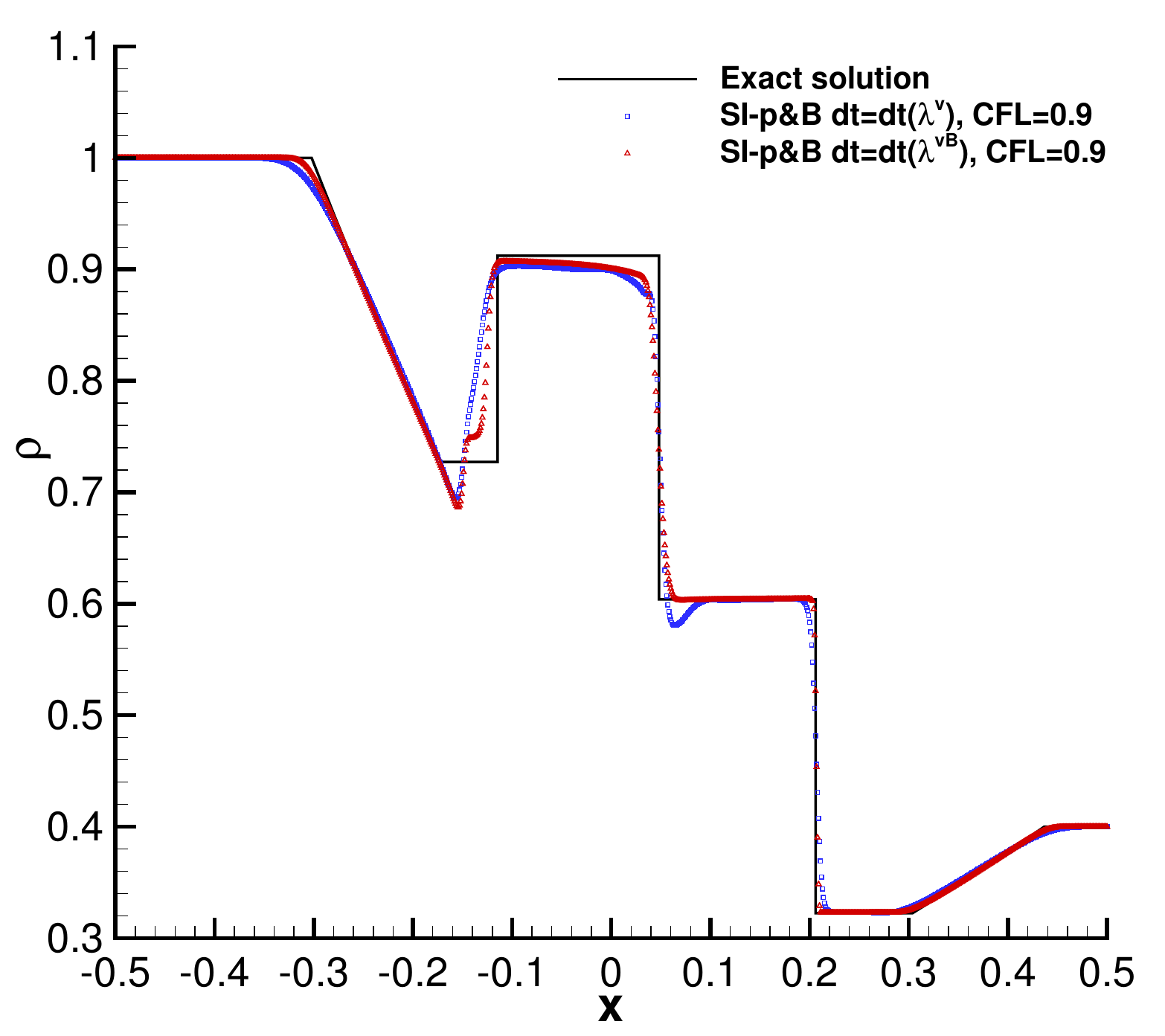}  &
\includegraphics[width=0.45\textwidth]{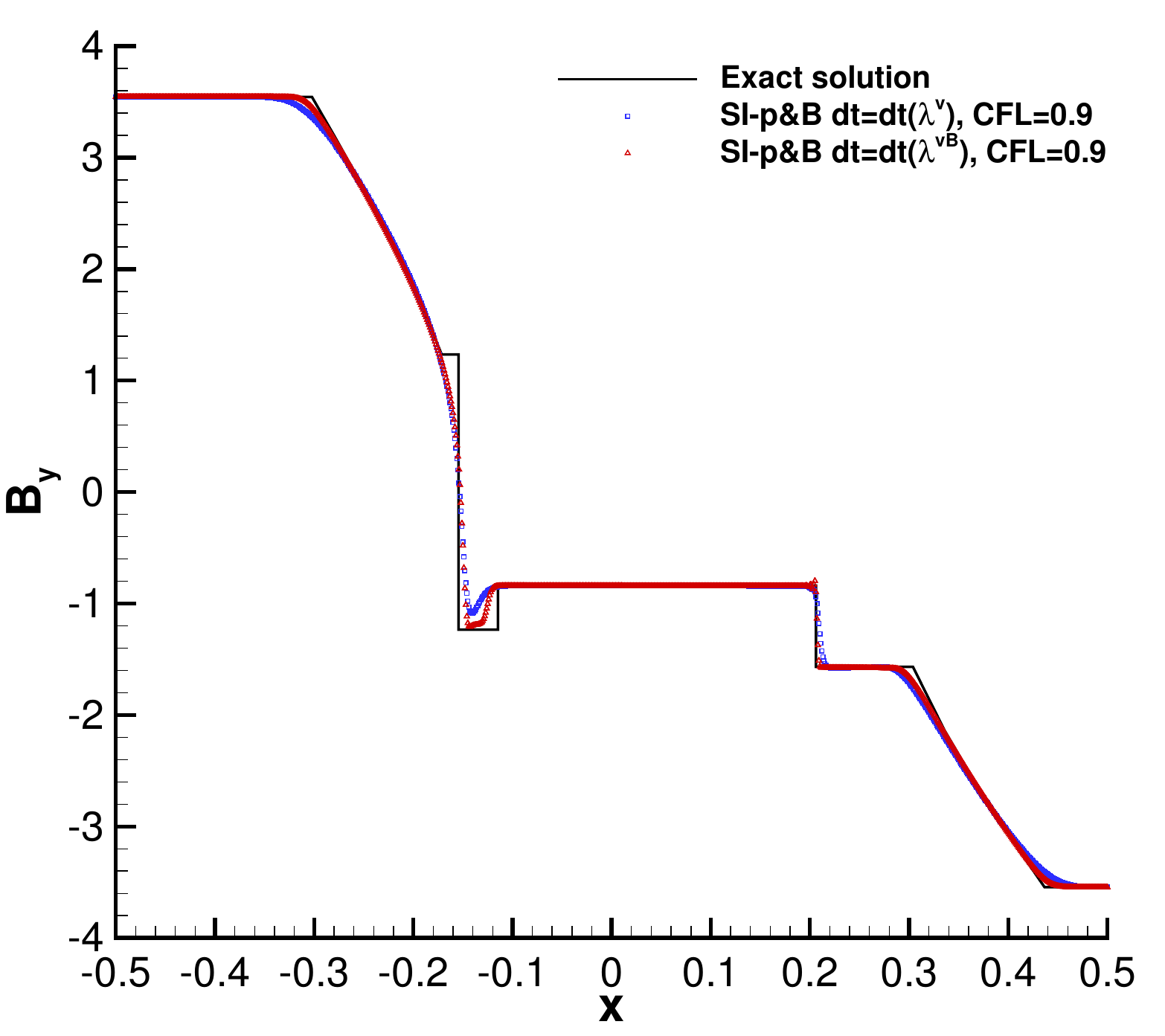}      
\end{tabular} 
\caption{Exact and numerical solution for the Riemann problem 
RP3 at $t=0.15$ (first row) and RP4 at $t=0.16$ (second row)
solving the ideal MHD equations with the here presented semi-implicit scheme. Density (left column) and magnetic field component $B_y$ (right column) are plotted, comparing the numerical solutions obtained after choosing a CFL time condition depending on the eigenvalues of the pure (convective) explicit terms, i.e. $\Delta t(\lambda^{v})$ (blue squares), and the eigenvalues $\Delta t(\lambda^{v\B})$ (red triangles). 
} 
\label{fig:rp2}
\end{center}
\end{figure}

\begin{figure}[!t]
\begin{center}
\begin{tabular}{cc}  
\includegraphics[width=0.45\textwidth]{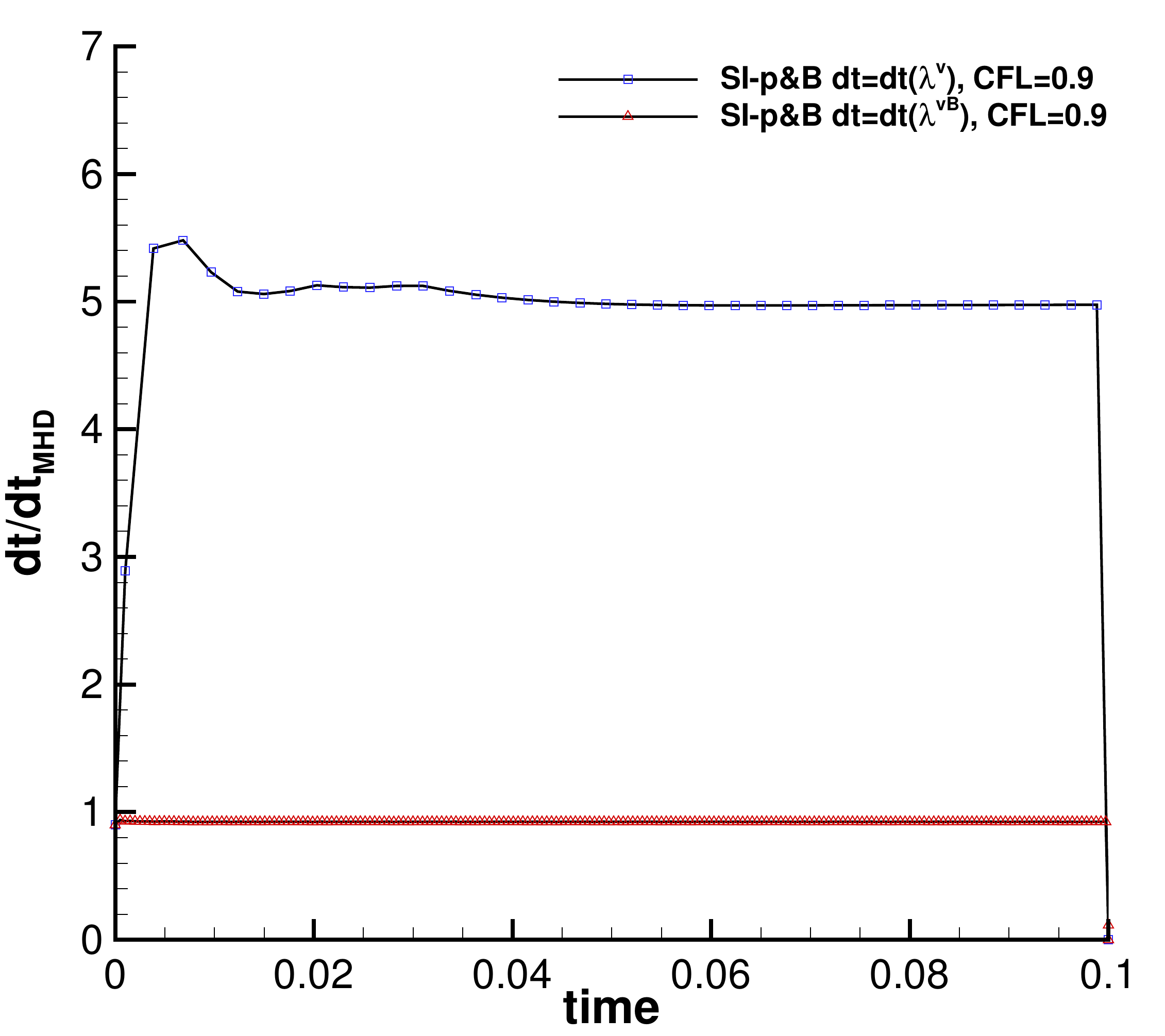} & 
\includegraphics[width=0.45\textwidth]{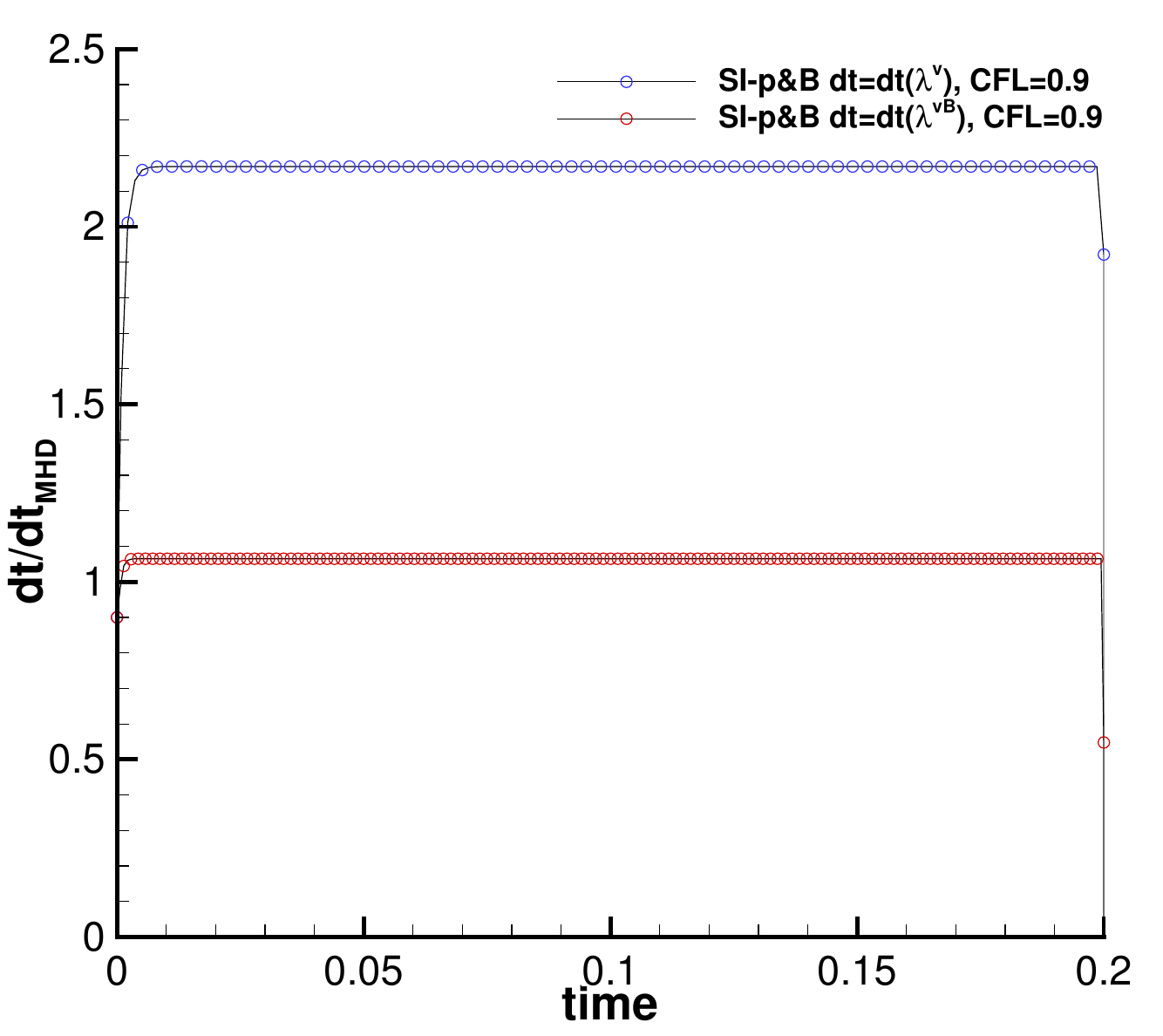} \\ 
\includegraphics[width=0.45\textwidth]{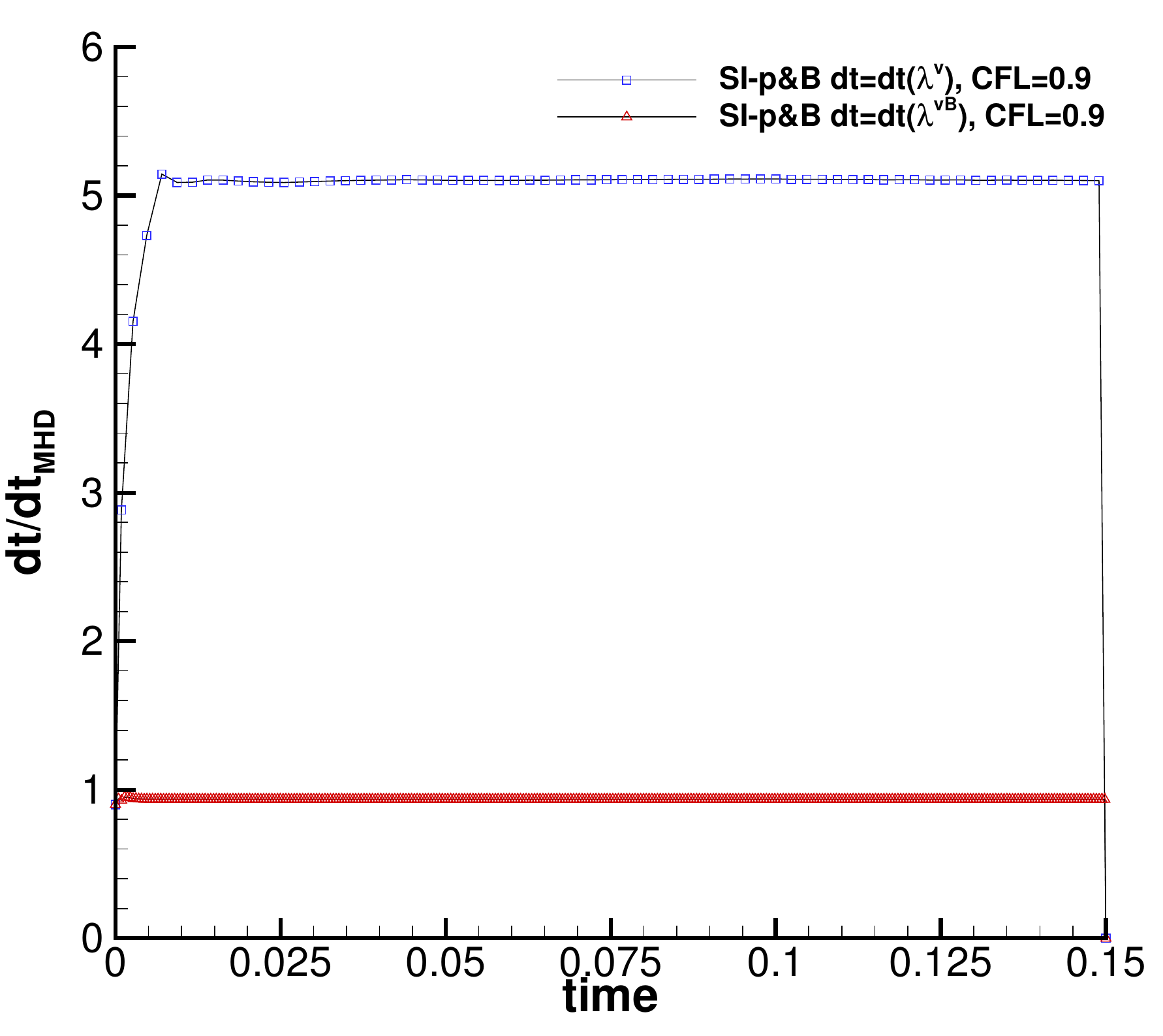}   &
\includegraphics[width=0.45\textwidth]{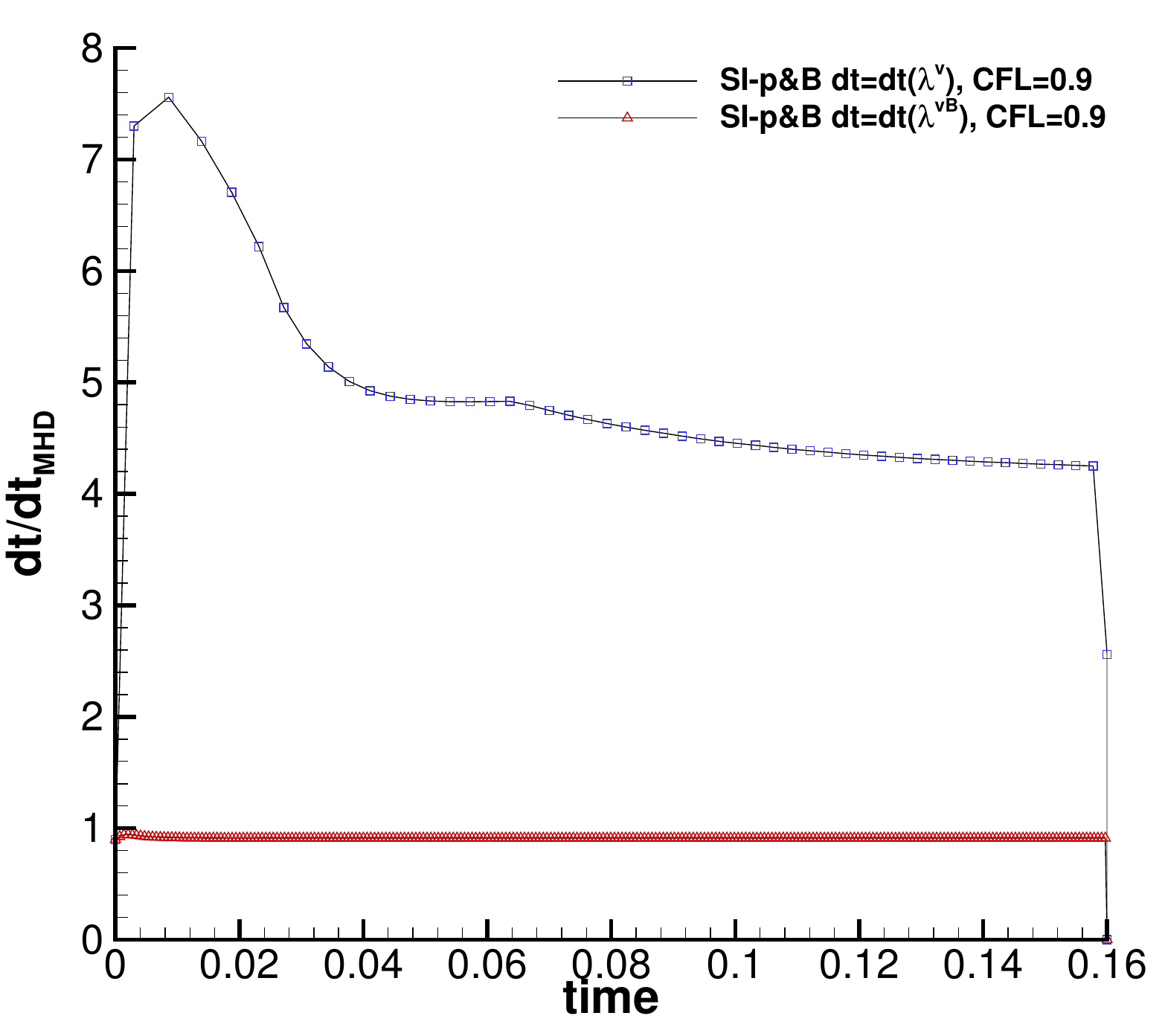}        
\end{tabular} 
\caption{The time evolution of the effective Courant number evaluated with respect to the eigenvalues $\lambda$ of the full MHD system, i.e.  $\Delta t / \Delta t (\lambda)$, for the numerical solution of the Riemann problems RP1 and RP2 (first row), RP3 and RP4 (second row), from left to right, respectively. } 
\label{fig:CFL}
\end{center}
\end{figure}

The computational domain is $x\in[-0.5,0.5]$, with mesh resolution $\Delta x = 10^{-3}$. In all our simulations we set the constant parameter $\CFL=0.9$ and $\theta_{\text{p}}=1$, while one had $\theta_{\text{B}}=0.55$ for RP1, RP2 and RP3 and $\theta_{\text{B}}=0.65$ for RP4.
The flow parameters are $\gamma =5/3$, $\mu=\eta=0$. The discontinuity is at $x_d=0$ for RP0, RP1 and RP4, while it is at $x_d=-0.1$ for RP2 and RP3. 
In Fig. \ref{fig:rp1}-\ref{fig:rp2} we report the scatter-plot of the results obtained by two different configurations of the same numerical scheme: in one case, the computational time-step is computed accordingly to the CFL restriction (\ref{eq:CFL}) based on the maximum eigenvalues of the pure convective subsystem, i.e. $dt=dt(\lambda^v)$, in the other based on the maximum eigenvalues of the convective and Alfvénic subsystems, i.e. $dt=dt(\lambda^{vB})$. To show the gain in terms of CFL numbers, in Fig.  \ref{fig:CFL} we plot the time evolution of the Courant number evaluated with respect to the full MHD eigen-system, in terms of the ratio $dt/dt(\lambda)$.

A reference solution for the proposed tests have been kindly provided  by S.A.E.G. Falle \cite{fallemhd,falle2}. 
The computed solution is shown to be in good agreement with others published works, see \cite{SIMHD,BrioWu,fallemhd,OsherUniversal,HPRmodelMHD}. Moreover, one can notice that when a convective and Alfvénic Courant number $dt=dt(\lambda^{vB})$ is used, the results fits very well those obtained in \cite{SIMHD}. On the other hand, the use of a purely convective Courant number $dt=dt(\lambda^v)$ may lead, for these tests, a gain of $2-5$ times in terms of Courant number, at the aid of numerical diffusion, due to the lower time-accuracy, and a slightly more oscillatory behavior due to the unresolved time-scales of the flows. Since a very basic Rusanov (or Local-Lax-Friedrich) numerical flux is used, the results are in line with what was expected. One could notice the artificial compound wave that is present in the plot of the density in Fig. \ref{fig:rp1}. This is a well known effect that is typically generated by any standard finite-volume solver, where a rotational-wave and a slow shock-wave are apparently approximated by an overcompressive intermediate-wave and a slow rarefaction-wave. The existence of such a solution is reported to be related to the non-convexity of the MHD equations, see references \cite{fallemhd,BrioWu}.

Then, the results confirm the desired robustness and shock-capturing capabilities even at the higher admissible Courant numbers. It is clear that the numerical error in the time resolution can be consistently reduced by lowering the Courant number, as shown in Fig. \ref{fig:rp1}-\ref{fig:rp2}. Moreover, in the definition of the computational time-step, one can choose the desired set of eigenvalues according to the time-scale of the physics of interest, after adding a proper amount of numerical stabilization terms both in the explicit and in the implicit solvers.

\begin{table}[!t]
 \caption{Initial states left (L) and right (R) for the primitive variables for the selected  Riemann problems of the ideal classical MHD equations. 
 In all cases $\gamma=5/3$. The discontinuity is at $x_d=0$ for RP0, RP1 and RP4, while it is 
at $x_d=-0.1$ for RP2 and RP3. } 
\begin{center} 
 \begin{tabular}{rcccccccc}
 \hline
 Case & $\rho$ & $u$ & $v$ & $w$ & $p$ & $B_x$ & $B_y$ & $B_z$        \\ 
 \hline   
 RP1 L: &  1.0    &  0.0     & 0.0    & 0.0      &  1.0     & $\frac{3}{4} \sqrt{4 \pi}$ &  $\sqrt{4 \pi}$  & 0.0       \\
     R: &  0.125  &  0.0     & 0.0    & 0.0      &  0.1     & $\frac{3}{4} \sqrt{4 \pi}$ & $-\sqrt{4 \pi}$  & 0.0       \\
 RP2 L: &  1.08   &  1.2     & 0.01   & 0.5      &  0.95    & 2.0 &  3.6     & 2.0            \\
     R: &  0.9891 &  -0.0131 & 0.0269 & 0.010037 &  0.97159 & 2.0 &  4.0244  & 2.0026         \\
 RP3 L: &  1.7    &  0.0     & 0.0    & 0.0      &  1.7     & 3.899398 &  3.544908  & 0.0              \\
     R: &  0.2    &  0.0     & 0.0    & -1.496891  &  0.2   & 3.899398 &  2.785898  & 2.192064         \\
 RP4 L: &  1.0    &  0.0     & 0.0    & 0.0      &  1.0     & $1.3 \sqrt{4 \pi}$ &  $\sqrt{4 \pi}$   & 0.0            \\
     R: &  0.4    &  0.0     & 0.0    & 0.0      &  0.4     & $1.3 \sqrt{4 \pi}$ &  $-\sqrt{4 \pi}$  & 0.0             \\
 \hline
 \end{tabular}
\end{center} 
 \label{tab.ic.mhd}
\end{table}

\subsection{Convergence rate}
In this section one smooth MHD flows is selected to test experimentally the convergence rate of the proposed scheme.
For this section, resistivity $\eta$ and artificial stabilizing terms $\mathbf{D}$ (see equation \ref{eq:dBx}) have been set to zero. 
In addition, the \texttt{minmod} flux limiter has been switched off in order to give a better estimate of the order of accuracy of the novel three-split semi-implicit scheme proposed in this paper.

\paragraph{Alfvén wave}

\begin{figure}
\centering  
  \includegraphics[width=0.7\textwidth]{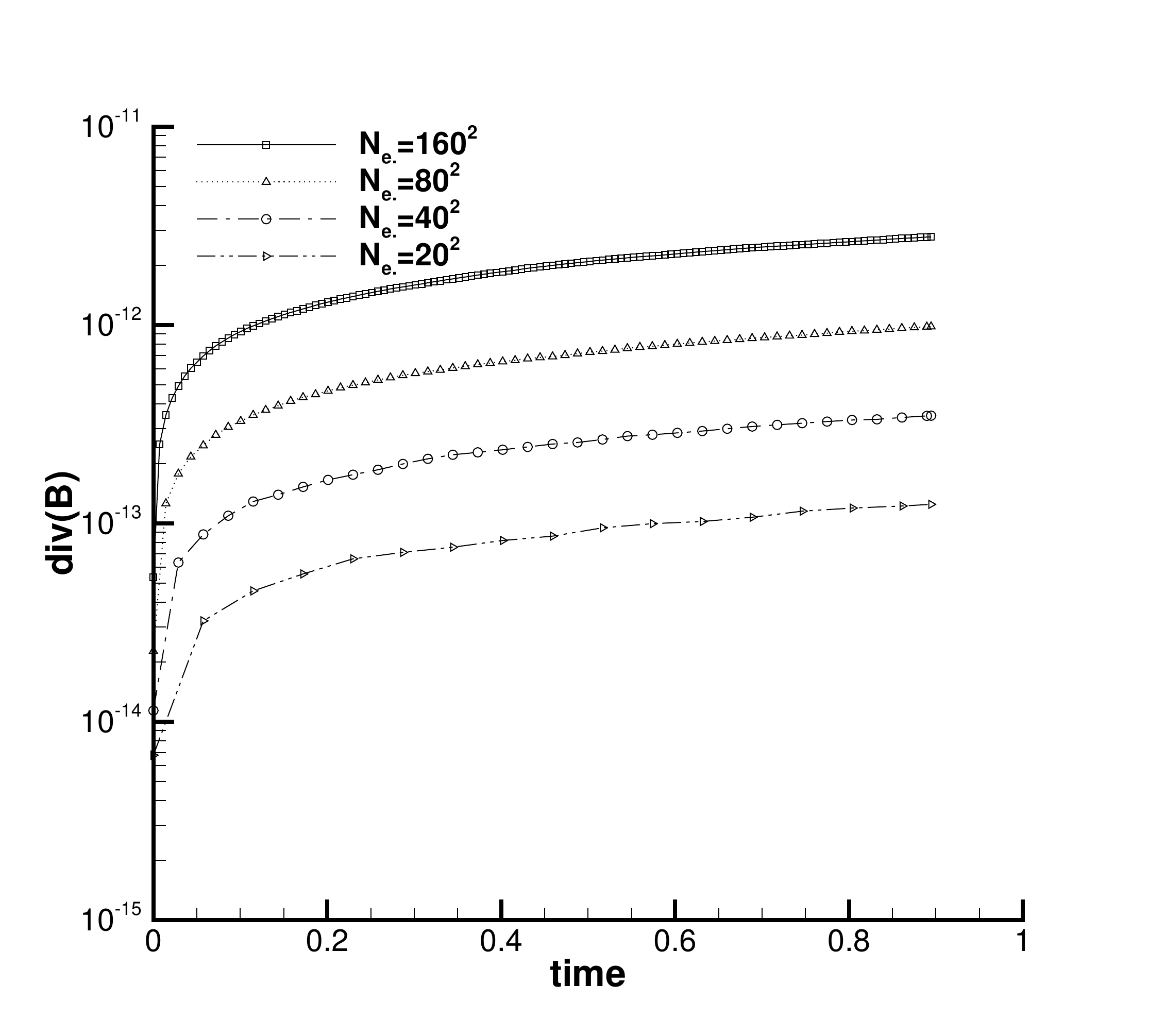} 
\caption{Time evolution of the $L_1$ norm of the divergence of the magnetic field  $\nabla\cdot \B$ for the Alfv\'en-wave test evaluated on the cell-nodes for different grid refinements.}
\label{fig:AW_divB}       
\end{figure}

The so-called Alfvén wave test is proposed for the convergence study. 
Being a smooth and time-dependent solution of the MHD equations, the Alfvén wave test is widely used to give an experimental estimate of the order of accuracy of any MHD solver. A low-Mach regime has been chosen after selecting a constant background pressure of $p=10^2$. Although the solution is three-dimensional, the problem is reduced to two dimensions by assuming the sinusoidal wave propagating on the $x-y$ plane with  periodic boundary conditions in the third space-dimension. 
In particular, the initial conditions are chosen to be 
\myequation{l}
{
\rho = 1, \\
\v = \alpha \left( - n_y \cos(\varphi), n_x \cos(\varphi), \sin(\varphi) \right),  \\
p = 10^2\\
\B =  \sqrt{2\pi}\left( n_x + n_y \alpha  \cos(\varphi), n_y - n_x \alpha  \cos(\varphi), -\alpha \sin(\varphi) \right),
}{\label{eq:Alfven}}
where
$$  \varphi = 2 \pi \frac{1}{n_y} \left[ n_x\left( x - n_x t\right) +  n_y \left( y - n_y t\right)\right].  $$
The direction of propagation is designed in order to be non-aligned with the grid, i.e.
$$  \mathbf{n} = ( n_x, n_y, n_z ) = \left( 1 , 2 , 0 \right)/\sqrt{5}. $$ 
A series of simulations have been performed after refining recursively the mesh-resolution. A second-order MUSCL scheme is used for the conservative fluxes (\ref{eq:PDEex}) and (\ref{eq:CVel}), while for the implicit part the theta-parameters are chosen as $\theta_p=\theta_B = 1/2$. Since a \texttt{minmod} flux limiter would have clipped the extrema of the sinusoidal wave, this test is simulated without the use of any flux-limiter for the evaluation of the conservative fluxes. The experimental results show that the expected \textbf{second-order} of accuracy is well reproduced, see Tab. \ref{tab:Alfven}, and the divergence-free condition is well preserved thanks to the proposed structure preserving formulation on staggered grids and despite the adopted nested recursion algorithm, see Fig. \ref{fig:AW_divB}. 
 \begin{table}[!t] 
 \centering
 \begin{tabular}{!{\extracolsep{-3pt}}ccccccccc!{}}
   \multicolumn{9}{c}{\textbf{Low Mach Alfv\'en-wave test}} \\
   \hline
& $N_{\text{element}}$ &  $L_1$ error &  $L_2$ error & $L_{\infty}$ error & $L_1$ or. & $L_2$ or. & $L_\infty$ or. &    Th. \\
   \cline{2-8}
   \multirow{4}{*}{$\rho$} 
& $ 20^2$ &   0.11476E-02	& 0.65322E-03 	& 0.62836E-03   	&   --- 		&   ---  		&  ---   & \multirow{4}{*}{\textbf{2}} \\ 
& $ 40^2$ &   0.19712E-03	& 0.11066E-03 	& 0.82484E-04   	&   2.54  	&   2.56  	&  2.93  &  									\\ 
& $ 80^2$ &   0.35965E-04	& 0.19978E-04 	& 0.14670E-04   	&   2.45  	&   2.47  	&  2.49  &  									\\ 
& $160^2$ &   0.77192E-05	& 0.42868E-05 	& 0.31000E-05   	&   2.22  	&   2.22  	&  2.24  &  									\\ 
  \cline{2-8}
  \multirow{4}{*}{$u$}  
& $ 20^2$ &   0.11875E\phantom{-}00	& 0.66122E-01		& 0.47556E-01 		& 	 --- 		&   ---  		&  ---   & 	\multirow{4}{*}{\textbf{2}}						   		\\ 
& $ 40^2$ &   0.28583E-01	& 0.16107E-01		& 0.12199E-01 		& 	2.05  	&   2.04  	&  1.96	 &                   \\
& $ 80^2$ &   0.86526E-02	& 0.48721E-02		& 0.36488E-02 		& 	1.72  	&   1.73  	&  1.74	 &                   \\
& $160^2$ &   0.27124E-02	& 0.15211E-02		& 0.11169E-02 		& 	1.67  	&   1.68  	&  1.71	 &                   \\
   \cline{2-8}                                                                      									
   \multirow{4}{*}{$p$}                                                             									 
& $ 20^2$ &   0.95619E-01	& 0.53641E-01	& 0.42323E-01		  &  --- 		&   ---  		&  ---  &   \multirow{4}{*}{\textbf{2}} 									\\ 
& $ 40^2$ &   0.26776E-01	& 0.14928E-01	& 0.11356E-01		  &  1.84 	&  1.85	&  1.90 		 & 						\\ 
& $ 80^2$ &   0.52476E-02	& 0.29212E-02	& 0.22013E-02		  &  2.35 	&  2.35	&  2.37  		 & 						\\ 
& $160^2$ &   0.11950E-02	& 0.66452E-03	& 0.49253E-03		  &  2.13 	&  2.14	&  2.16 		 & 						\\ 
    \cline{2-8}                                                                      									
    \multirow{4}{*}{$B_x$}                                                           									 
& $ 20^2$  & 0.40761E\phantom{-}00 	& 0.22747E\phantom{-}00 	& 0.16592E\phantom{-}00  	&   --- 		&   ---  		&  ---   &  \multirow{4}{*}{\textbf{2}} 									\\ 
& $ 40^2$  & 0.95749E-01 	& 0.54061E-01 	& 0.41267E-01  	&   2.09  	&   2.07  	&  2.01  &  									\\ 
& $ 80^2$  & 0.23237E-01 	& 0.13194E-01 	& 0.10252E-01  	&   2.04   	&   2.03  	&  2.01  &  									\\ 
& $160^2$  & 0.57346E-02 	& 0.32679E-02 	& 0.25407E-02  	&   2.02   	&   2.01  	&  2.01  &  									\\ 
   \hline
 \end{tabular}
 \caption{ \label{tab:Alfven} $L_1$, $L_2$ and $L_\infty$ errors and convergence rates for the 
   low-Mach Alfv\'en wave test obtained 
	with $\CFL=0.75$. Every norm is computed with respect to the grid where the corresponding variable is located, i.e. on the edges for the magnetic field, the faces for the velocity field, and the cell barycenters for pressure and matter-density.}
 \end{table}

\subsection{Two dimensional tests}

\paragraph{Stationary MHD vortex.} 

In this section we present a long-time simulation of the stationary isodensity vortex problem of \cite{Balsara2004}. A background state is given as 
\myequation{l}{
\left( \rho, \v, p, \B \right)  = \left( \rho, (0,0,0), 1, (0,0,0)\right)
}{}
while a fluctuation of velocity and magnetic field is added in the form

\myequation{l}{
\delta \v  =  \frac{v_0}{2\pi} e^{\frac{1}{2} (1-r^2)} \left( -y,x,0\right), \;\;
\delta \B  =  \frac{A_0}{2\pi} e^{\frac{1}{2} (1-r^2)} \left( -y,x,0\right) 
}{}
corresponding to a vector potential oriented in the $z$-direction $\delta A_z = \frac{A_0}{2\pi} e^{\frac{1}{2}(1-r^2)}$. An isodensity MHD-static solution can be obtained by setting the pressure
$$\delta p = \frac{1}{8\pi}\left(\frac{A_0}{2\pi}\right)^2 \left(1-r^2\right) e^{\left(1-r^2\right)} - \frac{1}{2} \left(\frac{v_0}{2\pi}\right)^2  e^{\left(1-r^2\right)},$$  see \cite{Balsara2004} for other details.

With this test one can further show some of the main benefit of using an implicit solver for the magnetic field relative to a standard explicit solver, even if both  make use of a divergence-free scheme.
The physical domain is a periodic box $(x,y)\in[-10,10]$ that is discretized with a mesh resolution of $\Delta x = \Delta y =0.1$. The parameters of the problem were $v_0=1$ and $A_0=\sqrt{4\pi}$.

\begin{figure} 
\centering 
\begin{tabular}{lr}
			\includegraphics[width=0.42\textwidth]{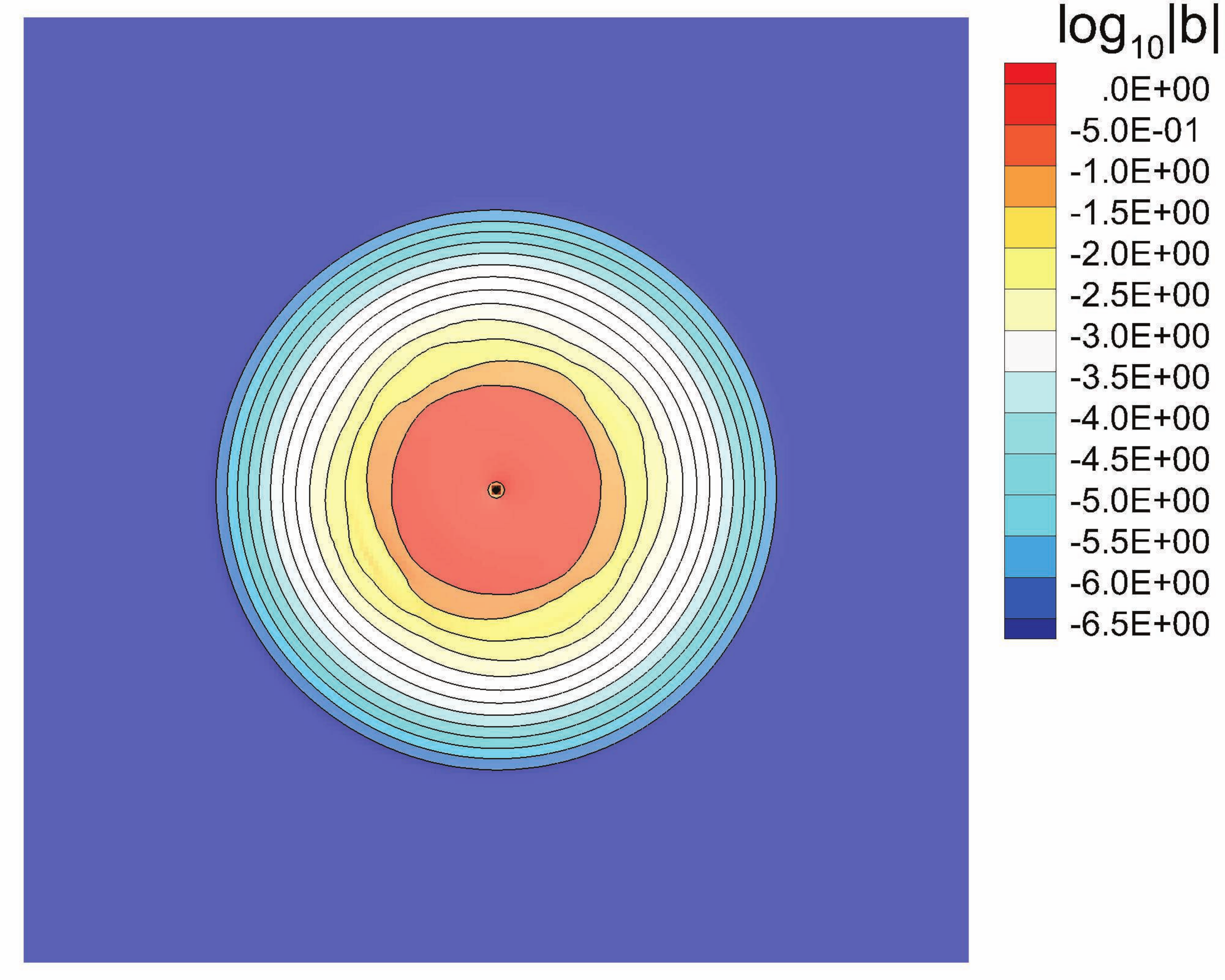} &
			\includegraphics[width=0.42\textwidth]{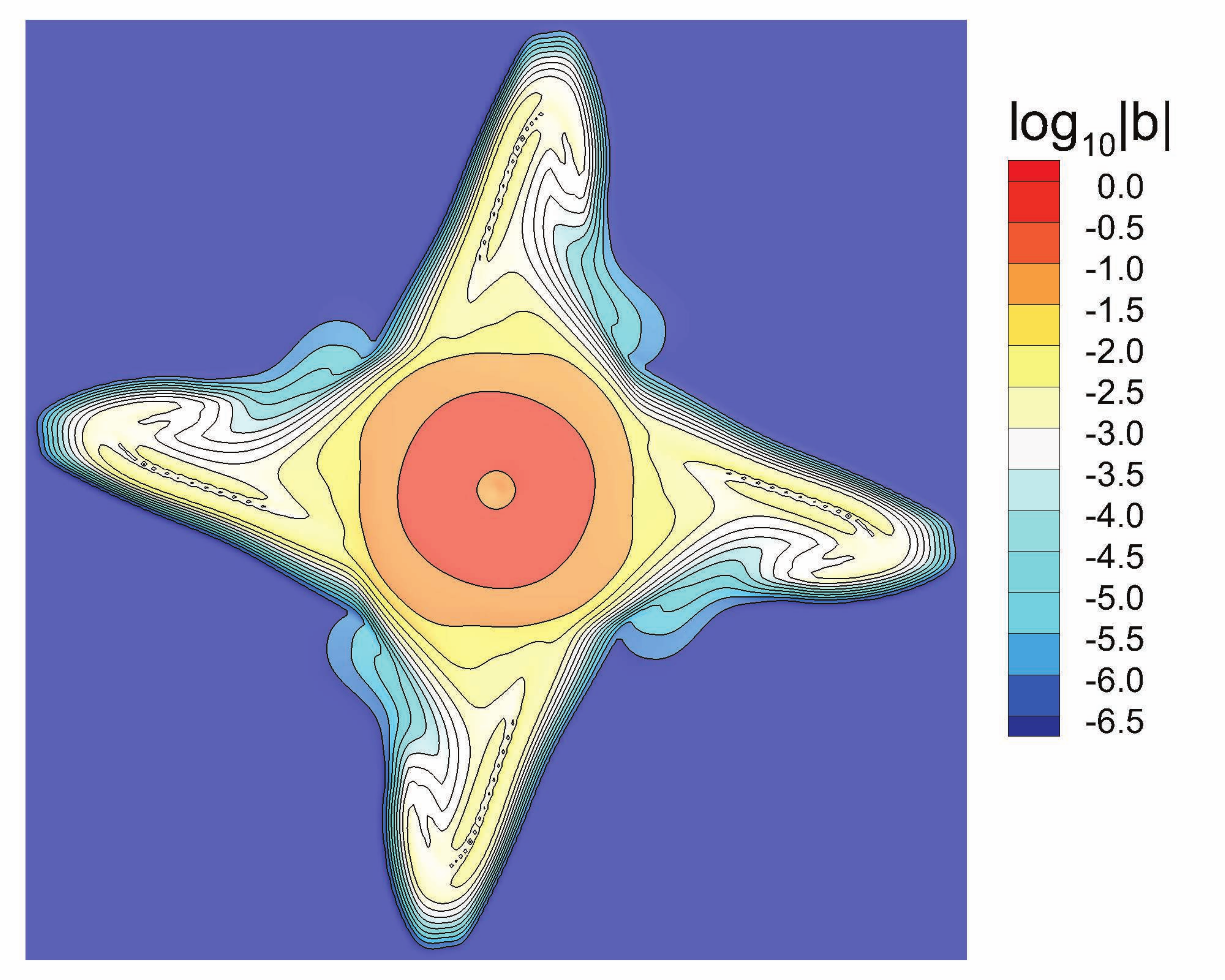}			\\
			\includegraphics[width=0.42\textwidth]{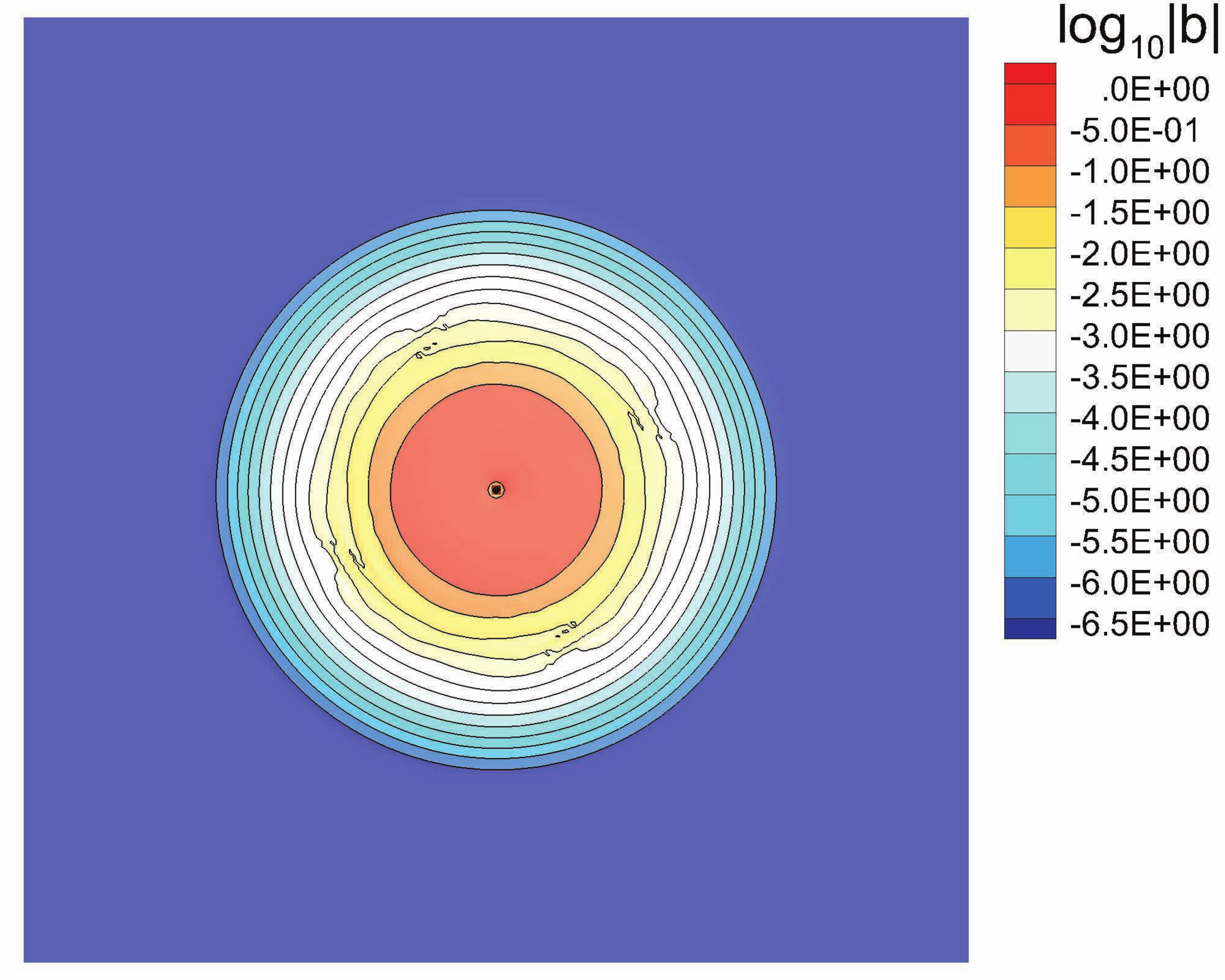} &
			\includegraphics[width=0.42\textwidth]{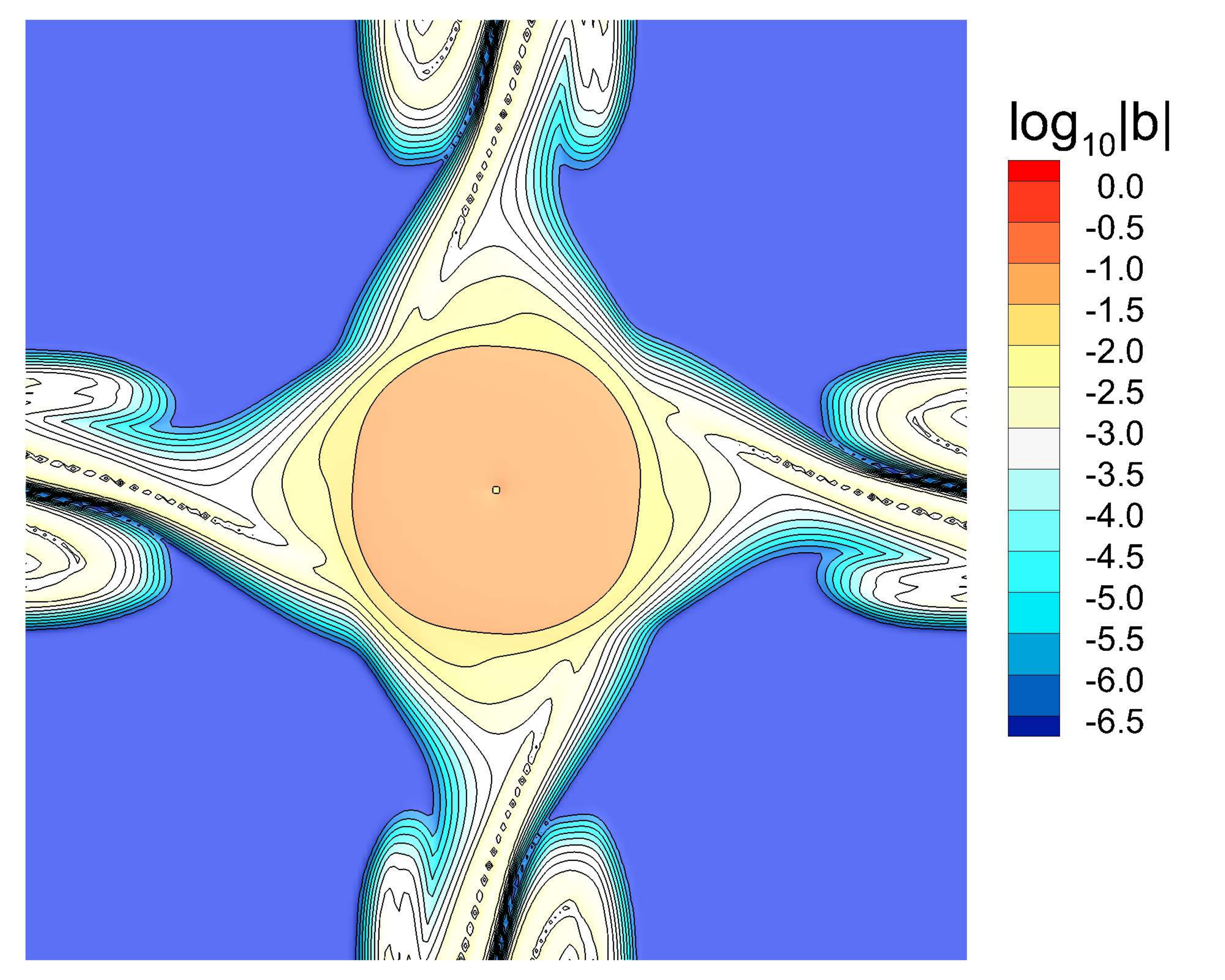} 
			\end{tabular}
\caption{Comparison of the numerical results of the normalized magnetic field $\mathbf{b}:=\B/\max(|\B_{t=0}|)$ for the long time $t_{\text{end}}=1000$ simulation of the MHD stationary vortex test. 
The contour lines of the $\log|\mathbf{b}|$ of the numerical solution obtained with our new semi-implicit structure-preserving algorithm SI-p\&B ({\textbf{left}}) and the pressure-implicit algorithm \cite{SIMHD} SI-p\ ({\textbf{right}}) are plotted for the time slice t=500 (\textbf{first row}) and t=1000 (\textbf{second row}). 
} \label{fig:IsodensityVortex}
\end{figure}

\begin{figure*} 
\centering 
\begin{tabular}{lr}
			\includegraphics[width=0.42\textwidth]{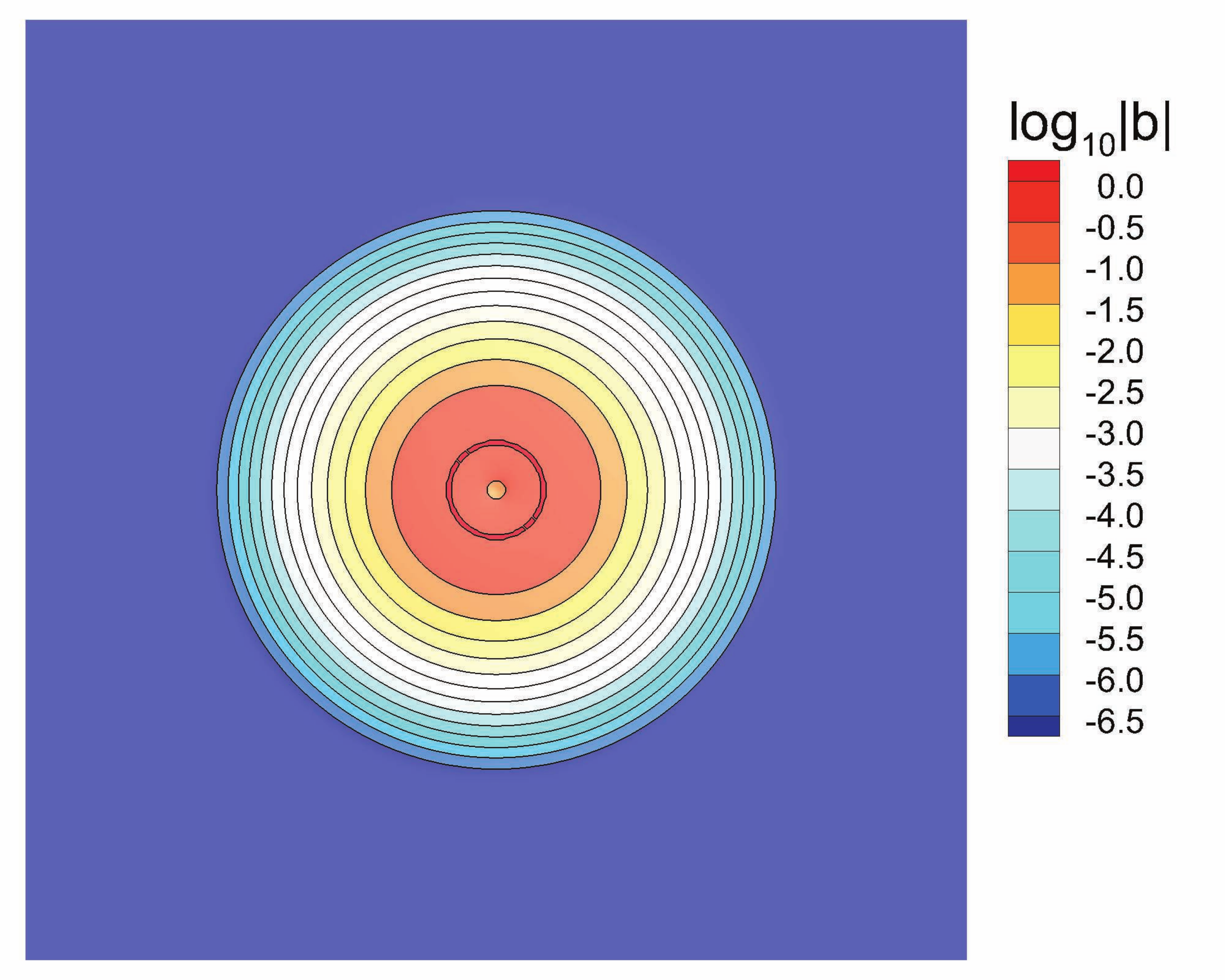} &
			\includegraphics[width=0.42\textwidth]{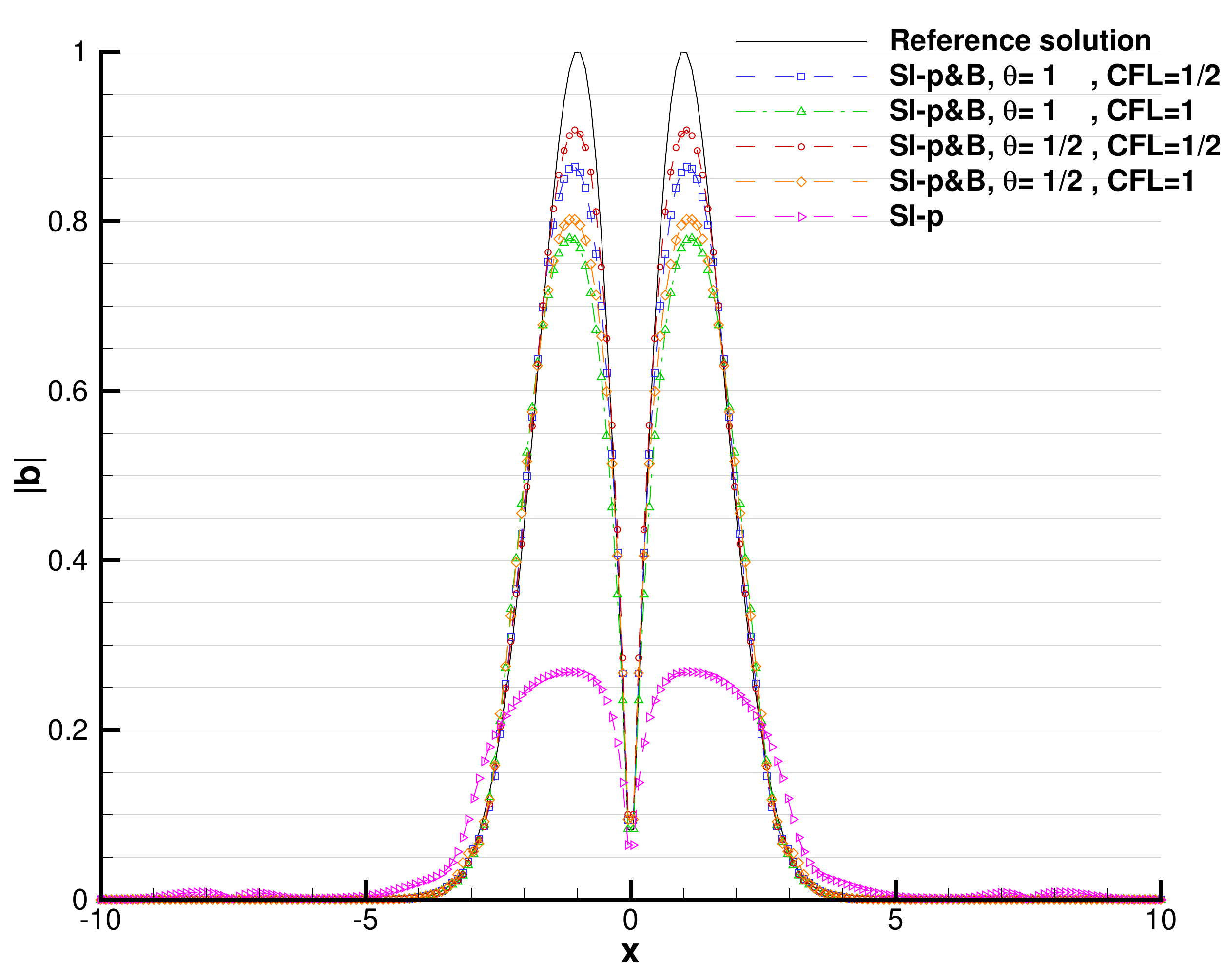}	
			\end{tabular}
\caption{Numerical results for the normalized magnetic field $\mathbf{b}:=\B/\max(|\B_{t=0}|)$ of the long time $t_{\text{end}}=1000$ simulation of the MHD stationary vortex test obtained with our new semi-implicit structure-preserving algorithm SI-p\&B at different $\CFL$ and $\theta_{\text{B}}$ parameters, compared to the reference solution and the results obtained by a pressure-implicit algorithm  SI-p  \cite{SIMHD}. 
At the left one has the corresponding reference solution in the two-dimensional domain, while at the right the one-dimensional cut $y=0$ of $\mathbf{b}$ is plotted for the two different schemes next to the reference solution. 
} \label{fig:IsodensityVortex2}
\end{figure*}

\begin{figure*} 
\centering 
\begin{tabular}{lr}
			\includegraphics[width=0.6\textwidth]{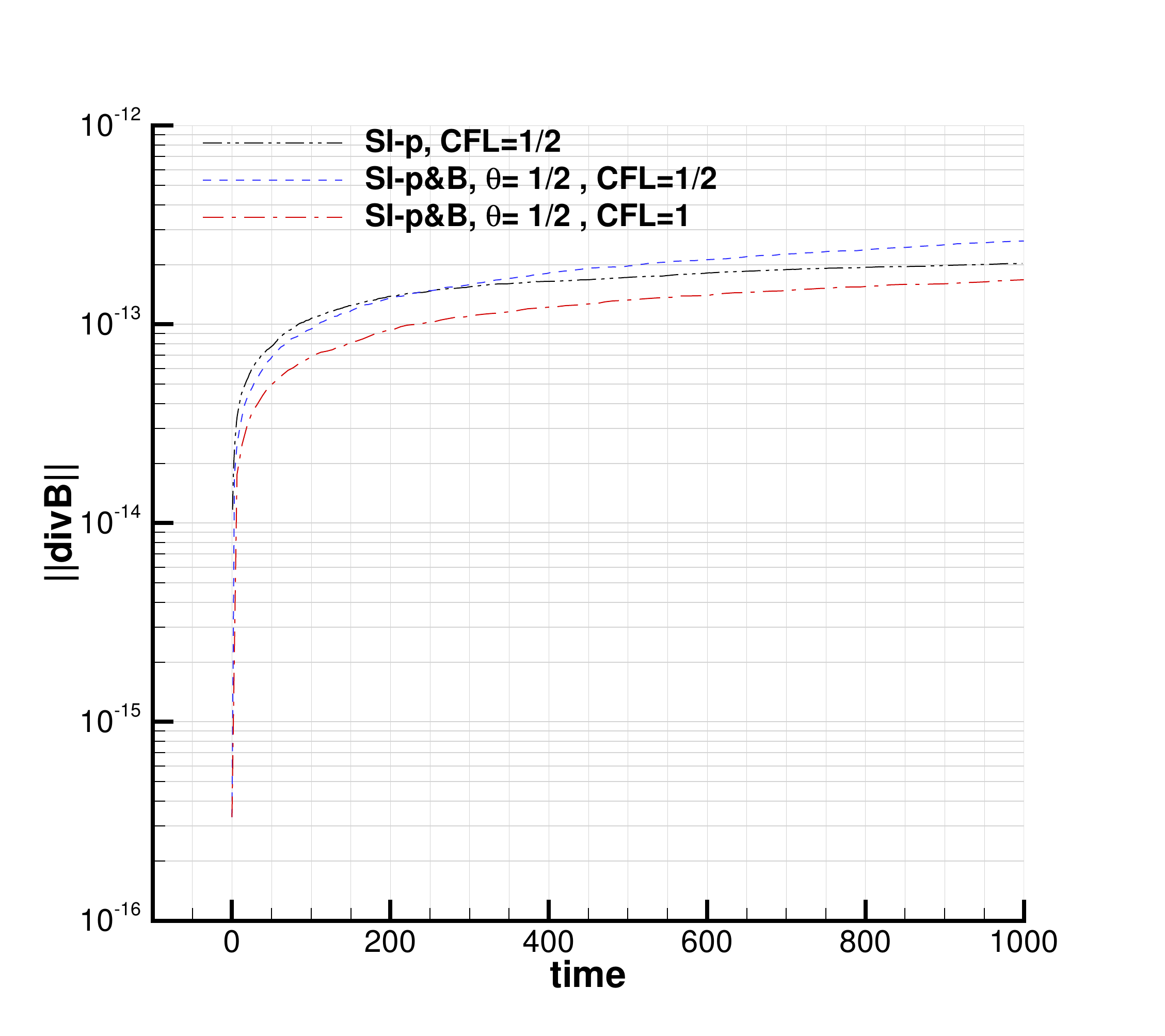} \\
			\includegraphics[width=0.6\textwidth]{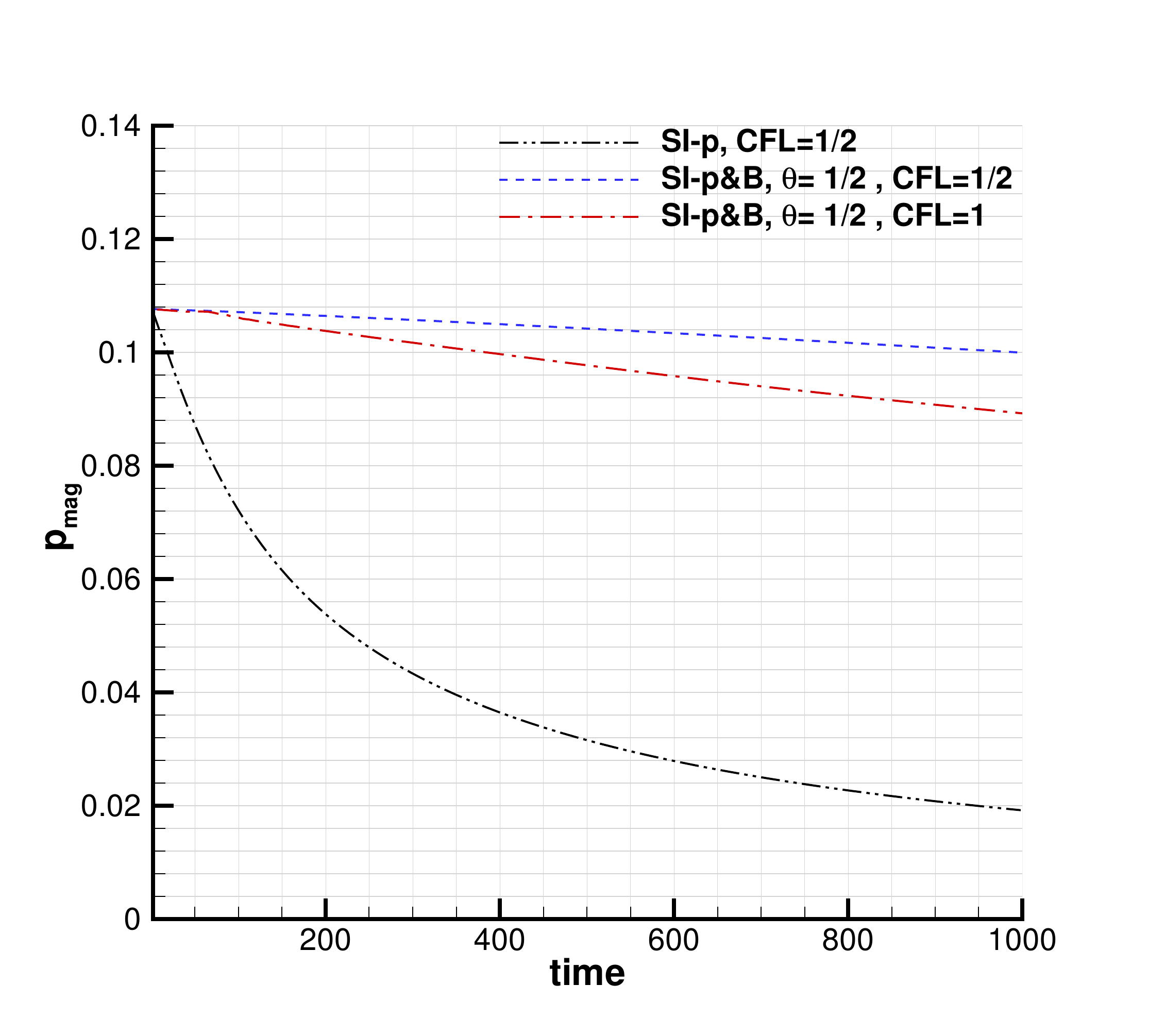} 
	\end{tabular}
\caption{Comparison of the $L_2$-norm of $\div\B$ (\textbf{left}) and $p_{\text{mag}}=\B^2/(8\pi)$ (\textbf{right}) for the long time $t=1000$ simulation of the MHD stationary vortex test obtained with the pressure-based semi-implicit, \cite{SIMHD}, and our new semi-implicit structure-preserving algorithm, with $\CFL=1/2$  and $\CFL=1$ on a $200\times 200$ grid.} \label{fig:IsodensityVortex_divB}
\end{figure*}

Here, the results of our novel three-split semi-implicit scheme are compared with those obtained by our previous two-split semi-implicit scheme \cite{SIMHD}.
Fig. \ref{fig:IsodensityVortex}-\ref{fig:IsodensityVortex2} show the numerical results obtained at two different times $t=500$ and $1000$ for the magnetic field, compared to the reference (initial) solution, together with a one-dimensional cut at $y=0$. 
Here one can easily verify that the computed solution significantly improved in quality with our new divergence-free semi-implicit three-split solver. Indeed, it is shown to be much less dissipative in the magnetic energy and it is able to preserve the stationarity of the initial state, while the numerical solution of the explicit solver is first strongly diffused, and secondly, it degenerate to a non-symmetric unstable solution. 
Fig.  \ref{fig:IsodensityVortex_divB} shows the time evolution of the magnetic energy and the $L_2$ norm of the divergence $\div{\B}$. 

\paragraph{Ideal Orszag-Tang vortex system.}

\begin{figure} 
\centering 
\begin{tabular}{lr}
			 \includegraphics[width=0.35\textwidth]{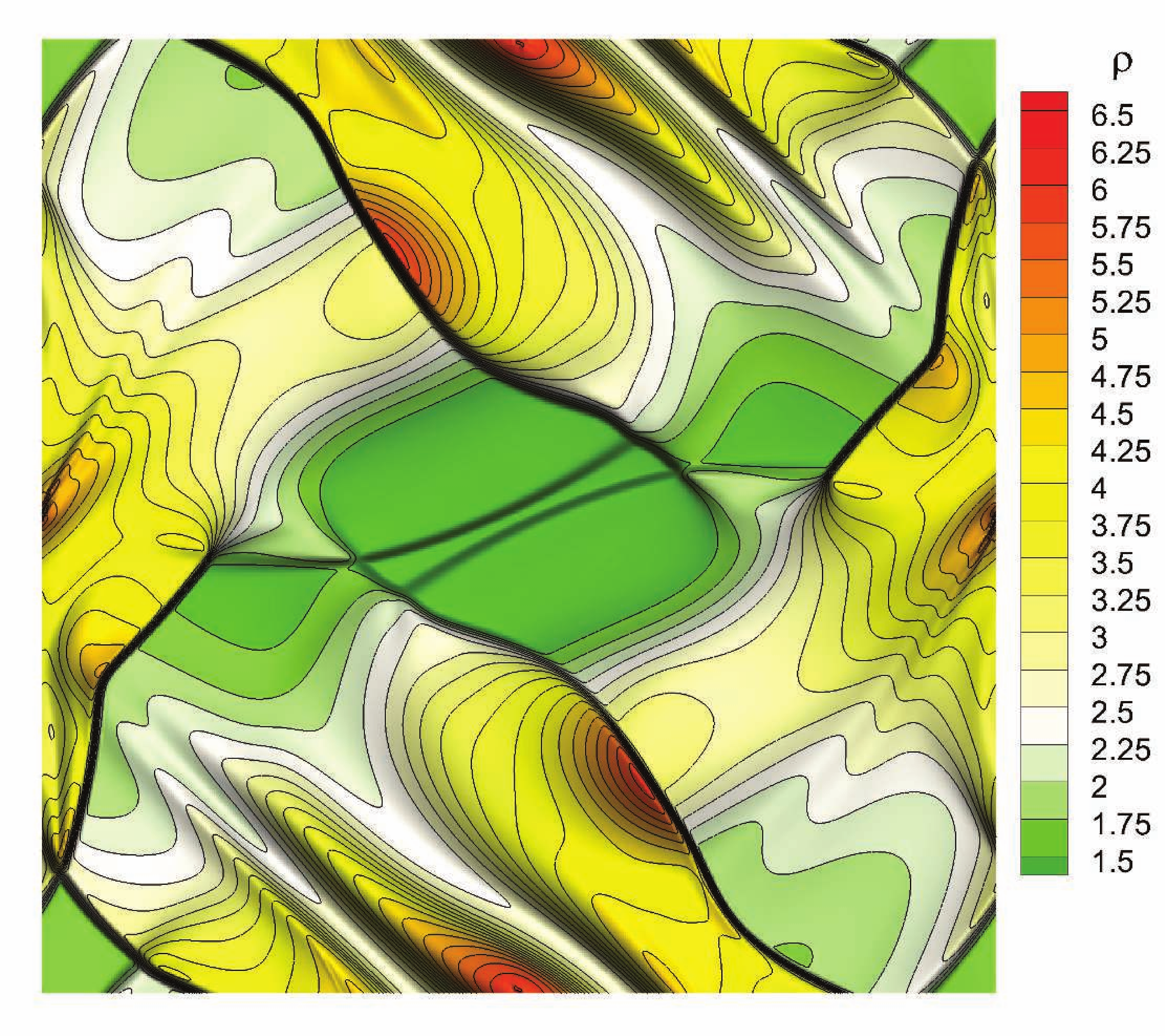} & 
													 \includegraphics[width=0.35\textwidth]{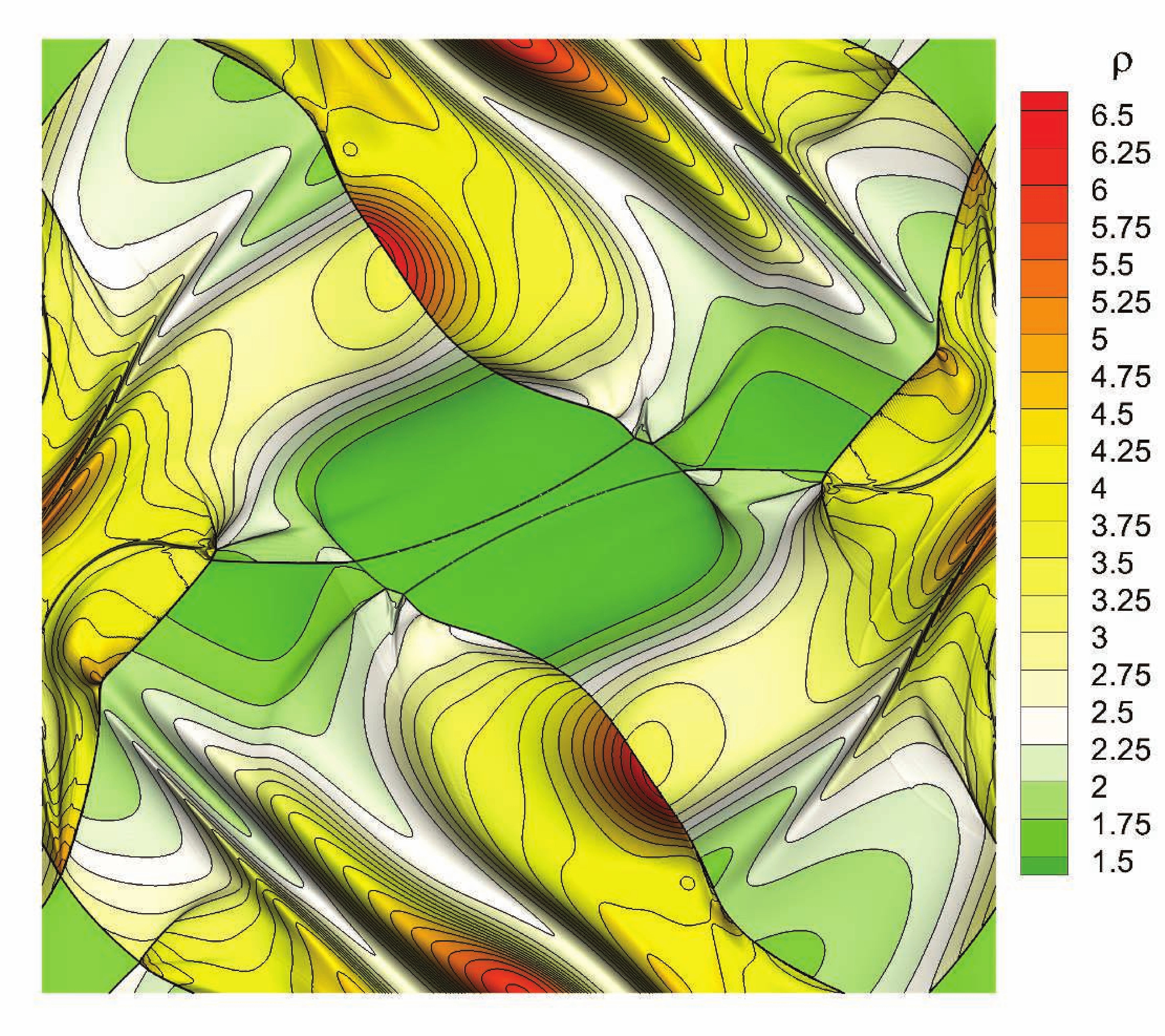} \\ 
			 \includegraphics[width=0.35\textwidth]{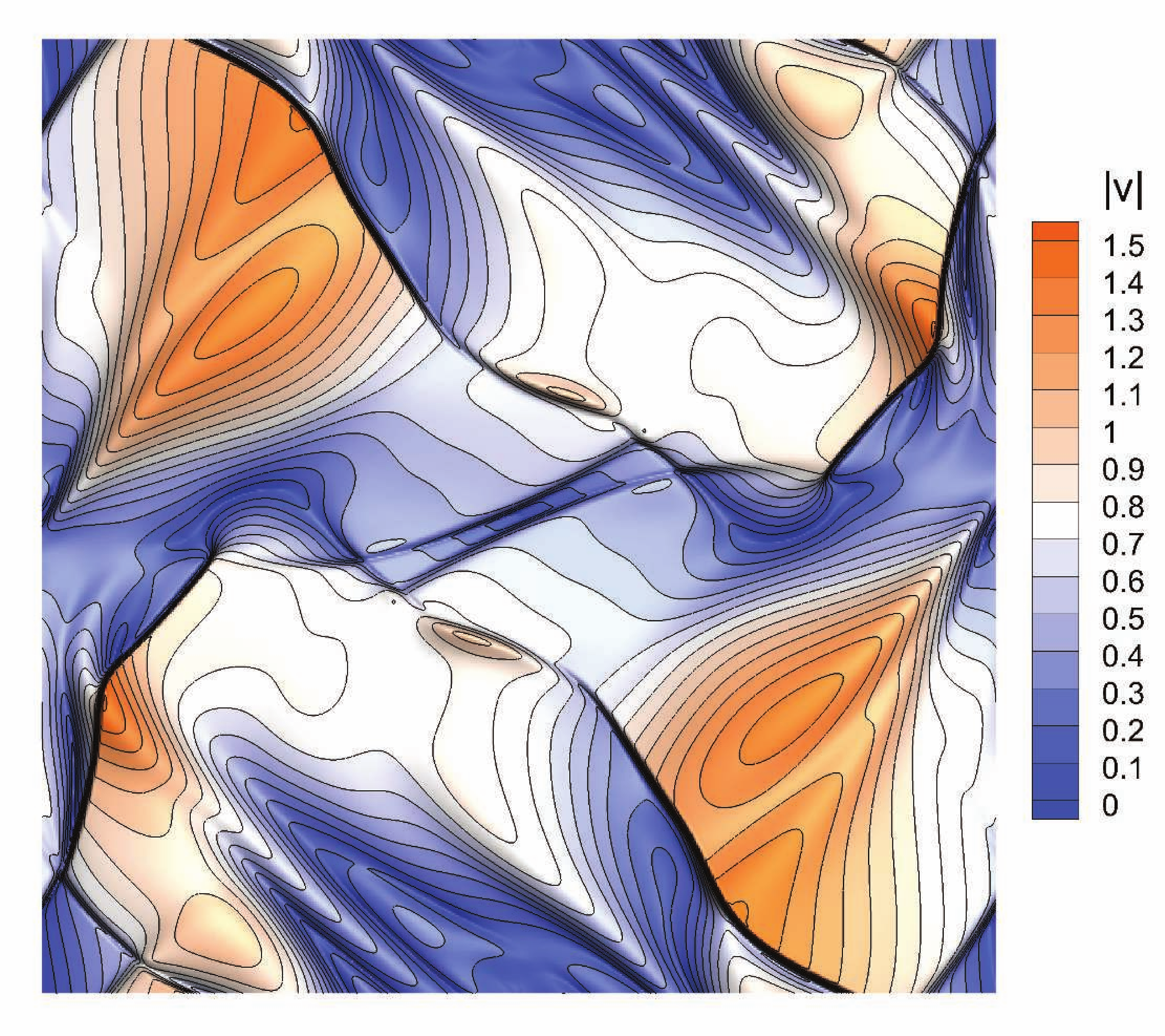} & 
													 \includegraphics[width=0.35\textwidth]{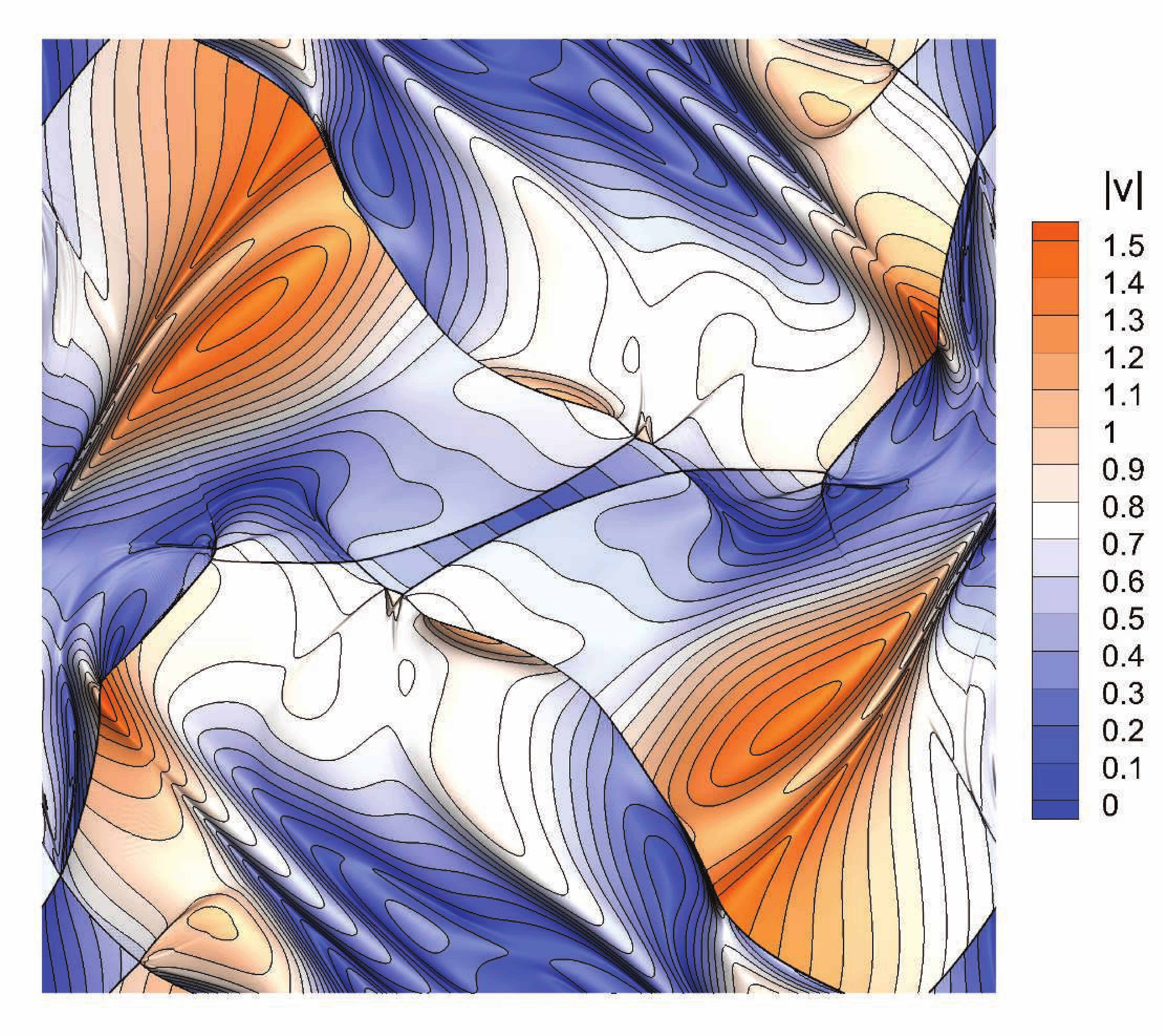} \\ 
			 \includegraphics[width=0.35\textwidth]{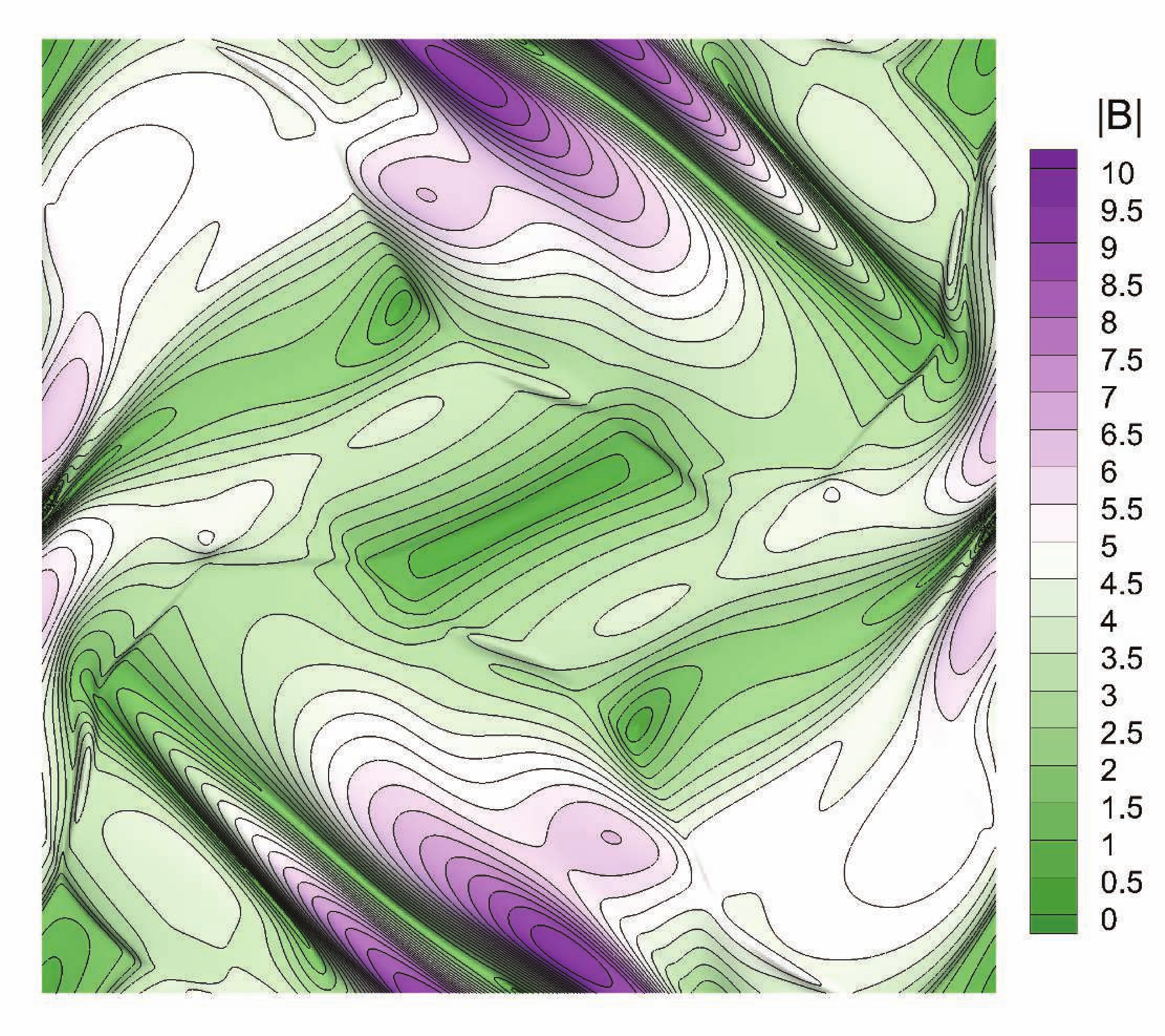} & 
													 \includegraphics[width=0.35\textwidth]{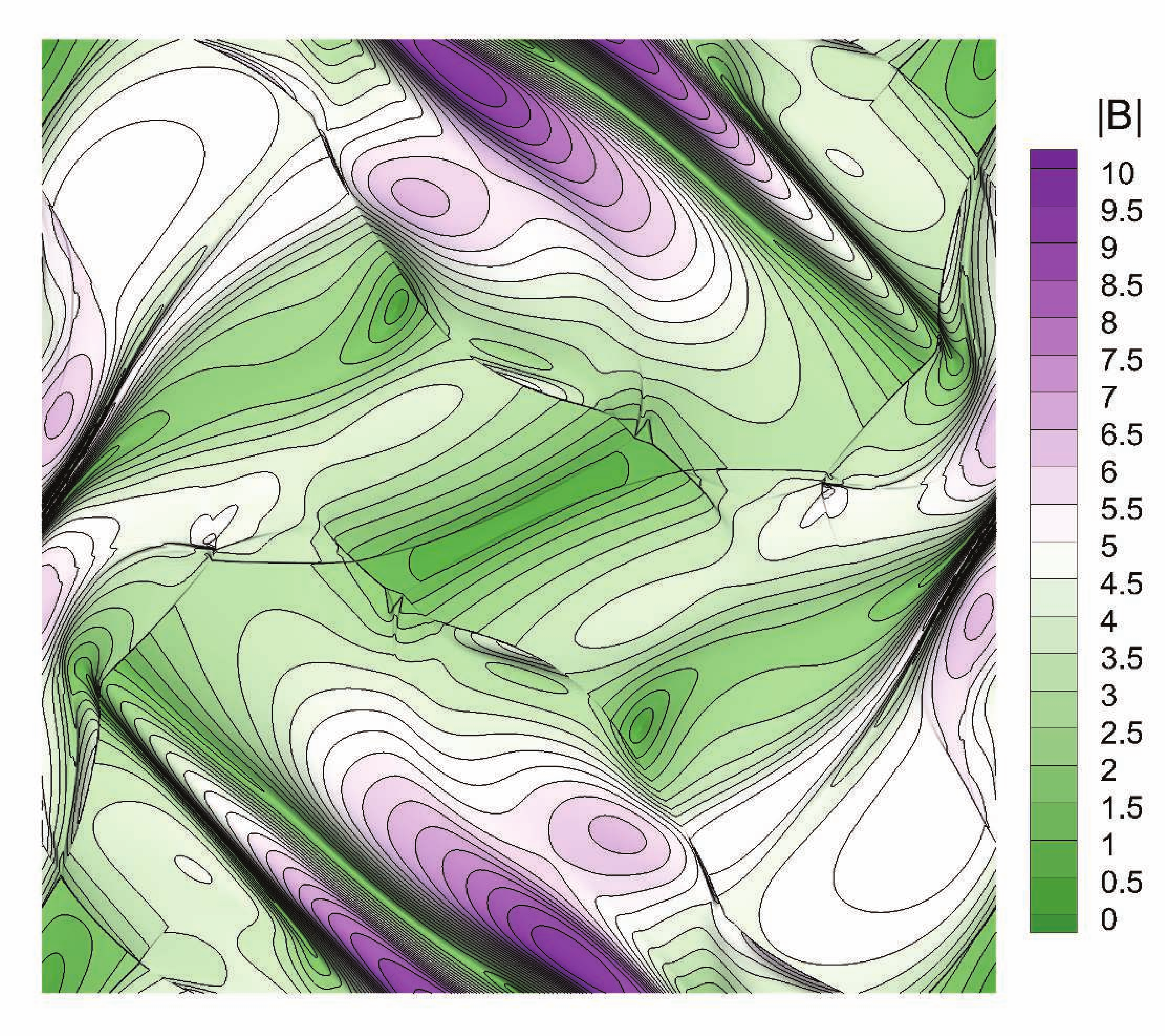} \\ 
			 \includegraphics[width=0.35\textwidth]{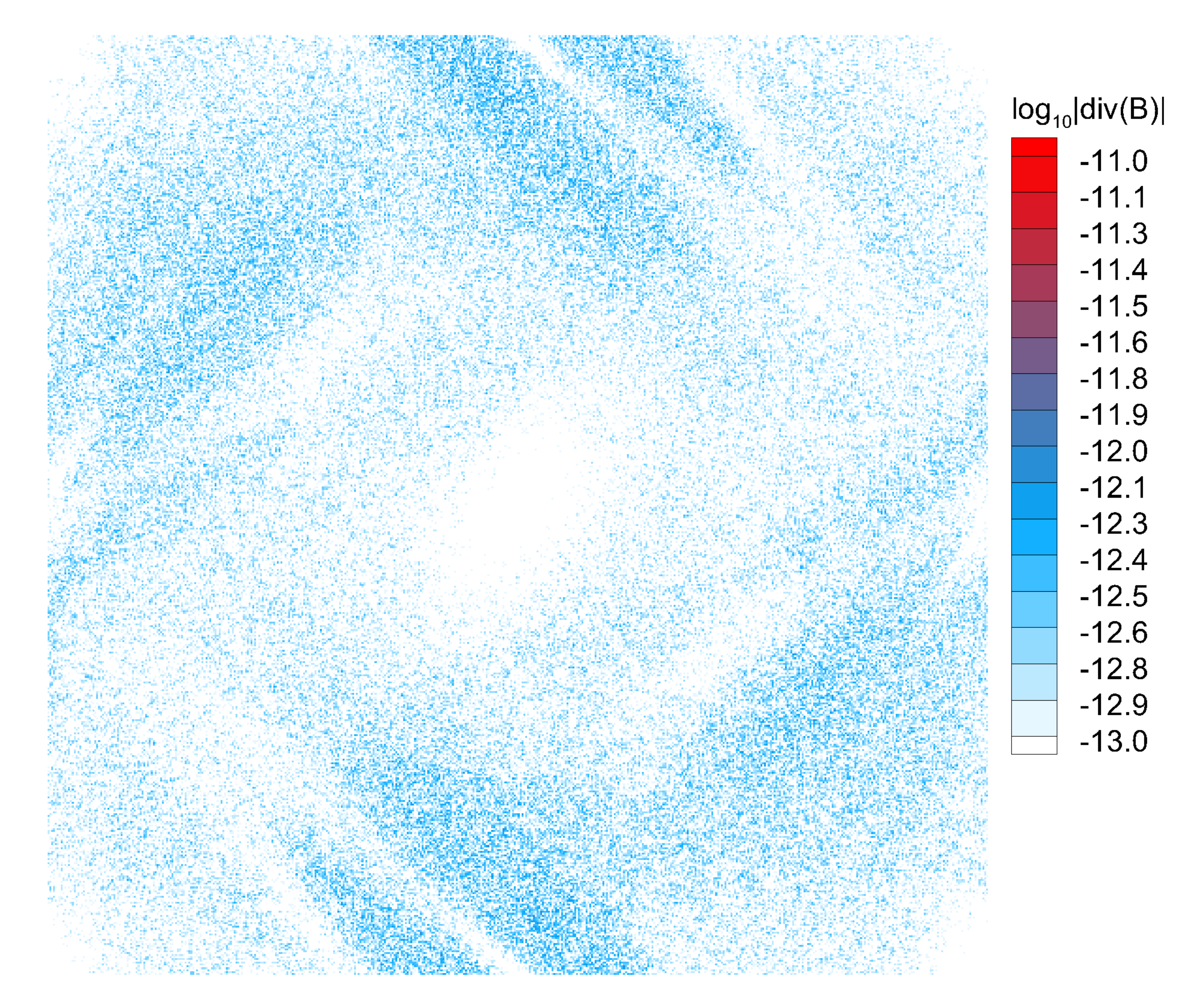} & 
													 \includegraphics[trim={-1.8cm 1cm 1.8cm 1cm},clip,width=0.35\textwidth]{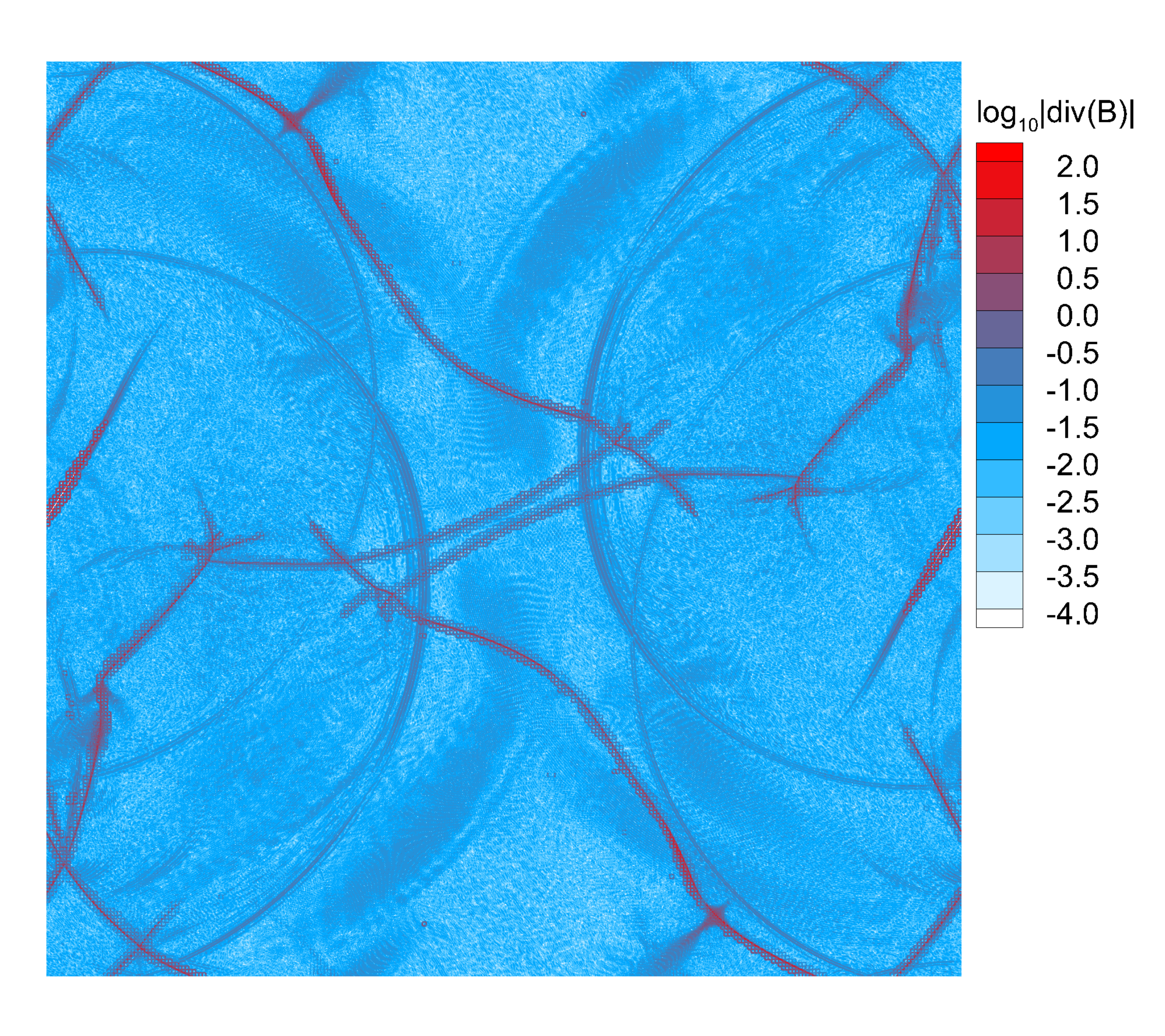}  
	\end{tabular}
\caption{Numerical results for the Orszag-Tang vortex system at time $t=2$ obtained with a $500\times 500$ grid with our new semi-implicit structure-preserving algorithm, after choosing $\CFL=0.9$, from the first row to the last: iso-contour lines for matter-density, absolute value of velocity and magnetic field, and divergence of the magnetic field in log-scale. 
On the right, as a reference solution, we choose a high-order ADER-DG-$\mathbb{P}_5$ scheme with an \emph{a-posteriori} sub-cell limiter with a standard hyperbolic divergence-cleaning, see \cite{Zanotti2015c}.} \label{fig:OT_t20}
\end{figure}

\begin{figure} 
\centering 
\begin{tabular}{lr}
			\includegraphics[width=0.35\textwidth]{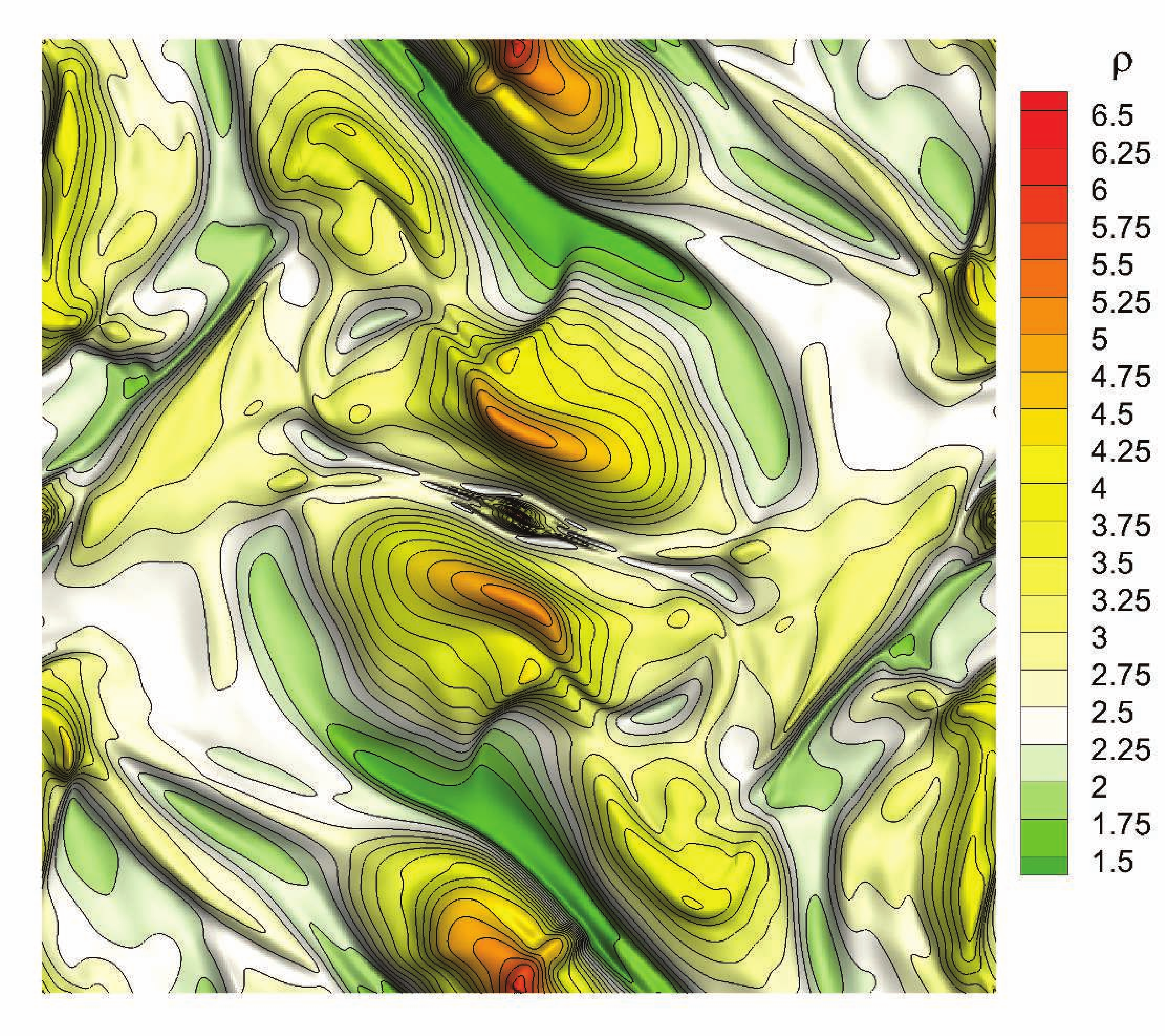} &
													 \includegraphics[width=0.35\textwidth]{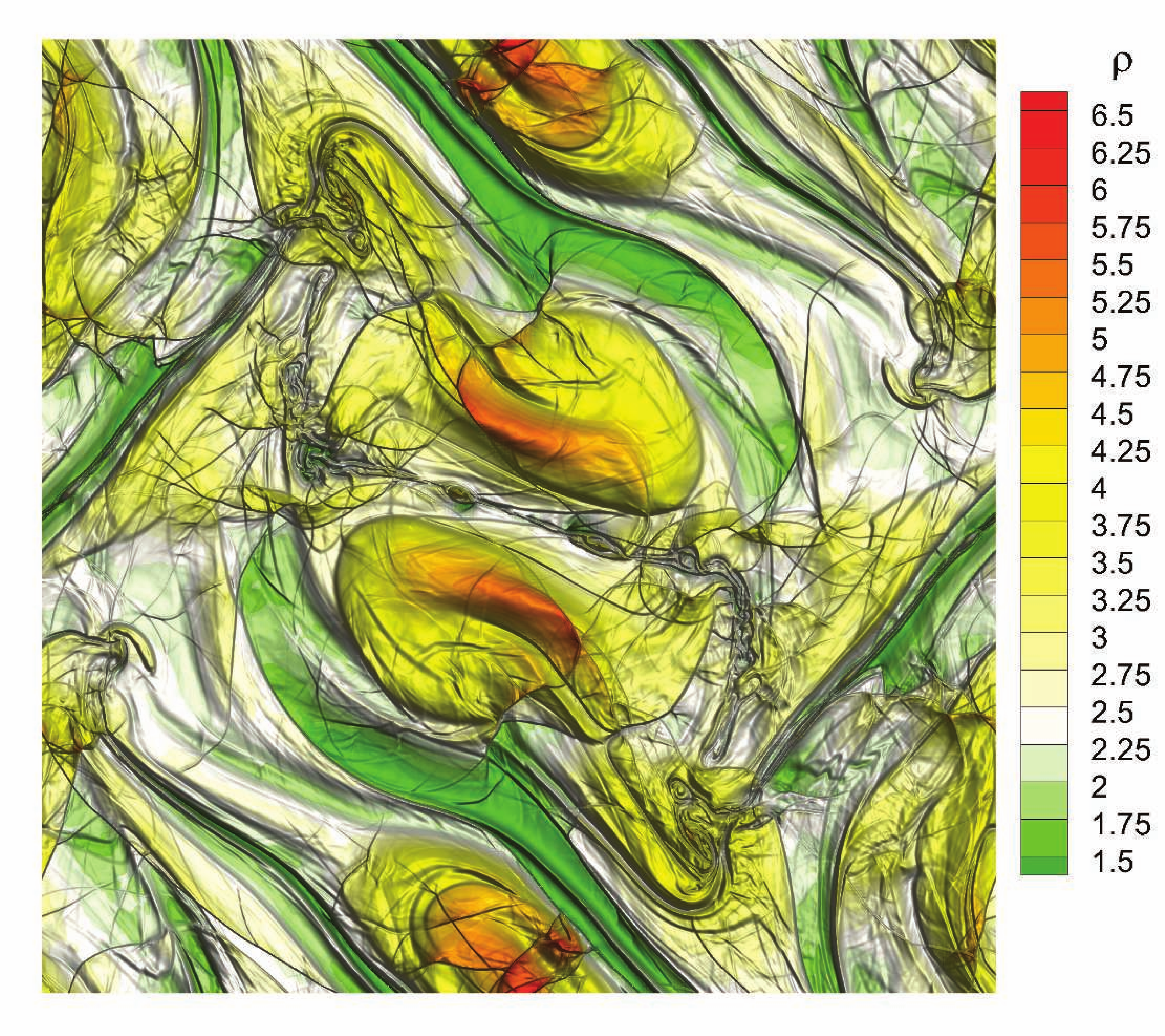} \\ 
			\includegraphics[width=0.35\textwidth]{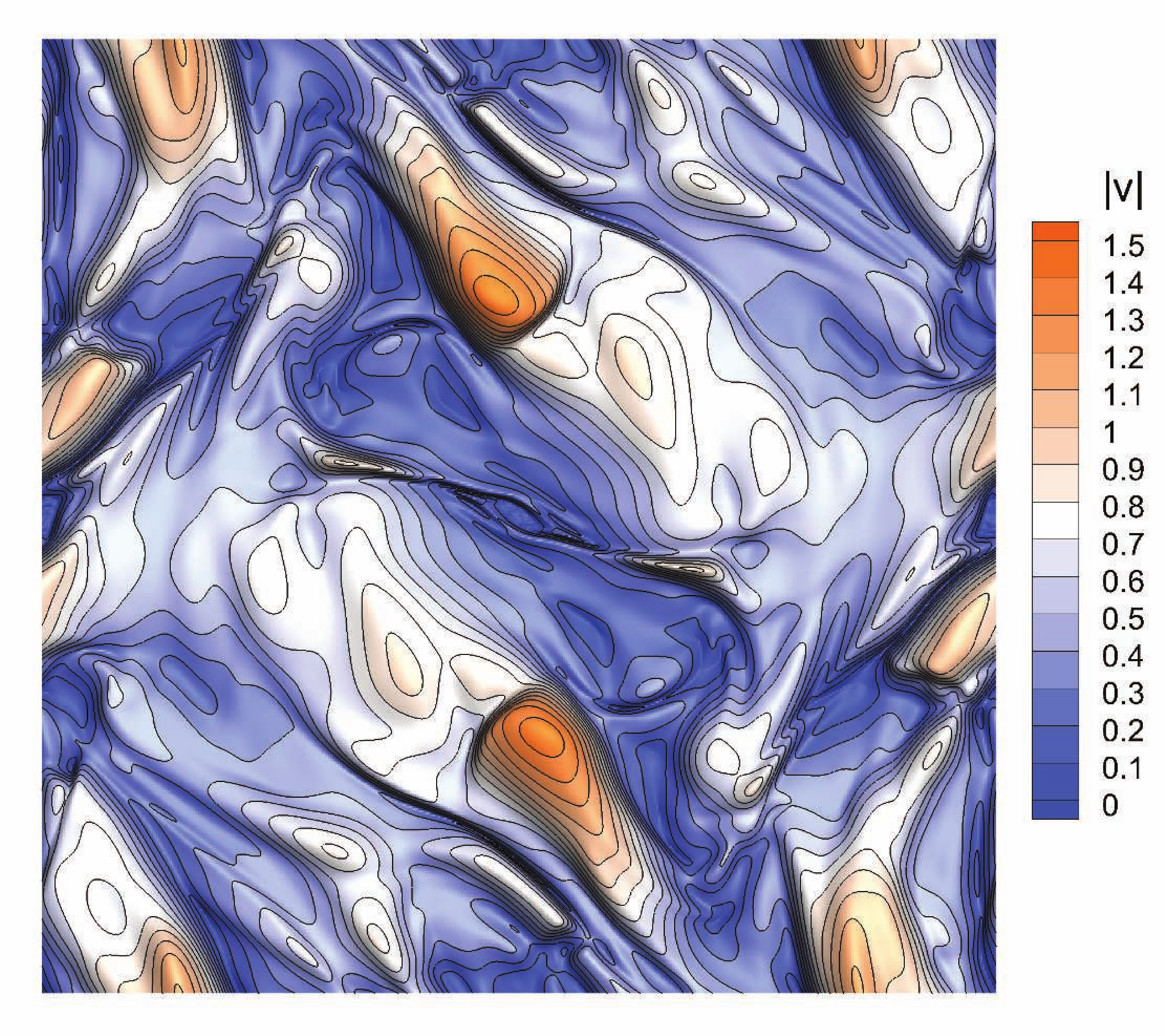} &
													 \includegraphics[width=0.35\textwidth]{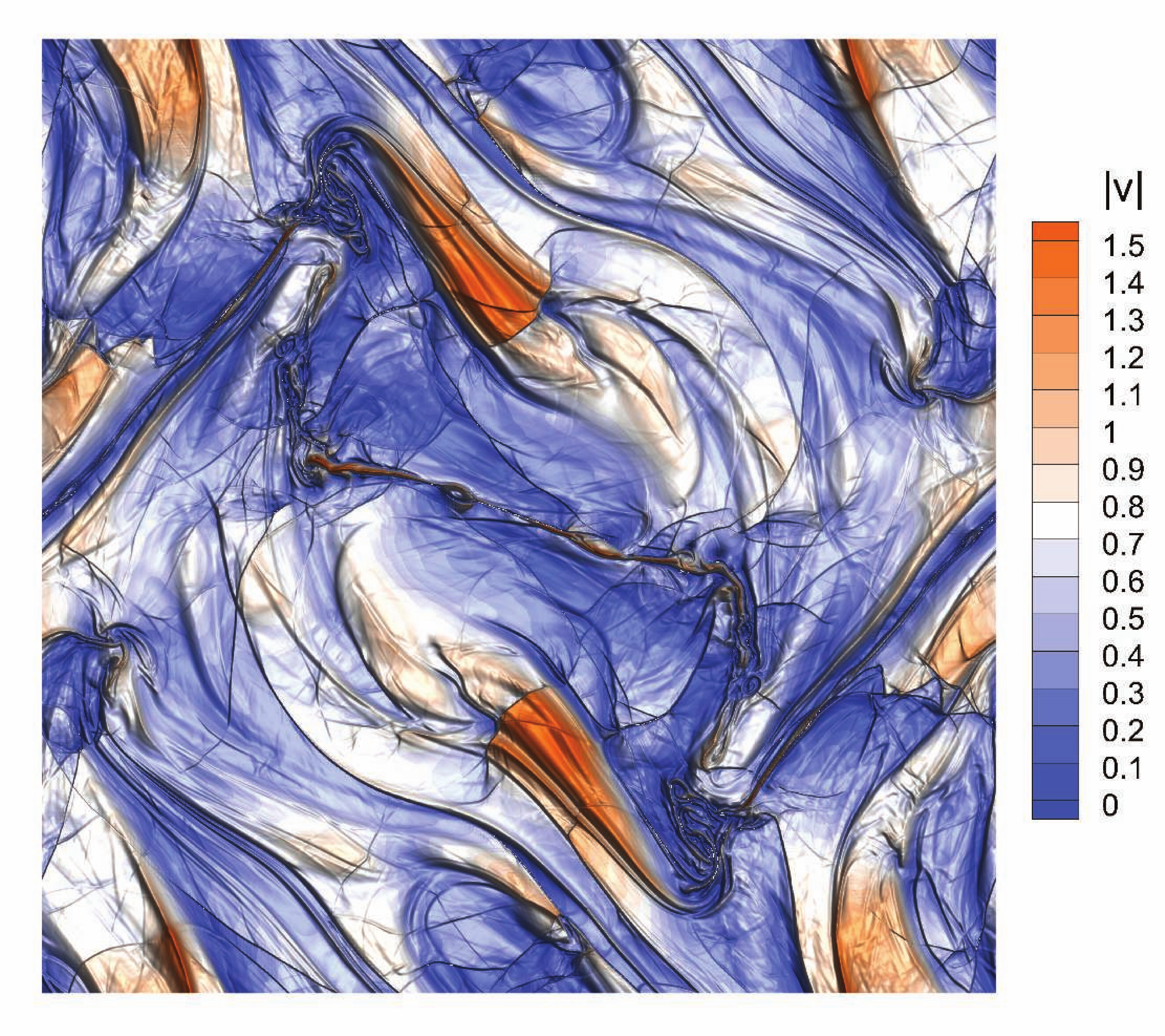}\\ 
			\includegraphics[width=0.35\textwidth]{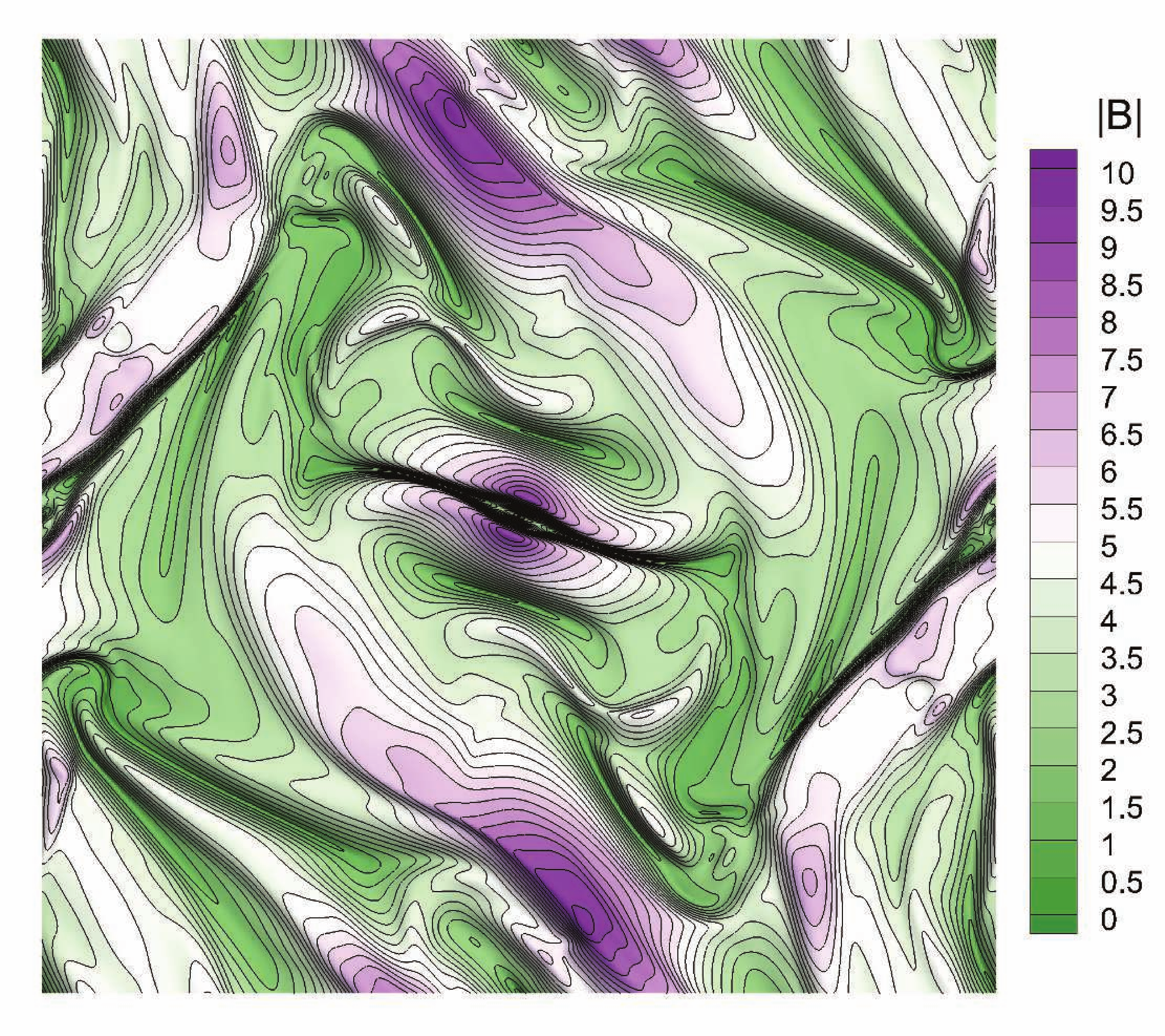} &
													 \includegraphics[width=0.35\textwidth]{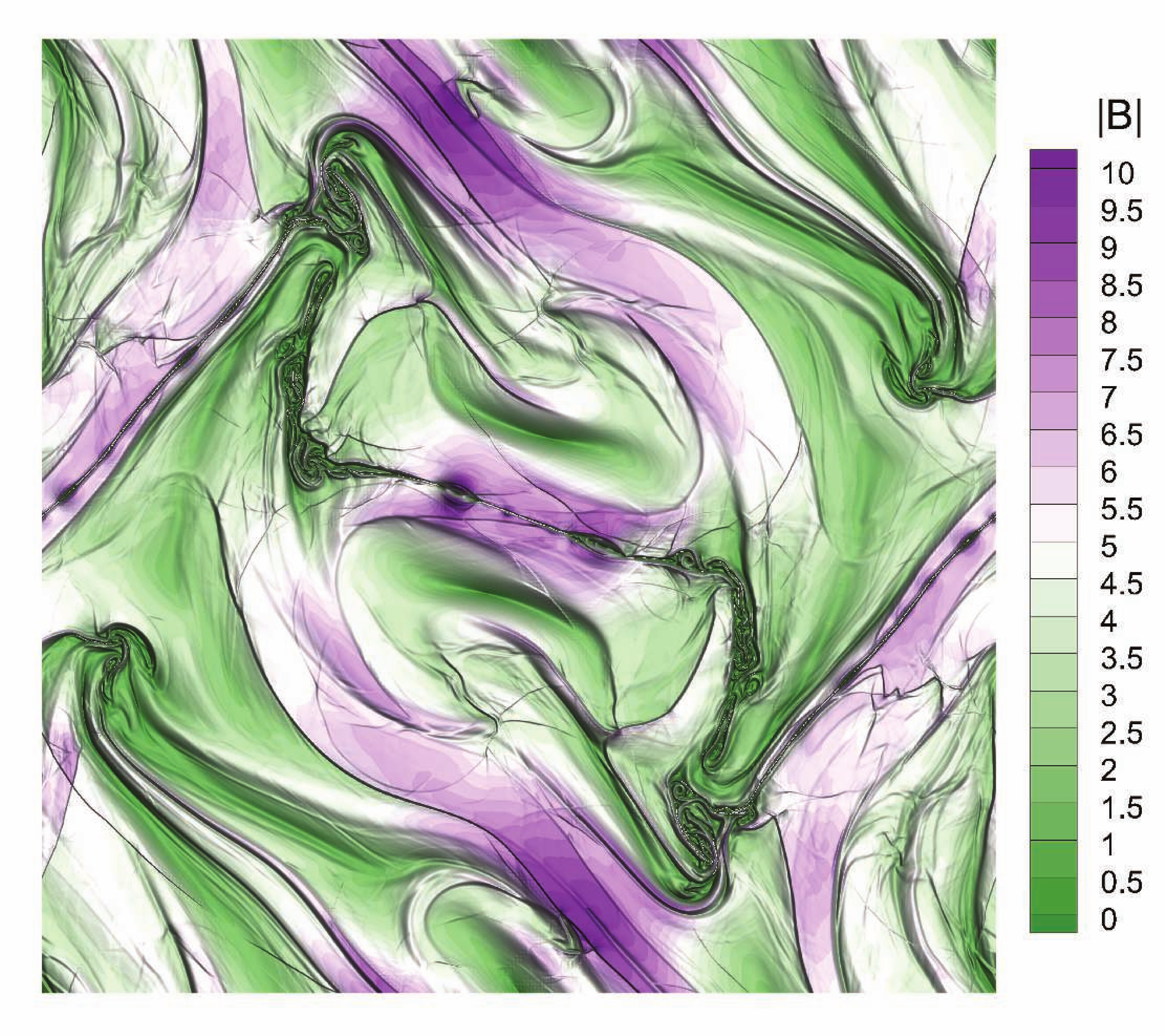} \\ 
			\includegraphics[width=0.35\textwidth]{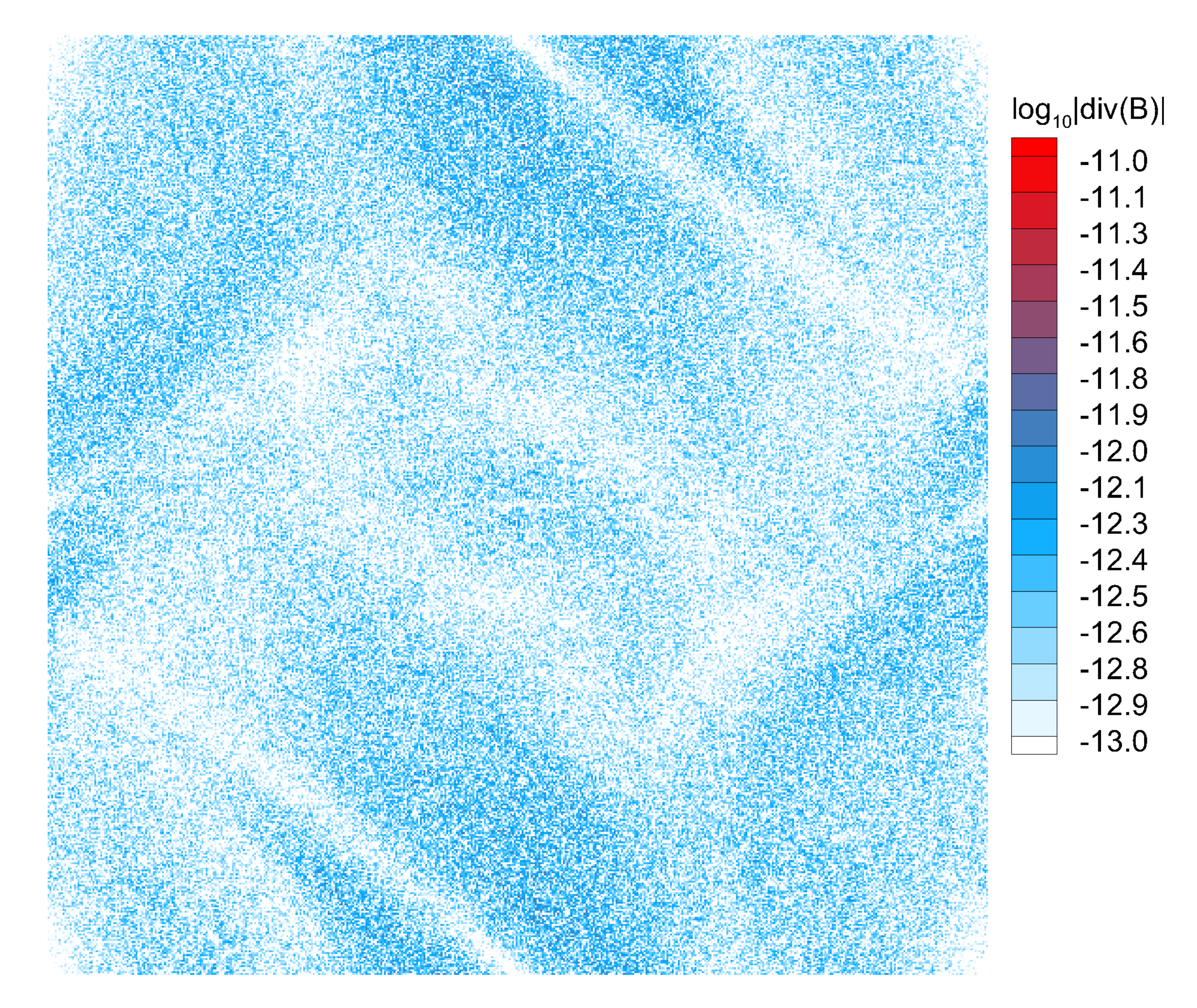} &
													 \includegraphics[trim={-1.8cm 1cm 1.8cm 1cm},clip,width=0.35\textwidth]{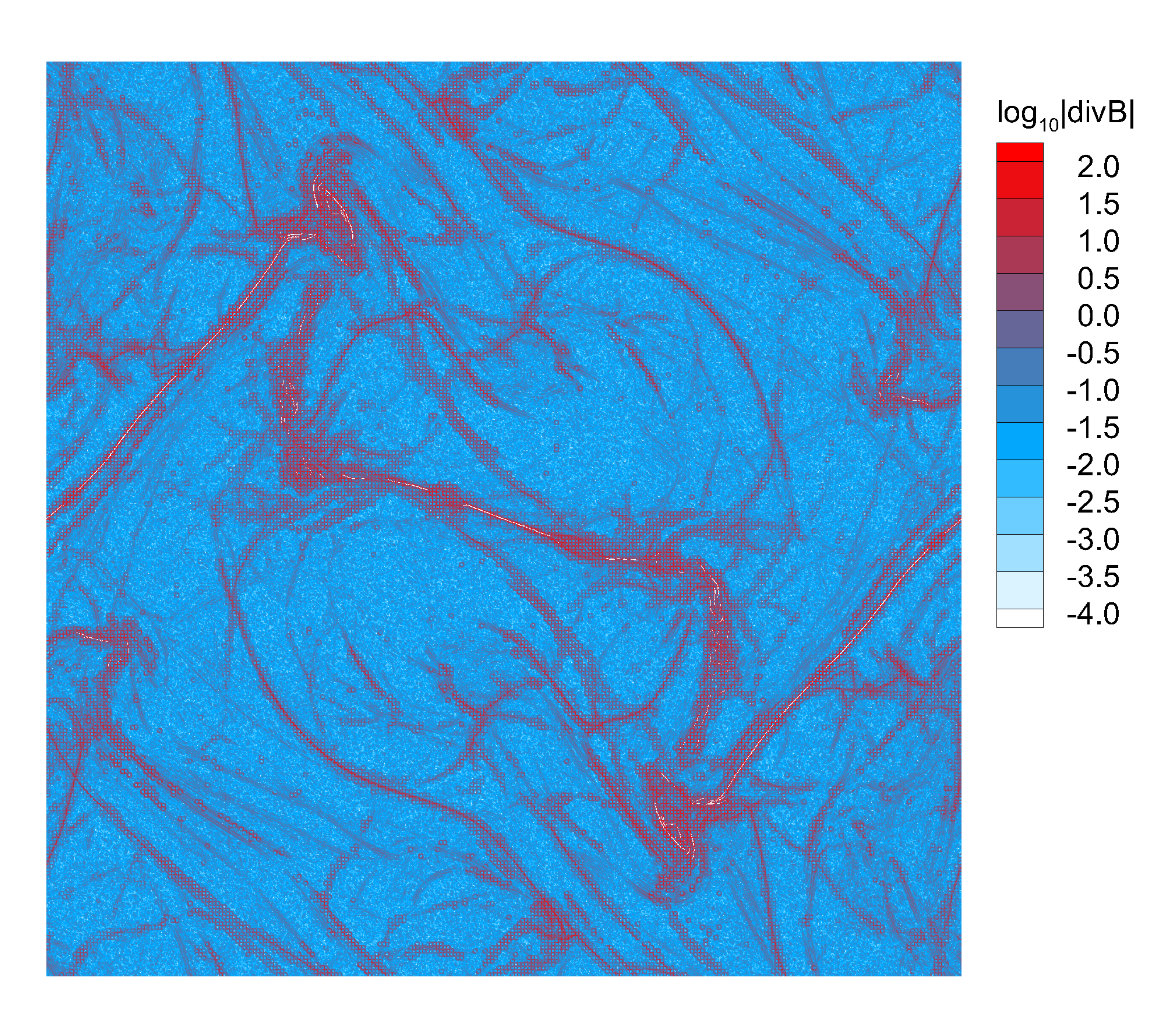}  
	\end{tabular}
\caption{Numerical results for the Orszag-Tang vortex system at time $t=5$ obtained with a $500\times 500$ grid with our new semi-implicit structure-preserving algorithm, after choosing $\CFL=0.9$, from the first row to the last: iso-contour lines for matter-density, absolute value of velocity and magnetic field, and divergence of the magnetic field in log-scale. 
 On the right, as a reference solution, we choose a high-order ADER-DG-$\mathbb{P}_5$ scheme with an a-posteriori sub-cell limiter with a standard hyperbolic divergence-cleaning, see \cite{Zanotti2015c}. Notice that the results differs because some turbulent instabilities are started and the kinematic viscosity is $\mu=0$.} \label{fig:OT_t50}
\end{figure}

Here, we present again the well-known Orszag-Tang vortex problem, originally proposed in \cite{OrszagTang}, and further investigated by \cite{PiconeDahlburg} and \cite{DahlburgPicone}. 
In this work, we refer to the configuration set by \cite{JiangWu}. 
%
The computational domain is a two-dimensional periodic box $\Omega=[0,2\pi]^2$, resolved with a $500\times 500$ grid, leading to a mesh resolution of about $\Delta x =\Delta y \sim 10^{-2}$.
The initial conditions are given by 
\myequation{lcl}
{
\rho = \gamma^2, & &\v  =  \left( - \sin\left(y\right), \sin \left(x \right), 0 \right), 						 		\nonumber\\
p = \gamma, 			& &	\B =  \sqrt{4\pi}\left( -  \sin\left(y\right), \sin \left(2x \right), 0\right).
}
{\label{eq:OrszagTang_ic}
}
where the ratio of specific heats is set to $\gamma=5/3$. In Fig. \ref{fig:OT_t20}-\ref{fig:OT_t50} some physical quantities related to the computed solution at times $t=2$ and $5$, respectively, are plotted next to a numerical reference solution.
The implicit solver is set up with $\theta_{\text{B}}=0.65$ and $\theta_{\text{p}}=1$

 For this test, the reference solution has been computed by means of a space-time ADER-DG-$\mathbb{P}_5$ scheme, supplemented by an \emph{a-posteriori} sub-cell WENO3 limiter, applied to the ideal MHD equations and hyperbolic divergence cleaning, and it has been published in \cite{Zanotti2015c}. For the history and an overview of high order ADER schemes, see \cite{toro4,titarevtoro,Toro:2006a,BTVC2016,FrontierADERGPR}.  
 
\begin{figure} 
\centering 
\begin{tabular}{lr}  
			 \includegraphics[width=0.5\textwidth]{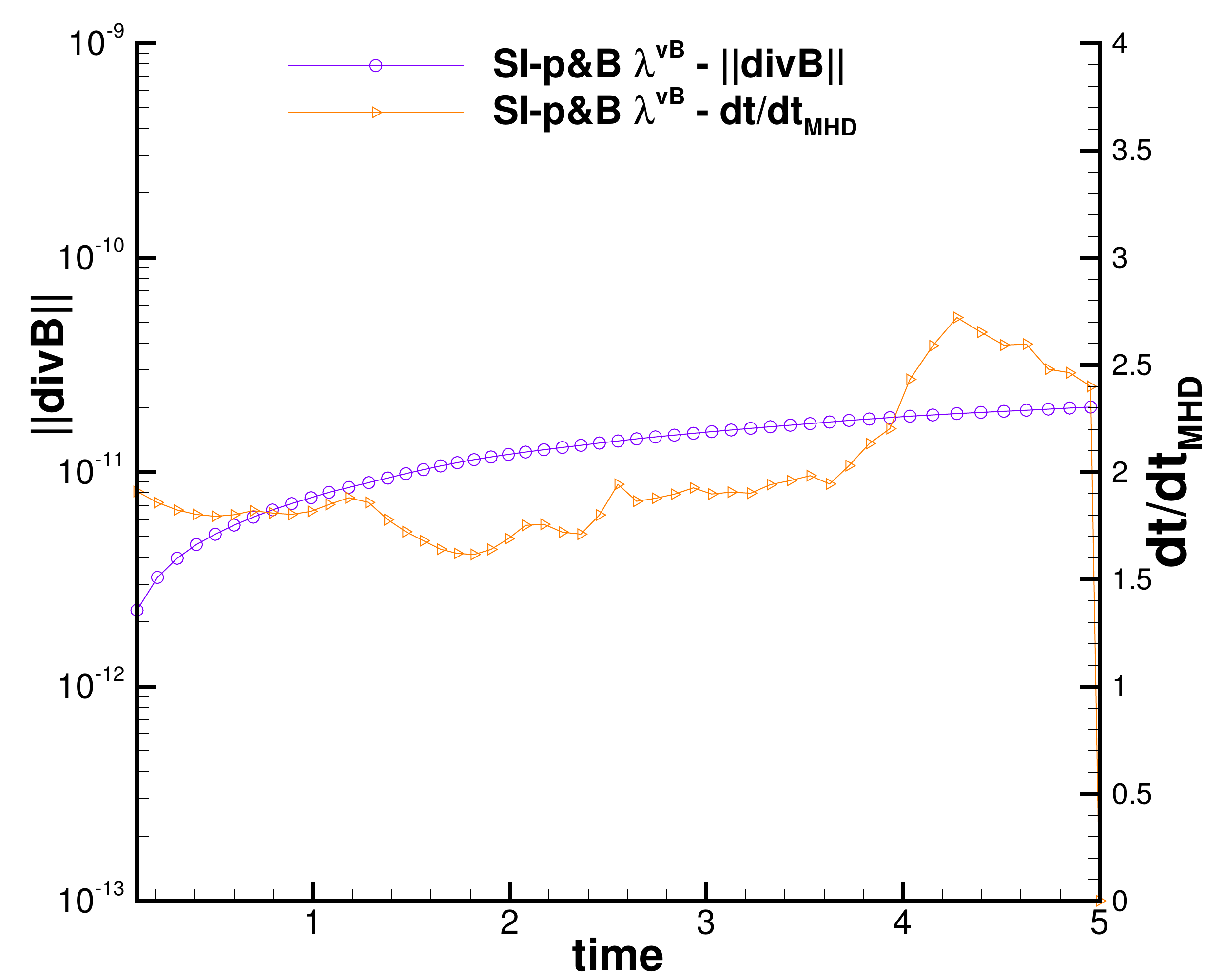} 
	\end{tabular}
\caption{Numerical results for the Orszag-Tang vortex system at time $t=0.5$ obtained with a $500\times 500$ grid with our new semi-implicit structure-preserving algorithm, after choosing $\CFL=0.9$: $L_2$ norm of the divergence of the magnetic field $||\div \B||$ and the ratio between the used computational time-step and the time-step as it would be obtained using the full set of eigenvalues of MHD, i.e. the time evolution of the Courant-Friedrichs-Levy number evaluated with respect to the full MHD system. These data are sampled every 20 time-steps.} \label{fig:OT_divB}
\end{figure}

Considering the numerical resolution, the computed results can be considered to be well compatible with the selected high-order reference solution. At time $t=2.0$, when the first flow-discontinuities are already generated, the two numerical solutions still show a very good agreement; while at time $t=5$, when some flow-instabilities becomes more dominant in the computational domain, the different mesh resolution becomes evident.
On the other hand, it turns interesting to look at the $\div\B$ errors obtained by the two different algorithms. In this case, our novel divergence-free three-split semi-implicit scheme shows a maximum of about $\div\B \sim 10^{-12}$, see again  the last rows in Fig. \ref{fig:OT_t20}-\ref{fig:OT_t50} and Fig.  \ref{fig:OT_divB} where the time evolution of the $L_2$ norm $||\div\B||$ is plotted together with the estimated Courant number, i.e. the ratio $dt(\lambda^{v})/dt(\lambda^{\text{MHD}})$.

\paragraph{Viscous and resistive Orszag-Tang vortex system.}

Here, the Orszag-Tang vortex system is set up again for the full viscous and resistive MHD equations, with flow parameters $\gamma = \frac{5}{3}$, $\mu = \eta = 10^{-2}$, $c_v=1$ and a Prandtl number of $Pr=1$, see \cite{WarburtonVRMHD,ADERVRMHD}.

\begin{figure} 
\centering 
\begin{tabular}{lcr}
			 \includegraphics[width=0.45\textwidth]{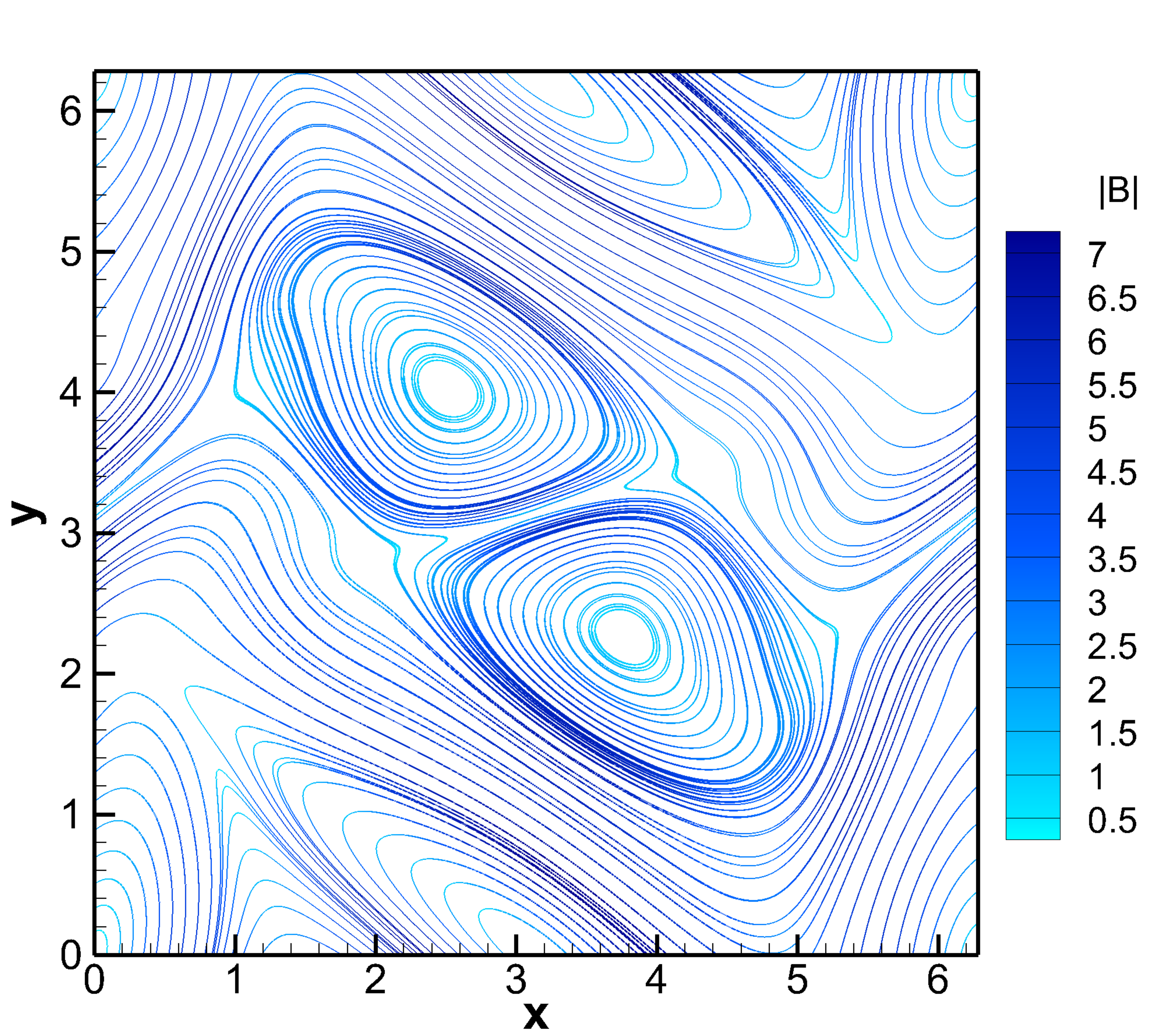} & 
			\includegraphics[width=0.45\textwidth]{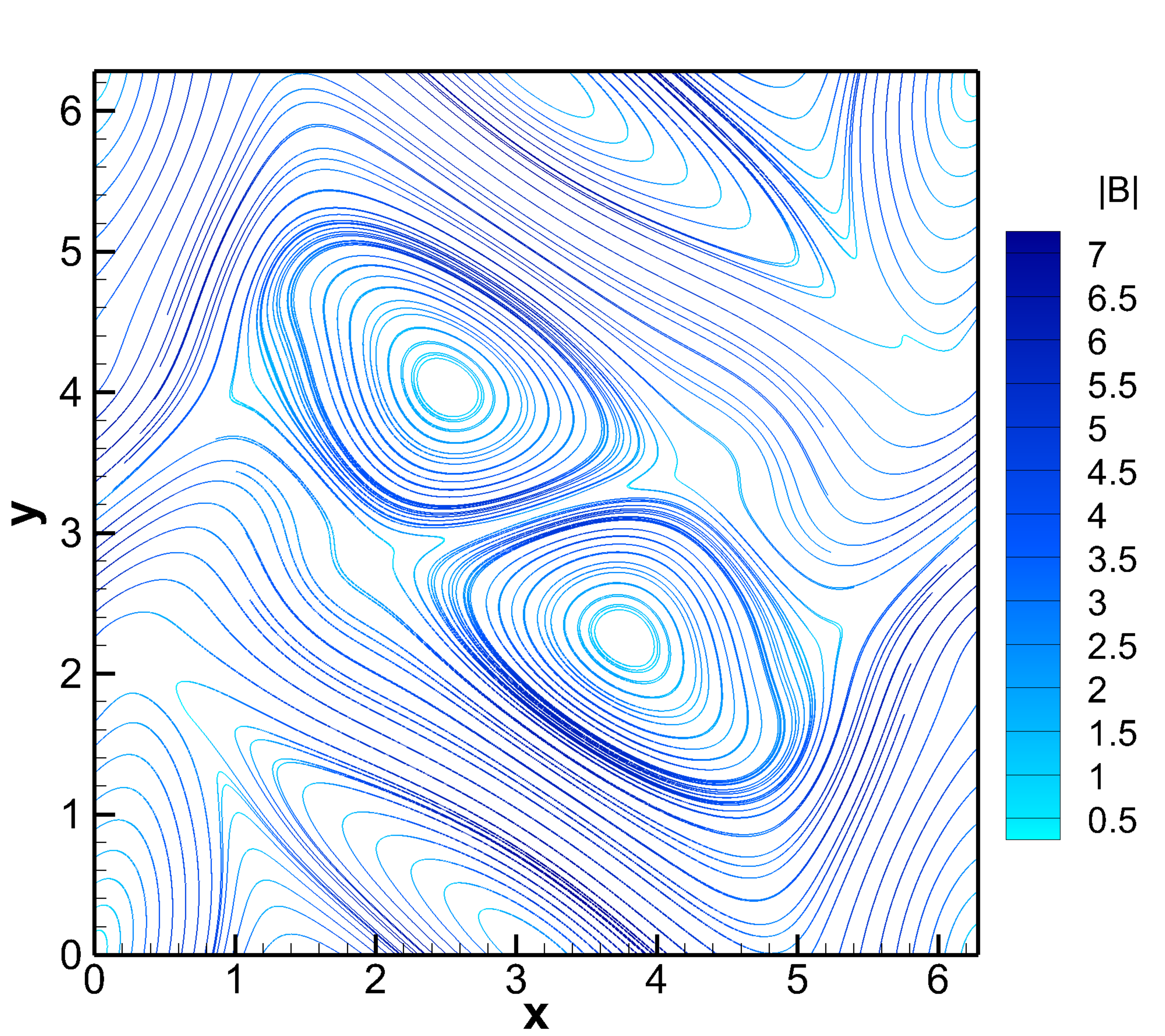}   \\
			 \includegraphics[width=0.45\textwidth]{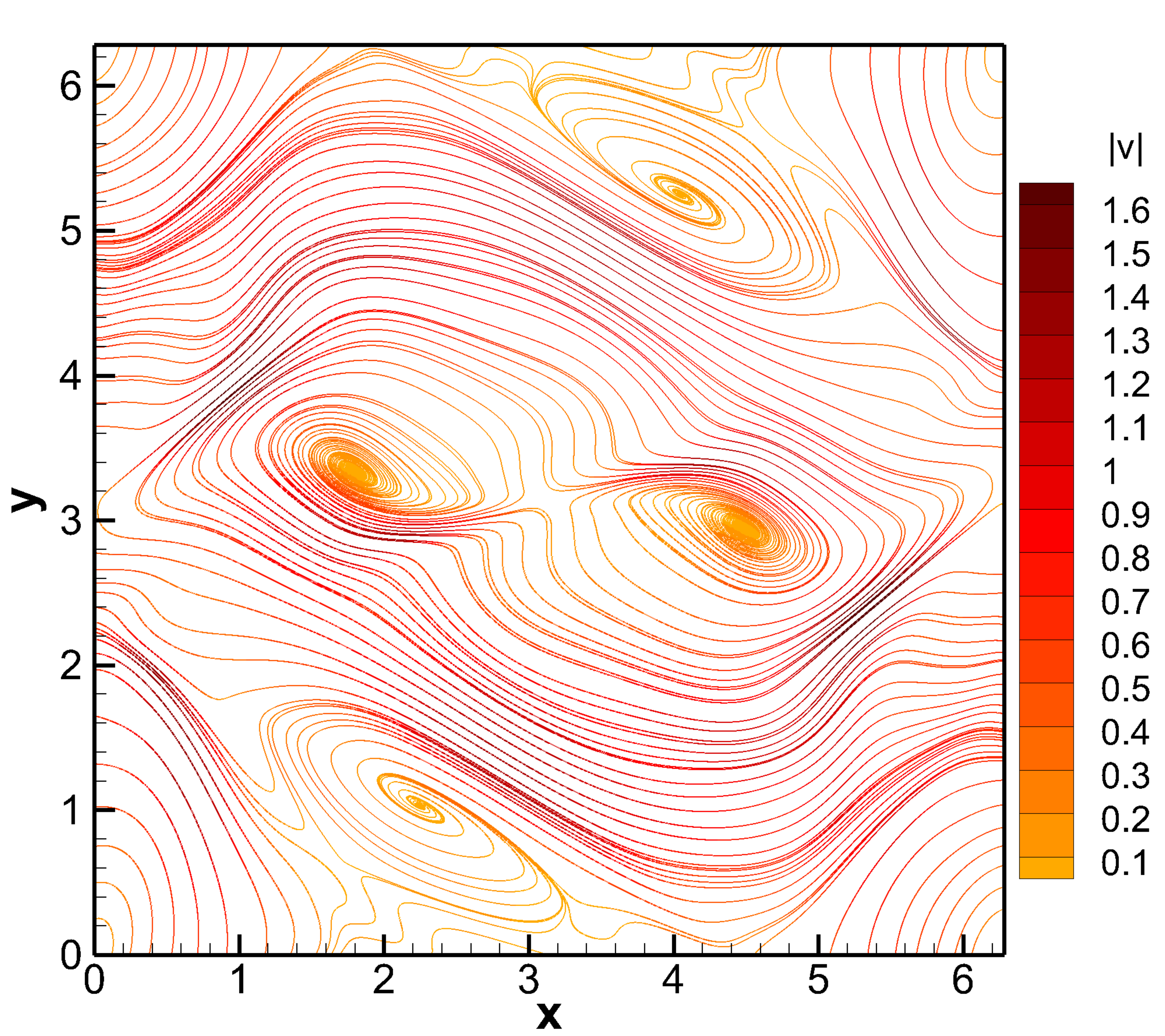} & 
			\includegraphics[width=0.45\textwidth]{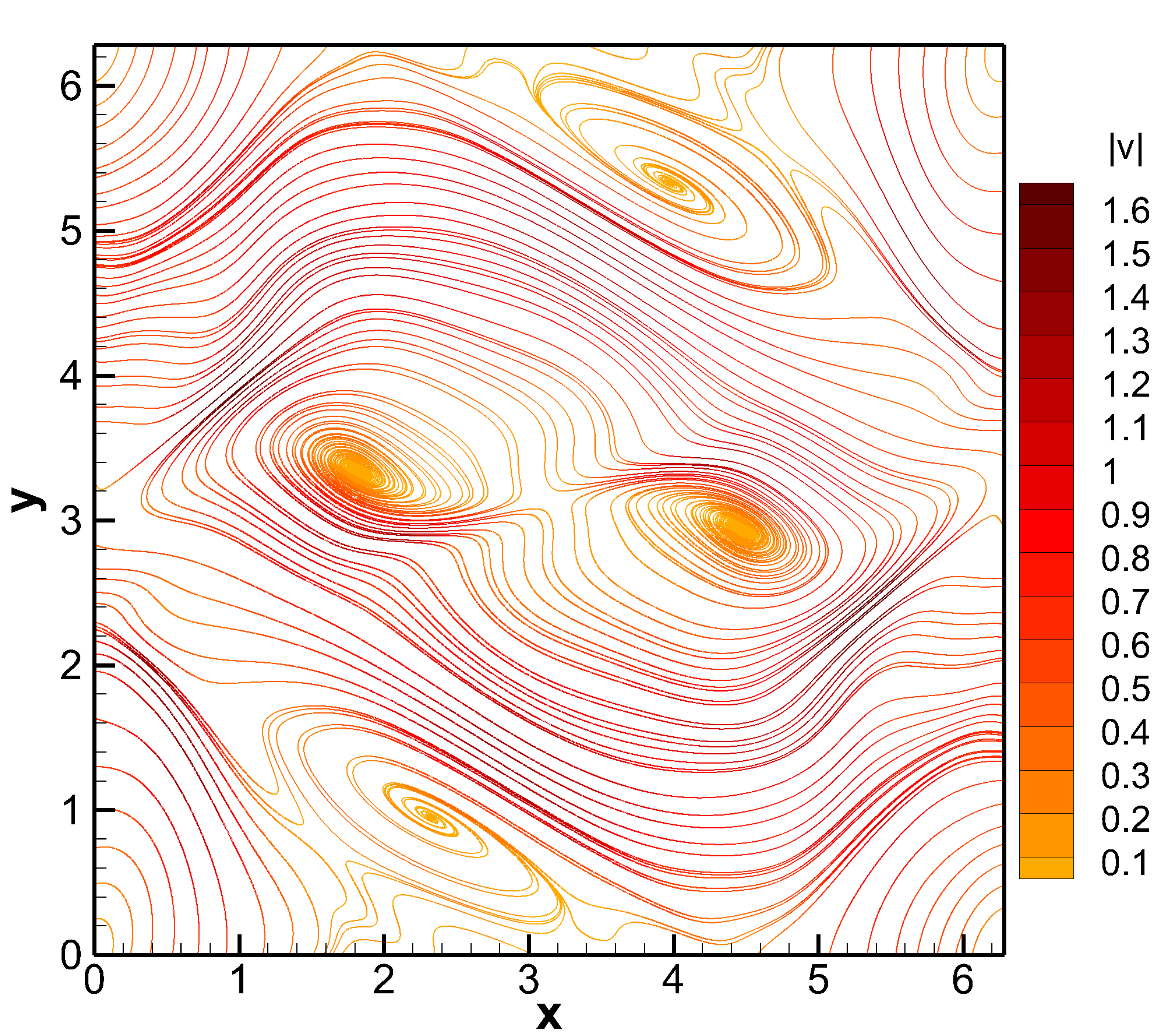}   
	\end{tabular}
\caption{
Numerical results for the viscous and resistive Orszag-Tang vortex system at time $t=2$ obtained with a $500\times 500$ grid with our new semi-implicit structure-preserving algorithm (right figures), compared to a reference solution obtained in \cite{ADERVRMHD} at the aid of a very high order accurate $P_NP_M$ scheme.
At the first row, the stream-lines of the magnetic field $\B$ are plotted, colored by the amplitude $|\B|$; at the bottom, one has the stream-lines of the velocity field $\v$, colored by its amplitude $|\v|$.
 } \label{fig:VROT}
\end{figure}

The physical and numerical domain is the same as before, a $500\times 500$ grid built on a two-dimensional periodic-box $(x,y)\in[0,2\pi]^2$, and a time-step evaluated with $\CFL=0.9$.
\myequation{l}
{
\rho = 1, \\
\v  =  \sqrt{4 \pi} \left( - \sin\left(y\right), \sin \left(x \right), 0 \right),\\
p =\frac{15}{4} + \frac{1}{4} \cos(4x) + \frac{4}{5} \cos(2x) \cos(y) - \cos(x) \cos(y) + \frac{1}{4} \cos(2y),\\
\B =  \sqrt{4\pi}\left( -  \sin\left(y\right), \sin \left(2x \right), 0\right).  
}
{\label{eq:VROrszagTang_ic}
}
Fig.  \ref{fig:VROT} shows the streamlines of the magnetic field and the velocity vector of the computed solution at time $t=2$ obtained with our structure preserving sonic-Alfvénic implicit solver, coupled with a MUSCL scheme based second-order TVD reconstruction for the explicit terms. Since the viscous terms are still integrated within a explicit conservative scheme, for this test there is essentially no special reason for using our Alfvénic implicit solver. In any case, this test is also widely used to check if all the terms in the PDE system are correctly integrated in space and time and many numerical reference solution are available in literature. In Fig.  \ref{fig:VROT}, a numerical reference solution computed with a very high order accurate $P_NP_M$ scheme, published by \cite{ADERVRMHD}, confirms a good agreement with our results. For this simulation, the $\theta$-parameters of the implicit solver are set to $\theta_{\text{B}} = \theta_{\text{p}} = 0.6$. 

\paragraph{MHD blast wave problem.}

A standard two-dimensional explosion problem for the ideal MHD equations has been introduced by \cite{BalsaraSpicer1999} and it is commonly used as robustness-check of many MHD solvers.

\begin{figure}[!t]
\centering 
\begin{tabular}{lcr}
			 \includegraphics[width=0.5\textwidth]{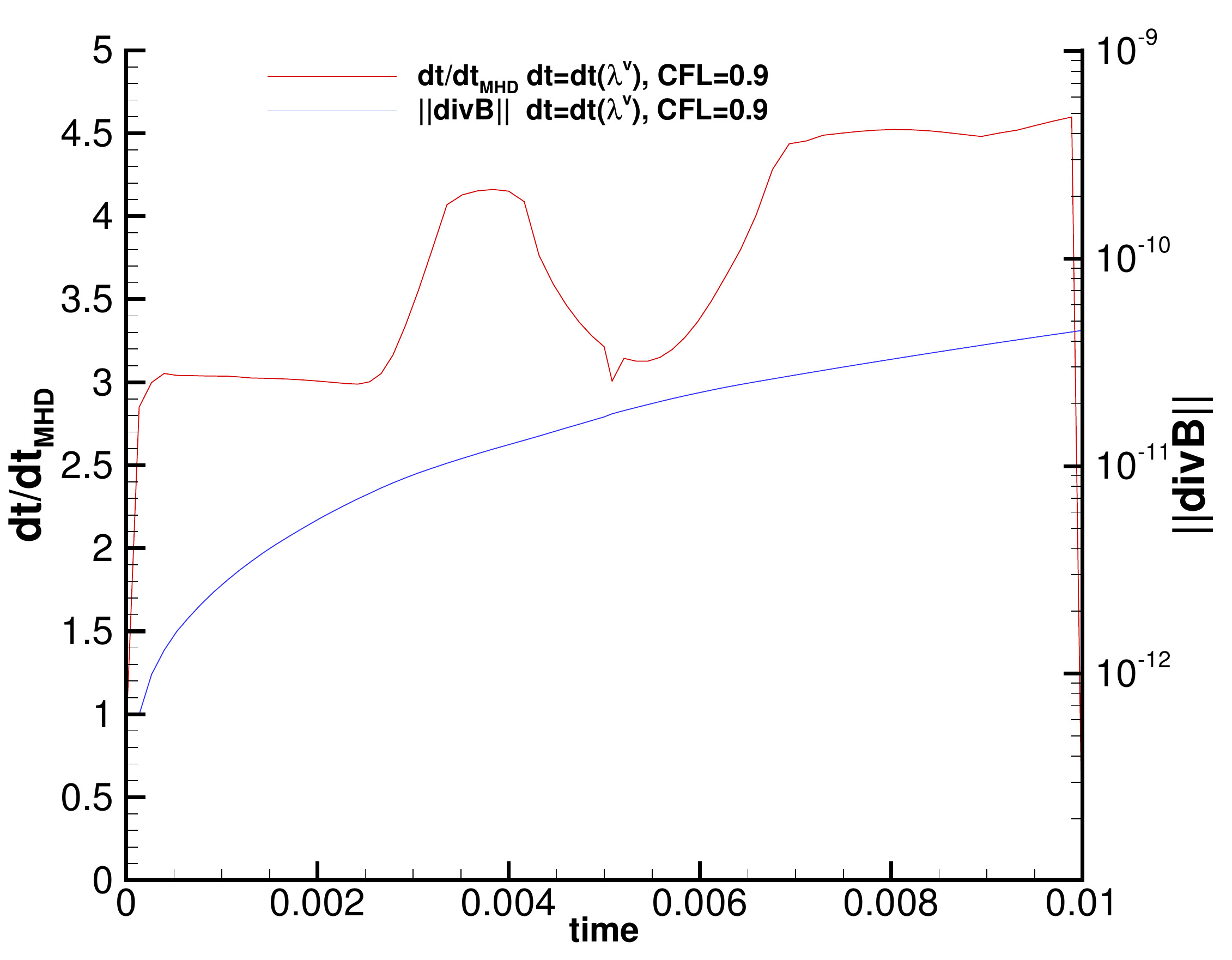}  
	\end{tabular}
\caption{Numerical results for the ideal Blast-Wave test at time $t=0.01$ obtained with a mesh resolution $\Delta x=\Delta y=10^{-3}$ with our new semi-implicit structure-preserving algorithm. The time evolution of the effective computational Courant number, i.e. $\Delta t/ \Delta t(\lambda)$, is plotted, together with the $L_2$ norm of the divergence of the magnetic field $||\div\B||$. 
 } \label{fig:BW_data}
\end{figure}

In this problem, a background constant magnetic field is set to $\B=(B_x, 0 ,0)$ with, in this work, $B_x=10^2$ while a strong discontinuity in  pressure is chosen in the form
\myequation{l}{
  p  = \left\{ 10^3 \;\; \text{if} \  r<R \,; \;\; 
  10^{-1} \;\;  \text{otherwise}\,,
\right.
	}{}
to trigger the explosion, where $r=\sqrt{x^2+y^2}$ is the radial position and $R=0.1$ is the initial radius of the high-pressure plasma.
. In this work, the discontinuity has been solved '\emph{as is}' without using any kind of linear smoother in space.
The computational domain is again a two-dimensional periodic box $(x,y)\in[-0.55,0.55]^2$ that is discretized with a mesh resolution of $\Delta x= \Delta y = 10^{-3}$. The time parameters of the new  semi-implicit three-split solver are set to $\theta_{\text{B}} = \theta_{\text{p}} = 0.6$, with Courant parameter $\CFL=0.9$.
With the current set up of the three-splitting scheme the numerical solutions, plotted in Fig.  \ref{fig:VROT} for the time slice $t=0.01$, is shown to be comparable to other published results in literature, see e.g. \cite{BalsaraSpicer1999}. It is also interesting also to look at the gain in terms of Courant number that is achieved tanks to the chosen implicit time discretization, see Fig.  \ref{fig:BW_data}, plotted together with the $L_2$ norm of the magnetic field $||\div\B||$.

\begin{figure}[!t]
\centering 
\begin{tabular}{lr}
			 \includegraphics[width=0.45\textwidth]{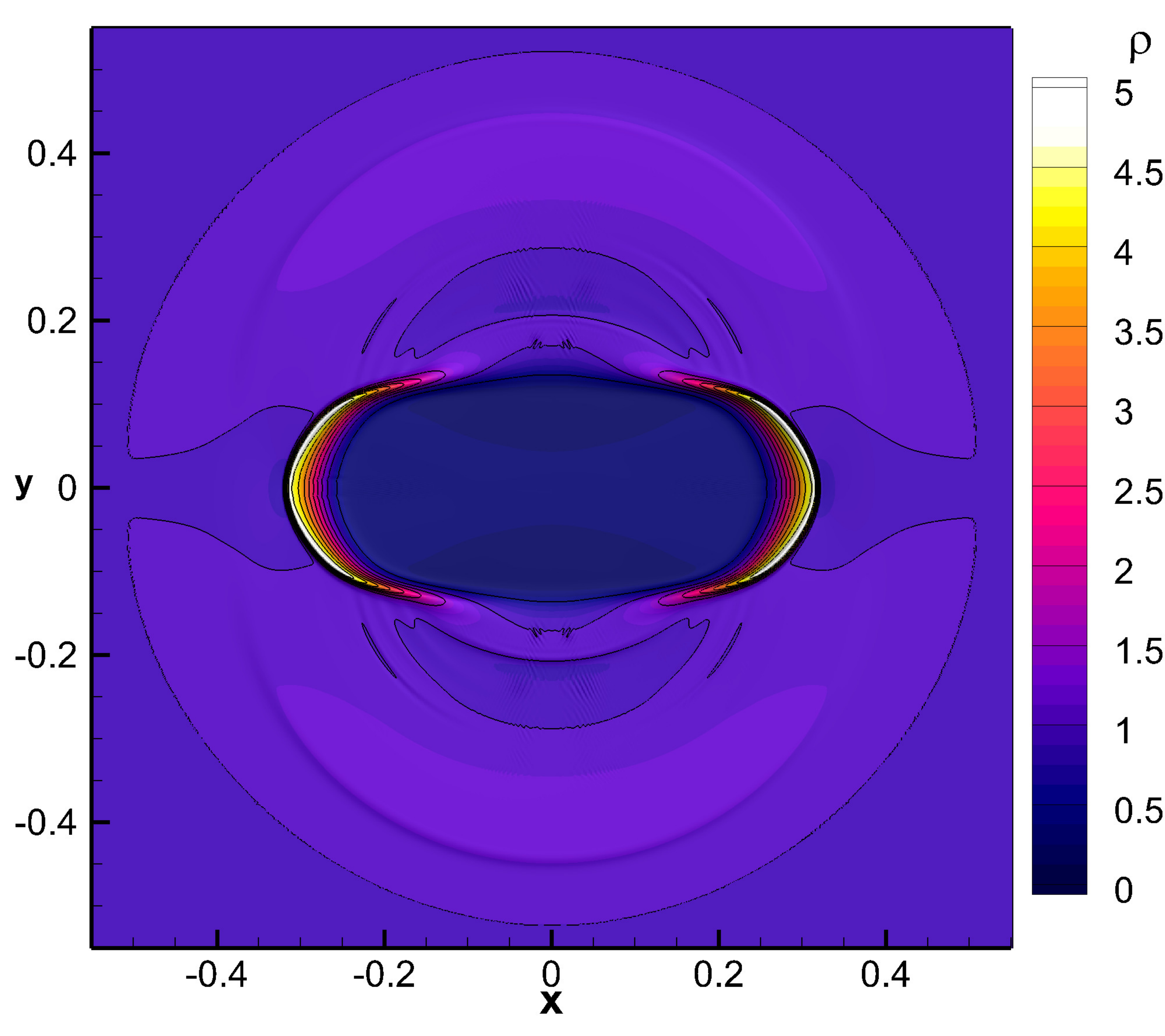} & 
			\includegraphics[width=0.45\textwidth]{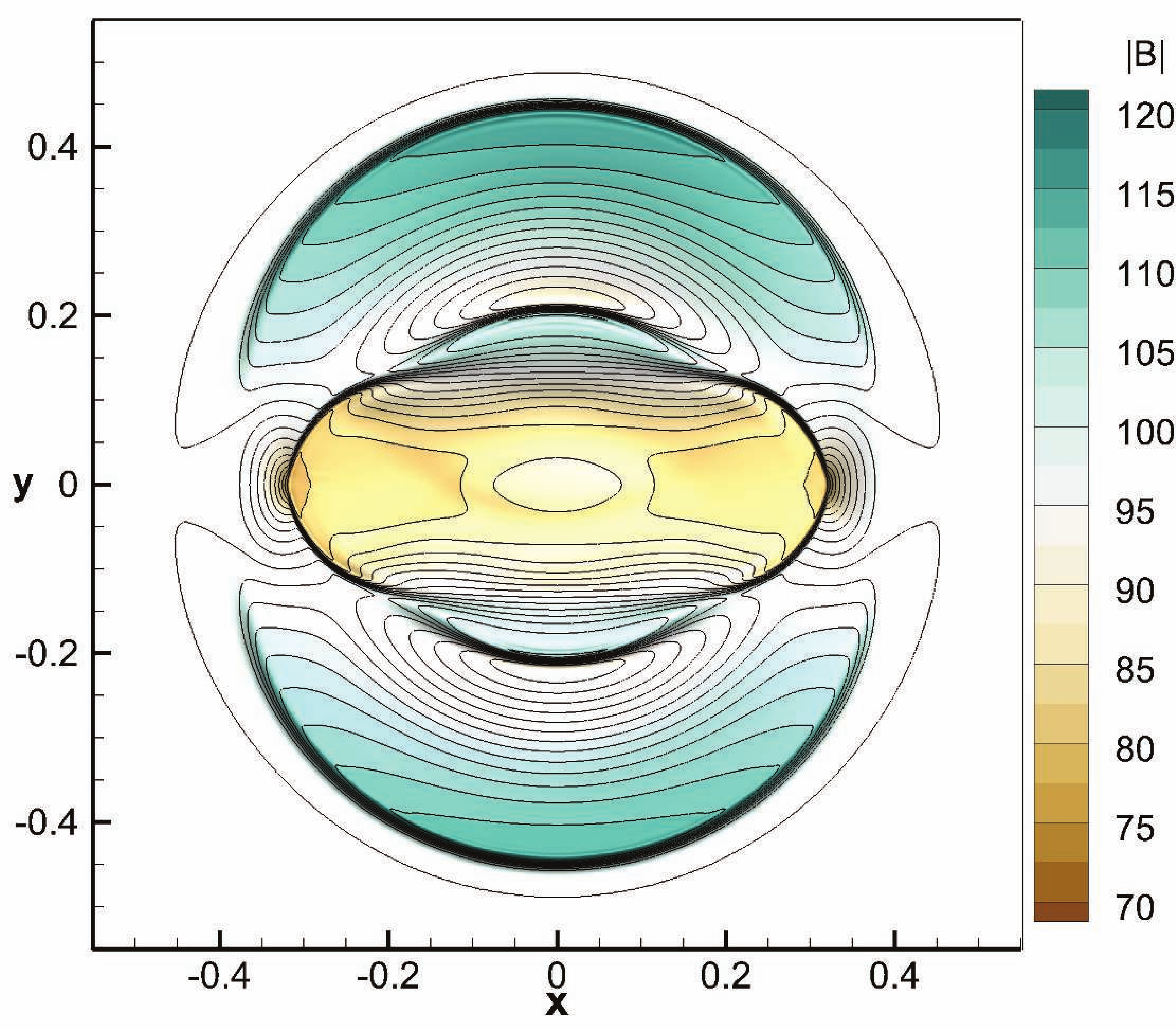}   
	\end{tabular}
\caption{Numerical results for the ideal Blast-Wave test at time $t=0.01$ obtained with a mesh resolution $\Delta x=\Delta y=10^{-3}$ with our new semi-implicit structure-preserving algorithm.
At the left, the density is plotted, while the amplitude $|\B|$ is at the right.
 } \label{fig:BW_01}
\end{figure}

\paragraph{MHD rotor problem.}

Another test that is widely used to check the robustness of an MHD solver, especially regarding the correct treatment of the divergence constraint $\div\B=0$, is the so called \emph{rotor}-test, introduced again by \cite{BalsaraSpicer1999}. Being configured for the ideal MHD equations, we set to zero $\mu=\eta=0$.
In this test, a radial discontinuity in the fluid density and velocity is configured in the form of a \emph{rotor} as
\myequation{l}{
  (\rho, \v )= \left\{ \begin{array}{ll}
						\left(10, \omega \times \mathbf{r} \right) & \text{if} \  r<R \,; \\
						\left(1 , (0,0,0) \right)  &  \text{otherwise}\,, \end{array}
\right.	}{}
	where $\mathbf{r}=(x,y)$ is the spatial vector position, and $r$ its module,
with a  background constant pressure and magnetic field $p=1$ and $\B=(B_x, 0,0)$ with $B_x=2.5$ and a constant angular velocity $\omega = 10 \hat{\mathbf{z}}$. 
In order to compare the numerical results, we follow \cite{Zanotti2015c} in introducing a narrow \emph{linear} smoother function to taper the discontinuity in the range $r\in[R,R+5\%R]$, so that we can also compare the final results.

\begin{figure} 
\centering 
\begin{tabular}{lr}
			\includegraphics[width=0.45\textwidth]{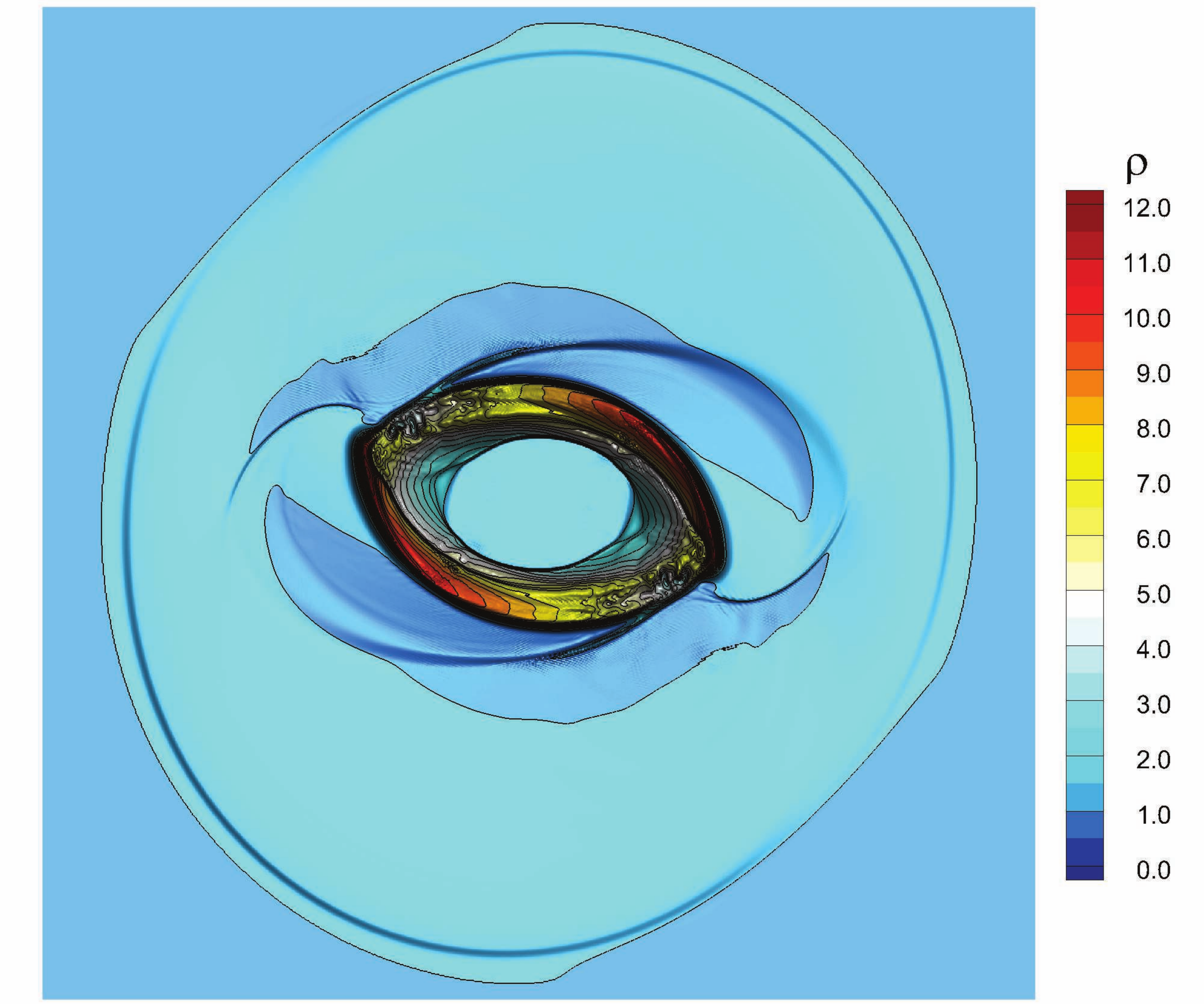} & 
			\includegraphics[width=0.45\textwidth]{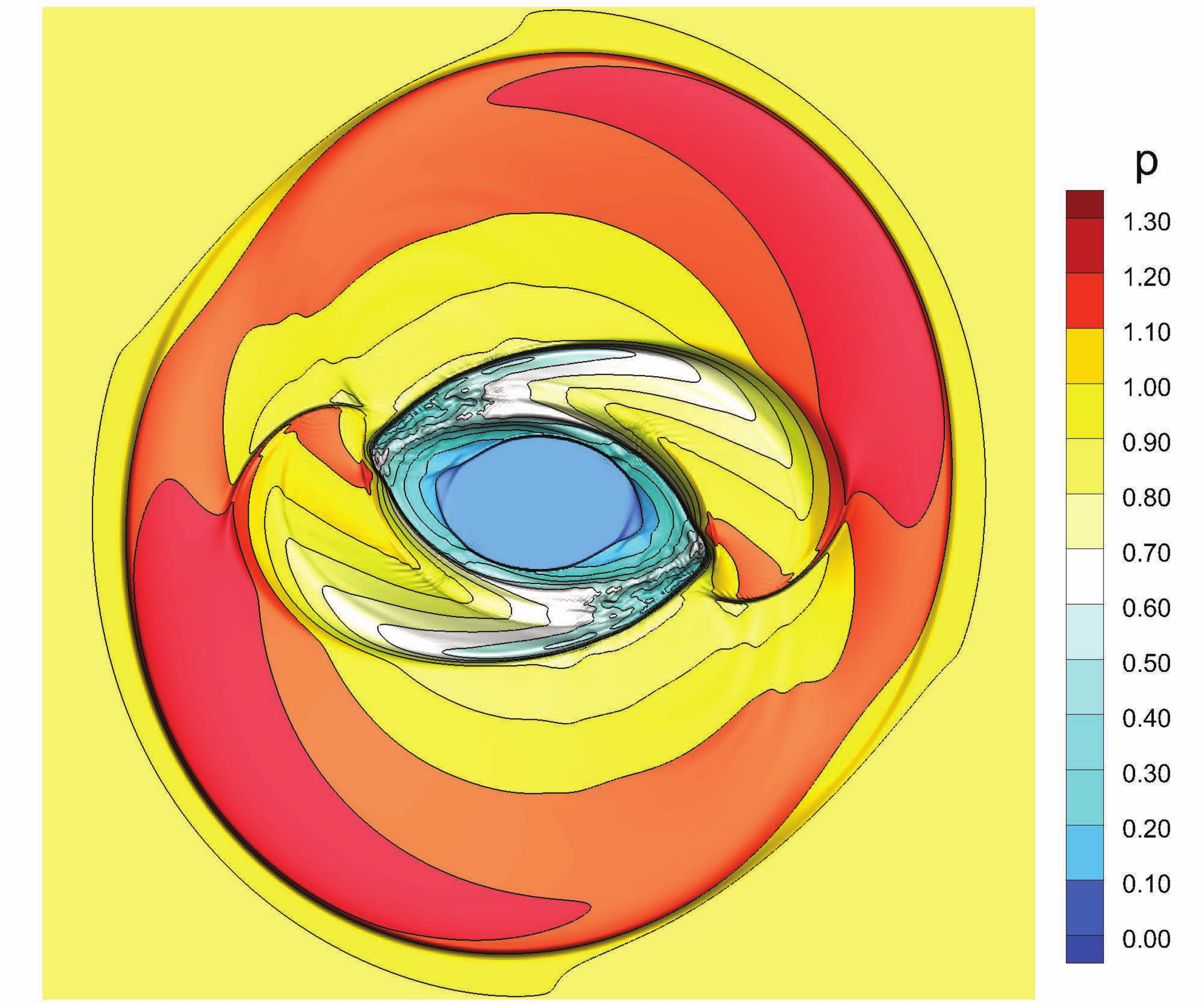}
			\\                                                          
			\includegraphics[width=0.45\textwidth]{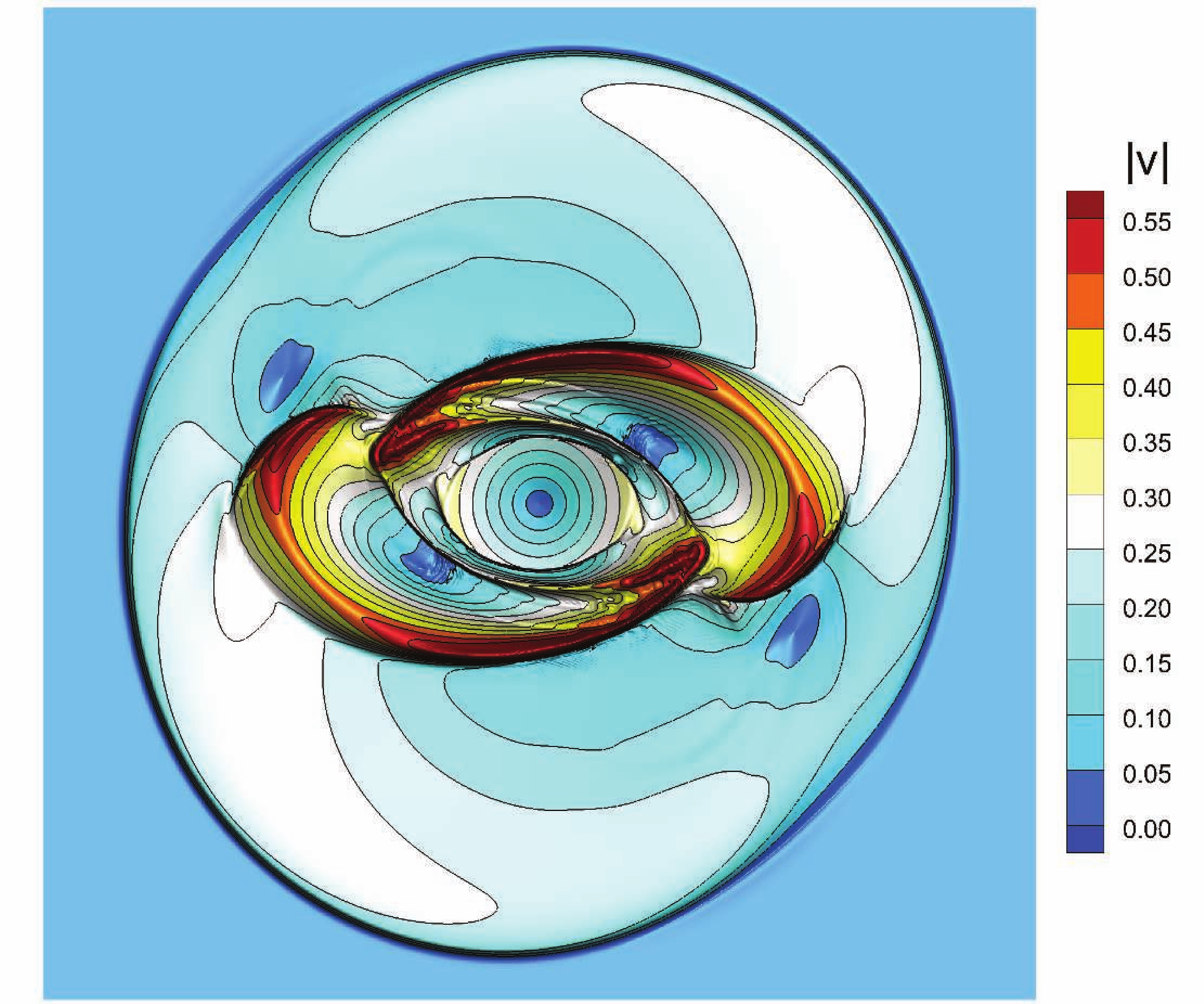} & 
			\includegraphics[width=0.45\textwidth]{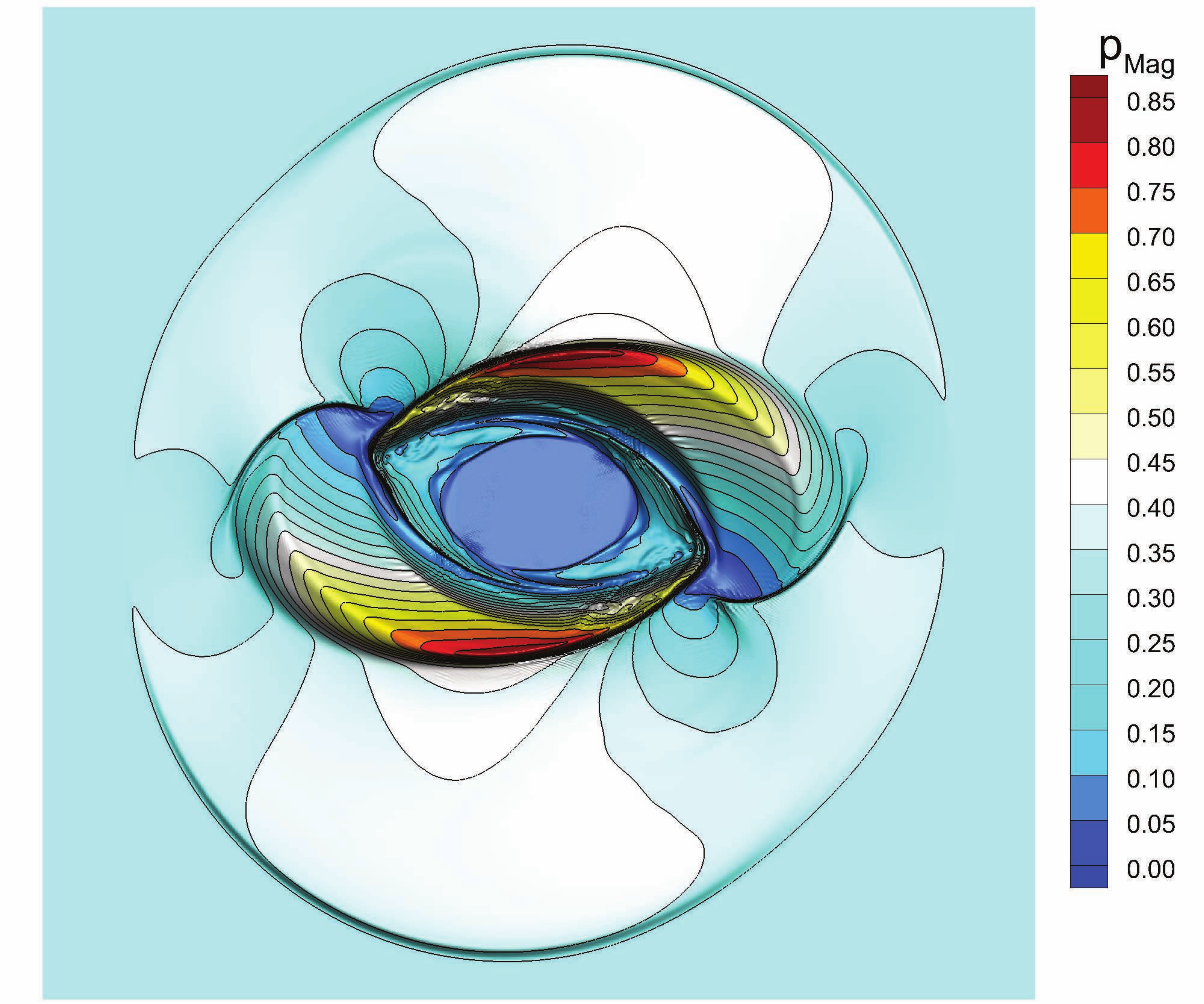} 
	\end{tabular}
\caption{Numerical results for the MHD rotor test at time $t=0.25$ obtained with a $500\times 500$ grid with our new semi-implicit structure-preserving algorithm, with $\CFL=1/4$: iso-countour lines for matter-density (top-left), pressure (top-right), absolute value of the velocity (bottom-left) and magnetic pressure (bottom-right) are plotted .} \label{fig:Rot}
\end{figure}

The numerical solution plotted in Fig.  \ref{fig:Rot} is obtained after setting the selected  MUSCL scheme based second-order TVD reconstructions for the nonlinear convective terms, and the implicitness parameters to $\theta_{\text{B}}=\frac{1}{2}$, $\theta_{\text{p}}=1$ for our new semi-implicit three-split solver, i.e. the method is formally second order accurate for the Alfv\'en waves. 

\begin{figure} 
\centering 
\begin{tabular}{cr}
\multicolumn{2}{c}{\includegraphics[width=0.5\textwidth]{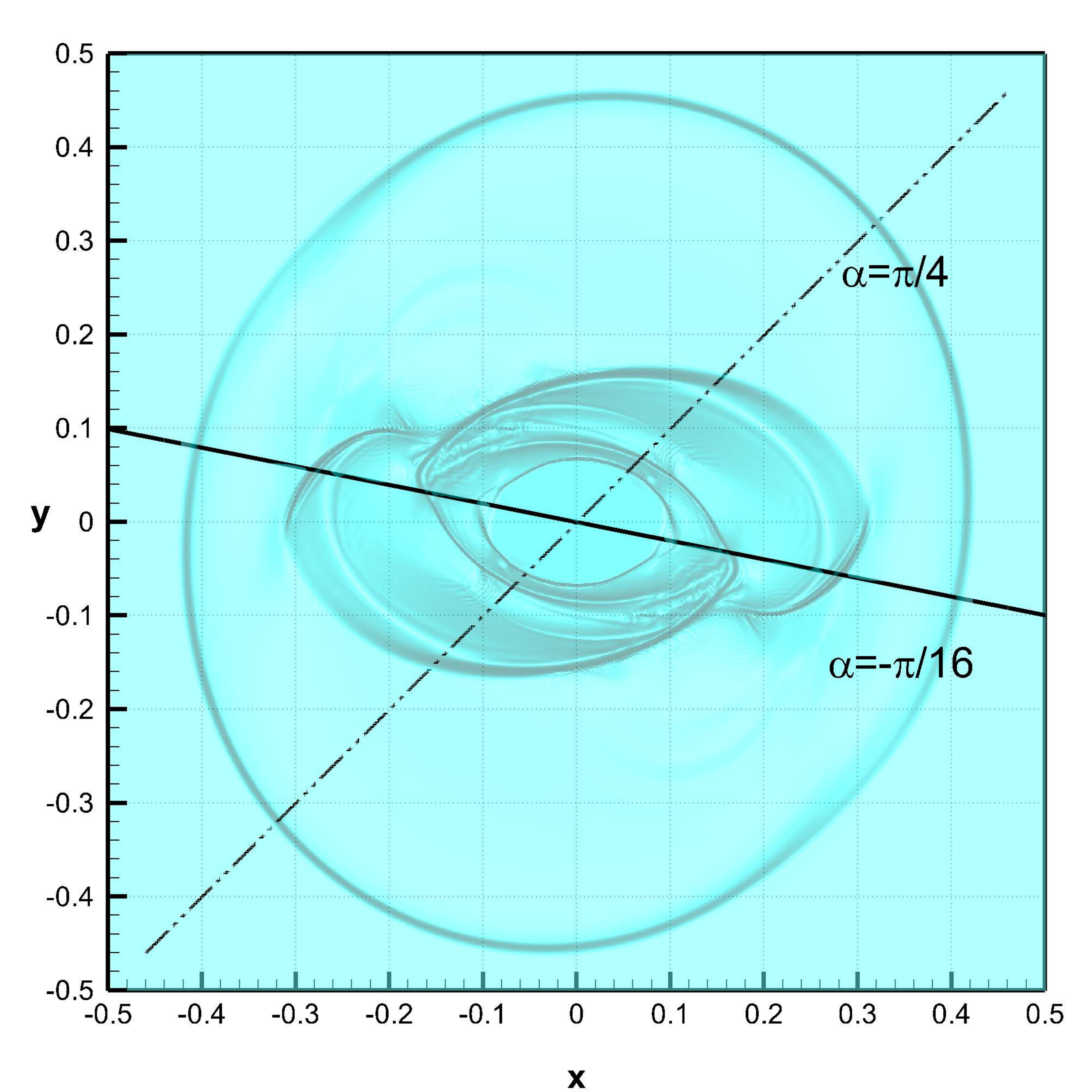} }\\
			\includegraphics[width=0.4\textwidth]{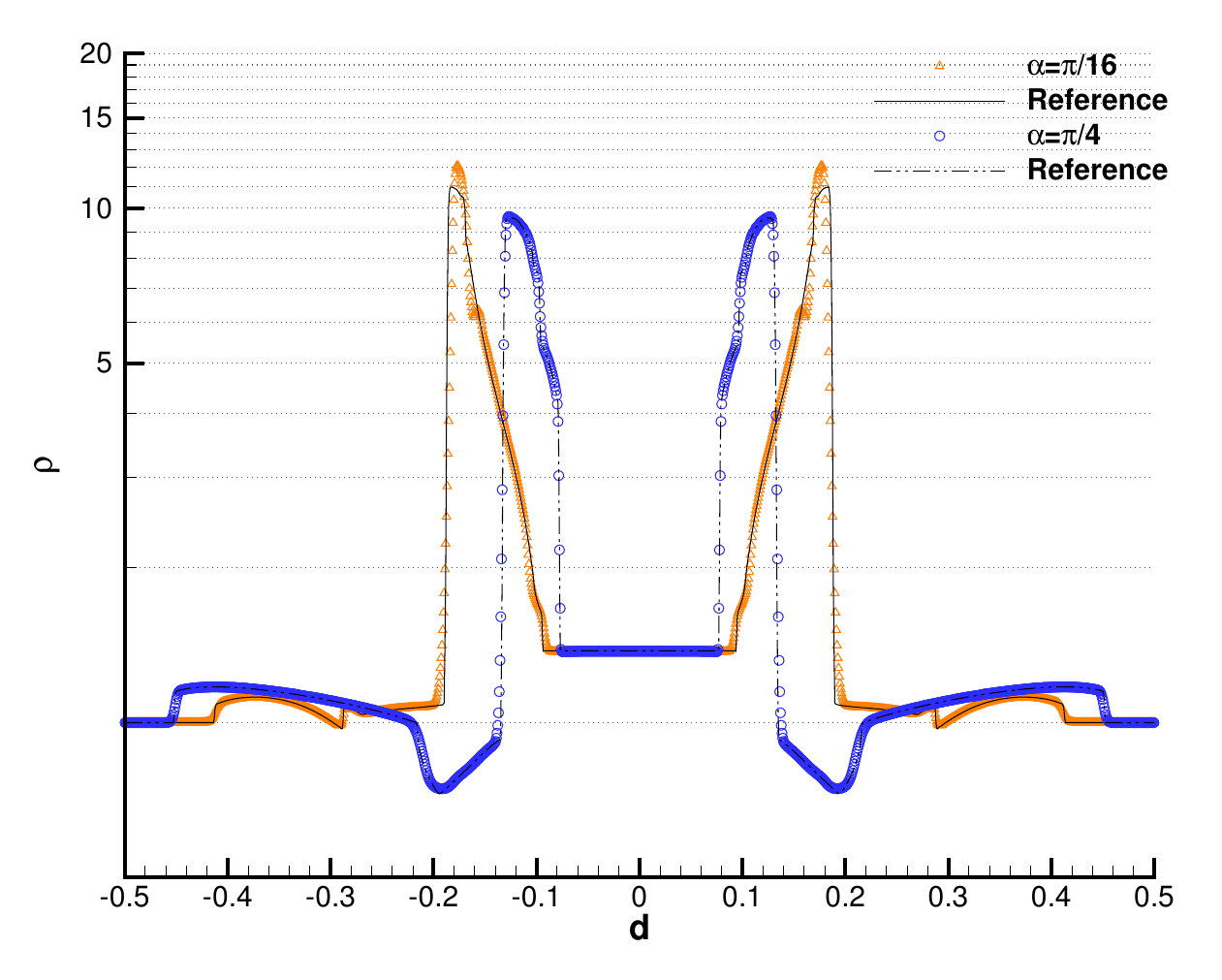} & 
			\includegraphics[width=0.4\textwidth]{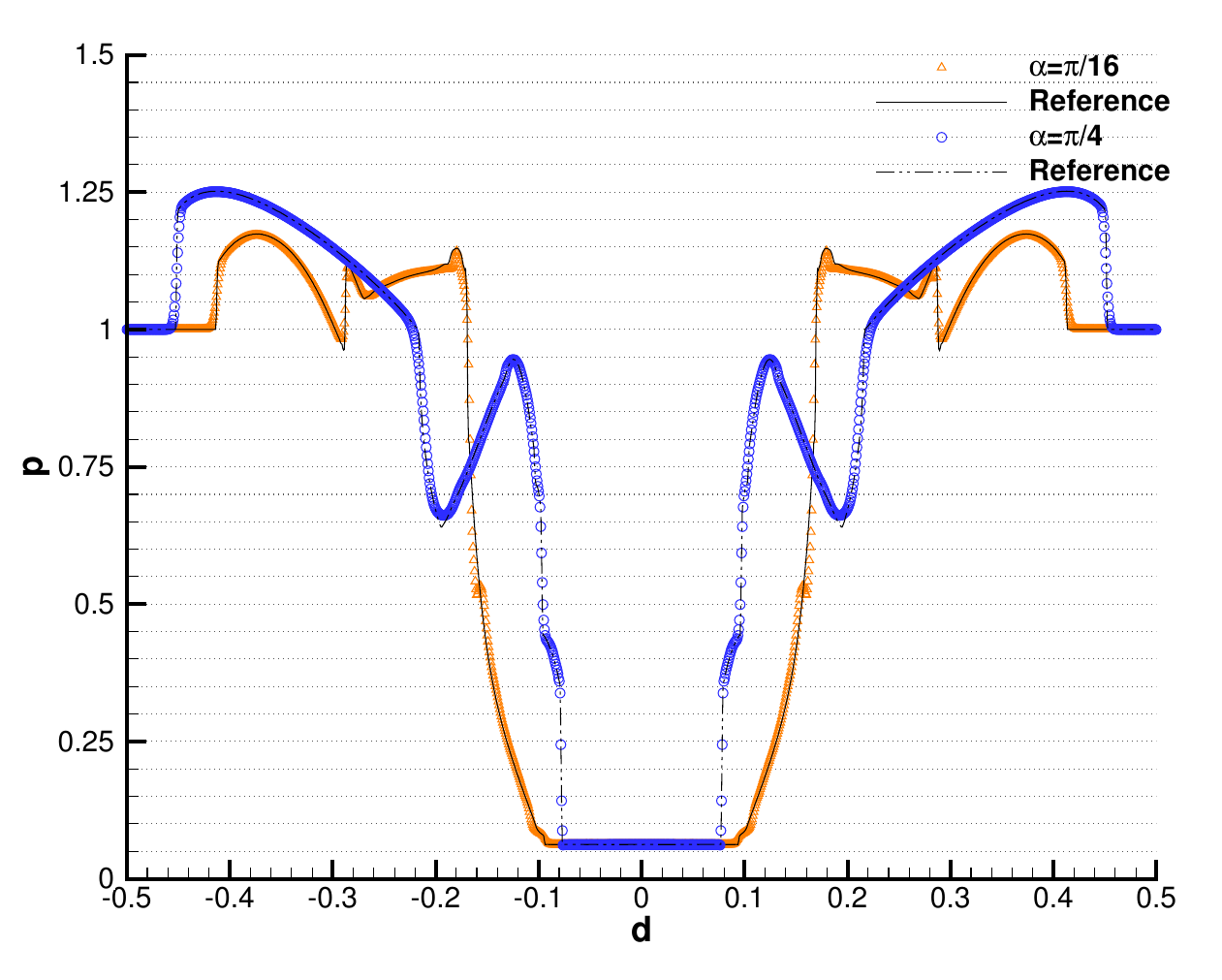}
			\\
			\includegraphics[width=0.4\textwidth]{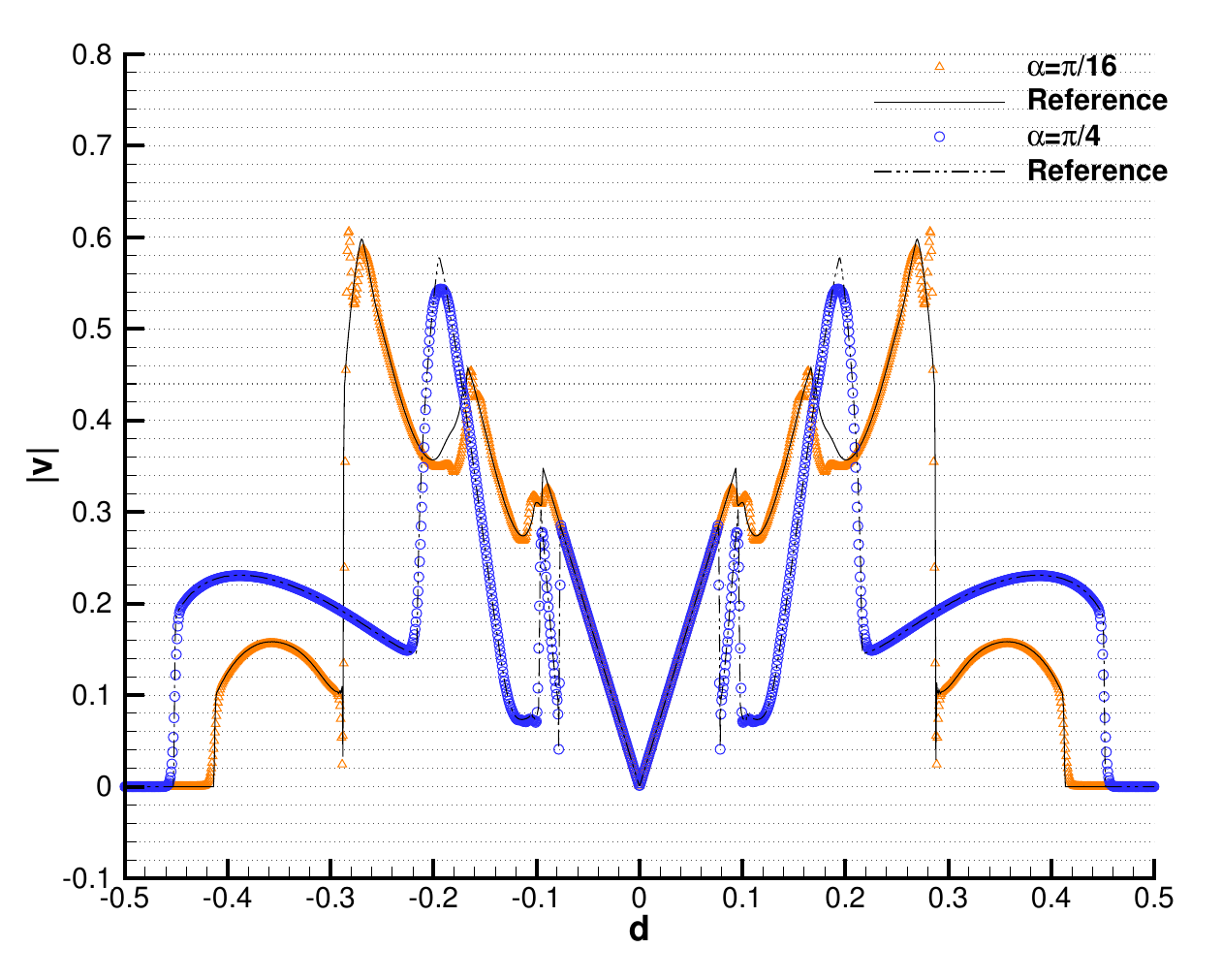} & 
			\includegraphics[width=0.4\textwidth]{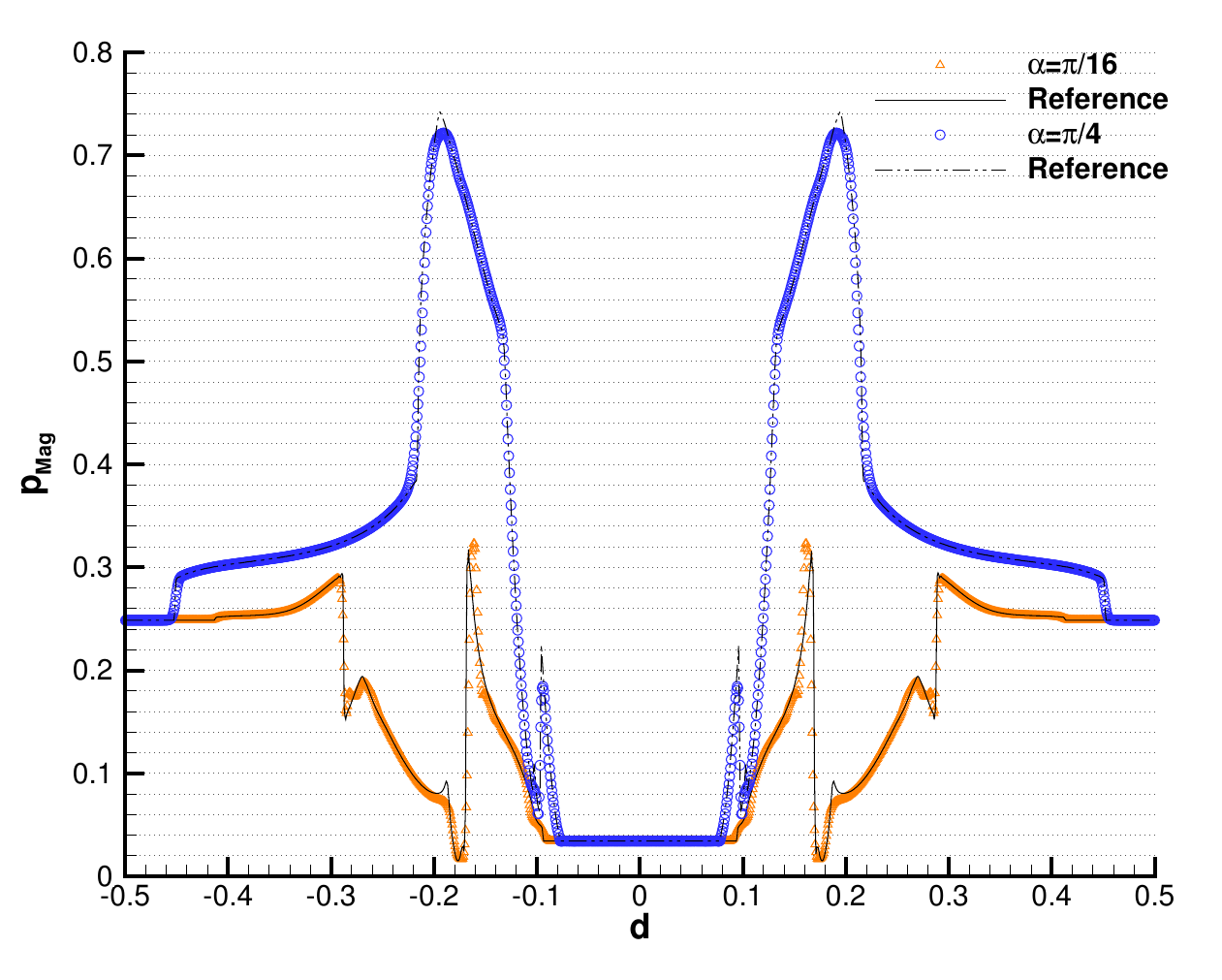} 
	\end{tabular}
\caption{Numerical results for the MHD rotor test at time $t=0.25$ obtained with a $500\times 500$ grid with our new semi-implicit structure-preserving algorithm, with $\CFL=1/4$: the one-dimensional cut along the lines $y/x=tg(\alpha)$, $\alpha=\pi/4,-\pi/16$, for some chosen physical quantities are plotted. At the top, the two-dimensional physical domain, together with the graphical representation of the interpolating lines are shown. Then, the numerical solution along the one-dimensional cut-lines are plotted next to a reference solution for matter-density (top-left), pressure (top-right), absolute value of the velocity (bottom-left) and magnetic pressure (bottom-right). As a reference solution, we choose a high-order ADER-DG-$\mathbb{P}_5$ scheme with an a-posteriori sub-cell limiter, see \cite{Zanotti2015c}.} \label{fig:Rot_1d}
\end{figure}

\begin{figure} 
\centering 
\begin{tabular}{cc}
			\includegraphics[width=0.4\textwidth]{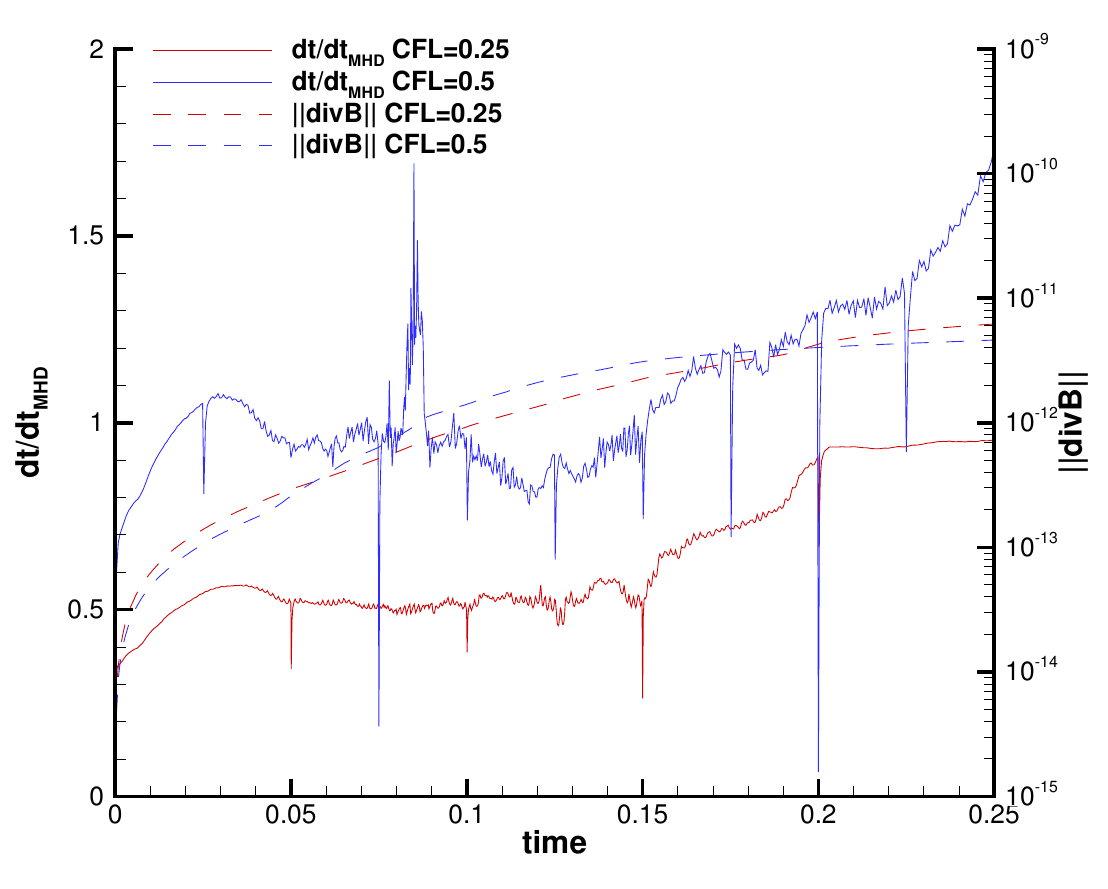} & 
			\includegraphics[width=0.4\textwidth]{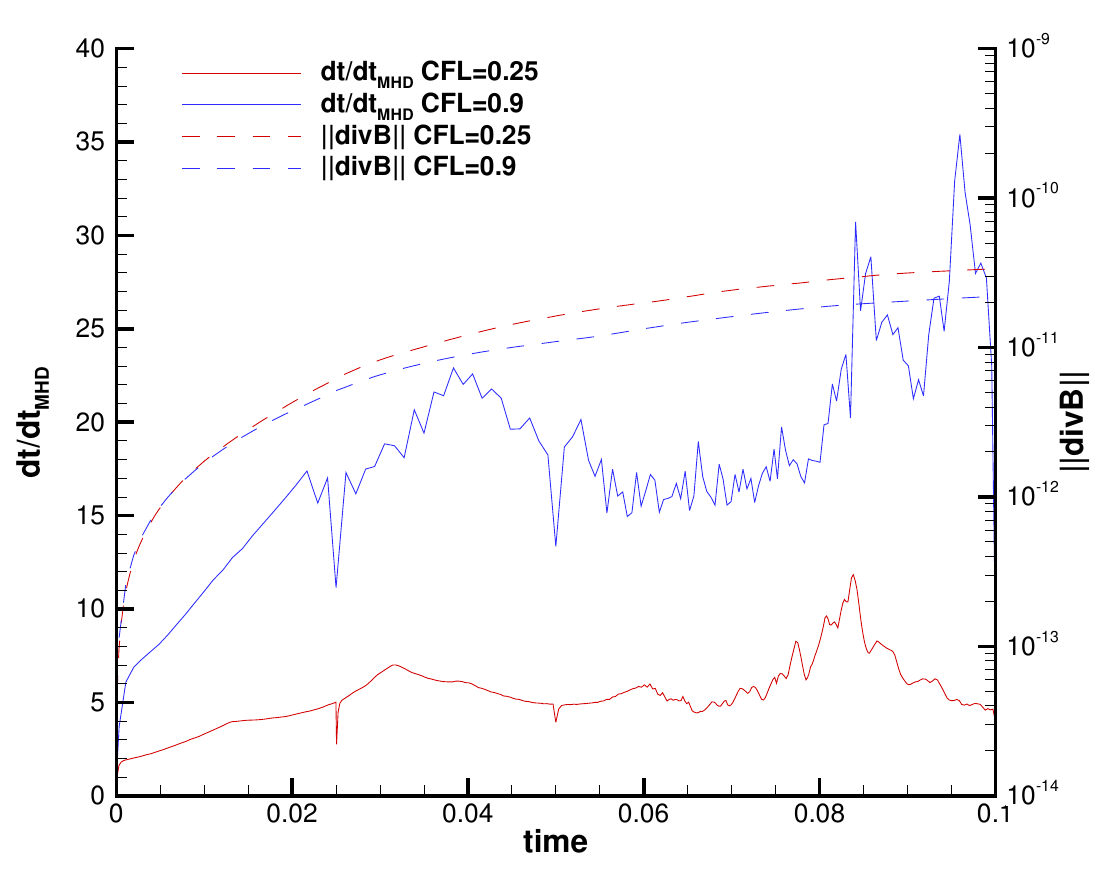} 
	\end{tabular}
\caption{Numerical results for the MHD rotor test with $(p_0,B_x)=(1,2.5)$ and $(p_0,B_x)=(1,25)$ 
 obtained with a $\Delta x=\Delta y=10^{-3}$ mesh resolution with our new semi-implicit structure-preserving algorithm, at different $\CFL$ parameters. Here it is plotted the time evolution of the effective Courant number, i.e. the ratio $dt/dt_{\text{max}}$, together with the L$_2$ norm of the divergence $\div \B$.  
Note that the eye-catching lower spikes in the plots correspond to the automatic time step adjustment used to adapt the printing process to the prescribed output time.} \label{fig:RotCFL}
\end{figure}

In Fig.  \ref{fig:Rot} the two-dimensional plots of some physical variables are shown for the time-slice $t=0.25$, and one can eventually verify that the computed results are well compatible with other published results, see \cite{Dumbser2008,balsarahlle2d,AMR3DCL,ADERdivB,HPRmodelMHD}. In order to give a quantitative comparison  in Fig.  \ref{fig:Rot_1d} the numerical solution is interpolated at the one-dimensional slices $\alpha=\pi/4$, $-\pi/16$ and compared with a reference solution computed at the aid of a high-order ADER-DG-$\mathbb{P}_5$ scheme with an a-posteriori sub-cell limiter and adaptive mesh-refinement, see \cite{Zanotti2015c,Fambri2020}. 
These plots shows a general good agreement between the results obtained with our novel structure-preserving semi-implicit three-split solver. Since the speed of the  Alfvén waves are of the same order of the other eigenvalues of the PDE system, one cannot take benefit in terms of the effective Courant number, see \ref{fig:RotCFL}. In any case, this was a good way to verify the robustness of the presented algorithm against a well-known test problem.

Then we also tested a lower beta and Alfvénic Mach rotor test by setting the constant magnetic field to $B_x=25$.
Fig.  \ref{fig:RotCFL} show the time evolution of the effective Courant number after choosing $\CFL=0.9$, $\frac{1}{2}$ or $\frac{1}{4}$. 

Results show that, also for this test, the presented implicit solver may become beneficial. Effective Courant numbers higher than $10$ are reached for this second test with lower Mach number with respect to the Alfvén wave speed, see Fig.  \ref{fig:RotCFL}.

\subsection{Three dimensional tests}

This last section is dedicated to the numerical validation of the proposed scheme for the full three-dimensional nonlinear MHD equations.
Here, we propose two nontrivial three-dimensional tests that are the 3D extension of the above presented blast-wave and viscous and resistive Orszag-Tang vortex.

\paragraph{3D MHD blast wave.}

\begin{figure}[!htbp]
\centering 
\begin{tabular}{lr}
			 \includegraphics[width=0.45\textwidth]{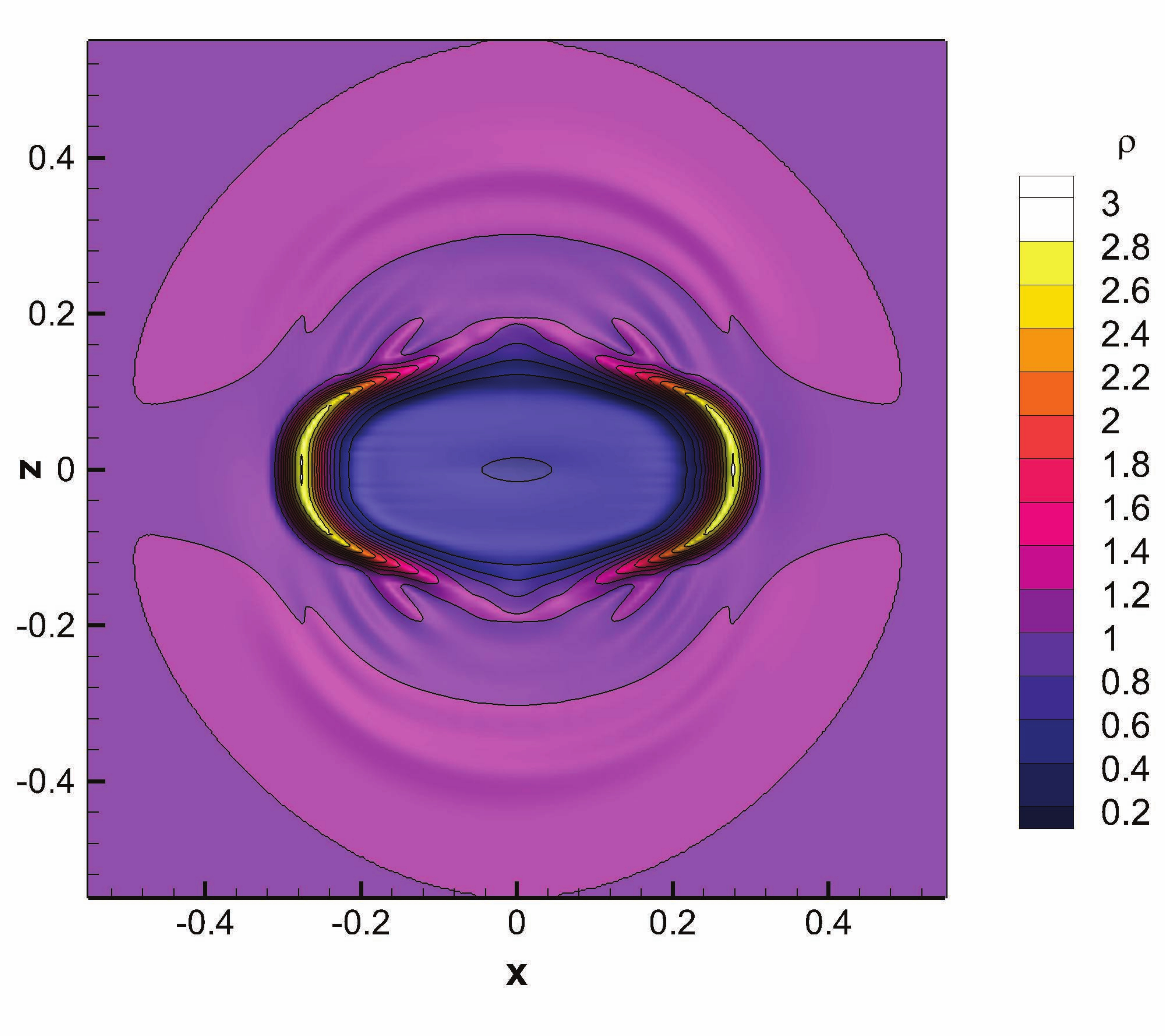}  &
			 \includegraphics[width=0.45\textwidth]{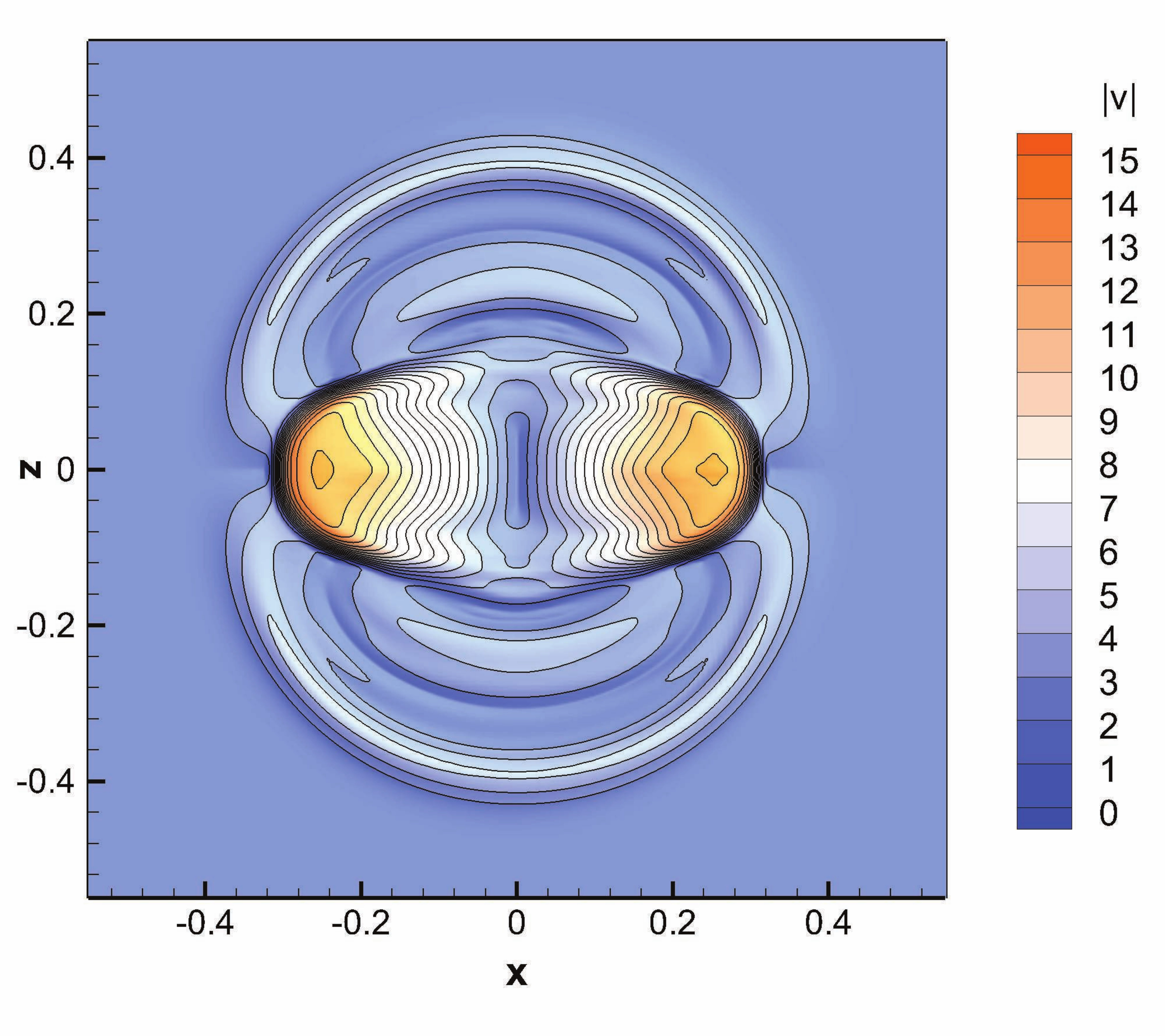}  \\
			 \includegraphics[width=0.45\textwidth]{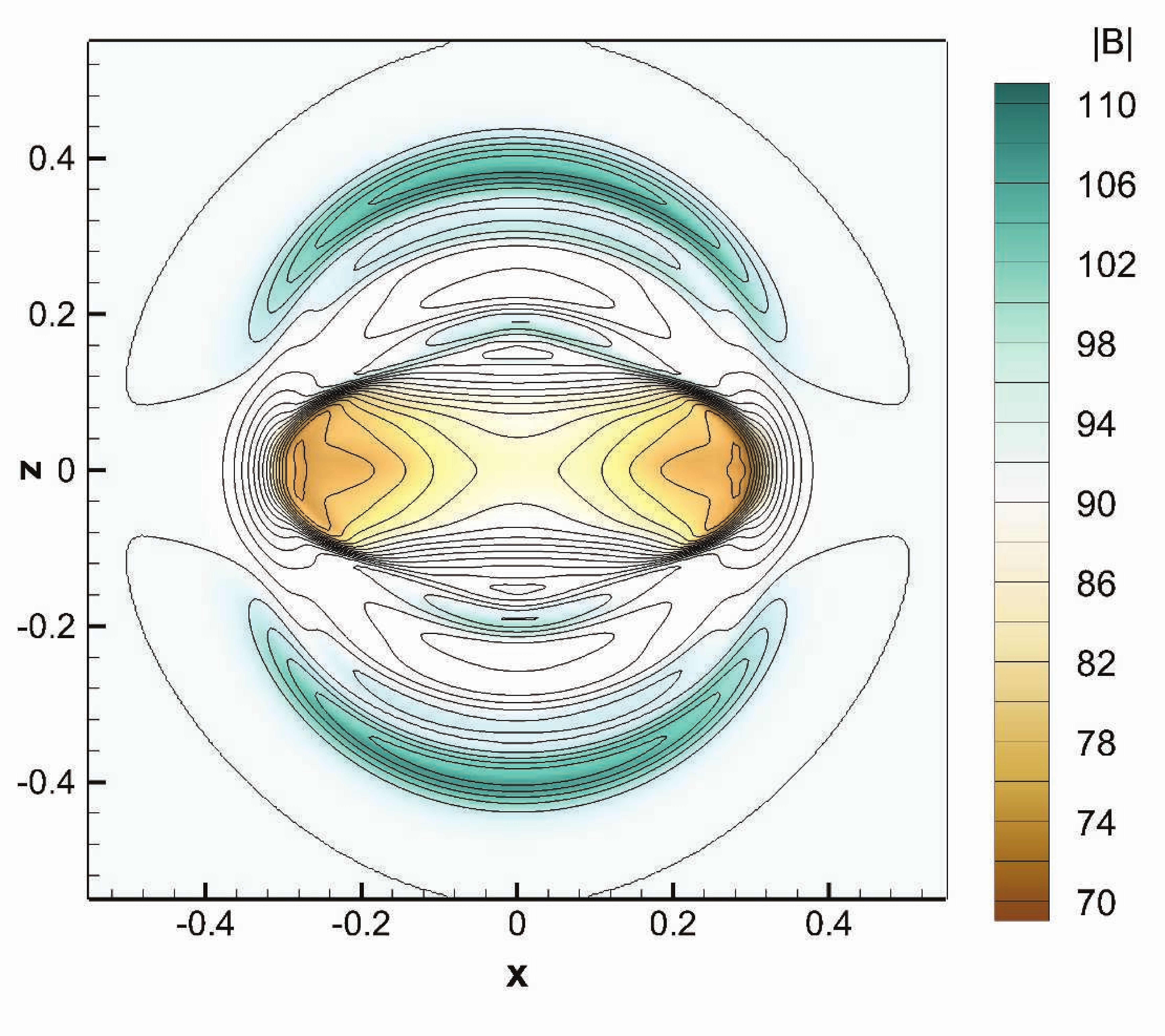}  &
			 \includegraphics[width=0.45\textwidth]{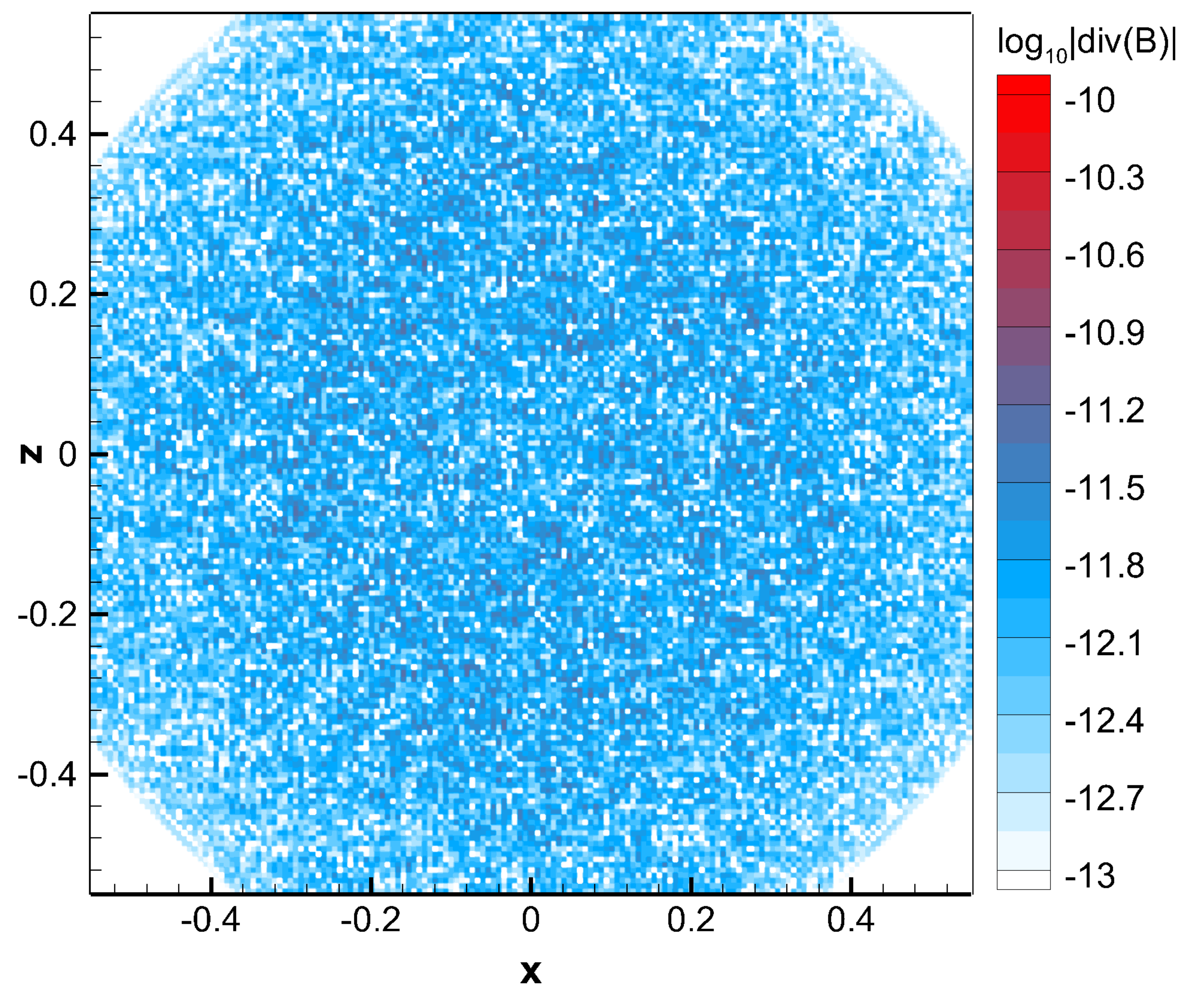}  \\
			 \includegraphics[width=0.45\textwidth]{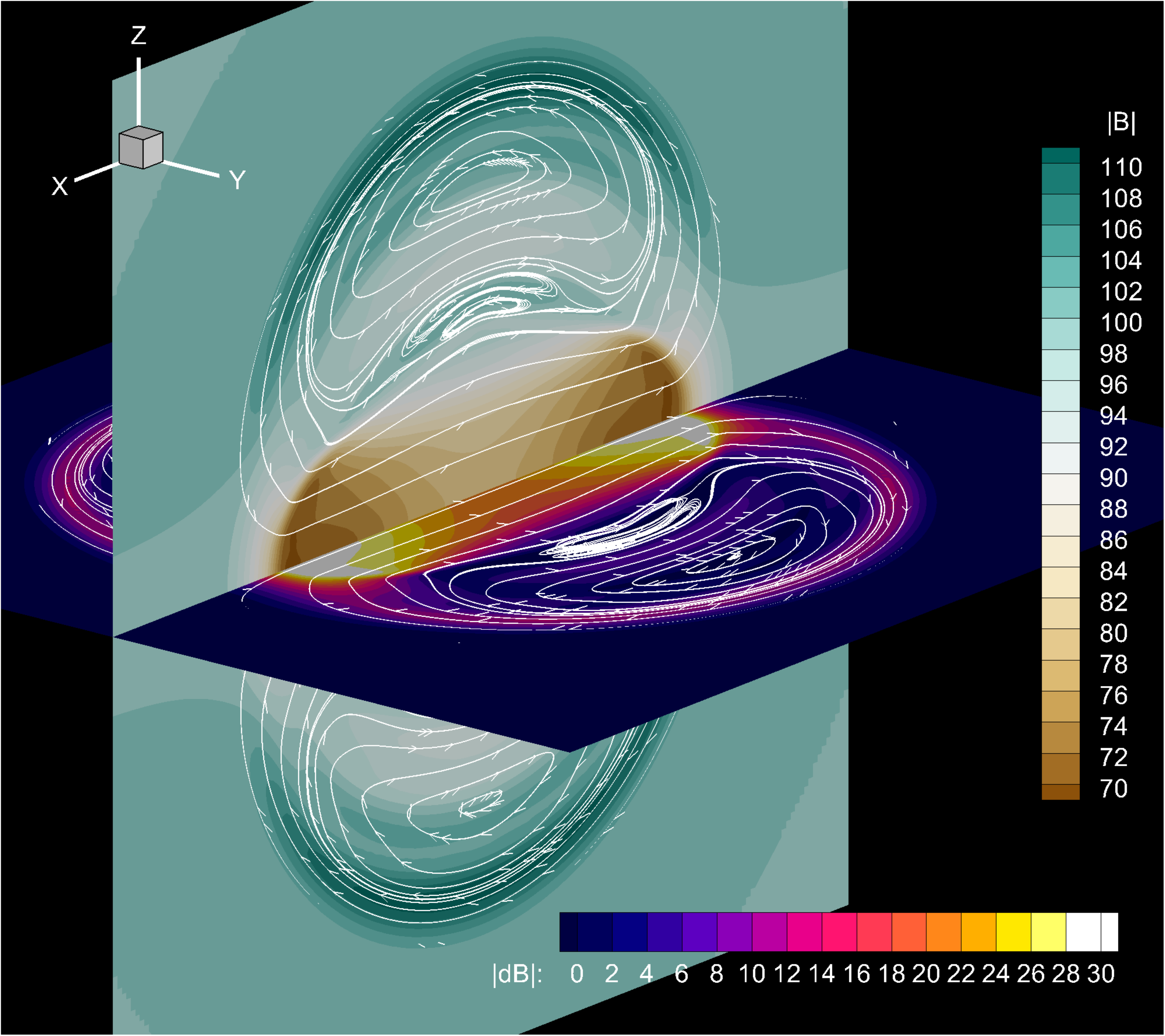} & 
		 	 \includegraphics[width=0.45\textwidth]{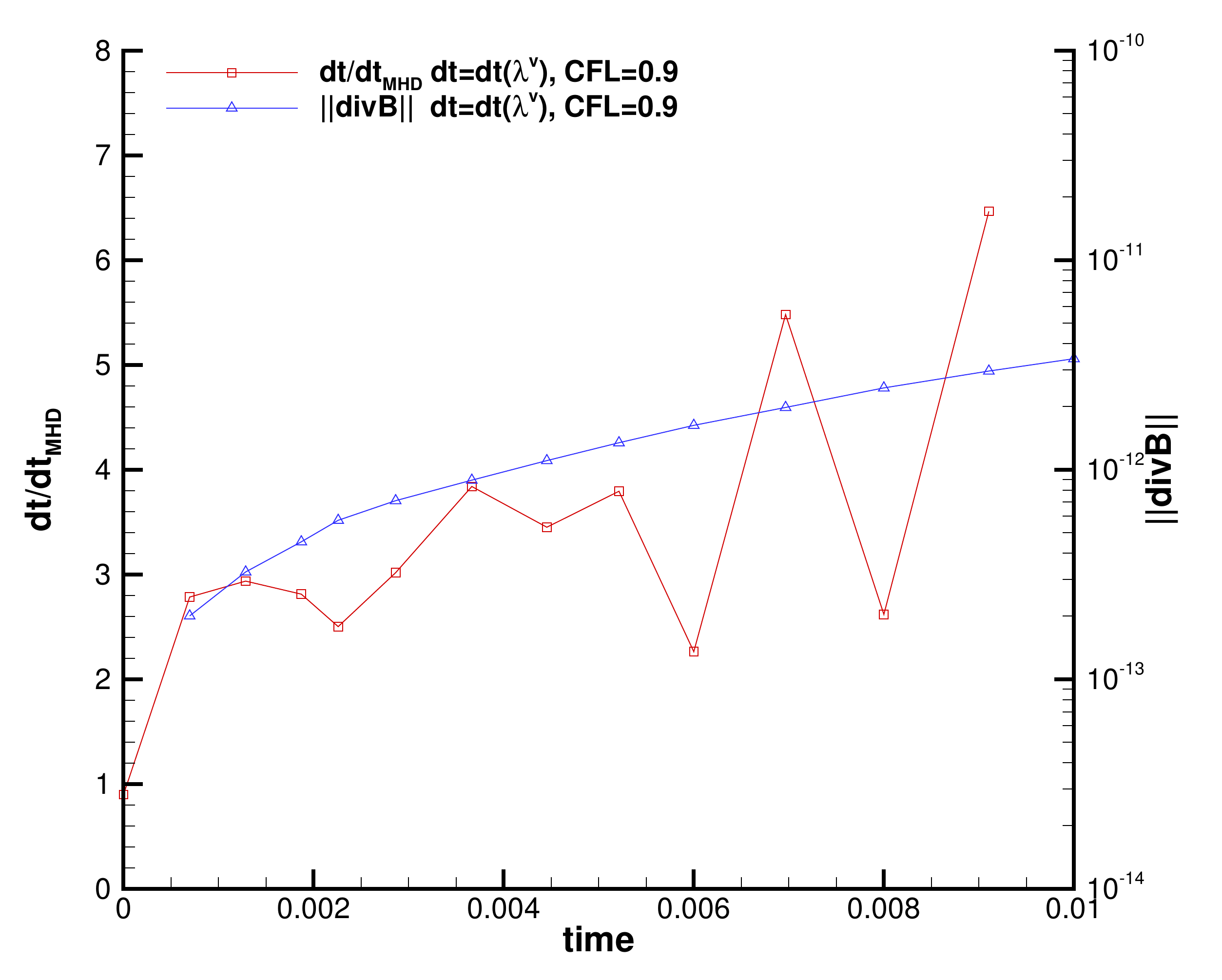}   
	\end{tabular}
\caption{Numerical results for the three-dimensional ideal Blast-Wave test at time $t=0.01$ obtained 
 with our new semi-implicit structure-preserving algorithm, interpolated along the 2D $\hat{xOz}$ plane. From left to right, on the first row the fluid density and the magnitude of the velocity field $|\v|$ are shown, while on the second row the magnitude of the magnetic field $|\B|$ and the logarithm of the absolute value of the divergence of the magnetic field $\log_{10}|\div\B|$ are plotted. On the third row,  a 3d plot shows the magnitude of the magnetic field $|\B|$ along the 2D  $\hat{xOz}$  plane, the magnitude of the variation of the magnetic field $|d\B|$ with respect to the initial configuration, i.e. $d\B:=\B(t,x)-\B(0,x)$, along the 2D $\hat{xOy}$ plane, as far as the field-lines of $d\B$ evaluated along these cut-planes; on the right, the time evolution of the effective computational Courant number, i.e. $\Delta t/ \Delta t(\lambda)$, is plotted, together with the $L_2$ norm of the divergence of the magnetic field $||\div\B||$, $\lambda$ being the full set of the MHD eigenvalues.
 } \label{fig:BW3d2}
\end{figure}

The three-dimensional blast wave problem can be built as the direct extension of what presented in the previous paragraph, see also \cite{PopovElizarova2015}, after re-defining the three-dimensional radial position $r=\sqrt{x^2+y^2+z^2}$.
The computational domain is a three-dimensional periodic box $(x,y,z)\in[-0.55,0.55]^3$ that is discretized with a mesh resolution of $\Delta x= \Delta y = \Delta z = 1/150$. The time parameters of the new  sonic-Alfvénic implicit solver are set to $\theta_{\text{B}} = 0.55$ and $\theta_{\text{p}} = 1$, with Courant parameter $\CFL=0.9$. In Fig.  \ref{fig:BW3d2} the numerical solution at the time slice $t=0.01$ is shown. The three-dimensional computed solution is interpolated along the 2D $\hat{xOz}$  plane, the plasma density, velocity and magnetic field amplitudes are shown, together with the $\div\B$ error in log-scale. As it is shown, in this three-dimensional test, the maximum absolute error in the divergence of the magnetic field was about $10^{-11}$.

A 3D view of the numerical solution is also shown, next to time-evolution of the estimated effective Courant number and the $L_2$ norm of the $\div\B$ error. Also in this case, the lower peaks in the effective Courant number $dt/dt_{\text{MHD}}$ are due to the automatic time step adjustment used to adapt the printing process to the prescribed output time.

Also in this case the computed solution agrees well with other published results, see \cite{PopovElizarova2015}.

\paragraph{3D Viscous and resistive Orszag-Tang vortex system.}

\begin{figure}[!htbp]
\centering 
\begin{tabular}{lr}
			 \includegraphics[width=0.45\textwidth]{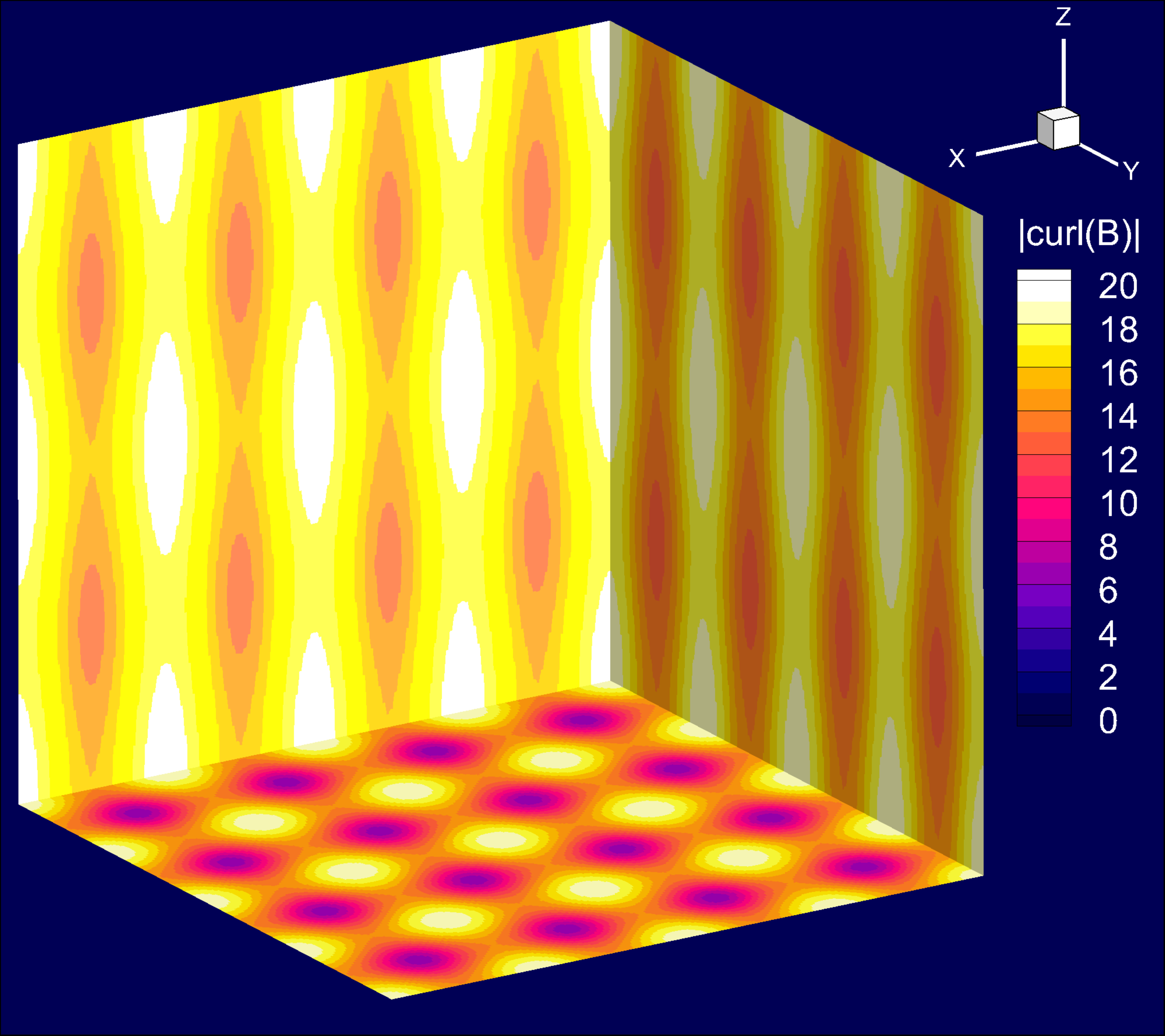}  &
			 \includegraphics[width=0.45\textwidth]{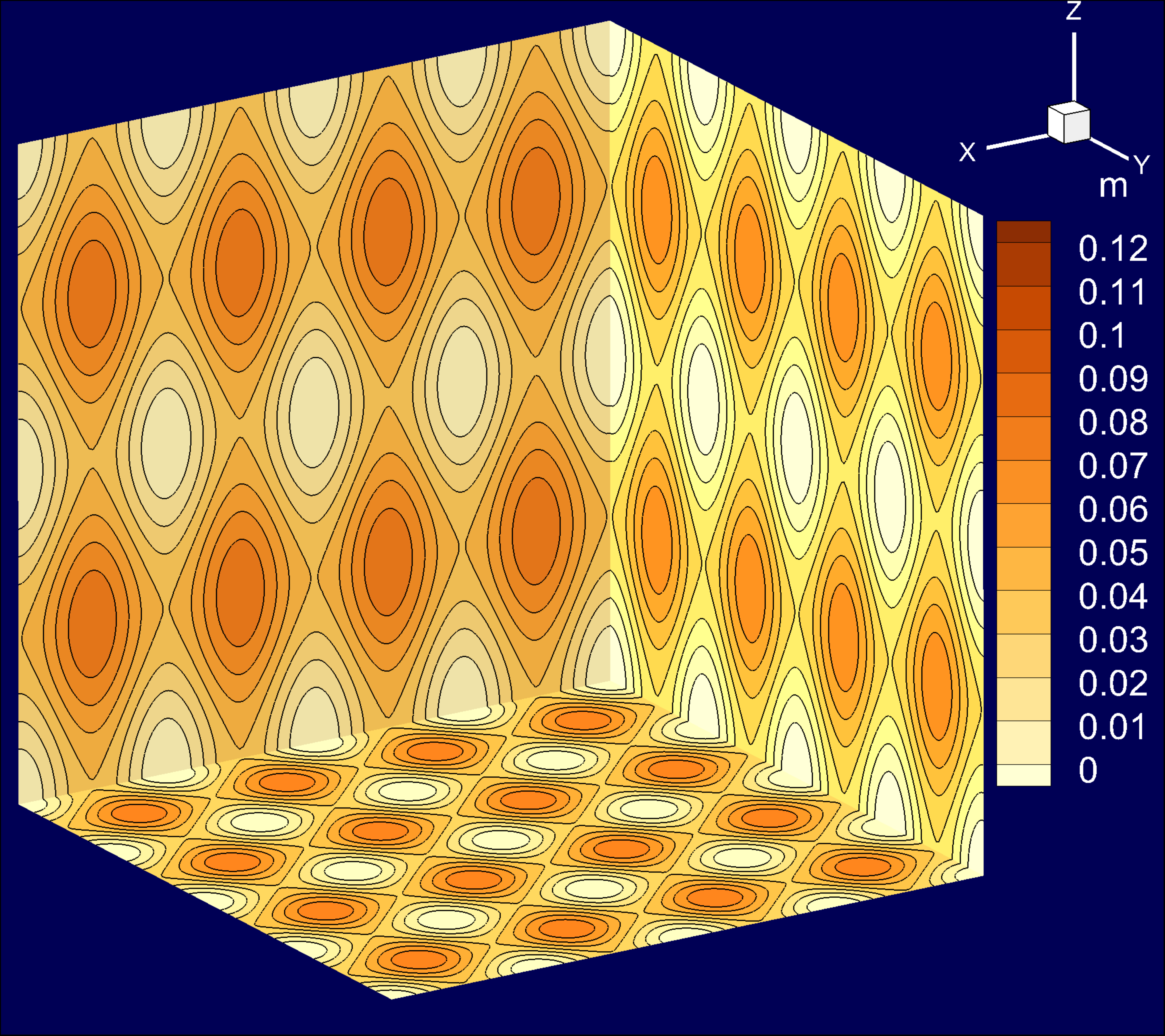}  \\
			 \includegraphics[width=0.45\textwidth]{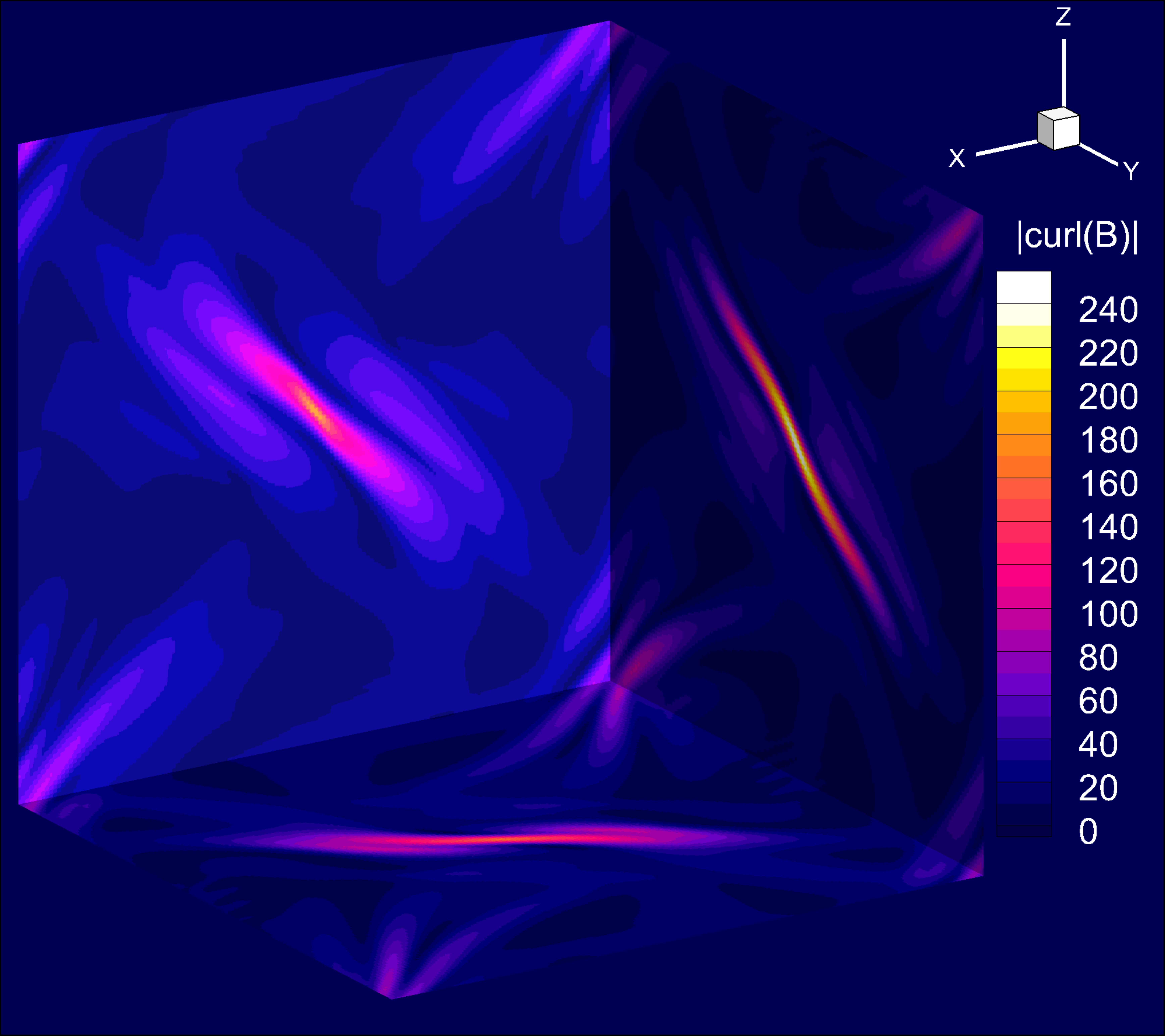}&
			 \includegraphics[width=0.45\textwidth]{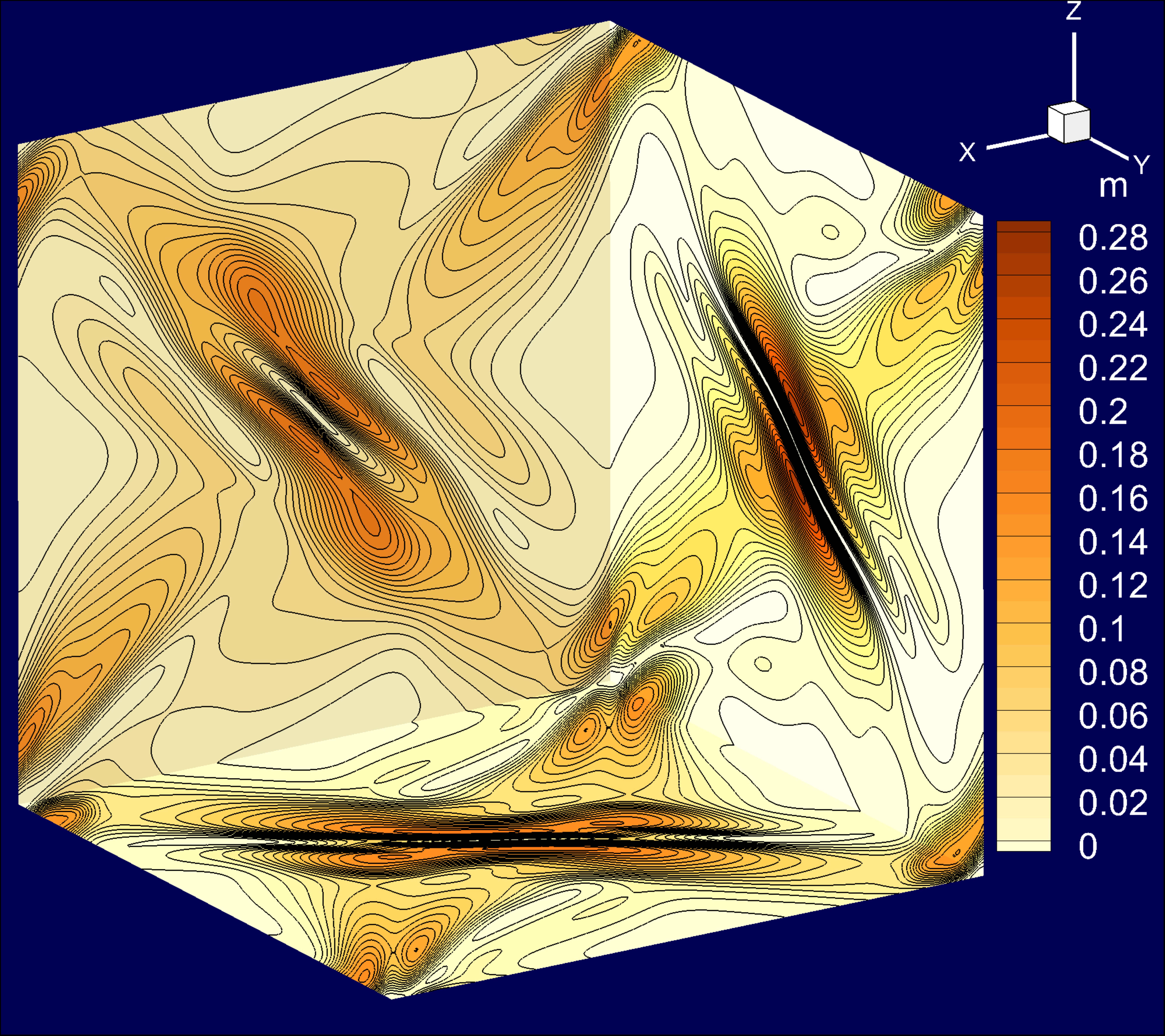}\\
			 \includegraphics[width=0.45\textwidth]{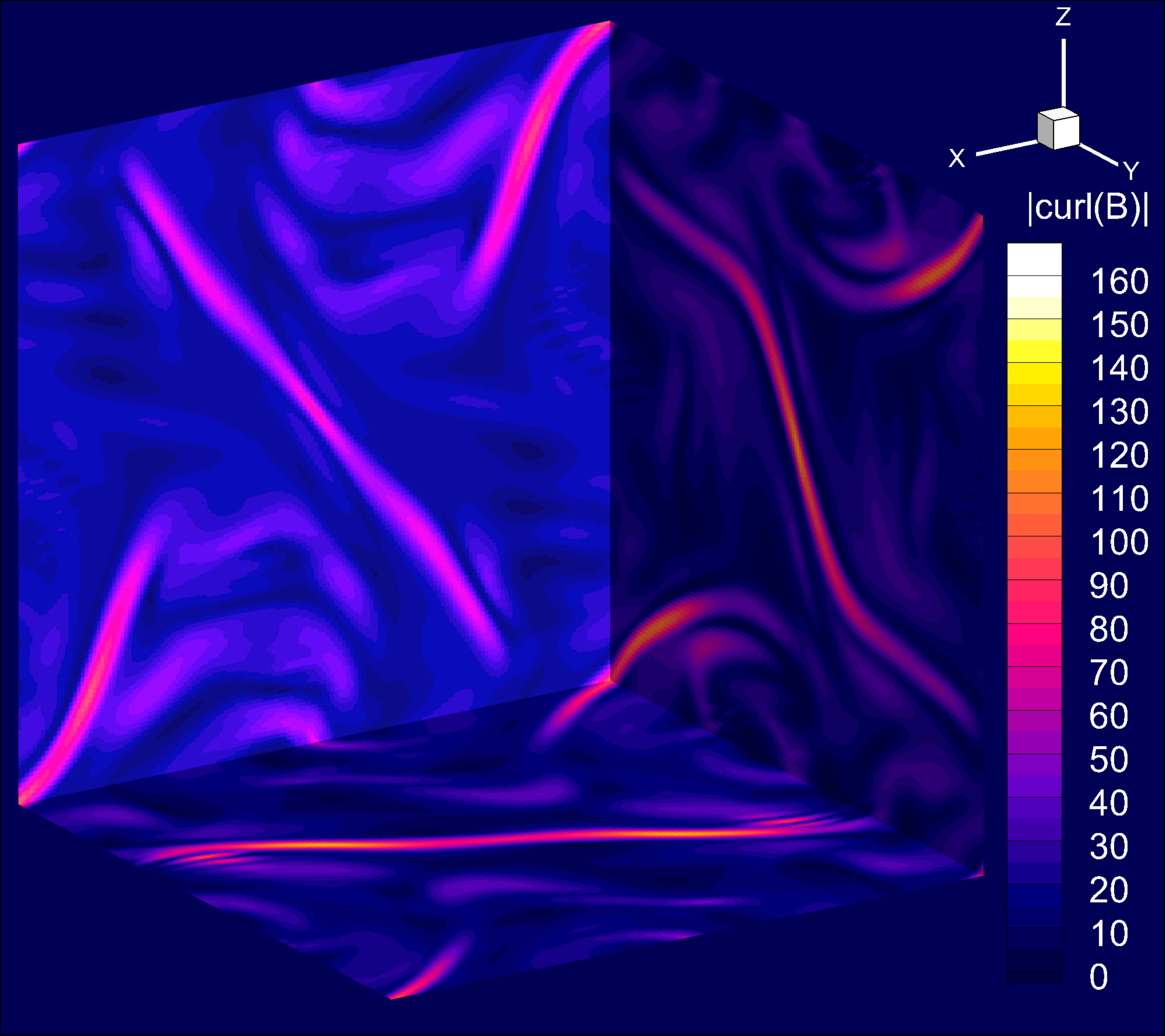}& 
			 \includegraphics[width=0.45\textwidth]{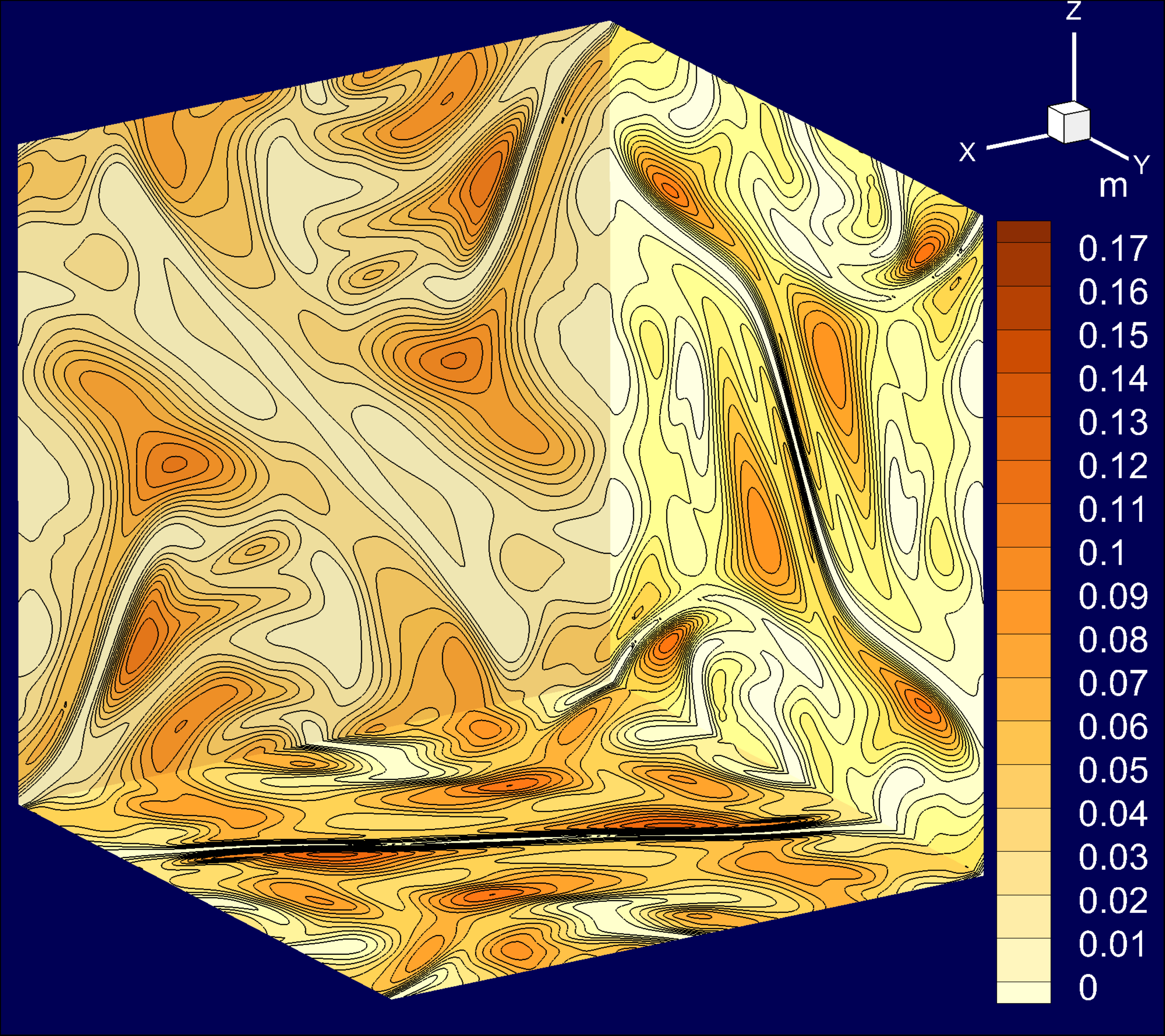}
	\end{tabular}
\caption{Numerical results for the three-dimensional viscous and resistive Orszag-Tang vortex system at times $t=0$, $0.25$ and $0.5$ obtained with a mesh resolution $\Delta x=\Delta y=\Delta z =1/150$ with our new semi-implicit structure-preserving algorithm,  from the top row to the bottom, respectively. At the left, the amplitude of the magnetic current $|\nabla \times \B|= |\mathbf{J}|$ is shown, while at the left the magnetic pressure $m=|\B|^2/(4\pi)$ is plotted. Here, the physical variables are interpolated along the three planes passing through the origin of the coordinate-axis, $\hat{xOy}$, $\hat{xOz}$ and $\hat{yOz}$.  
 } \label{fig:VROT3D}
\end{figure}

\begin{figure}[!htbp]
\centering 
\begin{tabular}{lr}
			 \includegraphics[width=0.45\textwidth]{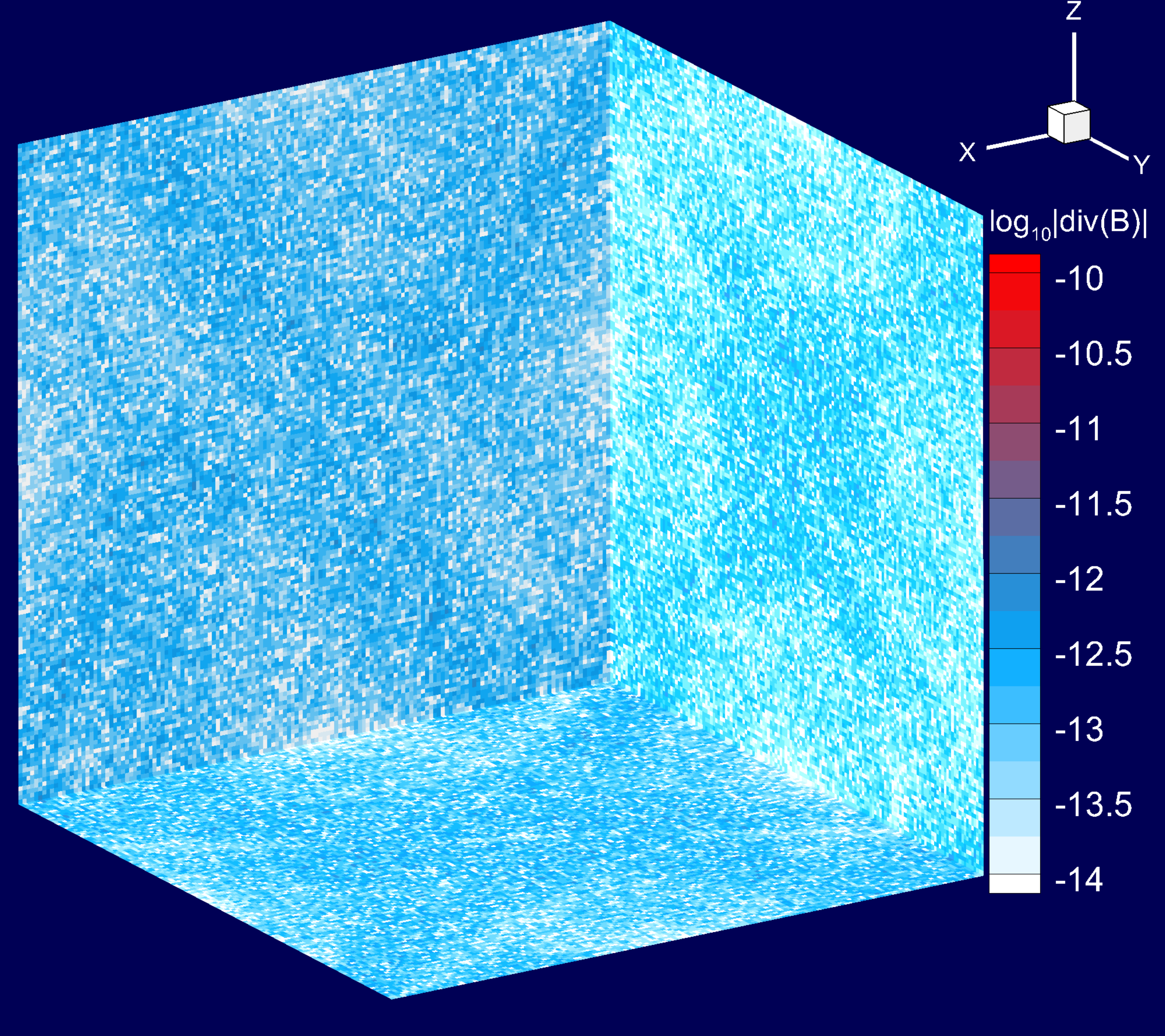}  &
			 \includegraphics[width=0.45\textwidth]{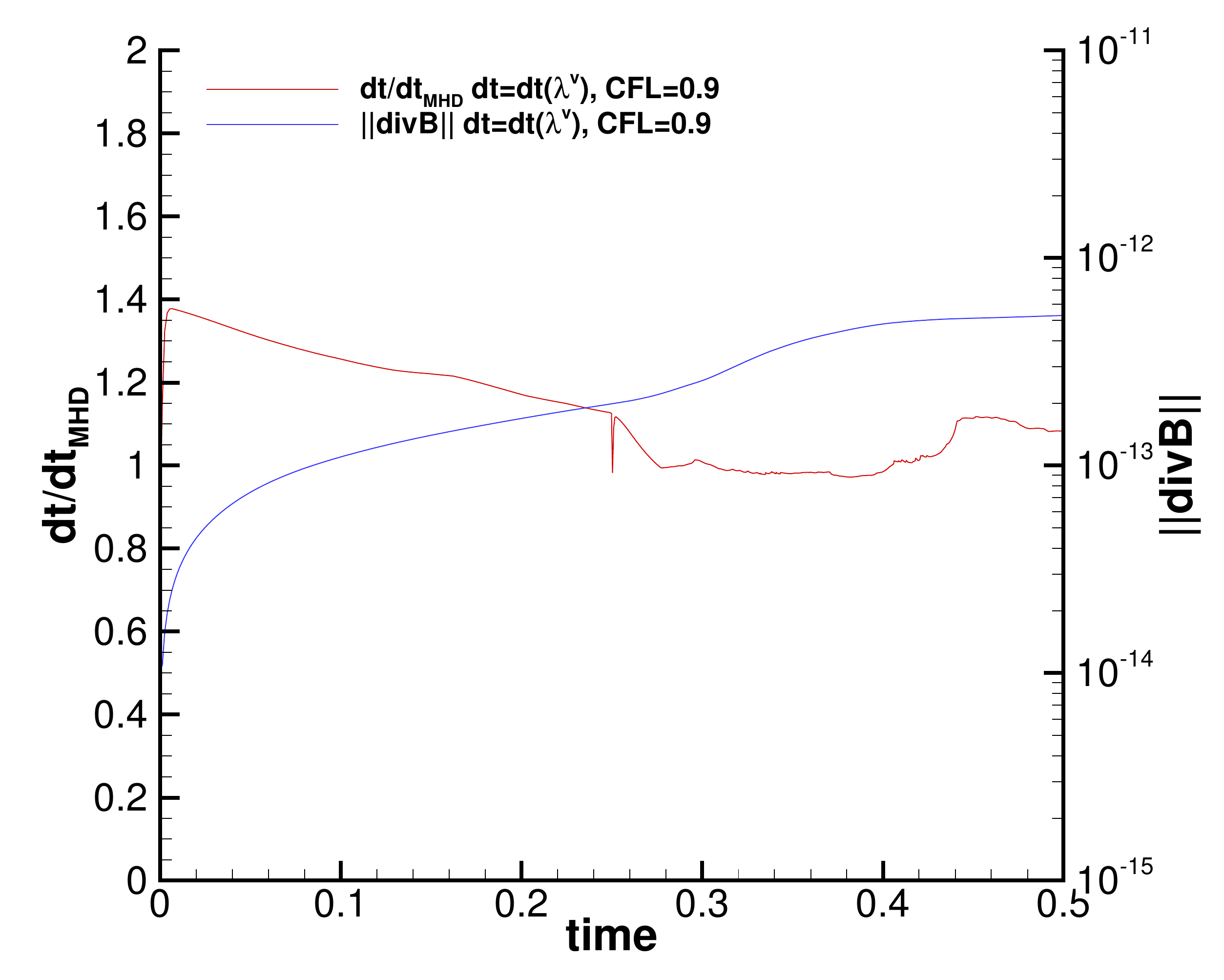}  
	\end{tabular}
\caption{Numerical results for the three-dimensional viscous and resistive Orszag-Tang vortex system obtained with a mesh resolution $\Delta x=\Delta y=\Delta z =1/150$ with our new semi-implicit structure-preserving algorithm. At the left, the divergence of the magnetic field $\div\B$ at time $t=0.5$ is interpolated along the three orthogonal 2D planes, $\hat{xOy}$, $\hat{xOz}$ and $\hat{yOz}$, while at the right the time evolution of the effective computational Courant number, i.e. $\Delta t/ \Delta t(\lambda)$, is plotted, together with the $L_2$ norm of the divergence of the magnetic field $||\div\B||$.
 } \label{fig:VROT3D2}
\end{figure}

Here, the viscous and resistive Orszag-Tang vortex system is set up in three-space dimensions, following \cite{BOHM2018}, and \cite{PopovElizarova2015} for ideal MHD. The plasma parameters are chosen to be $\gamma = \frac{5}{3}$, $\mu = 10^{-6}$, $\eta = 10^{-3}$, $c_v=1$ and a Prandtl number of $Pr=0.72$.
The physical and numerical domain is is a $150^3$ grid built on a two-dimensional unitary periodic-box $(x,y)\in[0,1]^2$, and a computational time-step evaluated with $\CFL=0.9$. The sonic-Alfvénic implicit solver is configured with $\theta_{\text{B}}=\theta_{\text{p}}=0.65$. The initial condition is
\myequation{l}
{
\rho = \frac{25}{36 \pi}, \\
\v  =   \left( - \sin\left(2\pi z\right), \sin \left(2\pi x \right), \sin \left(2\pi y \right)\right),\\
p =\frac{5}{12 \pi } ,\\
\B =   \left( -  \sin\left(2\pi z\right), \sin \left(4\pi x \right), \sin \left(4\pi y \right)\right).  
}
{\label{eq:VROT3D_ic}
} 
Similarly to the two dimensional dissipative or ideal version of the Orszag-Tang vortex test, this test is a valid choice  to validate the accuracy and robustness of the numerical algorithm in solving a highly nonlinear flow in three-space dimensions where steep gradients may easily arise. In such conditions, artificial negative pressures may be generated, leading to crash the simulations. In this work, the results have been successfully computed, and are compatible with other published results in literature. Fig.  \ref{fig:VROT3D} shows the computed solution for the magnitude of the current density $| \nabla \times \B|$ and the magnetic pressure $m= |\B|^2/16 \pi$ interpolated along the  2D $\hat{xOy}$, $\hat{xOy}$ and $\hat{yOz}$ planes. In Fig.  \ref{fig:VROT3D2} the divergence error $\div\B$ interpolated along the 2D planes is plotted, together with the time-evolution of the corresponding $L_2$ norm and the effective Courant number.


\section{Conclusion}
\label{sec:conclusion}

In this manuscript, a novel conservative second-order structure-preserving semi-implicit three-split  finite-volume/finite-difference scheme for the nonlinear viscous and resistive MHD equations has been presented. 
 A novel three-splitting method has been derived so that the purely hydrodynamic convective and  acoustic waves can be separated from the magnetic counterpart.  
In this way, an implicit time-discretization for the magneto-acoustic terms was possible, leading to a weaker CFL stability condition for the computational time-step.  
The here proposed semi-implicit solver is built in the form of a nested- iterative algorithm which leads to symmetric positive definite algebraic subsystems that can be efficiently solved with a standard matrix-free conjugate gradient method even without using any kind of preconditioner.
Another main novelty of this paper is the use of a completely different staggering of the electro-magnetic field quantities, with the electric field defined on the faces and the magnetic field defined on the edges of the primal control volumes. 
%
Numerical results are provided to show the main features of the algorithm: linear stability, in the sense of Lyapunov, has been shown numerically for a given constant equilibrium solution; thanks to the introduction of the implicit $\theta$ method, it is shown to be of up to second-order of accuracy for smooth solutions; the shock-capturing capabilities have been verified against a set of stringent MHD shock-tube problems; it has succeeded in robustness against the numerical simulation of non-trivial two- and three-dimensional test for both the ideal MHD and VRMHD; and finally, the effective gain in terms of accessible Courant number has been estimated for each test.   
These results encourage the author to find possible extensions of the here presented numerical method to higher orders of accuracy. In this direction, an extension to the high-order finite-element or finite-volume frameworks will be investigated. 
The author also believes that the here presented structure-preserving semi-implicit three-split  algorithm could become a good candidate for the construction of a structure-preserving hybrid FV/FE method, see e.g. \cite{BUSTO2018,BERMUDEZ2020,PODFVFE}, as well as a divergence-free limiter scheme for high-order accurate DG methods \cite{ADERDGVisc,Zanotti2015c,Fambri2020,Ioriatti2019,Ioriatti2020}. Last but not least, in the future we also plan to extend the new structure-preserving three-split method presented in this paper to the unified GPR model of continuum mechanics, see \cite{PeshRom2014,HPRmodel,HPRmodelMHD,FrontierADERGPR,SIGPR}.

\section*{Acknowledgments}
The author would like to thank Eric Sonnendrücker and Florian Hindenlang for the inspiring discussions on the topic.

%
 \section*{Funding}

The author did not receive funding from any external organization for the submitted work.

 \section*{Conflict of interest}

The author declares that he has no conflict of interest.

%
%
%


\bibliography{SIMHD}{}
\bibliographystyle{plain}


\end{document}